\documentclass[11pt, a4paper]{article}
\usepackage{latexsym}
\usepackage{amsmath}
\usepackage{color,soul}
\renewcommand{\epsilon}{\varepsilon}

\newcommand{\Var}{\mbox{\rm Var}}

\usepackage{multirow}

\renewcommand{\emph}[1]{#1}

\usepackage[margin=1.2in]{geometry}

\usepackage{amsmath,bm}
\usepackage{amssymb}

\usepackage{ntheorem}
\usepackage[ansinew]{inputenc}
\usepackage{graphicx}
\usepackage{epsfig}
\usepackage{dsfont}
\usepackage{bbm}
\usepackage{subfig}
\usepackage{url}

\usepackage{tikz}
\usetikzlibrary{shapes,arrows}
\usetikzlibrary{calc}
\usetikzlibrary{positioning,decorations.pathreplacing}

\usepackage{mdframed}

\usepackage{titling}
\thanksmarkseries{alph}

\usepackage{placeins}

\usepackage{float}
\newfloat{Algorithm}{t}{lop}

\usepackage[authoryear]{natbib}
\newtheorem{satz}{Theorem}[section]

\newtheorem{lem}[satz]{Lemma}

\newtheorem{kor}[satz]{Corollary}
\newtheorem{rem}[satz]{Remark}
\newtheorem{assumption}{Assumption} 

\def\3{\ss}

\newcommand{\E}{\mathbbm{E}}

\newcommand{\IN}{\mathbbm{N}}
\newcommand{\IZ}{\mathbbm{Z}}
\newcommand{\IR}{\mathbbm{R}}
\newcommand{\IP}{\mathbbm{P}}

\newcommand{\bea}{\begin{eqnarray*}}
\newcommand{\eea}{\end{eqnarray*}}
\newcommand{\be}{\begin{eqnarray}}
\newcommand{\ee}{\end{eqnarray}}

\newcommand{\ba}{\begin{array}}
\newcommand{\ea}{\end{array}}
\newcommand{\cum}{\text{\rm cum}}

\newcommand{\Cov}{\text{\rm Cov}}

\newcommand{\lstat}{\text{\rm loc.}}
\newcommand{\stat}{\text{\rm stat.}}

\def\3{\ss}

\newcommand{\pkg}[1]{{\normalfont\fontseries{b}\selectfont #1}}
\let\proglang=\textsf

\usepackage{hyperref}

\begin{document}

\title{Predictive, finite-sample model choice for time series\\ under stationarity and non-stationarity}
\author{
Tobias {\sc Kley}\thanks{Address: School of Mathematics (Faculty of Science), University Walk, Bristol BS8 1TW, UK. Email: \url{tobias.kley@bristol.ac.uk}. Partially supported by the Engineering and Physical Sciences Research Council grant no. EP/L014246/1.},
Philip {\sc Preu}\ss\thanks{Address: Department of Mathematics, Institute of Statistics, 44780 Bochum, Germany. Email: \url{philip.preuss@rub.de}. Supported by the Sonderforschungsbereich ``Statistical modelling of nonlinear dynamic processes" (SFB~823, Teilprojekt C1) of the Deutsche Forschungsgemeinschaft.},
Piotr {\sc Fryzlewicz}\thanks{Address:  Department of Statistics, Columbia House, Houghton Street, London, WC2A 2AE, UK. Email: \url{p.fryzlewicz@lse.ac.uk}. Supported by the Engineering and Physical Sciences Research Council grant no. EP/L014246/1.},\\
\vspace{-2mm} \\
$^{\rm a}$University of Bristol \vspace{1mm}\\
$^{\rm b}$Ruhr-Universit\" at Bochum \vspace{1mm}\\
$^{\rm c}$London School of Economics and Political Science \vspace{1mm}\\
}

\maketitle

\vspace{-10mm}
\begin{abstract}
	In statistical research there usually exists a choice between structurally simpler or more complex models. We argue that, even if a more complex, locally stationary time series model were true, then a simple, stationary time series model may be advantageous to work with under parameter uncertainty. We present a new model choice methodology, where one of two competing approaches is chosen based on its empirical, finite-sample performance with respect to prediction, in a manner that ensures interpretability.
	A rigorous, theoretical analysis of the procedure is provided. As an important side result we prove, for possibly diverging model order, that the localised Yule-Walker estimator is strongly, uniformly consistent under local stationarity.
	An \proglang{R} package, \pkg{forecastSNSTS}, is provided and used to apply the methodology to financial and meteorological data in empirical examples. We further provide an extensive simulation study and discuss when it is preferable to base forecasts on the more volatile time-varying estimates and when it is advantageous to forecast as if the data were from a stationary process, even though they might not be. 
\end{abstract}
covariance stationarity.

\textbf{MSC2010 subject classifications:} Primary 62M20; secondary 62M10.

\section{Introduction}

A well-trodden path in applied statistical research is to propose a model
believed to be a good approximation to the data-generating process,
and then to estimate the model parameters with a view to performing a specific
task, for example, prediction. However, even if the analyst were `lucky' and
chose the right model family, thereby reducing modelling bias, the resulting 
parameter estimators could be so variable that the selected model might well 
be sub-optimal from the point of view of the task in question. Choosing a slightly 
wrong model but with less variable parameter estimates
may well lead to superior performance in, for example, prediction. 
This effect is usually referred to as the bias-variance trade-off and it has frequently been discussed in the literature.
In this paper we explore how this unsurprising but interesting phenomenon could and should affect
model choice in the analysis of non-stationary time series.

Choosing between stationary and non-stationary modelling is, typically, an 
important step in the analysis of time series data. Stationarity, which
assumes that certain probabilistic properties of the time series model do
not evolve over time, is a key assumption in time series analysis, and several excellent
monographs focus on stationary modelling; see, e.\,g., \cite{Brillinger1975}, \cite{brodav1991} or \cite{priestley1981}.
However, in practice, many time series are deemed to be better-suited for
non-stationary modelling; this judgement can be based on diverse factors,
such as, for example, visual inspection, formal tests against stationarity, or the
observation that the data have been collected in a time-evolving environment
and therefore are unlikely to have come from a stationary model.

Early contributions to the literature of non-stationary time series are \cite{SubbaRao1970}, where the tvAR model was introduced, and \cite{Hallin1978}, who defined the tvARMA model. A general non-stationary time series framework was provided by \cite{Priestley1965}, who defined the evolutionary spectrum. A now particularly popular framework for the rigorous description of non-stationary time series models is that of local stationarity, in which the data are modelled locally as approximately stationary \citep{dah97,dahlhaus2012}.
We now illustrate the main idea of the paper using a simple example of a locally stationary time series model, the time-varying autoregressive model (of order~1)
\begin{equation*} \label{timevaryingar}
X_{t,T}=a(t/T)X_{t-1,T}+Z_t, \quad t=1,...,T,
\end{equation*}
with $T$ denoting the sample size, $a: [0,1] \rightarrow (-1,1)$ being some suitable function and $Z_t$ being an i.\,i.\,d.\ sequence with mean zero and variance one. Typically, to forecast future observations, one would require an estimate of $a(1)$, see e.\,g.\ \cite{chen2010}. Before constructing a suitable estimator, some analysts would wish to test if $a$ was indeed time-varying, and there exist a vast amount of techniques to validate the assumption of a constant second-order structure in this framework; see \cite{sacneu00}, \cite{pap09}, \cite{dwisub10}, \cite{pap10}, \cite{detprevet10}, \cite{Nason2013}, \cite{prevetdet12} or \cite{vogtdet2013}. If the process was found to be non-stationary, it would be tempting to estimate $a(1)$ by a localised estimate based on the most recent observations of $X_{t,T}$. This localisation would most likely reduce the bias of the estimator if the true dependency structure was indeed time-varying, but also increase its variance. However, if, for example, the function $a$ was varying only slowly over time, this 
estimation procedure might result in sub-optimal estimation from the point of view of the mean squared prediction error, yielding inferior forecasts compared to the classical stationary AR(1) model. This would be particularly likely if the test of stationarity employed at the start was not constructed with the same performance measure in mind (i.\,e., mean squared prediction error) and was therefore `detached' from the task in question (i.\,e., prediction). One of the findings of this paper is that even if the function $a$ varied over time, one should in some cases treat it as constant in order to obtain smaller prediction errors, or in other words, `prefer the wrong model' from the point of view of prediction.

The main aim of this paper is to propose an alternative model choice methodology in time series analysis that avoids the pitfalls of the above-mentioned process of testing followed by model choice.
More precisely, our work has the following objectives:
\begin{itemize}
\item
To propose a generic procedure for \emph{finite-sample} model choice which avoids the path of hypothesis testing but instead chooses the model that offers better empirical finite-sample performance in terms of prediction on a validation set, with associated performance guarantees for the test set of yet unobserved data. Although the procedure is proposed and analysed theoretically in the framework of choice between stationarity and local stationarity and in the context of prediction, the procedure is applicable more generally whenever a decision needs to be made between two competing approaches, and can therefore be viewed as model- and problem-free. At the end of Section~\ref{sec:remarksprec}, we provide two examples of other situations in which the general principle of our procedure can be applied.
\item
To suggest `rules of thumb' indicating when the (wrong) stationary model may be preferred in a time-varying, locally stationary situation from the point of view of forecasting; and when a time-varying model should be preferred.
\end{itemize}

Our procedure validates and puts on a solid footing the possibly counter-intuitive observation that it is sometimes beneficial to choose the `wrong' (but possibly simpler) model in time series analysis, if that model relies on more reliable estimators of its parameters than the right (but possibly more complex) model. While we stop short of conveying the message that simplicity in time series should always be preferred, part of our aim is to draw time series analysts' attention to the fact that particularly complex time series models may well appear attractive on first glance as they have the potential to capture features of the data well, but on the other hand can be so hard to estimate that this makes them inferior to simple and easy-to-estimate alternative models, even if the latter are wrong.

We now briefly describe related recent literature. The work of \cite{xiatong}, who, while discussing time series
prediction,
select the model based on the minimisation of up to $m$-step ahead prediction errors (rather than the usual 
1-step ahead ones) also appears to carry the general message that different models may be preferred for the same 
dataset depending on the task in question, or, in the language of the authors, on the `features to be matched'.
Besides similarities in this general outlook, our model-fitting methodology and the context in which it is proposed
are entirely different. Forecasting in the presence of structural changes is a widely studied topic in the
econometrics literature, see e.\,g.\ the comprehensive review by \cite{Rossi2013} and the references therein.
In particular, \cite{giraitis2012} also use the minimisation of the 1-step ahead prediction error as a basis
for model choice under non-stationarity, but, unlike us, do not consider the question of how this may lead to
the preference for the `wrong' model in finite samples.
\cite{DasPolitis2017} apply the model-free prediction principle of \cite{Politis2015} in the context of locally stationary time series and construct 1-step-ahead point and interval predictors.

Instead of pursuing the cross-validation approach, \cite{McDonaldEtAl2012} evaluate the upper bound on the generalisation 
error in time series forecasting, and use its heuristically estimated version to guide model choice. We note,
however, that this approach requires the estimation of some possibly difficult to estimate parameters, unlike cross-validation-based
approaches. 
The empirical mean squared prediction error (MSPE) which we will employ in our method is closely related to the population MSPE under parameter uncertainty. The strand of literature discussing this population quantity includes \cite{bai79} and \cite{rei80}, where approximating expressions were derived for stationary VAR time series.
For locally stationary tvMA($\infty$) processes, \cite{pal13} discuss optimal $h$-step ahead forecasting, in terms of the true model characteristics. Yet, they do not take parameter uncertainty into account.

While the main question we are concerned with is whether a stationary or a time-varying autoregressive model should be used for prediction, a nested question is what order the stationary or non-stationary model should have. 
Traditionally, order selection is done via minimisation of an information criterion, see, e.\,g., \cite{brodav1991}, p.\,301.
\cite{zhakor15} develop an adaptive criterion for model selection based on predictive risk.
\cite{Akaike1969} introduced the Final Prediction Error (FPE) as a figure of merit for a potential predictor and adopts a decision theoretic approach, called the minimum FPE procedure, where the predictor with the best FPE is chosen. In practice, the decision is then based on an estimate of the FPE. In \cite{Akaike1970b} a theoretical basis of the procedure is provided. \cite{pensan07} derive and compare MSPE for univariate and multivariate predictors when the parameters are known. They then define and estimate a criterion (a measure of predictability) to choose between these two prediction options. Their approach is similar to ours in spirit, but, firstly, it chooses between univariate and multivariate models while we consider stationary and non-stationary models and, secondly, their methodology works with the population MSPE (which moves the focus away from the observed data to the postulated model), while we work with the corresponding empirical quantity directly. This difference in approaching the problem also holds for another, more general class of special-purpose-criteria: the focused information criteria (FIC), which were introduced in \cite{clahjo2003}. The FIC methodology with the focus on choosing the model best suited for prediction was then applied in the field of time series analysis in \cite{clae2007}, where the best AR($p$) model for prediction is chosen, in \cite{roha2011}, where the best ARMA($p$,$q$) model for this purpose is chosen, and in \cite{brogal08}, where models for volatility forecasting are chosen. The idea of the FIC is that the model which minimises the \emph{asymptotic} MSPE is the best one and the FIC is then based on an estimator of that asymptotic MSPE. Contrary to this, our approach is based on the \emph{empirical} MSPE directly, which we believe to be the more relevant quantity in many applications. Contrary to the FIC which is based on the large-sample theory of the estimators involved, we provide finite-sample exponential bounds that imply a performance guarantee for our method. This approach can be advantageous, when it is preferred that the model choice also depends on the size of the sample, which in our view should be a natural requirement.

Our paper is organised as follows. In Section~\ref{sec:mexample} we provide a simple motivating example.
In Sections~\ref{sec:prec} and~\ref{sec:remarksprec} we introduce and comment on our
new time series model choice methodology.
The statistical properties of our procedure are discussed in Section~\ref{theorysectionar1}, where also the performance guarantee (Theorem~\ref{thm:main3}) is provided. The results of a simulation study and the analysis of three empirical examples can be found in Sections~\ref{simulationsection1} and~\ref{data_example}.
In Section~\ref{sec:propertiesYW} we discuss statistical properties of the local Yule-Walker estimator and prove its strong uniform consistency under local stationarity (Corollary~\ref{kor:propertiesYW}). We conclude with a summary in Section~\ref{sec:conclusion}. Proofs, technical details, additional tables and figures from the simulations section are gathered in Appendices~\ref{app:main}--\ref{a:extraSim}. Note that Appendices~\ref{app:proofAuxLem}--\ref{a:extraSim} are only available in the arXiv'ed version of the manuscript \citep{KleyEtAl2019}.

\section{Motivating example}\label{sec:mexample}

We consider the time-varying autoregressive (tvAR) model of order 2:
\begin{equation*} 
X_{t,T} = a_1(t/T) X_{t-1,T} + a_2(t/T) X_{t-2,T} + Z_t, \quad t=1,...,T,
\end{equation*}
where $a_1(u) := 0.15 + 0.15 u$, $a_2(u) := 0.25 - 0.15u$, and $Z_t$ is Gaussian white noise.
$X_{t,T}$ is a non-stationary process which lies in the locally stationary class
of \cite{dah97}.
We will now compare different forecasting procedures for $X_{0.9T, T}$, where $T \in \{50, 500, 5000\}$. The predictor that minimises the mean squared prediction error is given by
\[\hat X_{0.9 T,T}^{\rm oracle} 
		= 0.285 X_{0.9 T-1, T} + 0.115 X_{0.9 T-2, T}.\]
Yet, since in practice the underlying model is unknown, the analyst needs to
\vspace*{-0.3cm}\begin{itemize}
	\item[(1)] make assumptions regarding the model, and
	\item[(2)] estimate the assumed model's parameters.
\end{itemize}
\vspace*{-0.3cm}For the purpose of this illustration, we discuss four possible models. In the first two models we falsely assume that the data were stationary and model $X_{t,T}$ to satisfy a traditional, autoregressive (AR) equation.
\vspace*{-0.3cm}\begin{itemize}
	\item In the first of the two cases we assume an AR(1) model and
	\item in the second case we assume the model to be an AR(2) model.
\end{itemize}
\vspace*{-0.3cm}We further, discuss cases 3--4, where the correct class of models (tvAR) is assumed. Yet,
\vspace*{-0.3cm}\begin{itemize}
	\item in case three, we falsely assume a tvAR(1) model, before
	\item in case four, we correctly assume the model to be a tvAR(2) model.
\end{itemize}
\vspace*{-0.3cm}Note that the true model, the tvAR(2) model, is the most complex one of the four choices. In each of the models we estimate the parameters by solving the empirical Yule-Walker equations. In the case of the tvAR models we localise by using the segment $X_{0.9T - N, T}, \ldots X_{0.9 T-1, T}$. In the case of the traditional, stationary AR models we use all available observations $X_{1, T}, \ldots X_{0.9 T-1, T}$. Details on the estimation are deferred to Section~\ref{sec:prec}.

Denoting the localised Yule-Walker estimates of order $1$ by $\hat a_{N,T}^{(1)}(0.9T-1)$ and the ones of order $2$ by $\hat a_{1;N,T}^{(2)}(0.9T-1)$ and $\hat a_{2;N,T}^{(2)}(0.9T-1)$ we obtain the predictors
\begin{equation*}
\begin{split}
	\hat X_{0.9 T,T}^{1,N} & := \hat a_{N,T}^{(1)}(0.9T-1) X_{0.9T-1,T}, \\
	\hat X_{0.9 T,T}^{2,N} & := \hat a_{1;N,T}^{(2)}(0.9T-1) X_{0.9T-1,T} + \hat a_{2;N,T}^{(2)}(0.9T-1) X_{0.9T-2,T},
\end{split}
\end{equation*}
where $\hat X_{0.9 T,T}^{1,N}$ corresponds to the models of order $1$ and $\hat X_{0.9 T,T}^{2,N}$ corresponds to the models of order $2$. The segment length $N$ will be chosen as $0.9 T-1$ in the AR models and strictly smaller than this in the tvAR models.

\begin{figure}
    \begin{center}
     \includegraphics[width=\linewidth]{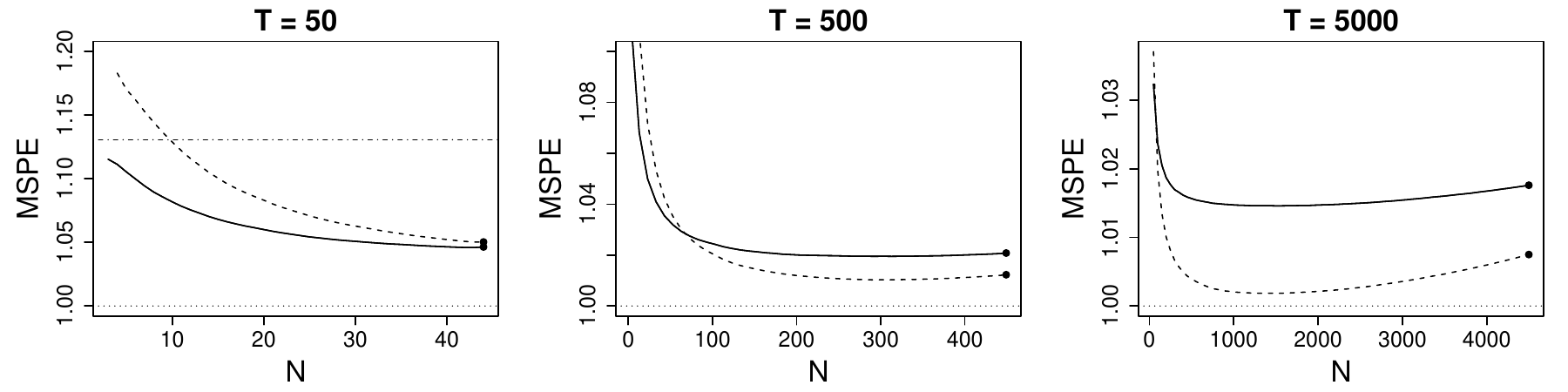}
    \end{center}
    \caption{Mean squared prediction errors (MSPEs) for forecasting $X_{0.9T, T}$ with predictors $\hat X_{0.9 T,T}^{1,N}$ and $\hat X_{0.9 T,T}^{2,N}$ associated with tvAR(1) and tvAR(2) modelling of the data, where $N$ varies. Left, middle and right column correspond to $T=50$, $T=500$ or $T=5000$, respectively. The solid lines corresponds to $\E (\hat X_{0.9 T,T}^{1,N} - X_{0.9T,T})^2$, the dashed line corresponds to $\E (\hat X_{0.9 T,T}^{2,N} - X_{0.9T,T})^2$. The endpoints of each line indicate the MSPEs of the predictor associated with the stationary AR(1) and AR(2) models. The dotted horizontal lines (at level 1.00) indicate the MSPE of the oracle predictor. The dashed-dotted line (approximately at level 1.13) indicates the variance of $X_{0.9T,T}$.
    \label{AR1example_mspe}}
\end{figure}

In Figure~\ref{AR1example_mspe}, we observe that the predictors associated with the simpler, stationary AR model perform better than or similarly well as the predictors associated with the more complex, locally stationary tvAR model if $T=50$ or $T=500$. If $T=5000$ the predictor associated with the locally stationary tvAR model performs visibly better in terms of its MSPE when the segment size $N$ is chosen appropriately. In conclusion, this example illustrates how it can be advantageous to assume a wrong, but structurally simpler model when only a short time series is available. In particular, the model chosen should depend on the task at hand (here: prediction) and on the amount of data available. For $T=50$ the best result is obtained by assuming the AR(1) model which is the simplest of the four candidates. When $T=500$ the more complex AR(2) model becomes advantageous. Note that this model is more complex than the AR(1) model and thus provides a better approximation to the true tvAR(2) mechanism, but is still simplifying, because it does not take the time-varying characteristics into account at all. Only when even more data (here: $T=5000$) are available, then the variability of the parameter estimates of the tvAR(2) model is small enough not to overshadow the modelling bias, which in this example is rather small.

Obviously, the bias-variance trade-off is at work here, which is well-known but interestingly, to our knowledge, has previously been unexplored in the important context of stationary versus non-stationary modelling for prediction. The observation to be made here, thus, is that finding the `right' model may not always be a suitable way of proceeding when it comes to the prediction of future observations. We point out that this observation was made in other contexts of time series analysis. For example, basic exponential smoothing is a widely used forecasting and trend extrapolation technique, and although it is well-known that it corresponds to standard Box-Jenkins forecasting in the ARIMA$(0,1,1)$ model, it is also frequently used for
data that does not follow it.

This paper investigates the question of what is the best model in terms of forecasting performance in the context of the choice between stationarity and non-stationarity.
To ask this question explicitly instead of applying a test for stationarity is 
important since the smallest sample size $T$ needed to reject the 
null hypothesis of stationarity may be smaller than the sample size needed to 
obtain improvement in terms of our task of interest, namely forecasting. In the 
following section, we will elaborate more on this question. Further,
in Section~\ref{simulationsection1}, we see, as results of a simulation study,
under which conditions using the true model is advantageous and when it can become disadvantageous.

\section{When (not) to use locally stationary models under local stationarity: the new model choice methodology} \label{ignore.stat}

\subsection{Precise description of the procedure}
\label{sec:prec}

We work in the framework of general locally stationary time series (a rigorous definition is deferred to Section~\ref{theorysectionar1}), in which the available data is a finite stretch $X_{1,T}, \ldots, X_{T,T}$ from an array $(X_{t,T})_{t \in \IZ, T \in \IN^{*}}$ of random variables with mean zero and finite variances.
Our aim is to determine a linear predictor for the unobserved $X_{T+h,T}$ from the observed $X_{1,T}, \ldots, X_{T,T}$.

Our proposal is to compare candidate $h$-step ahead predictors in terms of their empirical mean squared prediction error and choose the predictor with the best forecasting performance. To this end, we proceed as follows:

\textbf{Step 1.} Separate the final $2 m$ observations from the $T$ available observations. The observations with indices $M_0 := \{1,\ldots,T-2m\}$, $M_1 := \{T-2m+1,\ldots,T-m\}$ and $M_2 := \{T-m+1, \ldots, T\}$ will be referred to as the \emph{training set}, \emph{first validation set} and \emph{second validation set}, respectively. The set of unobserved data with the indices $M_3 := \{T+1, \ldots, T+m\}$ will be referred to as the \emph{test set}. The size $m$ of the separated sets will be small in comparison to the sample size $T$ (and hence also to the training set). Comments on why we require two distinct validation set are deferred to Section~\ref{sec:remarksprec}.

\textbf{Step 2.} Compute the linear $1$-step ahead prediction coefficients
	\begin{equation}
	\begin{split}
		\label{eqn:YW}
			\hat a_{N,T}^{(p)}(t) & := \big( \hat\Gamma_{N,T}^{(p)}(t) \big)^{-1} \hat\gamma_{N,T}^{(p)}(t)
				= \big( \hat a_{1;N,T}^{(p)}(t), \ldots, \hat a_{p;N,T}^{(p)}(t) \big)',
			\end{split}
	\end{equation}
	($a'$ denotes the transposed vector $a$) for $t+h \in M_1 \cup M_2$, $p = 1, \ldots, \max\mathcal{P}$, and $N \in \mathcal{N}$,
		\begin{equation}\label{eqn:hatGamma_t}
			\hat\Gamma_{N,T}^{(p)}(t) := \big[ \hat \gamma_{i-j;N,T}(t) \big]_{i,j = 1, \ldots, p}, \quad
			\hat \gamma_{N,T}^{(p)}(t) := \big( \hat \gamma_{1;N,T}(t), \ldots, \hat \gamma_{p;N,T}(t) \big)' \end{equation}
	and
	\begin{equation}\label{def:acf}
		\hat \gamma_{k;N,T}(t) := \frac{1}{N} \sum_{\ell=t-N+|k|+1}^{t} X_{\ell-|k|,T} X_{\ell,T}, \quad k = 0, \ldots, \max\mathcal{P}.
	\end{equation}
	 The set of possible model orders $\mathcal{P} \subset \{0,1,\ldots,\min\mathcal{N}-1\}$, with $\mathcal{P} \neq \emptyset$ and $\max\mathcal{P} \geq 1$, and the set of possible segment lengths $\mathcal{N} \subset \{\max\mathcal{P}+1, \ldots, T-2m-h+1\}$, with $\mathcal{N} \neq \emptyset$, are parameters to be specified by the user. Further comments on how they are to be chosen are deferred to Section~\ref{sec:remarksprec}.
	
\begin{figure}[t]
	\begin{center}
		\begin{tikzpicture}[
    every node/.style = {align=center},
          Line/.style = {-angle 90, shorten >=2pt},
          short_Line/.style = {-angle 90, shorten >=5pt},
    Brace/.style args = {#1}{semithick, decorate, decoration={brace,#1,raise=2pt,
                             pre=moveto,pre length=2pt,post=moveto,post length=2pt,}},
            ys/.style = {yshift=#1}
                    ]
\linespread{0.8}                    
\coordinate (a) at (0,0);
\coordinate[right=70mm of a]    (b);
\coordinate[right=20mm of b]    (c);
\coordinate[right=20mm of c]    (d);
\coordinate[right=20mm of d]    (e);
\coordinate[right=10mm of e]    (f);
\coordinate[left=3mm of b]    (a3);
\coordinate[left=15mm of a3]    (a2);
\coordinate[left=20mm of a2]    (a1);
\coordinate[below=12mm of a3]    (a3_L1);
\coordinate[below=12mm of a2]    (a2_L1);
\coordinate[below=20mm of a3]    (a3_L2);
\coordinate[below=20mm of a1]    (a1_L2);
\coordinate[below=28mm of a3]    (a3_L3);
\coordinate[below=28mm of a]    (a_L3);
\coordinate[right=4mm of b]    (b1);
\coordinate[left=7mm of c]    (c1);

\draw[Line] (a) -- (f) node[below left] {time};
\draw[Brace] (a) -- node[above=3pt] {Training\\ set $M_0$} (b);
\draw[Brace] (b) -- node[above=3pt] {Validation\\ set $M_1$} (c);
\draw[Brace] (c) -- node[above=3pt] {Validation\\ set $M_2$} (d);
\draw[Brace] (d) -- node[above=3pt] {Test\\ set $M_3$} (e);
\draw[Line] ([ys=-5mm] a1) node[below] {$t-N+1$} -- (a1);
\draw[Line] ([ys=-5mm] a2) node[below] {$t-p+1$} -- (a2);
\draw[Line] ([ys=-5.2mm] a3) node[below] {$t$} -- (a3);
\draw[Line] ([ys=-5mm] b1) node[below] {$t+h$} -- (b1);
\draw[Line] ([ys=-5.8mm] c1) node[below] {$s_1$} -- (c1);
\draw[short_Line] ([ys=12mm] a) node[above] {$1$} -- (a);
\draw[short_Line] ([ys=12mm] b) node[above] {$T-2m$} -- (b);
\draw[short_Line] ([ys=12mm] c) node[above] {$T-m$} -- (c);
\draw[short_Line] ([ys=12mm] d) node[above] {$T$} -- (d);
\draw[short_Line] ([ys=12mm] e) node[above] {$T+m$} -- (e);
\draw[Brace=mirror] (a2_L1) -- node[below=3pt] {for forecasts} (a3_L1);
\draw[Brace=mirror] (a1_L2) -- node[below=3pt] {for locally stationary coefficients} (a3_L2);
\draw[Brace=mirror] (a_L3) -- node[below=3pt] {for stationary coefficients} (a3_L3);
\end{tikzpicture}
	\end{center}
	\caption{Time line to illustrate the sets $M_j$, $j=0,1,2,3$ and relations of $t$, $p$, $h$, $m$ and $N$. Downward pointing arrows indicate first or last indices of the four sets. The three upward pointing arrows from the left and braces indicate the indices of the observations used to compute the forecasting coefficients and the observations that are weighted by the coefficients to constitute the forecasts. The upward pointing arrow second from the right indicates the index of an observation for which the forecast is computed. The rightmost upward pointing arrow indicates $s_1 := T-m-h$, the observation up to which the MSPEs can be evaluated; cf. eq.~\eqref{def:f}.}
	\label{fig:time_line}
\end{figure}
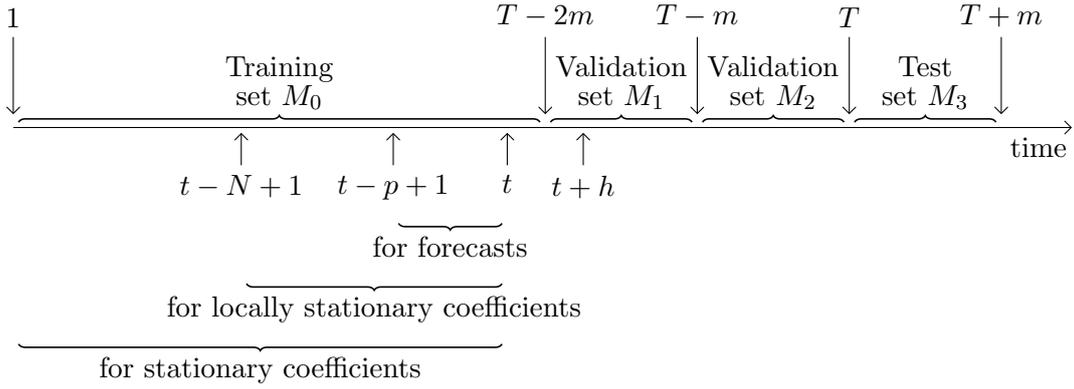

\textbf{Step 3.} Compute the linear $h$-step ahead prediction coefficients
\begin{equation}\label{def:vhat}
\begin{split}
	\big( \hat v_{N,T}^{(p,h)}(t) \big)'
				&  := \big(
				\hat v_{1;N,T}^{(p,h)}(t), \cdots, \hat v_{p;N,T}^{(p,h)}(t) \big)
				:= e'_1 \big( \hat A_{N,T}^{(p)}(t) \big)^h
				:= e'_1 \big( e_1 \big( \hat a_{N,T}^{(p)}(t) \big)' + H \big)^h,
\end{split}
\end{equation}
where $\hat a_{N,T}^{(p)}(t)$ is defined in~\eqref{eqn:YW}, $e_1$ denotes the first canonical unity vector of dimension $p$ and $H$ denotes a $p \times p$ Jordan block with all eigenvalues equal to zero; cf. equation~\eqref{eqn:def_e_H}, in the appendix. Comments on an equivalent, recursive definition are provided in Section~\ref{sec:remarksprec}.
Next, define $f^{{\lstat}}_{t,h;0,N} := 0$, $f^{\stat}_{t,h;0} := 0$ and, for $p \in \mathcal{P}\setminus\{0\}$ and $N \in \mathcal{N}$, compute
\begin{equation}
\label{eqn:lsEstim}
\begin{split}
	f^{\lstat}_{t,h;p,N}
		& := e'_1 \big( \hat A_{N,T}^{(p)}(t) \big)^h (X_t, X_{t-1}, \ldots, X_{t-p+1})'
		  := \sum_{i=1}^p \hat v^{(p,h)}_{i;N,T}(t) X_{t-i+1,T},
\end{split}
\end{equation}
\begin{equation}
\label{eqn:sEstim}
	\begin{split}
	f^{{\stat}}_{t,h;p}
		& := e'_1 \big( \hat A_{t,T}^{(p)}(t) \big)^h (X_t, X_{t-1}, \ldots, X_{t-p+1})'
		  := \sum_{i=1}^p \hat v^{(p,h)}_{i;t,T}(t) X_{t-i+1,T}
	\end{split}
\end{equation}
In Figure~\ref{fig:time_line}, a time line is shown that illustrates the relation of the sets $M_j$, $j=0,1,2,3$ and the quantities $t$, $p$, and $N$.

\textbf{Step 4.} Amongst predictors~\eqref{eqn:lsEstim} select $f^{\lstat}_{t,h} :=  f^{\lstat}_{t,h; \hat p_{\lstat}, \hat N_{\lstat}}$, with
\[(\hat p_{\lstat}, \hat N_{\lstat}) := \arg\min_{\substack{p \in \mathcal{P}\\ N \in \mathcal{N}}} \ \sum_{t+h \in M_1} \big( X_{t+h,T} - f^{\lstat}_{t,h;p,N} \big)^2,\]
and, amongst predictors~\eqref{eqn:sEstim} select $f^{{\stat}}_{t,h} :=  f^{{\stat}}_{t,h;\hat p_{\stat}}$, with
\[\hat p_{\stat} := \arg\min_{p \in \mathcal{P}} \ \sum_{t+h \in M_1} \big( X_{t+h,T} - f^{{\stat}}_{t,h;p} \big)^2.\]
Note that $f^{\lstat}_{t,h}$ and $f^{{\stat}}_{t,h}$ are the forecasts of type~\eqref{eqn:lsEstim} and~\eqref{eqn:sEstim} that minimise the empirical MSPE (on $M_1$) within the classes of tvAR and AR models of orders $p \in \mathcal{P}$, respectively.

\textbf{Step 5.} Use $f^{\lstat}_{t,h}$ as $h$-step ahead forecast of $X_{t+h}$, with $t+h > T$, if
\begin{equation} \label{decisionrule}
\hat R_{T,j}(h) := \frac{{\rm MSPE}_{T,j}^{{\stat}}(h)}{{\rm MSPE}_{T,j}^{\lstat}(h)} \geq 1 + \delta
\end{equation}
holds for $j=2$, and $f^{{\stat}}_{t,h}$ otherwise, where
\begin{equation}
\label{eqn:empMSPE}
{\rm MSPE}_{T,j}^{*}(h) := \frac{1}{m}\sum_{t+h \in M_j} (X_{t+h,T}- f_{t,h}^{*})^2,
\end{equation}
with $*$ indicating the corresponding model (we write `loc.' for the locally stationary approach and `stat.' for the stationary model) and $\delta \geq 0$ is a parameter by which the user of the procedure specifies which degree of superiority of the more complex procedure is required before it is preferred over the simpler alternative (cf. the end of Section~\ref{sec:remarksprec}). 

By Theorem~\ref{thm:main3} we have that, with an appropriately chosen $\delta$, the decision rule of type~\eqref{decisionrule} will, with high probability, prefer the same models for forecasting observations from the second validation set ($j=2$) and the test set ($j=3$).

\subsection{Remarks on the procedure}
\label{sec:remarksprec}

Some further explanations regarding the procedure are in order now. Our comments are organised according to the steps of the previous section.

\textbf{Step 1.} While it is common practice to separate one validation set when tuning the model parameters to avoid over-fitting, we require two such sets. This is necessary, because we would otherwise compare candidates in an unbalanced situation where $|\mathcal{P}|$ stationary predictors compete with $|\mathcal{N}| \times |\mathcal{P}|$ locally stationary ones. In our procedure, where we first choose the hyper-parameters by minimising the mean squared error on the first validation set and then choose between the two model classes by minimisation of the mean squared error on the second validation set, we achieve a fairer competition of the two model classes.

\textbf{Step 2.} The coefficients~\eqref{eqn:YW} are estimates for the coefficient functions $a_1(t/T), \ldots, a_p(t/T)$ if the data follows the tvAR($p$) model
\be \label{timevaryingarp}
X_{t,T}=\sum_{j=1}^p a_j(t/T) X_{t-j,T}+\sigma(t/T) Z_t, \quad t=1,...,T,
\ee
(see, for example, \cite{dahlgir1998}).
Recall that $Z_t$ is usually assumed to be white noise and that $X_{t,T}$ is non-stationary if at least one of the functions~$a_j$, $j=1,\ldots,p$, or~$\sigma$ is non-constant. A recursive algorithm to estimate the parameters was described and analysed in~\cite{MoulinesEtAl2005}.

We are interested in \emph{linear forecasts} that will perform well for time series possessing a \emph{general dependency structure}. The tvAR($p$) model~\eqref{timevaryingarp} is a natural choice to approximate the linear dynamics of the observed, non-stationary time series, because in this model the coefficient functions at time $t/T$ coincide with the $1$-step ahead prediction coefficients (of order $p$) which define the best linear predictor. In Section~\ref{sec:propertiesYW}, we show that $\hat a_{N,T}^{(p)}(t)$ from Step~2 can be used as estimates for the 1-step ahead linear prediction coefficients
\[ \tilde a_{T}^{(p)}(t) := \arg\min_{a = (a_1, \ldots, a_p)' \in \IR^p}
		\E \Big[\Big( X_{t,T} - \sum_{j=1}^p a_j X_{t-j,T} \Big)^2\Big],\]
also when the observations do not satisfy~\eqref{timevaryingarp}. A forecasting procedure derived within the tvAR($p$) model can therefore be expected to behave reasonably, irrespective of whether the tvAR($p$) model is true or \emph{just} an approximation to the truth. Note that we use the tvAR($p$) model to approximate the dynamic structure of the data in Section~\ref{sec:remarksprec} and most of our examples in Section~\ref{simulationsection1} are of this kind, but we do \emph{not} assume that the data actually satisfies it.

\textbf{Step 3.} Linear $h$-step ahead predictors can either be obtained by iterating model equation~\eqref{timevaryingarp} or by using a separate model for each $h$ in which the indices of the sum on the right hand side run from $j=h, \ldots, p+h-1$. These approaches have been referred to as the \emph{plug-in method} and the \emph{direct method}, respectively. A comparison of the two approaches can, for example, be found in \cite{bha96}, where results for a class of linear, stationary processes were derived. We employ the plug-in method.

The coefficients $\hat v_{N,T}^{(p,h)}(t)$ defined in~\eqref{def:vhat} can be computed efficiently via the recursion:
\begin{equation*}
	\begin{split}
		\hat v_{i;N,T}^{(p,1)}(t) & := \hat a_{i;N,T}^{(p)}(t), \quad i=1, \ldots, p, \\
		\hat v_{i;N,T}^{(p,\eta)}(t)
			& := \hat a_{i;N,T}^{(p)}(t) \hat v_{1;N,T}^{(p,\eta-1)}(t)
				+ \hat v_{i+1;N,T}^{(p,\eta-1)}(t) I\{i \leq p-1\}, \quad \eta = 2, 3, \ldots, h.
	\end{split}
\end{equation*}

From the previous comments it can be seen how the predictors $f_{t,h;p,N}^{\lstat}$ and $f_{t,h;p}^{\stat}$ relate to the choice of modelling the time series' dynamics by a tvAR($p$) or AR($p$) model, respectively. In each of these model classes, increasing the order $p$ will give a better approximation of the dynamics, but increase the complexity of the model, and make it more difficult to deal with under parameter uncertainty.

The parameters $\mathcal{P}$ and $\mathcal{N}$ are sets of integers to be chosen by the user. The choice should depend on~$T$. $p \in \mathcal{P}$ determines the order of the tvAR($p$) model that is used to approximate the dynamics.
$N \in \mathcal{N}$ determines the degree of locality in the estimation of the coefficients. The parameters $p \in \mathcal{P}$ and $N \in \mathcal{N}$ will influence the degree of bias and variance of the predictor. Our selection mechanism will balance them implicitly.

\textit{Traditional choice of $N$.} It is obvious that the variance of the estimator can decrease when a larger segment is used, but that the non-stationarity will potentially inflict an additional bias that increases with $N$. Under the condition that $N/T + T/N^2 = o(1)$, \cite{dahlgir1998} derive asymptotic expansions for the local Yule-Walker estimator's bias and variance for a centred sample. It follows from their results, that for the one-sided sample we require for forecasting, $N$ should be chosen at the order of $T^{2/3}$, with the constant depending on the second derivatives of the true model quantities, which are unknown and difficult to estimate. The choice of $\mathcal{N}$ should thus, ideally, be such that $N \asymp T^{2/3}$, for all $N \in \mathcal{N}$. In practice, since the true model parameters are unknown, this rate provides very little guidance to the user of the method. We recommend, though, to adhere to two facts: the upper and lower bound of $\mathcal{N}$ should be bounded away from $0$ and $T$, respectively. In other words, we recommend to choose $\mathcal{N}$ with $\min \mathcal{N}$ \emph{large enough}, for the performance guarantee to be valid (cf. Theorem~\ref{thm:main3}) and $\max \mathcal{N}$ being \emph{substantially smaller} than $T$, to ensure that there is a clear boundary between the locally stationary and the stationary approach. \cite{RichterDahlhaus2017} propose to adaptively choose a bandwidth for local M-estimators by minimising a cross-validation functional.

\textit{Traditional choice of $p$.} As described in the beginning of this section we use the tvAR($p$) model to \emph{approximate} the dynamic structure of the data. Intuitively, we have that the larger the order $p$ the better the approximation to the true dynamic structure. In opposition to the previously discussed question of how to choose the segment length $N$, we here have that a smaller $p$ will inflict a modelling bias, while a larger $p$ will typically be accompanied by an inflation of the variance of the estimation, because it implies that more parameters need to be estimated. Traditionally, the model order is chosen by minimizing information criteria as for example AIC or BIC. \cite{clae2007} propose to use a version of the focused information criterion (FIC, see \cite{clahjo2003}) to select the model order of a stationary AR($p$) model optimal with respect to forecasting when the true model is known to be AR($\infty$). However, as mentioned in the introduction, the FIC-based methods employs an estimator of the asymptotic MSPE, while our approach is based on the empirical MSPE, which facilitates our focus on the finite sample performance.

\textbf{Step 4 and 5.} Our procedure performs two stages of selections. Firstly (in Step 4), it selects the model order $p$ and, for the locally stationary approach, the segment length~$N$ by comparing predictors within each class of models under consideration (i.\,e., time-varying or non-time-varying autoregressive models). The parameters $p$ and $N$ are chosen such that the empirical MSPE (predicting observations from $M_1$) is minimised. Secondly (in Step 5), a final competition of the winners is performed to select among the two classes of models. The procedure that minimises the empirical MSPE (predicting $M_2$) is selected and used for forecasting of the test set ($M_3$).
In our theoretical analysis of the next section (see, in particular, Theorem~\ref{thm:main3}) we show that the proposed procedure will, with high probability, choose the same class of models on the validation as on the test set, implying that the procedure with the best empirical performance will be selected.

\textit{The parameter $\delta$.} By introduction of the parameter $\delta \geq 0$ the user is given additional control over which model the procedure prefers. In the simplest case, $\delta = 0$, this reduces to a straight choice between the two model classes, whereby the time-varying model is chosen if it performs better or equally well. Choosing $\delta > 0$ introduces penalisation against the choice of more complex models. In this case, the predictor derived from the more complex, locally stationary model is only chosen if it performs at least $\delta \cdot 100\%$ better than the one derived from the simpler, stationary model.

\textbf{Generalisations.} Besides linear predictions for stationary or locally stationary time series models, the general principle of our method can also be applied in many other situation. To illustrate this, we outline two examples below.

\textit{Non-linear predictions with neural networks.} In this scenario we either choose a neural network trained with the $N$ most recent observations (i.\,e., loc.) or with all available observations (i.\,e., stat.). To this end proceed as follows: with the available data partitioned as described in Step 1, consider a range of candidate networks (with different network topologies) suitable for forecasting the observations from the first validation set. Train them either with the $N$ most recent observations (loc.) or with all available observations (stat.). After first choosing the network for which we see the smallest MSPE on the first validation set within each class (loc. or stat.), we then choose that class for which the winner from the previous step obtains the best performance on the second validation set.

\textit{Predictors obtained from locally stationary or long-memory time series models.} Long-range dependence and non-stationarity can lead to the same stylised facts in financial time series; cf. \cite{mikosch2004}. Choosing a test, for example from \cite{PreussVetter2013}, to distinguish between the two model classes seems to be a sensible approach, but this might not lead to the best choice if the aim is to choose a model for the purpose of prediction. With the model choice methodology in this paper, one proceeds as follows: first, fit a set of long-range dependence models and a range of locally stationary models, use the implied predictors to forecast the observations from the first validation set and choose the model with the best forecasting performance within each model class. Then, choose between the two model classes by comparing the winning models within each class with respect to their empirical performance in predicting the observations from the second validation set.

\subsection{Performance guarantee: theoretical result for the general case} \label{theorysectionar1}

In this section, we establish theoretical results
that will facilitate our analysis of the model choice suggested
by decision rule~\eqref{decisionrule}. We show 
that the probability of choosing different models on the 
validation and the test set decays to zero at an exponential 
rate, which can be viewed as a \emph{performance guarantee} of
our model choice methodology.

To rigorously prove the results, some definitions and assumptions are in order. Throughout this paper, we work with the doubly indexed process $(X_{t,T})_{t \in \IZ, T \in \IN^*}$. The first index (i.\,e., $t$) refers to the time. The second index (i.\,e., $T$) indicates how well the covariance structure of $(X_{t,T})_{t \in \IZ}$ can be approximated locally by the autocovariance function of a stationary process. We will assume that, for large $T$, segments of observations $X_{t,T}$ with their indices $t \approx u T$ are approximately weakly stationary. The parameter $u$ is continuous and often referred to as the rescaled time. If the index $T$ coincides with the number of observations in a time series, then $u \in [0,1]$ (cf. \cite{Dahlhaus1996}). This restriction is not necessary and in fact, because we will consider $m$ unobservables (to be forecast) in addition to the $T$ observations (available at the time when the forecasting is done) it is more convenient to allow $u > 1$, as was also done by \cite{rousan18}.

The following definitions from \cite{rousan18} are required for our assumptions.
For an array $(X_{t,T})_{t \in \IZ, T \in \IN^*}$ with finite variances, the \emph{time-varying covariance function} is defined for all $t \in \IZ, T\in\IN^{*}$ and $k \in \IZ$ as
	\begin{equation}\label{eqn:time_varying_covariance_function}
		\tilde\gamma_{k,T}(t) = \Cov\left(X_{t,T}, X_{t-k,T}\right). 
	\end{equation}

A \emph{local spectral density}
$f$ is a $\IR^2\to\IR_{+}$ function, $(2\pi)$-periodic and locally integrable with
respect to the second variable. The \emph{local covariance function}~$\gamma$
associated with  the time-varying spectral density $f$ is defined for $(k,u) \in \IZ \times \IR$ by
	\begin{equation}\label{equation:local_covariance_function}
		\gamma_{k}(u) = \int\limits_{-\pi}^{\pi}\exp\left({\rm i} k \lambda\right)f(u,\lambda){\rm d}\lambda. 
	\end{equation}

The first five assumptions are specific to the kind of data we may apply our result to.

\begin{assumption}[Local stationarity, \cite{rousan18}]
\label{a:loc_stat}
	Let the array of random variables $(X_{t,T})_{t\in\IZ, T\in\IN^{*}}$ have finite
	variances. We say that $(X_{t,T})_{t\in\IZ,T\in\IN^{*}}$ is
	\emph{locally stationary with local spectral density $f$} if the time-varying
	covariance function of $(X_{t,T})_{t\in\IZ, T\in\IN^{*}}$ and the
	local covariance function associated with $f$ satisfy
	\begin{equation} \label{eqn:def:locally_stationary}
		\left|\tilde\gamma_{k,T}\left(t\right)-\gamma_k\left(\frac{t}{T}\right)\right| \leq \frac{C}{T},
	\end{equation}
	where $C \geq 0$  is a constant.
\end{assumption}

\begin{assumption}[Geometrically $\alpha$-mixing]
\label{a:mixing}
There exist constants $K > 0$ and $\rho > 1$ such that, for every $n \in \IN$,
\begin{equation}
	\label{eqn:mixing}
	\alpha(n) := \sup_{T \in \IN^*} \sup_{t \in \IZ} \sup_{A \in \sigma(X_{s,T} : s \leq t)} \sup_{B \in \sigma(X_{s,T} : s \geq t+n)}
	\big| \IP(A \cap B) - \IP(A) \IP(B) \big| \leq K \rho^{-n}.
\end{equation}
\end{assumption}

\begin{assumption} \label{a:sd_bounded}
The local spectral density $f$ is bounded from above and below:
\begin{equation}
	\label{eqn:sd_bounded}
	0 < m_f \leq f(u,\lambda) \leq M_f.
\end{equation}
\end{assumption}

\begin{assumption} \label{a:sd_der_bounded}
The local spectral density $f$ is continuously differentiable with respect to the first argument and the partial derivative is uniformly bounded. More precisely, assume the existence of $M'_f \geq 0$ such that
\begin{equation}
	\label{eqn:sd_der_bounded}
	\Big| \frac{\partial}{\partial u} f(u,\lambda) \Big| \leq M'_f.
\end{equation}
\end{assumption}

\begin{assumption}[Bernstein-type condition] \label{a:MomentCond}
There exist $c > 0$ and $d \geq 1/2$, such that
\begin{equation*}
	\E |X_{t,T}|^k \leq c^{k-2} (k!)^{d} \Var(X_{t,T}) \quad t \in \IZ; \ k=2,3,\ldots
\end{equation*}
\end{assumption}

The assumptions are reasonable and non-restrictive in the sense that many popular and widely used time series models (e.\,g., a wide range of tvARMA models) satisfy the full set of assumptions.
The notion of local stationarity we impose (Assumption~\ref{a:loc_stat}) goes beyond that of locally stationary linear processes and, in particular, we do not require the data to be tvAR. Assumption~\ref{a:loc_stat} is satisfied for second order stationary process (then we have $C=0$), the general (linear) locally stationary process introduced by \cite{Dahlhaus1996}, but also non-linear processes as elaborated by~\cite{rousan18}. Assumption~\ref{a:mixing} is satisfied for a broad class of (linear and non-linear) locally stationary time series models; see, for example, \cite{fryzlewicz2011} or \cite{vogt2012}. Assumptions~\ref{a:sd_bounded} and~\ref{a:sd_der_bounded} can be verified by considering the local spectral density when it is given explicitly. For example, in the tvAR model that we used to motivate our prediction approach in Section~\ref{sec:remarksprec}, see~\eqref{timevaryingarp}, and as examples in Section~\ref{simulationsection1}, the local spectral density and local covariances are naturally those of the stationary AR process when the parameter $u$ of the coefficient functions is chosen as $t/T$. We will refer to these AR processes as the \emph{tangent processes} of the tvAR process. Similar assumptions with respect to the local spectral density are common in the literature; cf.~\cite{dah97}. Processes with sub-Gaussian marginal distributions satisfy Assumption~\ref{a:MomentCond}; cf. Lemma~\ref{lem:MomentCond} in the appendix.
We recall, from Section~\ref{sec:remarksprec}, that the tvAR($p$) model is used to approximate the linear dynamic structure of the data, but that we do \emph{not} assume that the data actually satisfies it. Thus our results apply in a more general context. We require Assumptions~\ref{a:mixing} and~\ref{a:MomentCond} to prove that the probabilities in our results decay at an exponential rate. 

As a consequence of Assumptions~\ref{a:loc_stat} and~\ref{a:sd_bounded}, we have
\begin{equation}
\label{eqn:cons:sd_bounded}
	\pi m_f \leq \sigma_{t,T}^2 := \Var(X_{t,T}) \leq 3 \pi M_f, \quad \text{for all $T \geq \frac{C}{\pi m_f}$.}
\end{equation}
Further, by Assumption~\ref{a:sd_der_bounded} and Leibniz's integral rule, we have that $\gamma'_k(u)$ exists and has the following form
\begin{equation}\label{eqn:gammaPrime}
	\gamma'_k(u) := \frac{\partial}{\partial u} \gamma_k(u)
	= \frac{\partial}{\partial u} \int\limits_{-\pi}^{\pi}\exp\left({\rm i}\ell\lambda\right)f(u,\lambda){\rm d}\lambda
	= \int\limits_{-\pi}^{\pi}\exp\left({\rm i}\ell\lambda\right) \frac{\partial}{\partial u} f(u,\lambda){\rm d}\lambda,
	\end{equation}
which, in particular, implies that $|\gamma'_k(u)|\leq 2\pi M'_f$.

Assumptions~\ref{a:m_and_N} and~\ref{a:T}, which we state below, are more specific to our procedure. They concern minimum requirements for the size $m$ of the validation sets and the minimum segment size $\min\mathcal{N}$ which are used to compute the forecast as well as the number $T$ of observations required to be available at the time the forecasts are to be determined. To precisely state the final two assumptions, we will define $q(\delta)$ that quantifies the difference between the two approaches in terms of their expected empirical mean square prediction error forecasting performance.

To make the definition of $q(\delta)$ precise in an accessible manner we now present it from the inside outwards. At the core we have the local covariance function defined in~\eqref{equation:local_covariance_function} and averaged versions
\begin{equation}\label{eqn:gammaDelta_u}
\begin{aligned}
		\gamma^{(p)}_{\Delta}(u) & := \int_0^1 \gamma^{(p)}(u+\Delta (x-1)) {\rm d}x,
		& \gamma^{(p)}(u) & := [\gamma_1(u)\;\ldots\;\gamma_p(u)]',\\
		\Gamma^{(p)}_{\Delta}(u) & := \int_0^1 \Gamma^{(p)}(u+\Delta (x-1)) {\rm d}x,
		& \Gamma^{(p)}(u) & := (\gamma_{i-j}(u);\,i,j=1,\ldots,p).
\end{aligned}
\end{equation}
If $\Delta := (N-|k|)/T$ or $N/T$ and $u=t/T$, then the entries $\int_0^1 \gamma_k(u+\Delta (x-1)) {\rm d}x$ in $\gamma^{(p)}_{\Delta}(u)$ and $\Gamma^{(p)}_{\Delta}(u)$ are approximations for the expectation $\E \hat \gamma_{k;N,T}(t)$ of the lag $k$ autocovariance estimate $\hat \gamma_{k;N,T}(t)$ computed from $X_{t-N+1,T}, \ldots, X_{t,T}$; cf. Lemma~\ref{lem:exp_gamma_k}. This seemingly complicated construction is necessary, because we do not require that $N/T$ is negligible. By allowing $\Delta > 0$ we can capture the evolving second moments of the processes. Further note that, for every $u \in \IR$ and $\Delta \geq 0$, the averaged local autocovariances form the autocovariance function of a stationary process that can be seen as an average of the stationary approximations $X_t(\cdot)$ over the local times in $[u-\Delta,u]$. Solving the Yule-Walker equations for this average process yields
\begin{equation}\label{def:a_Delta}
	a^{(p)}_{\Delta}(u) := \big( a^{(p)}_{1,\Delta}(u), \ldots, a^{(p)}_{p,\Delta}(u) \big)' := \Gamma^{(p)}_{\Delta}(u)^{-1} \gamma^{(p)}_{\Delta}(u).
\end{equation}
As can be seen from Theorem~\ref{thm:propertiesYW} and Lemma~\ref{lem:rel_a_bar}, $a^{(p)}_{\Delta}(u)$ is an approximation to the limit of the Yule-Walker estimate obtained from $X_{t-N+1,T}, \ldots, X_{t,T}$. It further is related to the 1-step ahead linear forecasting coefficients, as can be seen from Lemma~\ref{lem:rel_a_1}. The $h$-step ahead counterpart of $a^{(p)}_{\Delta}(u)$ is defined as
\begin{equation}\label{def:v_Delta}
\begin{split}
	\big( v^{(p,h)}_{\Delta}(u) \big)'
				&  := \big(
				v_{1;\Delta}^{(p,h)}(u), v_{2;\Delta}^{(p,h)}(u), \cdots, v_{p;\Delta}^{(p,h)}(u)	\big) \\
				& := e'_1 \big( A_{\Delta}^{(p)}(u) \big)^h
				:= e'_1 \big( e_1 \big( a_{\Delta}^{(p)}(t) \big)' + H \big)^h,
\end{split}
\end{equation}
where $e_1$ and $H$ are the same as in~\eqref{def:vhat}.
Then, for  $u \in \IR$, $\Delta_1, \Delta_2 \geq 0$, the functions ${\rm MSPE}_{\Delta_1, \Delta_2}^{(p,h)}(u)$ are defined as
\begin{equation}
\label{eqn:thMSPE}
{\rm MSPE}_{\Delta_1, \Delta_2}^{(p,h)}(u) :=
\int_0^1 g^{(p,h)}_{\Delta_1}\Big( u + \Delta_2 (1-x) \Big) {\rm d}x,
\end{equation}
where $g^{(0,h)}_{\Delta}(u) := \gamma_0(u)$ and, for $p \in \IN^*$, with $\gamma_0^{(p,h)}(u) := \big( \gamma_h(u), \ldots, \gamma_{h+p-1}(u) \big)'$,
\begin{equation}\label{eqn:defg}
		g^{(p,h)}_{\Delta}(u) := \gamma_0(u) - 2 \big( v_{\Delta}^{(p,h)}(u) \big)' \gamma_0^{(p,h)}(u) + \big( v_{\Delta}^{(p,h)}(u) \big)' \Gamma_0^{(p)}(u) v_{\Delta}^{(p,h)}(u).
\end{equation}
From Lemmas~\ref{lem:main} and~\ref{lem:rel_MSPE}, it can be seen that
\begin{equation*}
	{\rm MSPE}_{s,m,N,T}^{(p,h)} := \frac{1}{m} \sum_{t=s+1}^{s+m} \big( X_{t+h,T} - \sum_{i=1}^p \hat v^{(p,h)}_{i;N,T}(t) X_{t-i+1,T} \big)^2,
\end{equation*}
concentrates around ${\rm MSPE}_{\Delta_1, \Delta_2}^{(p,h)}(u)$ with $\Delta_1 = N/T$, $\Delta_2 = m/T$ and $u = s/T$. Note that two arguments $\Delta_1$ and $\Delta_2$ are required to allow for the averaging of possible effects due to non-stationarity originating from (a) either the computation of the forecasting coefficients or (b) the computation of the mean squared prediction errors. The quantity $g^{(p,h)}_{N/T}(t/T)$ approximates the MSPE of $f^{{\lstat}}_{t,h;p,N}$ defined in~\eqref{eqn:lsEstim}. In the case of 1-step ahead forecasts we can simplify the expression in~\eqref{eqn:defg} to
\begin{equation}\label{eqn:g1}
	g^{(p,1)}_{\Delta}(u) = \E[ (\hat X_{t}^{(p)}(u) - X_{t}(u))^2]
		+ \big\| a_{\Delta}^{(p)}(u) - a_0^{(p)}(u) \big\|_{\Gamma_0^{(p)}(u)}^2,
\end{equation}
where $\hat X_{t}^{(p)}(u) := \sum_{j=1}^p a^{(p)}_{j,0}(u) X_{t-j}(u)$ is the best linear 1-step ahead forecast for $X_t(u)$ and $\| x \|_{\Gamma}^2 := x' \Gamma x$ denotes the quadratic form associated with $\Gamma$. Decomposition~\eqref{eqn:g1} is into two non-negative quantities. The first term only depends on the characteristics of the stationary tangent process $X_t(u)$ and will be a decreasing sequence with index $p$ for any~$u$. The second term is the squared weighted difference of the forecasting coefficients obtained from the stationary approximation at time~$u$ and the forecasting coefficients obtained from the non-stationary data; more precisely from the stationary approximations $X_t(\cdot)$ ``averaged'' over $[u-\Delta, u]$.

The final two assumptions require that the size $m$ of the validation sets, the smallest segment size $\min\mathcal{N}$ from which locally stationary forecasting coefficients are computed, and the number of available observations $T$ are large enough in relation to the maximum model order $\max\mathcal{P}$, the forecasting horizon $h$, and the minimum possible difference of performance of stationary and locally stationary forecasts in terms of MSPE, which we measure by
\begin{equation}\label{def:f}
	q(\delta) := \min_{\substack{p_1, p_2 \in \mathcal{P} \\ N \in \mathcal{N}}} \ \Big| {\rm MSPE}_{s_1/T,m/T}^{(p_1,h)}(\frac{s_1}{T})
- (1+\delta) \cdot {\rm MSPE}_{N/T,m/T}^{(p_2,h)}(\frac{s_1}{T}) \Big|, \quad
\end{equation}
where $s_1 := T-m-h$.

Assumption~\ref{a:m_and_N} requires the size $m$ of the validation sets and the smallest possible segment lengths $N \in \mathcal{N}$ from which to estimate the forecasting coefficients to be `large enough'.
\begin{assumption}[Minimum size for $m$ and $\min\mathcal{N}$]\label{a:m_and_N}
Let $K_0 := 4 C_0 (2 C_0+1)$.
For $\delta \geq 0$, $m, T, h \in \IN^*$, 
$\mathcal{P} \subset \{0,1,\ldots,\min\mathcal{N}-1\}$, such that $\mathcal{P} \neq \emptyset$, $\max\mathcal{P} \geq 1$, and $\emptyset \neq \mathcal{N} \subset \{\max\mathcal{P}+1, \ldots, T-2m-h+1\}$, assume that
\begin{equation*}
	\min\mathcal{N} \geq 8 h 2^h \big( C_0 \big)^{2 h+1} (\max\mathcal{P})^2 \max\Big\{\frac{20 (1+\delta)}{q(\delta)}, 1\Big\} \Big[ 6 (2\pi M'_f + C) + 1 \Big] 
\end{equation*}
and
\begin{multline}\label{cond:f} 
	\max\Big\{ \Big(\frac{h+\max\mathcal{P}}{m}\Big)^{\frac{1+4d}{3+8d}} K_0^{h} (\max\mathcal{P})^2, \Big(\frac{\max\mathcal{P}}{\min\mathcal{N}-\max\mathcal{P}}\Big)^{\frac{1+2d}{3+4d}} K_0^{h} (\max\mathcal{P})^3 h \Big\} \\
	 < \frac{q(\delta)}{20 (1+\delta)}. 
\end{multline}
\end{assumption}

Assumption~\ref{a:T} requires the sample size $T$ to be `large enough'.
\begin{assumption}[Minimum sample size $T$]\label{a:T}
With $C_0$ and $C_1$ defined in~\eqref{eqn:def_C0}, in the appendix, and $C$ and $M'_f$ from Assumptions~\ref{a:loc_stat} and~\ref{a:sd_der_bounded}, respectively, 
	\begin{equation*}
		T \geq \max\Big\{ 6 h 2^{h} C_1 (\max\mathcal{P})^2, 4 m \big( 2 h + 1 \big) \big( C_0 \big)^{2 h+1} M'_f \frac{20 (1+\delta)}{q(\delta)} \Big\}.
	\end{equation*}
\end{assumption}

The intuition behind the final two assumptions is that if two forecasts exist, one stationary and one locally stationary, that behave similarly well in terms of approximations to their expected empirical mean squared errors, then $m$ and $\min\mathcal{N}$ need to be large enough (in relation to $q(\delta), h$, and $\max\mathcal{P}$). Further, we require that $T$ exceeds a specified level (depending on $q(\delta), h, \max\mathcal{P}$, and $m$) to be able to provide bounds of the error of approximation of the local stationary process with the tangent process.
The specific form of Assumptions~\ref{a:m_and_N} and~\ref{a:T} are due to technical reasons in our proof and, in fact, our simulation results in Section~\ref{simulationsection1} suggest that the probability bounded in Theorem~\ref{thm:main3} will also be large for $T$ smaller than the threshold, as long as $\delta$ is chosen appropriately.
The quantity $q(\delta)$ is constructed to measure the difference between the MSPEs of the stationary predictors for different $p_1$ and the MSPEs of the locally stationary predictors for different $(p_2, N)$ scaled by a factor of $1+\delta$. Assumptions~\ref{a:m_and_N} and~\ref{a:T} are slightly stronger than necessary, as we do not only require only those procedures to perform differently for the $p_1$ and $(p_2, N)$ that yield the best result, but we require it for any combination. This is due to our method of proof. On the other hand, it is obvious that some condition like this is required for consistency of the procedure, because if there is no difference in performance either approach may equally likely be chosen. It is important to note that in the situation where both approaches perform equally well we do not need the selection to be consistent.

The quantity $q(\delta)$ depends on the model under consideration and, as $|\mathcal{P}|$ and $|\mathcal{N}|$ get larger, may potentially tend to zero. Thus, to employ Theorem~\ref{thm:main3} in practice, one has to analyse $q(\delta)$ to determine the right bounds stated in Assumptions~\ref{a:m_and_N} and~\ref{a:T}.
In Section~\ref{theorySpecialCase} we show how this can be done in the special case where $\mathcal{P} = \{1\}$ and $h=1$. There we show that if $\delta$ is chosen large enough or, in the case where the true model is non-stationary, if $\delta$ is chosen small enough, then $q$ is bounded away from $0$. If $q(\delta) > \varepsilon_0 > 0$, then, even in an asymptotic framework where $h$ and $\max\mathcal{P}$ do not need to be bounded and $m, \min\mathcal{N} \rightarrow \infty$ as $T \rightarrow \infty$, then condition~\eqref{cond:f} will hold for $T$ large enough, if
\[(h+\max\mathcal{P}) \big( K_0^{h} (\max\mathcal{P})^2 \big)^{\frac{3+8d}{1+4d}} = o(m), \quad \text{and} \quad
(\max\mathcal{P})^{\frac{10+14d}{1+2d}} \big( K_0^{h} h \big)^{\frac{3+4d}{1+2d}} = o(\min\mathcal{N}),
\]

Note that, $(\max\mathcal{P})^{1+2\frac{3+8d}{1+4d}} \leq (\max\mathcal{P})^{17/3}$ and $(\max\mathcal{P})^{\frac{10+14d}{1+2d}} \leq (\max\mathcal{P})^{17/2}$.
Therefore, if $h = O(1)$, we have that condition~\eqref{cond:f} will hold for $T$ large enough, if $\max\mathcal{P} = O(m^{3/17})$ and $\max\mathcal{P} = O\big((\min\mathcal{N})^{2/17}\big)$.

For the finite sample case, the quantity $q(\delta)$ can easily be computed for any tvAR($p$) model. A function performing the necessary calculations is provided in our \proglang{R} package \pkg{forecastSNSTS}. Numerical illustrations are provided in Section~\ref{simulationsection1}.

We are now ready to state the main result that guarantees that our procedure will, with high probability, choose the predictor that achieves the best empirical performance on the test set.

\begin{satz}\label{thm:main3}
Let $(X_{t,T})_{t \in \IZ, T \in \IN^*}$ satisfy
Assumptions~\ref{a:loc_stat}--\ref{a:MomentCond} and $\E X_{t,T} = 0$. Further, let $\delta, m, T, h, \mathcal{P}$, and $\mathcal{N}$ be such that Assumptions~\ref{a:m_and_N}--\ref{a:T} are satisfied.
Then, with $\hat R_{T,j}(h)$, $j=2,3$, defined in~\eqref{decisionrule}, we have
\begin{equation*}
\begin{split}
	& \IP\Big( (\hat R_{T,2}(h) \geq 1 + \delta \text{ and }
		\hat R_{T,3}(h) \geq 1 + \delta)
	\text{ or }
		(\hat R_{T,2}(h) < 1 + \delta \text{ and }
		\hat R_{T,3}(h) < 1 + \delta) \Big) \\
	& \quad \geq 1 - 6 D_1 |\mathcal{P}|^2  | \mathcal{N} |  \Bigg[ (\max\mathcal{P})^2 \exp\Bigg(-D_2 \Big( \frac{m}{h+\max\mathcal{P}} \Big)^{1/(3+8d)} \Bigg) \\
	& \hspace{4cm} + m (\max\mathcal{P})^3 \exp\Bigg(-D_3 \Big( \frac{\min\mathcal{N}-\max\mathcal{P}}{\max\mathcal{P}} \Big)^{1/(3+4d)} \Bigg) \Bigg],
	\end{split}
\end{equation*}
where $D_1, D_2$ and $D_3$ are constants, defined in~\eqref{def:D}, in the appendix, that
only depend on $K$ and $\rho$, $m_f$ and $M_f$, and $c$ and $d$ the constants from Assumption~\ref{a:mixing}, \ref{a:sd_bounded} and \ref{a:MomentCond}, respectively.
\end{satz}

The proof of Theorem~\ref{thm:main3} is long and technical and therefore deferred to Section~\ref{app:main}.
The probability in Theorem~\ref{thm:main3} tends to one if $m \gg (h+\max\mathcal{P}) [\log( |\mathcal{P}|^2 |\mathcal{N}| (\max\mathcal{P})^2 )]^{3+8d}$ and $\min\mathcal{N} \gg \max\mathcal{P} [\log(|\mathcal{P}|^2 |\mathcal{N}| m (\max\mathcal{P})^3 )]^{3+4d}$, where we have used the notation $a_T \gg b_T$ for $a_T / b_T \rightarrow \infty$, as $T \rightarrow \infty$.
Thus, Theorem~\ref{thm:main3} provides a ``performance guarantee'' of
our model choice methodology in the sense that it asserts that, with high probability, the method which \emph{we have observed to perform better empirically} in forecasting the observations from the second validation set \emph{will also perform better empirically in forecasting the future}, not yet observed values of the test set.

\subsection{Theoretical results for a simple, special cases} \label{theorySpecialCase}

To illustrate the usefulness of Theorem~\ref{thm:main3} we now discuss the special case in which the model order is pre-determined to be 1, for both locally stationary and stationary forecasts, and the forecasting horizon is 1-step ahead; i.\,e., $\mathcal{P} := \{1\}$ and $h=1$. Though this special case is usually not of practical interest, restricting ourselves will allow to illustrate how the general conditions simplify and can more easily be understood. For the simplification we proceed by finding lower bounds for $q(\delta)$ (uniformly in $T$ and $N \in \mathcal{N}$) which in turn allows us to state more explicit conditions that imply Assumptions~\ref{a:m_and_N} and~\ref{a:T}.

To apply Theorem~\ref{thm:main3}, we require that the MSPE of the stationary predictors are not to close to $1+\delta$ times the MSPE of the locally stationary predictors (cf. Assumption~\ref{a:m_and_N}).
Therefore, we now consider the following two cases:
\begin{itemize}
	\item[(a)] The parameter $\delta$ is chosen large enough.
	\item[(b)] The parameter $\delta$ is chosen small enough and the true model is non-stationary.
\end{itemize}
To make the requirements precise, we define
\begin{equation}\label{def:rho}
	\rho := \sup_{1-m/T \leq u \leq 1} \Big| \frac{\gamma_1(u)}{\gamma_0(u)} \Big|, 
\end{equation}
\begin{equation}\label{def:D_sup}
	D_{\sup} := \sup_{1-m/T \leq u \leq 1} \Bigg| \frac{\int_0^1 \gamma_1\big(u + \frac{s_1}{T}(x-1)\big) {\rm d}x}{\int_0^1 \gamma_0\big(u + \frac{s_1}{T}(x-1)\big) {\rm d}x}
- \frac{\gamma_1(u)}{\gamma_0(u)} \Bigg|,
\end{equation}
and
\begin{equation}\label{def:D_inf}
	D_{\inf} := \inf_{1-m/T \leq u \leq 1} \Bigg| \frac{\int_0^1 \gamma_1\big(u + \frac{s_1}{T}(x-1)\big) {\rm d}x}{\int_0^1 \gamma_0\big(u + \frac{s_1}{T}(x-1)\big) {\rm d}x}
- \frac{\gamma_1(u)}{\gamma_0(u)} \Bigg|,
\end{equation}
where $\gamma_0(u)$ and $\gamma_1(u)$ are the local autocovariances from Assumption~\ref{a:loc_stat}. The suprema and the infimum are with respect to points $u$ of the second validation set. Averaging of autocovariances in the first terms of $D_{\inf}$ and $D_{\sup}$ is across the training set and first validation set. Note that $D_{\inf} \leq D_{\sup} \leq 2$ and that $D_{\inf}$ is a measure for the non-stationarity of the training set. In particular, it will vanish if the data stems from a stationary process. Further, note that $\rho$ is a measure for the strength of serial dependence.

The simplified conditions that imply Assumptions~\ref{a:m_and_N} and~\ref{a:T} for the special case, will be stated in terms of $\rho$, $D_{\sup}$ and $D_{\inf}$. Note that, also in the case where $\mathcal{P} = \{1\}$, the quantity $q(\delta)$ in Assumptions~\ref{a:m_and_N} and~\ref{a:T} depends on $\mathcal{N}$, but the $D_{\inf}$, $D_{\sup}$ and $\rho$ only depend on $m$, $T$, $\gamma_1(\cdot)$ and $\gamma_0(\cdot)$. Therefore, the conditions in Lemmas~\ref{thm:main3:kor1} and~\ref{thm:main3:kor2} are indeed simpler than Assumptions~\ref{a:m_and_N} and~\ref{a:T}. Further note that the local autocovariances $\gamma_k(\cdot)$ can be determined easily for many time series models. If, for example the data stems from a tvAR(1) process with coefficient function~$a$, then we have
$\gamma_k(u) = a(u)^{|k|} / (1-a(u)^2)$, $k \in \IZ$.

We now state two results about the special case of the procedure for $1$-step ahead forecasting. The first result illustrates that the modified procedure will be consistent if $\delta$ is chosen large enough:

\begin{lem}\label{thm:main3:kor1}
Let $(X_{t,T})_{t \in \IZ, T \in \IN^*}$ satisfy Assumptions~\ref{a:loc_stat}--\ref{a:MomentCond}, and $\E X_{t,T} = 0$. Assume that $\rho < 1$ and $\delta \geq 2 D_{\sup}^2 / \big( 1 - \rho^2 \big)$,
where $\rho$ and $D_{\sup}$ are defined in~\eqref{def:rho} and~\eqref{def:D_sup}. Then, $q(\delta) \geq \delta \pi m_f ( 1 - \rho^2)$, where $m_f$ is from Assumption~\ref{a:sd_bounded}. In particular, this implies that constants $K_1$, $K_2$ and $K_3$, defined in the proof, exist such that, if
$m > K_1$ and $\min\mathcal{N} > K_2$
then Assumption~\ref{a:m_and_N} holds. Further, if $T \geq K_3 m$,
then Assumption~\ref{a:T} holds.
\end{lem}

Further more, we have as a second result that if the true model is non-stationary in the sense that the quantity $D_{\inf}$ is large compared to $N/T$ for all $N \in \mathcal{N}$, then we also have consistency for $\delta$'s that are small enough:

\begin{lem}\label{thm:main3:kor2}
Let $(X_{t,T})_{t \in \IZ, T \in \IN^*}$ satisfy Assumptions~\ref{a:loc_stat}--\ref{a:MomentCond}, $\E X_{t,T} = 0$, and
\begin{equation}\label{thm:main3:kor2:assN}
	D_{\inf}^2 \geq 2 \Big( \frac{M'_f}{m_f} \frac{\max\mathcal{N}}{T}\Big)^2,
\end{equation}
with $D_{\inf}$ defined in~\eqref{def:D_inf}.
Assume that $\delta \leq \frac{1}{8} D_{\inf}^2$,
Then, $q(\delta) \geq \pi D_{\inf}^2 m_f / 2$, where $m_f$ is from Assumption~\ref{a:sd_bounded}.
In particular, this implies that constants $K_4$, $K_5$ and $K_6$, defined in the proof, exist such that, if
$m > K_4$ and $\min\mathcal{N} > K_5$
then Assumption~\ref{a:m_and_N} holds. Further, if $T \geq K_6 m$
then Assumption~\ref{a:T} holds.
\end{lem}

By Lemma~\ref{thm:main3:kor1} we have that, in the case where $\mathcal{P} = \{1\}$, $h=1$ and $\delta \geq 0$ have been fixed, Assumption~\ref{a:m_and_N} will hold if $m$ and $\min\mathcal{N}$ are chosen larger than some constant. This requirement is not restrictive, in the sense that we would typically consider $m$ and $\min\mathcal{N}$ to diverge as $T$ diverges, such that by Theorem~\ref{thm:main3} the probability for consistent model choice will tend to one. In Lemma~\ref{thm:main3:kor2} the restrictions on $m$ and $\min\mathcal{N}$ are even weaker, as in a typical application $\max\mathcal{N}/T$ will tend to 0. In both Lemmas~\ref{thm:main3:kor1} and~\ref{thm:main3:kor2} the condition that implies Assumption~\ref{a:T} to hold is that $T$ is chosen larger than a multiple of $m$, which is eventually satisfied if $m/T$ tends to zero.

\begin{rem}\label{rem:P_mN}
In Lemmas~\ref{thm:main3:kor1} and~\ref{thm:main3:kor2} a lower bound of the form
\begin{equation}\label{lowerbound_q}
	\frac{q(\delta)}{20(1+\delta)} \geq \varepsilon_0
\end{equation}
is proven, for the special case where $\mathcal{P} = \{1\}$ and $h=1$. This lower bound implies that Assumption~\ref{a:m_and_N} holds, but it is in fact stronger, as Assumption~\ref{a:m_and_N} allows for $q(\delta)$ tending to 0, as $|\mathcal{N}|$ and $|\mathcal{P}|$ increase, as long as $m$ and $\min\mathcal{N}$ are increasing fast enough. Under condition~\eqref{lowerbound_q} and the conditions of Theorem~3.1 we have the following, stronger result:
\begin{equation*}
\begin{split}
	& \IP\Big( (\hat R_{T,2}(h) \geq 1 + \delta \text{ and }
		\hat R_{T,3}(h) \geq 1 + \delta)
	\text{ or }
		(\hat R_{T,2}(h) < 1 + \delta \text{ and }
		\hat R_{T,3}(h) < 1 + \delta) \Big) \\
	& \geq 1 - 6 D_1 |\mathcal{P}|^2 | \mathcal{N} |  \Bigg[ (\max\mathcal{P})^2 \exp\Bigg(-D_2 \varepsilon_0^{\frac{1}{2+4d}} \Big( \frac{m}{K_0^h(\max\mathcal{P})^2 (h+\max\mathcal{P})} \Big)^{\frac{1}{2+4d}} \Bigg) \\
	& \hspace{5cm} + m (\max\mathcal{P})^3 \exp\Bigg(-D_3 \varepsilon_0^{\frac{1}{2+2d}} \Big( \frac{\min\mathcal{N}-\max\mathcal{P}}{h K_0^h (\max\mathcal{P})^4} \Big)^{\frac{1}{2+2d}} \Bigg) \Bigg],
	\end{split}
\end{equation*}
which can be proved along the same lines of the proof of Theorem~3.1, together with inequality~\eqref{eqn:nicer_bound} from the proof of Lemma~\ref{lem:nicerboundP}, which is available in the arXiv'ed version of the manuscript \citep{KleyEtAl2019}.

In particular, when the parameters $\max\mathcal{P} = (\max\mathcal{P})(T)$, $h = h(T)$ and $\varepsilon_0 = \varepsilon_0(T)$ are bounded sequences ($\varepsilon_0$ also bounded away from zero), we get the following bound:
\begin{align}
	& \IP\Big( (\hat R_{T,2}(h) \geq 1 + \delta \text{ and }
		\hat R_{T,3}(h) \geq 1 + \delta)
	\text{ or }
		(\hat R_{T,2}(h) < 1 + \delta \text{ and }
		\hat R_{T,3}(h) < 1 + \delta) \Big) \nonumber \\
	& \geq 1 - \kappa_1 | \mathcal{N} | \Big( \exp(-m^{1/(2+4d)} \kappa_2) + m \exp(-N^{1/(2+2d)} \kappa_3) \Big) \label{eqn:bound_h1P1}
	\end{align}
where $\kappa_1, \kappa_2, \kappa_3$ are constants that do not depend on $m$ or $N$ and $d$ is the constant from Assumption~\ref{a:MomentCond} (e.\,g., for sub-Gaussian processes: $d=1/2$).
\end{rem}

\section{Simulations} \label{simulationsection1}

In this section we discuss finite sample properties of the estimates $\hat R_{T,i}(h)$, defined in~\eqref{decisionrule}, and their population counterparts
$ R_{T,j}(h) := (\E({\rm MSPE}_{T,j}^{{\stat}}(h)))/(\E({\rm MSPE}_{T,j}^{\lstat}(h)))$.
The simulation study was conducted with the \proglang{R} package \texttt{forecastSNSTS} \citep{R,forecastSNSTS}, available from The Comprehensive R Archive Network (CRAN). In particular, we investigate the performance of decision rule \eqref{decisionrule}. To this end, we have considered 15 different tvAR models. Three of the models are stationary, the other 12 are non-stationary. Amongst the non-stationary processes we have some where the covariance structure changes quickly and some where the covariances change slowly. Further, we will have examples where the processes given by the parameters at some local time $u$ are almost unit root and some where they are not.

For each of the models we proceed as follows. We simulate sequences of length  $T+m = n \in \{100, 200, 500, 1000, 2000, 4000, 6000, 8000, 10000\}$. The $T+m$ observations, with $T$ and $m$ as in Section~\ref{ignore.stat}, contain the training, validation and test set. We separate the test and validation sets of length
$m := \lfloor n^{.85}/4 \rfloor$. Thus,
$n_{i} := n-(3-i) \lfloor n^{.85}/4 \rfloor$, $i=0,...,3$,
mark the end indices of the training set, the validation sets and the test set, respectively. We have chosen $m$ as a function of $n$ in such a way that $m = o(n)$ and $m \rightarrow \infty$, as $n \rightarrow \infty$.
The sizes of the three sets therefore are $12, 22, 49, 88, 159, 288, 406, 519$, and $627$ for the different sequence sizes, respectively.

As described in Section~\ref{sec:prec} we then, for any $h=1,\ldots,H := 10$, determine linear \hbox{$h$-step} ahead predictions for $X_{t+h,T}$ with $t+h \in \{n_{0}+1, \ldots, n_{1}\}$. We determine the `stationary predictions', with coefficients estimated for a given $p=0,\ldots,p_{\max} := 7$, from $X_{1,T}, \ldots, X_{t,T}$ by $\hat v_{t,T}^{(p,h)}(t)$ from Step~3 of the procedure. For simplicity, we have chosen the same $p_{\max}$ for every $T$. We further determine `locally stationary predictions' where the coefficients $\hat v_{N,T}^{(p,h)}(t)$ are used for $p=0,\ldots,p_{\max}$ and
\[\mathcal{N} = \{ N := N_{\min} + i \lceil (N_{\max} - N_{\min})/ 25 \rceil \, : \, i \in \IN, N \leq N_{\max} \},\]
where $N_{\min} := \lfloor (n/2)^{4/5} \rfloor$ and $N_{\max} := \lfloor n^{4/5} \rfloor$.
Instead of considering every integer between $N_{\min}$ and $N_{\max}$ as a possible segment size, we restrict the number $\#\mathcal{N}$ of possible values for $N$ to a maximum of 25 elements to reduce computation time. The results did not change significantly when a larger number of elements was used.
We then compare the predictors with respect to their empirical mean squared prediction error (MSPE) on the first validation set and, according to Step 4 of the procedure, choose the stationary predictor with $\hat p_{\stat}$ that minimises the MSPE on $M_1$ amongst all stationary predictors and the locally stationary predictor with $(\hat p_{\lstat}, \hat N_{\lstat})$ which minimizes the MSPE on $M_1$ amongst all the locally stationary predictors.

For those two predictors we then determine the empirical mean squared prediction errors ${\rm MSPE}_{T,2}^*$ and ${\rm MSPE}_{T,3}^*$, defined in~\eqref{eqn:empMSPE}, on the validation and test set, respectively.
We record seven pieces of information: $\hat p_{\stat}$, $\hat p_{\lstat}, \hat N_{\lstat}$, ${\rm MSPE}_{T,2}^{\stat}$, ${\rm MSPE}_{T,2}^{\lstat}$, ${\rm MSPE}_{T,3}^{\stat}$, and ${\rm MSPE}_{T,3}^{\lstat}$.
We replicate the experiment $10000$ times.

Now we define the first two models. Both are tvAR(1) models defined by two periodic coefficient functions, namely the models are
\begin{align}
\label{periodic1}
X_{t,T} & = (0.8+0.19\sin(4\pi\frac{t}{T})) X_{t-1,T}+Z_t, \\
\label{periodic2}
X_{t,T} & = (0.3+0.19\sin(4\pi\frac{t}{T})) X_{t-1,T}+Z_t.
\end{align}
The innovations $Z_t$ are i.\,i.\,d Gaussian white noise.
In this section we discuss the above two models in detail. The remaining processes are defined in Appendix~\ref{a:extraSim} \citep{KleyEtAl2019}, where to also the corresponding tables and figures for them are being deferred.

\begin{figure}
\centering 
\includegraphics[width=0.24\textwidth]{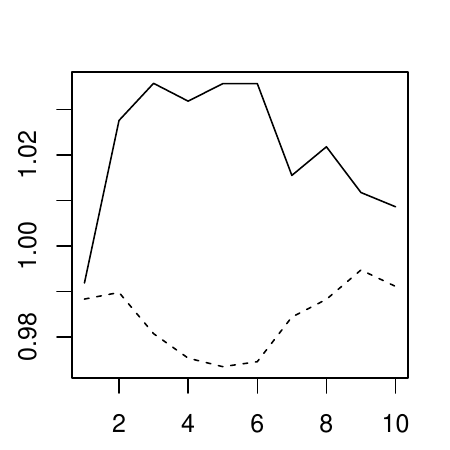}~
\includegraphics[width=0.24\textwidth]{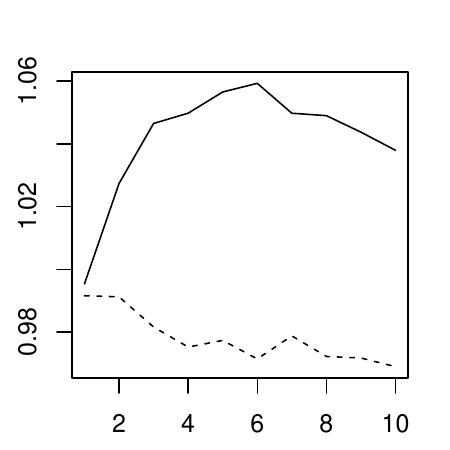}~
\includegraphics[width=0.24\textwidth]{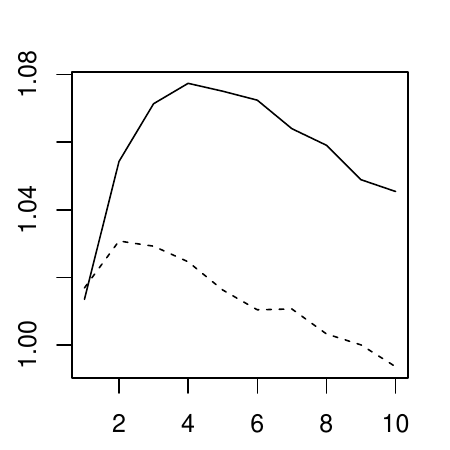}~
\includegraphics[width=0.24\textwidth]{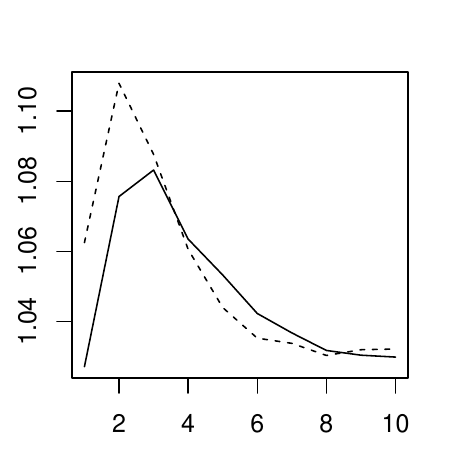} \\
\includegraphics[width=0.24\textwidth]{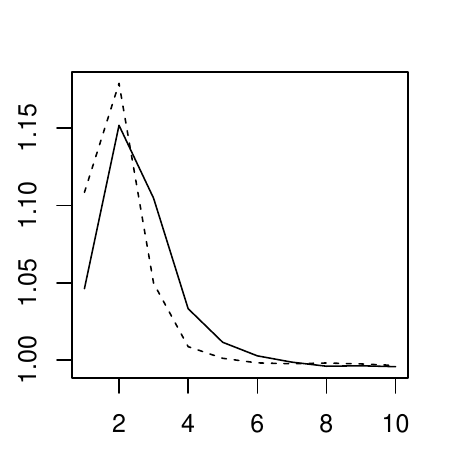}~
\includegraphics[width=0.24\textwidth]{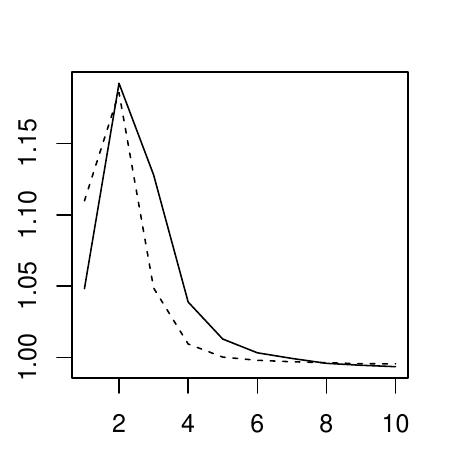}~
\includegraphics[width=0.24\textwidth]{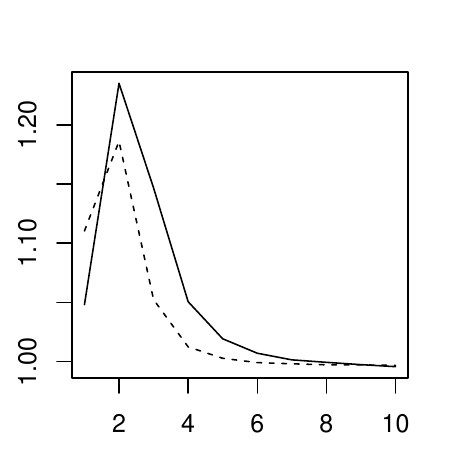}~
\includegraphics[width=0.24\textwidth]{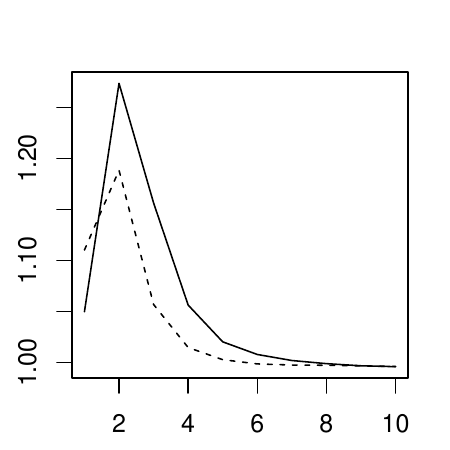}
   \caption{\it Plot of $h \mapsto R_{T,i}(h)$ for model \eqref{periodic1} and different values of $n$ (from left to right: n=100, n=200, n=500, n=1000 [first row], n=4000, n=6000, n=8000, n=10000 [second row]). Solid line: $i=3$ (test set), dashed line: $i=2$ (second validation set).} \label{MSPEperiodic1}
\end{figure}

\begin{table} \scriptsize
\begin{center}
\begin{tabular}{|c||c|c|c||c|c|c|}
\hline
 \multirow{2}{*}{$n$} & & \multicolumn{2}{|c||}{$R^{{\rm s}, {\rm ls}}_{T,2}(1)$}  & & \multicolumn{2}{|c|}{$R^{{\rm s}, {\rm ls}}_{T,2}(5)$} \\ 
 & & $\geq 1.01$  & $< 1.01$  & & $\geq 1.01$  & $< 1.01$ \\ 
 \hline 
 \multirow{2}{*}{$100$} & $R^{{\rm s}, {\rm ls}}_{T,3}(1) \geq 1.01$  & 0.1825 & 0.2777  & $R^{{\rm s}, {\rm ls}}_{T,3}(5) \geq 1.01$  & 0.1747 & 0.2479 \\ 
& $R^{{\rm s}, {\rm ls}}_{T,3}(1) < 1.01$  & 0.1888 & 0.351 & $R^{{\rm s}, {\rm ls}}_{T,3}(5) < 1.01$  & 0.1424 & 0.435 \\

\hline\hline
 \multirow{2}{*}{$n$} & & \multicolumn{2}{|c||}{$R^{{\rm s}, {\rm ls}}_{T,2}(1)$}  & & \multicolumn{2}{|c|}{$R^{{\rm s}, {\rm ls}}_{T,2}(5)$} \\ 
 & & $\geq 1.2$  & $< 1.2$  & & $\geq 1.2$  & $< 1.2$ \\ 
 \hline 
 \multirow{2}{*}{$1000$} & $R^{{\rm s}, {\rm ls}}_{T,3}(1) \geq 1.2$  & 5e-04 & 0.0055  & $R^{{\rm s}, {\rm ls}}_{T,3}(5) \geq 1.2$  & 0.0758 & 0.0636 \\ 
& $R^{{\rm s}, {\rm ls}}_{T,3}(1) < 1.2$  & 0.0063 & 0.9877 & $R^{{\rm s}, {\rm ls}}_{T,3}(5) < 1.2$  & 0.0699 & 0.7907 \\

\hline\hline
 \multirow{2}{*}{$n$} & & \multicolumn{2}{|c||}{$R^{{\rm s}, {\rm ls}}_{T,2}(1)$}  & & \multicolumn{2}{|c|}{$R^{{\rm s}, {\rm ls}}_{T,2}(5)$} \\ 
 & & $\geq 1$  & $< 1$  & & $\geq 1$  & $< 1$ \\ 
 \hline 
 \multirow{2}{*}{$10000$} & $R^{{\rm s}, {\rm ls}}_{T,3}(1) \geq 1$  & 0.9916 & 0  & $R^{{\rm s}, {\rm ls}}_{T,3}(5) \geq 1$  & 0.7567 & 0.2054 \\ 
& $R^{{\rm s}, {\rm ls}}_{T,3}(1) < 1$  & 0.0084 & 0 & $R^{{\rm s}, {\rm ls}}_{T,3}(5) < 1$  & 0.0251 & 0.0128 \\

\hline\hline
 \multirow{2}{*}{$n$} & & \multicolumn{2}{|c||}{$R^{{\rm s}, {\rm ls}}_{T,2}(1)$}  & & \multicolumn{2}{|c|}{$R^{{\rm s}, {\rm ls}}_{T,2}(5)$} \\ 
 & & $\geq 1.05$  & $< 1.05$  & & $\geq 1.05$  & $< 1.05$ \\ 
 \hline 
 \multirow{2}{*}{$10000$} & $R^{{\rm s}, {\rm ls}}_{T,3}(1) \geq 1.05$  & 0.4917 & 4e-04  & $R^{{\rm s}, {\rm ls}}_{T,3}(5) \geq 1.05$  & 0.0019 & 0.1698 \\ 
& $R^{{\rm s}, {\rm ls}}_{T,3}(1) < 1.05$  & 0.5077 & 2e-04 & $R^{{\rm s}, {\rm ls}}_{T,3}(5) < 1.05$  & 0.0025 & 0.8258 \\

\hline\hline
 \multirow{2}{*}{$n$} & & \multicolumn{2}{|c||}{$R^{{\rm s}, {\rm ls}}_{T,2}(1)$}  & & \multicolumn{2}{|c|}{$R^{{\rm s}, {\rm ls}}_{T,2}(5)$} \\ 
 & & $\geq 1.1$  & $< 1.1$  & & $\geq 1.1$  & $< 1.1$ \\ 
 \hline 
 \multirow{2}{*}{$10000$} & $R^{{\rm s}, {\rm ls}}_{T,3}(1) \geq 1.1$  & 0.0033 & 9e-04  & $R^{{\rm s}, {\rm ls}}_{T,3}(5) \geq 1.1$  & 1e-04 & 0.0119 \\ 
& $R^{{\rm s}, {\rm ls}}_{T,3}(1) < 1.1$  & 0.7025 & 0.2933 & $R^{{\rm s}, {\rm ls}}_{T,3}(5) < 1.1$  & 1e-04 & 0.9879 \\

\hline\hline
 \multirow{2}{*}{$n$} & & \multicolumn{2}{|c||}{$R^{{\rm s}, {\rm ls}}_{T,2}(1)$}  & & \multicolumn{2}{|c|}{$R^{{\rm s}, {\rm ls}}_{T,2}(5)$} \\ 
 & & $\geq 1.15$  & $< 1.15$  & & $\geq 1.15$  & $< 1.15$ \\ 
 \hline 
 \multirow{2}{*}{$10000$} & $R^{{\rm s}, {\rm ls}}_{T,3}(1) \geq 1.15$  & 0 & 0  & $R^{{\rm s}, {\rm ls}}_{T,3}(5) \geq 1.15$  & 0 & 7e-04 \\ 
& $R^{{\rm s}, {\rm ls}}_{T,3}(1) < 1.15$  & 0.0188 & 0.9812 & $R^{{\rm s}, {\rm ls}}_{T,3}(5) < 1.15$  & 2e-04 & 0.9991 \\

\hline
\end{tabular}
\caption{\textit{Proportions of the individual events in~\eqref{samedecision} for the process \eqref{periodic1} and selected combinations of $n$ and $\delta$.}}  \label{MSPEanalysisperiodic1c}
\end{center}
\end{table}

In Figure~\ref{MSPEperiodic1}, note that, since in the numerator we have the MSPE for the best stationary predictor and in the denominator the MSPE for the best locally stationary predictor, a ratio above 1 corresponds to the situation where the best locally stationary predictor outperforms the best stationary predictor. It can be seen whether this happens on average, while in Table~\ref{MSPEanalysisperiodic1} we can see the proportion of simulated cases in which this has happens. In Figure~\ref{MSPEperiodic1}, we thus observe that, for $n=100$, the stationary approach performs better on average across all values of $h$ on both the test and the second validation set. For $n = 200$ the locally stationary approach performs better for $3 \leq h \leq 6$ on the test set, while the stationary approach still excels for all $h$ on the second validation set. For $n \geq 500$ the locally stationary approach is better across all values of $h$ on the test set and for $2 \leq h \leq 4$ it outperforms the stationary approach on the second validation set. For $n \geq 1000$ the locally stationary approach is always as least as good as the stationary approach for all $h$. It is striking that, for this particular model and for the larger $n$'s we see that as $h$ gets larger the two approaches (stationary and locally stationary) perform almost equally well on average, which can be seen from the lines in Figure~\ref{MSPEperiodic1} being close to one. Another important observation is that, as $n$ gets larger and $m/n$ gets smaller, we see the lines for the validation and test set converging, which is in line with what Theorem~\ref{thm:main3} suggests should happen.

We now, briefly, compare the outcome of model~\eqref{periodic1} to that of model~\eqref{periodic2}; details are shown in Appendix~\ref{a:extraSim} \citep{KleyEtAl2019}. Note that in model~\eqref{periodic1} the coefficient function ranges from $0.61$ to $0.99$, placing some of its tangent processes close to the unit root. In model~\eqref{periodic2} the coefficient function ranges from $0.11$ to $0.49$. Thus, the two models have the same variation of the coefficient function, but in model~\eqref{periodic2} the tangent processes are further away from the unit root. In Figure~\ref{MSPEperiodic2}, it can be seen that the stationary approach is preferred over the locally stationary approach for sequences up to length $n=1000$. Further, we observe that the advantage of using the locally stationary approach for sequences of length $n \geq 4000$ is minuscule and visible only for $1$-step ahead forecasting. For the other models we can make similar observations:

\textbf{Rules of Thumb.} The locally stationary approach outperforms the stationary approach only if either the sequence is long, or the coefficient function exhibits considerable variation, or the tangent processes (cf. the comment after Assumptions~\ref{a:loc_stat}--\ref{a:MomentCond}) are close to the unit root. In any other case the stationary approach can be chosen without (a large) loss.

Our observation that the locally stationary forecast performs better when the stationary approximations are near unit root may possibly be explained by the fact that the coefficient of a near unit root AR(1) process can be estimated at a better rate than in the classical case where the rate is $T^{-1/2}$; cf. \cite{ChanWei1987,DzhaparidzeEtAl1994}. A rigorous analysis of the issue is beyond the scope of this paper and left for future research.

\begin{table} \scriptsize
\begin{center}
	\begin{tabular}{|c||c|c|c|c|c|c|c|c|c|}\hline 
 &\multicolumn{9}{|c|}{\eqref{samedecision} holds for $h=1$, $i=2$} \\ 
\hline 
$\delta$ & $n$=100 & $n$=200 & $n$=500 & $n$=1000 & $n$=2000 & $n$=4000 & $n$=6000 & $n$=8000 & $n$=10000 \\ 
\hline \hline 
0  & 0.4182 & 0.4414 & 0.5873 & 0.8844 & 0.9962 & 0.9999 & 1 & 1 & 1 \\ 
0.01  & 0.3713 & 0.3821 & 0.5152 & 0.8422 & 0.993 & 0.9999 & 1 & 1 & 1 \\ 
0.05  & 0.2431 & 0.2172 & 0.27 & 0.5986 & 0.9106 & 0.9855 & 0.9958 & 0.9985 & 0.9994 \\ 
0.1  & 0.1518 & 0.1029 & 0.0893 & 0.2347 & 0.4748 & 0.6212 & 0.67 & 0.6995 & 0.7058 \\ 
0.15  & 0.0983 & 0.0483 & 0.024 & 0.0475 & 0.071 & 0.0617 & 0.0449 & 0.0305 & 0.0188 \\ 
0.2  & 0.0622 & 0.0233 & 0.0072 & 0.0068 & 0.0046 & 6e-04 & 1e-04 & 0 & 0 \\ 
0.4  & 0.0138 & 0.0026 & 8e-04 & 2e-04 & 5e-04 & 1e-04 & 0 & 0 & 0 \\ 
0.6  & 0.0059 & 0.0015 & 4e-04 & 2e-04 & 4e-04 & 0 & 0 & 0 & 0 \\ 
\hline 
\end{tabular} 
	\begin{tabular}{|c||c|c|c|c|c|c|c|c|c|}\hline 
 &\multicolumn{9}{|c|}{\eqref{samedecision} holds for $h=1$, $i=3$} \\ 
\hline 
$\delta$ & $n$=100 & $n$=200 & $n$=500 & $n$=1000 & $n$=2000 & $n$=4000 & $n$=6000 & $n$=8000 & $n$=10000 \\ 
\hline \hline 
0  & 0.4968 & 0.5118 & 0.5989 & 0.6929 & 0.8099 & 0.9279 & 0.9681 & 0.9829 & 0.9916 \\ 
0.01  & 0.4602 & 0.4671 & 0.5436 & 0.6357 & 0.7485 & 0.8725 & 0.9276 & 0.9559 & 0.9751 \\ 
0.05  & 0.3357 & 0.3103 & 0.3372 & 0.3905 & 0.416 & 0.4577 & 0.4776 & 0.4732 & 0.4921 \\ 
0.1  & 0.2292 & 0.1745 & 0.1481 & 0.1427 & 0.0954 & 0.0458 & 0.0195 & 0.0113 & 0.0042 \\ 
0.15  & 0.1487 & 0.0942 & 0.0536 & 0.0313 & 0.0098 & 3e-04 & 2e-04 & 0 & 0 \\ 
0.2  & 0.0983 & 0.0491 & 0.0194 & 0.006 & 7e-04 & 1e-04 & 0 & 0 & 0 \\ 
0.4  & 0.0209 & 0.006 & 0.0013 & 3e-04 & 5e-04 & 1e-04 & 0 & 0 & 0 \\ 
0.6  & 0.0078 & 0.0023 & 7e-04 & 2e-04 & 5e-04 & 1e-04 & 0 & 0 & 0 \\ 
\hline 
\end{tabular} 
	\begin{tabular}{|c||c|c|c|c|c|c|c|c|c|}\hline 
 &\multicolumn{9}{|c|}{\eqref{samedecision} holds for $h=5$, $i=2$} \\ 
\hline 
$\delta$ & $n$=100 & $n$=200 & $n$=500 & $n$=1000 & $n$=2000 & $n$=4000 & $n$=6000 & $n$=8000 & $n$=10000 \\ 
\hline \hline 
0  & 0.3365 & 0.3525 & 0.3646 & 0.323 & 0.1831 & 0.1944 & 0.2508 & 0.3045 & 0.359 \\ 
0.01  & 0.3171 & 0.3246 & 0.3416 & 0.3137 & 0.1648 & 0.1409 & 0.1642 & 0.1853 & 0.2101 \\ 
0.05  & 0.2535 & 0.236 & 0.2609 & 0.2798 & 0.1216 & 0.0307 & 0.0121 & 0.0073 & 0.0044 \\ 
0.1  & 0.1968 & 0.159 & 0.18 & 0.2378 & 0.0985 & 0.0132 & 0.0013 & 3e-04 & 2e-04 \\ 
0.15  & 0.1541 & 0.1101 & 0.1235 & 0.1929 & 0.0844 & 0.012 & 0.0011 & 1e-04 & 2e-04 \\ 
0.2  & 0.1244 & 0.0757 & 0.0843 & 0.1457 & 0.0696 & 0.0101 & 9e-04 & 1e-04 & 2e-04 \\ 
0.4  & 0.0575 & 0.0207 & 0.0164 & 0.0263 & 0.0101 & 0.0013 & 1e-04 & 0 & 0 \\ 
0.6  & 0.0321 & 0.0068 & 0.0032 & 0.0033 & 4e-04 & 0 & 0 & 0 & 0 \\ 
\hline 
\end{tabular} 
	\begin{tabular}{|c||c|c|c|c|c|c|c|c|c|}\hline 
 &\multicolumn{9}{|c|}{\eqref{samedecision} holds for $h=5$, $i=3$} \\ 
\hline 
$\delta$ & $n$=100 & $n$=200 & $n$=500 & $n$=1000 & $n$=2000 & $n$=4000 & $n$=6000 & $n$=8000 & $n$=10000 \\ 
\hline \hline 
0  & 0.4423 & 0.5487 & 0.5253 & 0.3178 & 0.2155 & 0.2979 & 0.3955 & 0.4736 & 0.5393 \\ 
0.01  & 0.4226 & 0.5269 & 0.5054 & 0.3086 & 0.1943 & 0.2561 & 0.3474 & 0.4249 & 0.4931 \\ 
0.05  & 0.3563 & 0.4417 & 0.4281 & 0.2705 & 0.1312 & 0.0999 & 0.1168 & 0.1477 & 0.1717 \\ 
0.1  & 0.2823 & 0.3492 & 0.3378 & 0.2245 & 0.0859 & 0.0243 & 0.0153 & 0.0119 & 0.012 \\ 
0.15  & 0.2291 & 0.2779 & 0.2602 & 0.1806 & 0.0625 & 0.0098 & 0.0023 & 8e-04 & 7e-04 \\ 
0.2  & 0.1864 & 0.2192 & 0.1923 & 0.1394 & 0.0451 & 0.0065 & 8e-04 & 1e-04 & 0 \\ 
0.4  & 0.0815 & 0.0823 & 0.0561 & 0.0398 & 0.0081 & 5e-04 & 0 & 0 & 0 \\ 
0.6  & 0.0406 & 0.0303 & 0.0164 & 0.0071 & 4e-04 & 0 & 0 & 0 & 0 \\ 
\hline 
\end{tabular} 
\caption{\textit{Proportion of \eqref{decisionrule}  being fulfilled for the process \eqref{periodic1} and different values of $h$, $\delta$ and $n$.}}  \label{MSPEanalysisperiodic1}
\end{center}
\end{table}

\begin{table} \scriptsize
\begin{center}
	\begin{tabular}{|c||c|c|c|c|c|c|c|c|c|}\hline 
 &\multicolumn{9}{|c|}{\eqref{samedecision} holds for $h=1$} \\ 
\hline 
$\delta$ & $n$=100 & $n$=200 & $n$=500 & $n$=1000 & $n$=2000 & $n$=4000 & $n$=6000 & $n$=8000 & $n$=10000 \\ 
\hline \hline 
0  & 0.5254 & 0.5176 & 0.5088 & 0.6377 & 0.8071 & 0.9278 & 0.9681 & 0.9829 & 0.9916 \\ 
0.01  & 0.5335 & 0.5308 & 0.4984 & 0.5877 & 0.7441 & 0.8724 & 0.9276 & 0.9559 & 0.9751 \\ 
0.05  & 0.606 & 0.6265 & 0.5928 & 0.4899 & 0.4346 & 0.46 & 0.4776 & 0.4737 & 0.4919 \\ 
0.1  & 0.7114 & 0.7714 & 0.7972 & 0.7034 & 0.53 & 0.395 & 0.3381 & 0.3054 & 0.2966 \\ 
0.15  & 0.8028 & 0.8743 & 0.9284 & 0.9256 & 0.9214 & 0.9384 & 0.9549 & 0.9695 & 0.9812 \\ 
0.2  & 0.8667 & 0.9358 & 0.9756 & 0.9882 & 0.9957 & 0.9995 & 0.9999 & 1 & 1 \\ 
0.4  & 0.9727 & 0.9936 & 0.9995 & 0.9999 & 1 & 1 & 1 & 1 & 1 \\ 
0.6  & 0.9903 & 0.9982 & 0.9993 & 1 & 0.9999 & 0.9999 & 1 & 1 & 1 \\ 
\hline 
\end{tabular} 
	\begin{tabular}{|c||c|c|c|c|c|c|c|c|c|}\hline 
 &\multicolumn{9}{|c|}{\eqref{samedecision} holds for $h=5$} \\ 
\hline 
$\delta$ & $n$=100 & $n$=200 & $n$=500 & $n$=1000 & $n$=2000 & $n$=4000 & $n$=6000 & $n$=8000 & $n$=10000 \\ 
\hline \hline 
0  & 0.6054 & 0.5684 & 0.6755 & 0.8888 & 0.879 & 0.7975 & 0.7729 & 0.7699 & 0.7695 \\ 
0.01  & 0.6097 & 0.5711 & 0.673 & 0.8851 & 0.8819 & 0.783 & 0.7324 & 0.6858 & 0.6538 \\ 
0.05  & 0.633 & 0.5991 & 0.6752 & 0.8705 & 0.9246 & 0.9008 & 0.8801 & 0.8492 & 0.8277 \\ 
0.1  & 0.6817 & 0.6582 & 0.7082 & 0.8601 & 0.9546 & 0.9825 & 0.9854 & 0.988 & 0.988 \\ 
0.15  & 0.7252 & 0.715 & 0.7547 & 0.8583 & 0.9551 & 0.9938 & 0.998 & 0.9993 & 0.9991 \\ 
0.2  & 0.7664 & 0.7677 & 0.807 & 0.8665 & 0.9553 & 0.9944 & 0.9991 & 1 & 0.9998 \\ 
0.4  & 0.8816 & 0.9046 & 0.9373 & 0.9463 & 0.9872 & 0.9984 & 0.9999 & 1 & 1 \\ 
0.6  & 0.9337 & 0.9647 & 0.9814 & 0.99 & 0.9992 & 1 & 1 & 1 & 1 \\ 
\hline 
\end{tabular} 
\caption{\textit{Proportion of \eqref{samedecision}  being fulfilled for the process \eqref{periodic1} and different values of $h$, $\delta$ and $n$.}}  \label{MSPEanalysisperiodic1b}
\end{center}
\end{table}

\begin{table} \scriptsize
\begin{center}
	\begin{tabular}{|c||c|c|c|c|c|c|c|c|c|}\hline 
 &\multicolumn{9}{|p{12.5cm}|}{\centering Value of $q(\delta)$, defined in \eqref{cond:f}, for $h=1$} \\ 
\hline 
$\delta$ & $n$=100 & $n$=200 & $n$=500 & $n$=1000 & $n$=2000 & $n$=4000 & $n$=6000 & $n$=8000 & $n$=10000 \\ 
\hline \hline 
0  & 0 & 0 & 0 & 0 & 0 & 0 & 0 & 0 & 0 \\ 
0.01  & 6.2e-05 & 1.1e-05 & 7.6e-05 & 0.02 & 0.02 & 0.021 & 0.021 & 0.021 & 0.021 \\ 
0.05  & 0.00018 & 1.7e-05 & 2e-05 & 6.4e-05 & 0.012 & 0.0087 & 0.0054 & 0.0032 & 0.0014 \\ 
0.1  & 0.023 & 6.4e-05 & 2.6e-05 & 3.7e-06 & 0.0056 & 0.011 & 0.015 & 0.017 & 0.018 \\ 
0.15  & 0.076 & 0.051 & 0.044 & 0.049 & 0.056 & 0.062 & 0.066 & 0.067 & 0.069 \\ 
0.2  & 0.13 & 0.1 & 0.094 & 0.1 & 0.11 & 0.11 & 0.12 & 0.12 & 0.12 \\ 
0.4  & 0.2 & 0.28 & 0.3 & 0.3 & 0.31 & 0.31 & 0.32 & 0.32 & 0.32 \\ 
0.6  & 3e-04 & 0.054 & 0.19 & 0.33 & 0.42 & 0.47 & 0.49 & 0.51 & 0.52 \\ 
\hline 
\end{tabular} 
	\begin{tabular}{|c||c|c|c|c|c|c|c|c|c|}\hline 
 &\multicolumn{9}{|p{12.5cm}|}{\centering Value of $q(\delta)$, defined in \eqref{cond:f}, for $h=5$} \\ 
\hline 
$\delta$ & $n$=100 & $n$=200 & $n$=500 & $n$=1000 & $n$=2000 & $n$=4000 & $n$=6000 & $n$=8000 & $n$=10000 \\ 
\hline \hline 
0  & 0 & 0 & 0 & 0 & 0 & 0 & 0 & 0 & 0 \\ 
0.01  & 3.3e-06 & 0.00029 & 0.00041 & 6.5e-05 & 0.02 & 0.021 & 0.021 & 0.021 & 0.021 \\ 
0.05  & 0.00016 & 4.9e-06 & 0.001 & 0.042 & 0.034 & 0.026 & 0.021 & 0.018 & 0.015 \\ 
0.1  & 0.005 & 0.0015 & 3.9e-05 & 0.14 & 0.13 & 0.13 & 0.12 & 0.12 & 0.12 \\ 
0.15  & 0.094 & 5.5e-05 & 0.00022 & 0.055 & 0.16 & 0.13 & 0.12 & 0.11 & 0.1 \\ 
0.2  & 0.2 & 0.00088 & 0.00018 & 0.00091 & 0.058 & 0.031 & 0.015 & 0.0046 & 0.0035 \\ 
0.4  & 4e-04 & 0.05 & 0.099 & 0.14 & 0.17 & 0.2 & 0.22 & 0.23 & 0.24 \\ 
0.6  & 0.33 & 0.4 & 0.46 & 0.52 & 0.56 & 0.6 & 0.62 & 0.64 & 0.65 \\ 
\hline 
\end{tabular} 
\caption{\textit{Values of $q(\delta)$, defined in \eqref{cond:f}, for the process \eqref{periodic1} and different values of $h$, $\delta$ and $n$.}} \label{MSPEanalysisperiodic1d}
\end{center}
\end{table}

The proportions shown in Table~\ref{MSPEanalysisperiodic1b} provide information on the consistency of the procedure, as we see the proportion of cases in which the same procedure (stationary or locally stationary) is chosen on both the test set and the second validation set.
This validates Theorem~\ref{thm:main3} for the example.
It is interesting to compare the observed proportions with the corresponding value of $q(\delta)$, which we provide in Table~\ref{MSPEanalysisperiodic1d}. We see that a larger proportion typically goes along with a larger value of $q(\delta)$ indicating the relevance of condition~\eqref{cond:f}.
To make it more precise: the tables are concerned with the proportion for which the decision rule~\eqref{decisionrule} yields the same result no matter if we take $i=2$ or $i=3$, i.\,e. we count what proportion of runs satisfies
\begin{equation} \label{samedecision}
(\hat R_{T,2}(h) \geq 1 + \delta \text{ and }
		\hat R_{T,3}(h) \geq 1 + \delta)
	\text{ or }
		(\hat R_{T,2}(h) < 1 + \delta \text{ and }
		\hat R_{T,3}(h) < 1 + \delta).
\end{equation}
We see that if $\delta$ is chosen large enough then the probability for the event~\eqref{samedecision} approaches~$1$, as $T$ and $m$ increase. More precisely, this is the case, if $\delta$ is chosen smaller than the ratio of MSPEs depicted in Figure~\ref{MSPEperiodic1} on both the validation and test set or larger than both those ratios. This is as expected from Corollary~\ref{thm:main3:kor1} and~\ref{thm:main3:kor2}. A more detailed analysis is possible, employing the information provided in Table~\ref{MSPEanalysisperiodic1c}.
In the third row of tables we see, for example, that for $n=10000$ and $\delta=0$ the procedure will consistently choose the locally stationary approach on both the test set and the second validation set for $1$-step ahead forecasting. For $n=10000$ and $\delta=0.05$, on the other hand, we see that the procedure almost consistently chooses the locally stationary approach on the validation set while it is rather undecided (50\%-50\%) on the test set. For $\delta=0.1$ the procedure almost consistently chooses the stationary approach on the test set and is to some degree undecided (70\%-30\%) on the second validation set. Finally, if $\delta=0.15$, we see that the stationary approach gets chosen almost consistently on both validation and test sets. This is just what we would expect, as a smaller $\delta$ must lead to the locally stationary approach being preferred, as the more complex locally stationary approach only gets selected if the empirical MSPE of the stationary approach is at least $(1+\delta)$-times of the empirical MSPE of the locally stationary approach.

The remaining part of the simulation studied is deferred to Section~\ref{a:extraSim} \citep{KleyEtAl2019}.

\section{Data examples}
\label{data_example}

\subsection{London Housing Prices}
\label{data_example_ukhpi}

\begin{figure}
    \begin{center}
     \includegraphics[width=\linewidth]{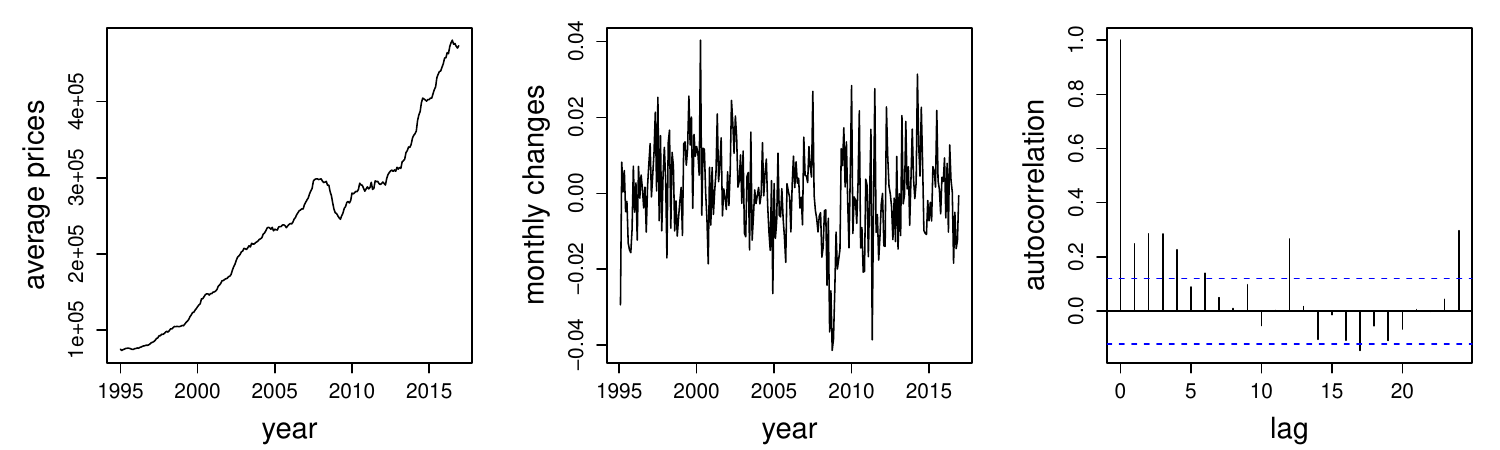}
    \end{center}
    \caption{Data for London from the UK House Price Index. Left: average monthly housing prices. Middle: monthly changes of average housing prices in percent, demeaned by subtracting arithmetic mean. Right: autocorrelation function, computed from the sequence in the middle.
    \label{hp_data}}
\end{figure}

We analyse average housing prices from the UK House Price Index (HPI). The HPI is updated monthly with data from the Land Registry, the Registers of Scotland, and the Land and Property Services Northern Ireland. The data is combined by the Office of National Statistics using hedonic regression; cf. \cite{HPI}.
The sequence we used for the analysis contains 264 monthly index values from 1995 to 2016. It was obtained as follows: In the `customise your search' part of the `search the UK house price index' form we have selected the `English region' London, the period from 01-1995 to 12-2016, and then obtained the `average price' for `all property types'.
The data is depicted in the left panel of Figure~\ref{hp_data}. For the analysis we consider $T+m=263$ monthly changes (in percent). The prices are centred by subtracting the arithmetic mean prior to the analysis. We clearly see autocorrelation at lags less or equal than 4 and at lag 12 in the right panel of Figure~\ref{hp_data}.

We then compute the 1-step to 6-step prediction coefficients, defined in~\eqref{def:vhat}, with which we can predict an observation $X_{t+h}$ from $X_t, \ldots, X_{t-p+1}$, where $X_{t+h}$ is an observation made either in 2014, 2015 or 2016, respectively. We choose $p=0,1,\ldots,18$, where $p=0$ shall mean that we are predicting with 0. Note that the maximum $p$ was chosen larger than 12, as we are dealing with monthly data and dependence at lag 12 can be seen from the autocorrelation function. We consider the stationary predictors as well as locally stationary predictors with $N=50, 51, \ldots, 87 = \lceil 263^{4/5} \rceil$.

\begin{figure}
    \begin{center}
     \includegraphics[width=\linewidth]{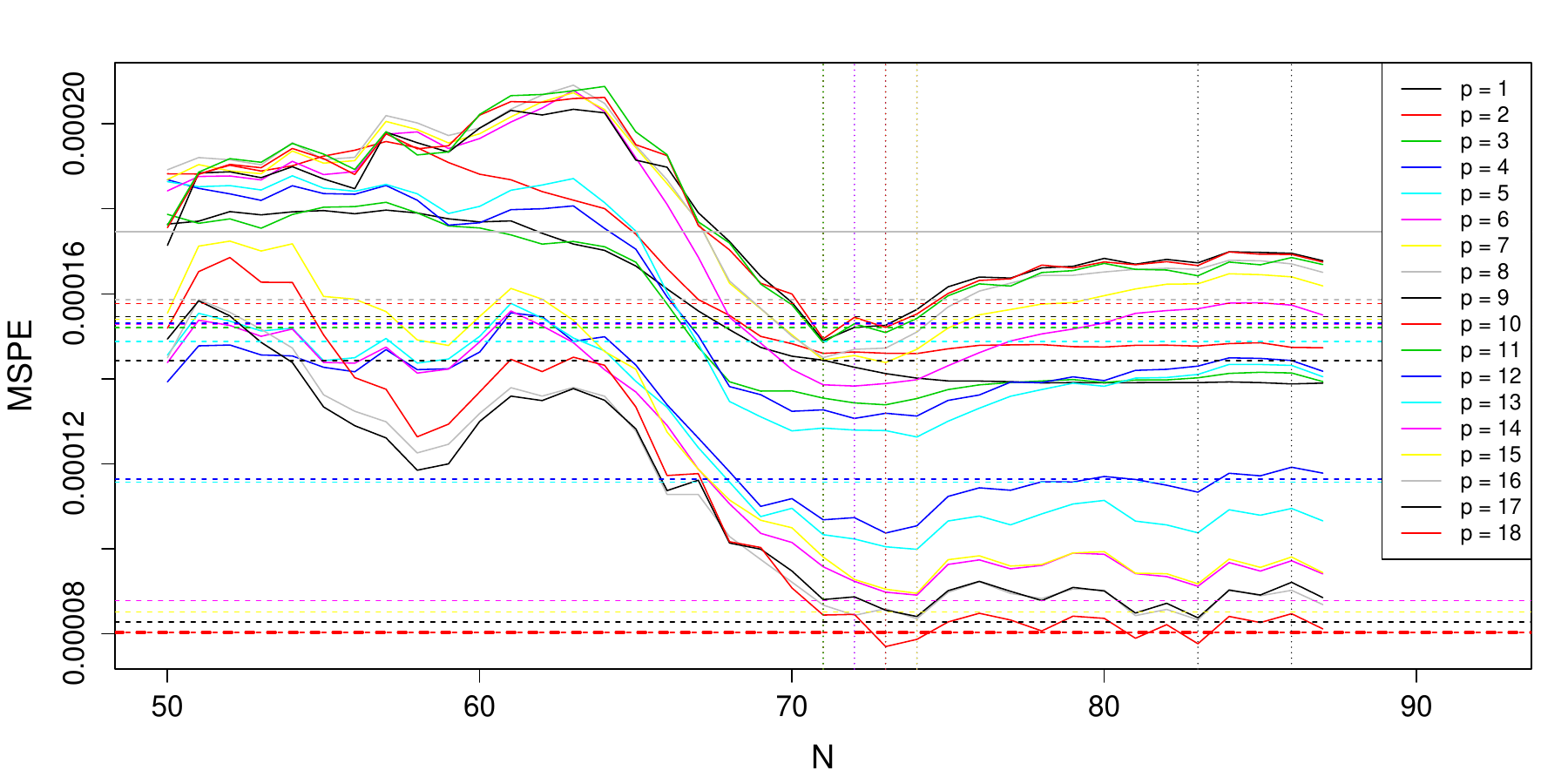}
     \includegraphics[width=\linewidth]{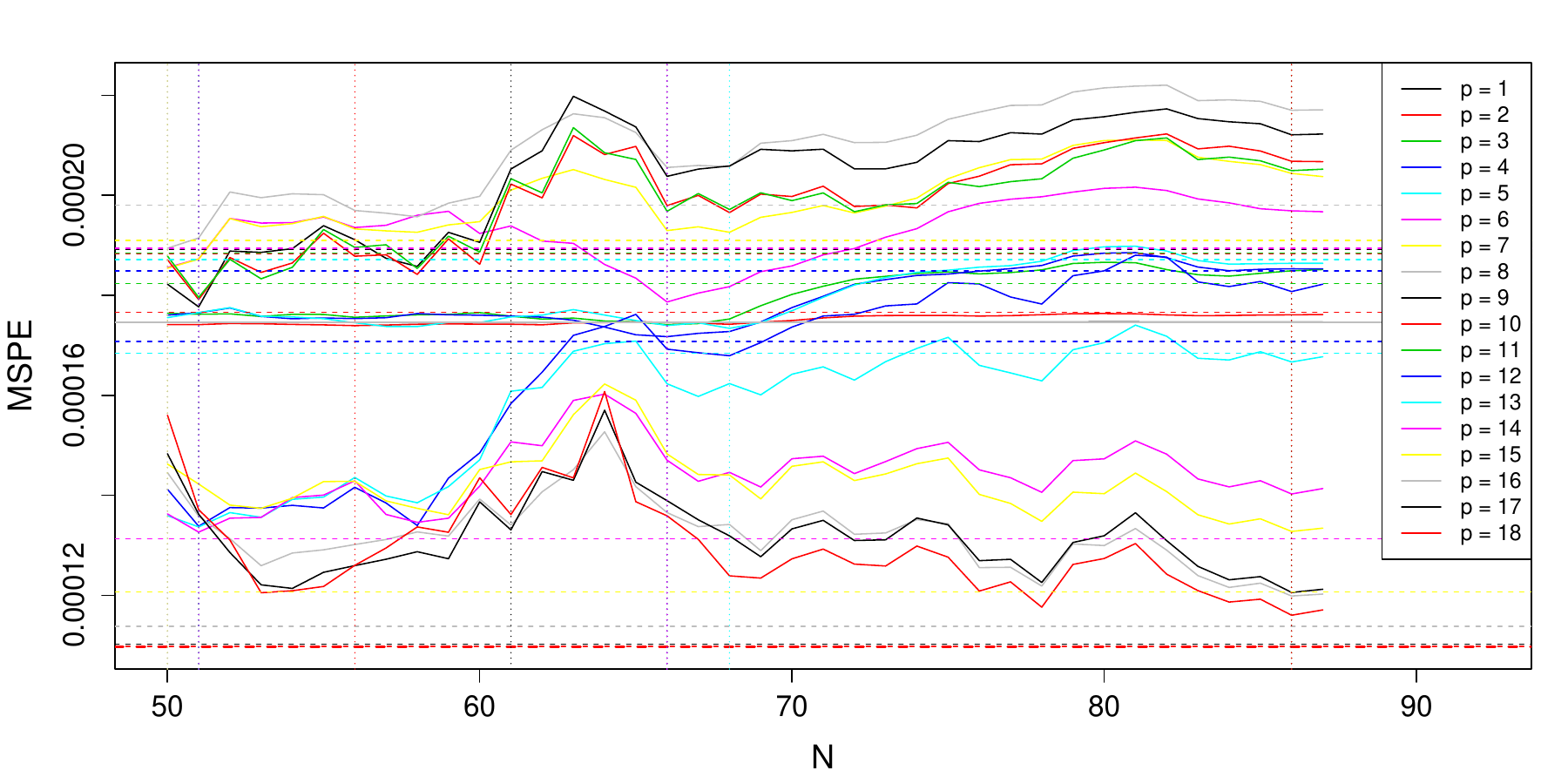}
    \end{center}
    \caption{Empirical mean squared prediction errors (MSPEs) computed on the frist validation set (predicting the 12 observations from 2014). Top panel shows MSPEs for 1-step ahead prediction. Bottom panel shows MSPEs for 6-step ahead prediction. The colours indicate which $p$ was used. The colour code is described in the plot's legend. The solid lines correspond to the MSPEs for different $N$ when the locally stationary approach is used. The dashed lines show the MSPE when the stationary approach is used. The horizontal grey line indicates the MSPE for the trivial forecasts ($f^{\lstat}_{t,h;0,N}$ and $f^{\stat}_{t,h;0}$). The MSPE in this case is $0.000175$.
    \label{hp_data_mspe}}
\end{figure}

Interestingly, in Figure~\ref{hp_data_mspe}, we observe that the MSPE of the locally stationary forecasts are typically larger than corresponding ones of the stationary forecasts.

As described in our procedure we now determine the $\hat p_{\stat}$, $\hat p_{\lstat}$, and $\hat N$ that minimise the MSPE within each class of predictors. For $1$-step ahead prediction we find
$\hat p_{\stat} = 18$, $\hat p_{\lstat} = 18$, and $\hat N = 73$.
For $6$-step ahead prediction we find
$\hat p_{\stat} = 18$, $\hat p_{\lstat} = 18$, and $\hat N = 86$. The numbers are summarised in Table~\ref{hp_numbers}.

\begin{table} \scriptsize
\begin{center}
\begin{tabular}{|c||c|c|c|c|c|}
\hline
 $h$ & $\hat p_{\stat}$ & ${\rm MSPE}_{T,1}^{\stat}(h)$ & $\hat p_{\lstat}$ & $\hat N_{\lstat}$ & ${\rm MSPE}_{T,1}^{\lstat}(h)$ \\
\hline \hline
1 &  18 & 8.033024e-05 &  18  & 73 & 7.701586e-05 \\
2 &  18 & 8.547987e-05 &  18  & 72 & 9.027318e-05 \\
3 &  18 & 9.362087e-05 &  18  & 71 & 9.512262e-05 \\
4 &  18 & 1.079008e-04 &  18  & 71 & 1.039368e-04 \\
5 &  18 & 1.164369e-04 &  18  & 87 & 1.291897e-04 \\
6 &  18 & 1.097551e-04 &  18  & 86 & 1.160201e-04 \\
\hline
\end{tabular}
\hspace*{-0.25cm}\begin{tabular}{|c||c|c|c||c|c|c|}
\hline
 $h$ & ${\rm MSPE}_{T,2}^{\stat}(h)$ & ${\rm MSPE}_{T,2}^{\lstat}(h)$ & $\hat R_{T,2}(h)$ & ${\rm MSPE}_{T,3}^{\stat}(h)$ & ${\rm MSPE}_{T,3}^{\lstat}(h)$ & $\hat R_{T,3}(h)$\\
\hline \hline
1 &  3.473298e-05 & 3.501655e-05 & 0.992 &  9.740925e-05 & 0.0001385059 & 0.703 \\
2 &  3.560845e-05 & 4.308688e-05 & 0.826 &  9.547598e-05 & 0.0001351634 & 0.706 \\
3 &  4.31916e-05  & 4.21518e-05  & 1.025 &  0.0001052688 & 0.0001309526 & 0.804 \\
4 &  4.57004e-05  & 4.429208e-05 & 1.032 &  0.0001053983 & 0.0001421635 & 0.741 \\
5 &  5.970928e-05 & 4.943228e-05 & 1.208 &  0.0001210628 & 0.0001195622 & 1.012 \\
6 &  6.412237e-05 & 5.234349e-05 & 1.225 &  0.0001152908 & 0.0001146555 & 1.006 \\
\hline
\end{tabular}
\caption{\textit{Minimum empirical mean squared prediction errors (MSPEs) for $h$-step ahead prediction, $h=1,\ldots,6$, of the house price data. Top table shows values computed on the first validation set. Bottom table shows values computed on the second validation set and on the test set.}}  \label{hp_numbers}
\end{center}
\end{table}

We then determine the MSPE for forecasting the observations from the second validation set (here: the year 2015) using these predictors. For $1$-step ahead prediction we find that ${\rm MSPE}_{251,2}^{\stat}(1) = 3.47 \cdot 10^{-5}$ and ${\rm MSPE}_{251,2}^{\lstat}(1) = 3.50 \cdot 10^{-5}$, with ${\rm MSPE}_{T,j}^{*}(h)$ defined in~\eqref{eqn:empMSPE}. For $6$-step ahead prediction we find that ${\rm MSPE}_{251,2}^{\stat}(6) = 6.41 \cdot 10^{-5}$ and ${\rm MSPE}_{251,2}^{\lstat}(6) = 5.23 \cdot 10^{-5}$. Consequently, we decide to use the stationary approach for $1$-step and the locally stationary approach for $6$-step ahead forecasting of the observations made in 2016.

The MSPEs computed from $1$-step ahead forecasting the observations from the test set (here: the year 2016) are ${\rm MSPE}_{251,3}^{\stat}(1) = 9.74 \cdot 10^{-5}$ and ${\rm MSPE}_{251,3}^{\lstat}(1) = 1.39 \cdot 10^{-4}$. The MSPEs computed from $6$-step ahead forecasting the observations from 2016 are ${\rm MSPE}_{251,3}^{\stat}(6) = 1.153 \cdot 10^{-4}$ and ${\rm MSPE}_{251,3}^{\lstat}(6) = 1.147 \cdot 10^{-4}$. We have thus chosen the better performing procedure for 1-step and 6-step ahead forecasting.

In conclusion, our analysis has revealed that, from the point of view of 1-month ahead prediction of the 2016 observations, treating the data as stationary does not have a negative effect.  We were able to see that using the estimates from the stationary AR(18) model gave us better predictions than using the (locally stationary) estimates of segments of 73 month (roughly 6 years). For the 6-month ahead prediction the local estimates are better, but only by a small margin. Contrary to what one might naively expect, the impact of, for example, the 2008-2009 financial crisis on the stationary estimates is not profound enough to substantially worsen the predictors' performance. 

\subsection{Temperatures Hohenpei{\ss}enberg}
\label{data_example_hpb}

\begin{figure}
    \begin{center}
     \includegraphics[width=\linewidth]{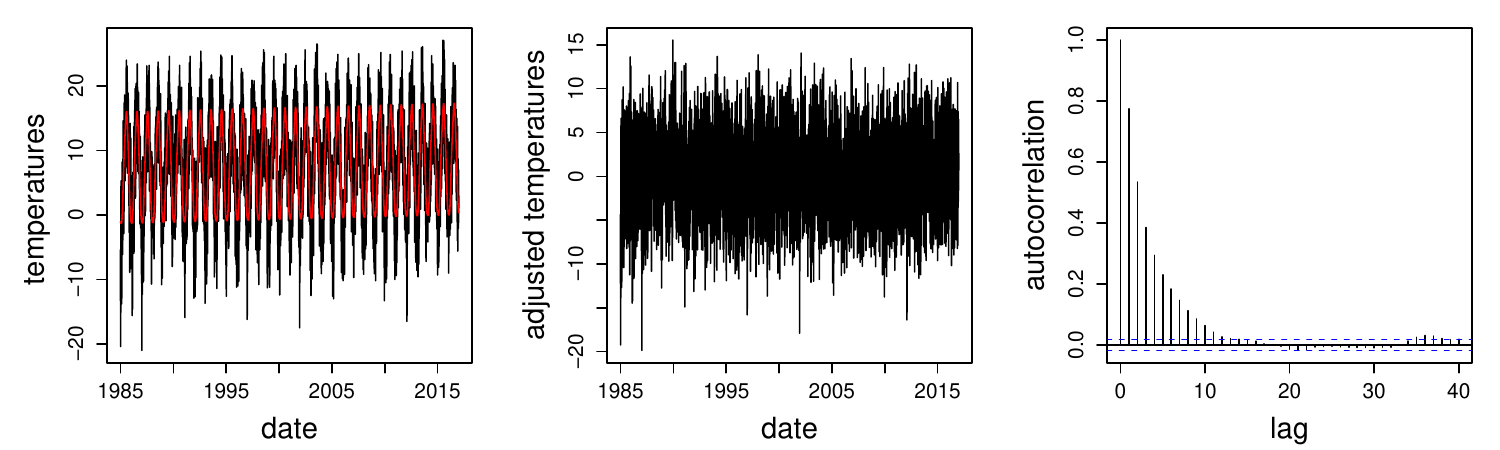}
    \end{center}
    \caption{Temperature data Hohenpei{\ss}enberg. Left: daily temperatures and fitted harmonic regression model. Middle: adjusted data (demeaned and detrended). Right: autocorrelation function, computed from the sequence in the middle.
    \label{hpb_data}}
\end{figure}

In this example, we analyse seasonally adjusted, daily temperature data collected at the meteorological observatory in Hohenpei{\ss}enberg (Germany).
More precisely, we use $n = T+m = 11680$ observations of daily mean temperatures that were recorded between 1985 and 2016.%
\footnote{The data was obtained from \url{http://www.dwd.de/DE/klimaumwelt/cdc/cdc_node.html}.}
The data are shown in the left panel of Figure~\ref{hpb_data}. To eliminate the clearly visible trend and seasonality, we have fitted a harmonic linear regression model of the form
\[ y_t = c + \alpha t + \sum_{i=1}^4 \Big( \beta_i \sin(2\pi t i/365) + \gamma_i \cos(2\pi t i/365) \Big), \]
to capture the trend and annual variation. The red curve in the left panel of Figure~\ref{hpb_data} is the prediction of the fitted model. We then consider the residuals of this model which are shown in the middle panel of Figure~\ref{hpb_data}. The right panel of Figure~\ref{hpb_data} shows the autocorrelation function, which clearly indicates that serial dependence is present. \cite{BirrEtAl2016} analyse the same data set and fit a stationary ARMA(3,1) model to capture the serial dependence.

\begin{figure}
    \begin{center}
     \includegraphics[width=\linewidth]{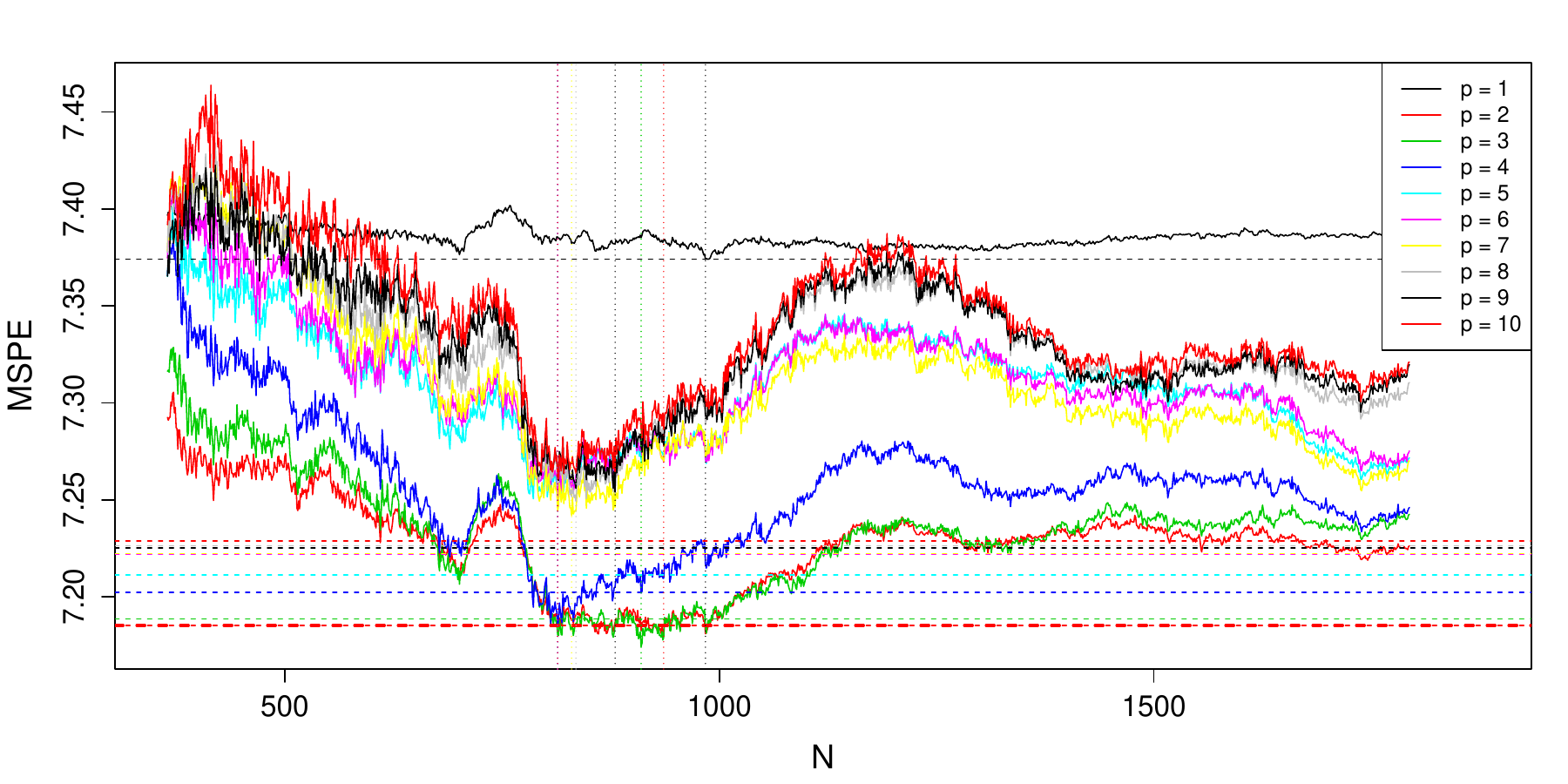}
     \includegraphics[width=\linewidth]{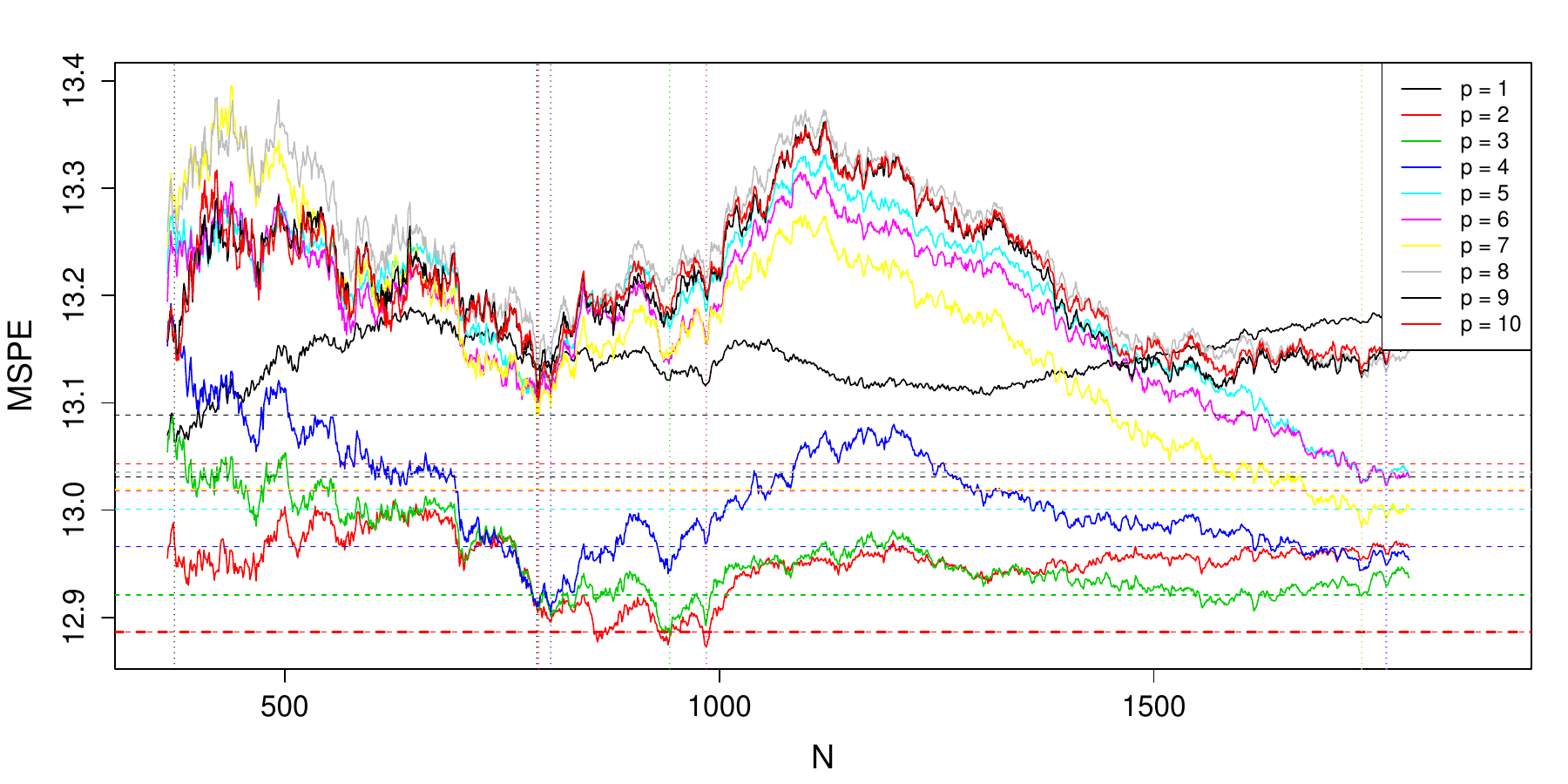}
    \end{center}
    \caption{Empirical mean squared prediction errors (MSPEs) computed on the first validation set (predicting the 365 observations from 2014) of the temperature data. Top panel shows MSPEs for 1-step ahead prediction. Bottom panel shows MSPEs for 2-step ahead prediction. The colours indicate which $p$ was used. The colour code is described in the plot's legend. The solid lines correspond to the MSPEs for different $N$ when the locally stationary approach is used. The dashed lines show the MSPE when the stationary approach is used.
    \label{hpb_data_mspe}}
\end{figure}

In Figure~\ref{hpb_data_mspe}, the MSPE are presented in the same manner as in Section~\ref{data_example_ukhpi}. In this example we have chosen $p_{\max} = 10$ and $\mathcal{N} := \{365, 366 \ldots, \lceil n^{4/5} \rceil\} = \{365, 366, \ldots, 1794, 1795\}$ and $m := 365$. The MSPE corresponding to $p=0$ is 110.2 in this example and therefore not visible in the plot.

By minimising the empirical MSPE on the first validation set the procedure chooses, for the stationary approach $\hat p_{\stat}=2$ for $h=1,2$. For the locally stationary approach the procedure chooses $(\hat p_{\lstat}, \hat N_{\lstat})=(3,910)$ and $(\hat p_{\lstat}, \hat N_{\lstat})=(2,985)$ for $h=1$ and $h=2$, respectively. Empirical MSPEs for other values of $p$ and $N$ are shown in Figure~\ref{hpb_data_mspe}. The numbers are summarised in Table~\ref{hpb_numbers}.

\begin{table} \scriptsize
\begin{center}
\begin{tabular}{|c||c|c|c|c|c|}
\hline
 $h$ & $\hat p_{\stat}$ & ${\rm MSPE}_{T,1}^{\stat}(h)$ & $\hat p_{\lstat}$ & $\hat N_{\lstat}$ & ${\rm MSPE}_{T,1}^{\lstat}(h)$ \\
\hline \hline
1 &   2 &  7.185208 &    3 &  910 &  7.173272 \\
2 &   2 & 12.886257 &    2 &  985 & 12.870544 \\
3 &   2 & 15.397509 &    2 &  870 & 15.343298 \\
4 &   2 & 16.605640 &    2 &  800 & 16.504915 \\
5 &   2 & 17.226943 &    2 &  800 & 17.093823 \\
\hline
\end{tabular}
\hspace*{-0.25cm}\begin{tabular}{|c||c|c|c||c|c|c|}
\hline
 $h$ & ${\rm MSPE}_{T,2}^{\stat}(h)$ & ${\rm MSPE}_{T,2}^{\lstat}(h)$ & $\hat R_{T,2}(h)$ & ${\rm MSPE}_{T,3}^{\stat}(h)$ & ${\rm MSPE}_{T,3}^{\lstat}(h)$ & $\hat R_{T,3}(h)$\\
\hline \hline
1 &  8.10974 & 8.058095  & 1.006 &  8.08899 & 7.967895 & 1.015 \\
2 &  14.86848 & 14.94354 & 0.995 &  15.42535 & 15.39907 & 1.001 \\
3 &  17.72551 & 17.92775  & 0.989 &  17.4254 & 17.3617 & 1.004 \\
4 &  19.63724 & 19.8143 & 0.991 &  17.68487 & 17.60241 & 1.005 \\
5 &  20.97236 & 21.0989 & 0.994 &  17.92979 & 17.83498  & 1.005 \\
\hline
\end{tabular}
\caption{\textit{Minimum empirical mean squared prediction errors (MSPEs) for $h$-step ahead prediction, $h=1,2,3,4,5$, of the temperature data Hohenpei{\ss}enberg. Top table shows values computed on the first validation set. Bottom table shows values computed on the second validation set and on the test set.}}  \label{hpb_numbers}
\end{center}
\end{table}

For $1$-step ahead forecasting and on validation 2 set this yields, given the $\hat p_{\stat}$ chosen by the procedure, that ${\rm MSPE}^{\stat}_{11315,2}(1) = 8.11$ and, given the  $\hat p_{\lstat}$ and $N$ chosen by the procedure, that ${\rm MSPE}^{\lstat}_{11315,2}(1) = 8.06$. Similarly, for $2$-step ahead forecasting, we have ${\rm MSPE}^{\stat}_{11315,2}(2) = 14.87$ and ${\rm MSPE}^{\lstat}_{11315,2}(2) = 14.94$. The respective ratios are both very close to $1$. The procedure thus chooses the stationary approach over the locally stationary approach if $\delta = 0.01$ is chosen and, obviously, this superiority will continue to hold if $\delta $ is chosen larger than that.
On the test set we have ${\rm MSPE}^{\stat}_{11315,3}(1) = 8.09$ and ${\rm MSPE}^{\lstat}_{11315,3}(1) = 7.97$ for $1$-step ahead forecasting. Likewise, for $2$-step ahead forecasting, we have ${\rm MSPE}^{\stat}_{11315,3}(2) = 15.43$ and ${\rm MSPE}^{\lstat}_{11315,3}(2) = 15.41$. Thus, again, both approaches for $1$-step and $2$-step ahead forecasting behave almost equally well and we see that had we chosen $\delta > 0.015$ our procedure chose the stationary approach, which performs almost equally well as the more complicated locally stationary approach.

In conclusion, in this example, we have provided clear evidence that the temperature data, after adjusting for trend and seasonality, collected in the Hohenpei{\ss}enberg observatory, from the point of view of prediction, can be treated as if they were stationary. We see that using the estimates related to a AR(2) [or AR(3)] model yielded forecasts that in all cases perform almost equally well as the estimates localised to the segment suggested by the procedure (using the past 2.2--2.7 years; 800--910 days). This observation is remarkable, in the sense that, in 30 years of data an analyst might typically expect non-stationarity (e.\,g., changes due to global warming) to worsen the predictions. Our conclusion indicates that the variation of covariance structure might be less substantial than the change in mean. Note that our procedure did not consistently chose the approach with the better performance on the test set, but that both approaches perform almost equally well on either set. It is thus legitimate to use the simpler, stationary approach.

\subsection{Volatility around the time of the EU referendum in the UK, 2016}
\label{data_example_brexit}

\begin{figure}
    \begin{center}
		 \includegraphics[width=\linewidth]{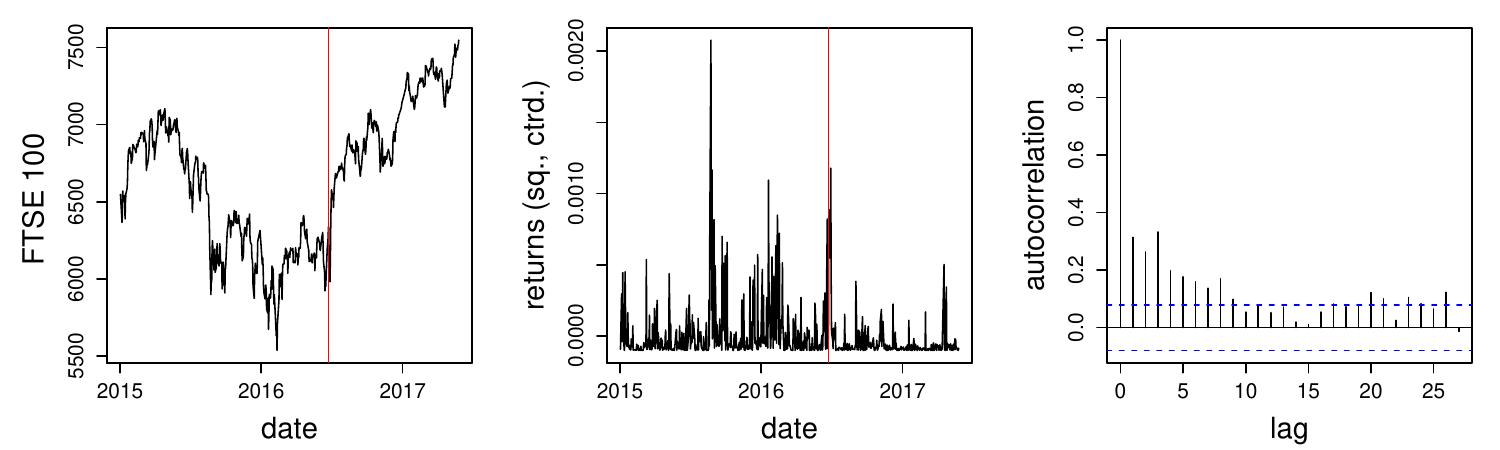}
    \end{center}
    \caption{Volatility of the FTSE 100 Index, for 2 January 2015 to 26 May 2017. Left: FTSE 100 closing price. Middle: squared and centred returns. Right: autocorrelation function, computed from the sequence in the middle. Red vertical line in the left and middle plot marks 23/06/2016, the day of the EU referendum in the UK.
    \label{ftse_data}}
\end{figure}

This example is about forecasting volatility of the FTSE 100 stock index. More precisely, we consider a sequence of $n=T+m=607$ (daily) opening prices $p_{\rm open}$ and closing prices $p_{\rm close}$, dated from 2 January 2015 to 26 May 2017.%
\footnote{The data was obtained from \url{http://www.finanzen.net/index/FTSE_100/Historisch}.}
The analysis is then based on the sequence $ ( (p_{\rm close} - p_{\rm open})/p_{\rm close})^2$, centred by subtracting the arithmetic mean of this sequence. The data are shown in Figure~\ref{ftse_data}.

We separate the final 60 observations of the data as test set, first validation set and second validation set (used for determining the model orders and segment sizes). Each set is of size $m := 20$. Visual inspection of these 60 observations suggested that some returns are unusually small or large. Indeed, the returns of 1~March, 18~April, and 24~April 2017 are either more than 1.5 times the interquartile range (IQR) smaller than the lower quartile or 1.5 times the IQR larger than the upper quartile. By Tukey's criterion they can thus be classified as outliers. To better deal with the outliers, we use a robustified measure of accuracy to compare the forecasts in this example. More precisely, instead of the MSPE in Steps 4 and 5 of our procedure, we now use an empirical trimmed mean of absolute prediction errors (trMAPE), where we trim the largest $25\%$, averaging only the remaining 15 out of 20 absolute errors. We have further chosen $p_{\max} = 8$ and $\mathcal{N} := \{40, 41, \ldots, 250\}$.

\begin{figure}
		\vspace{-2cm}
    \begin{center}
     \includegraphics[width=\linewidth]{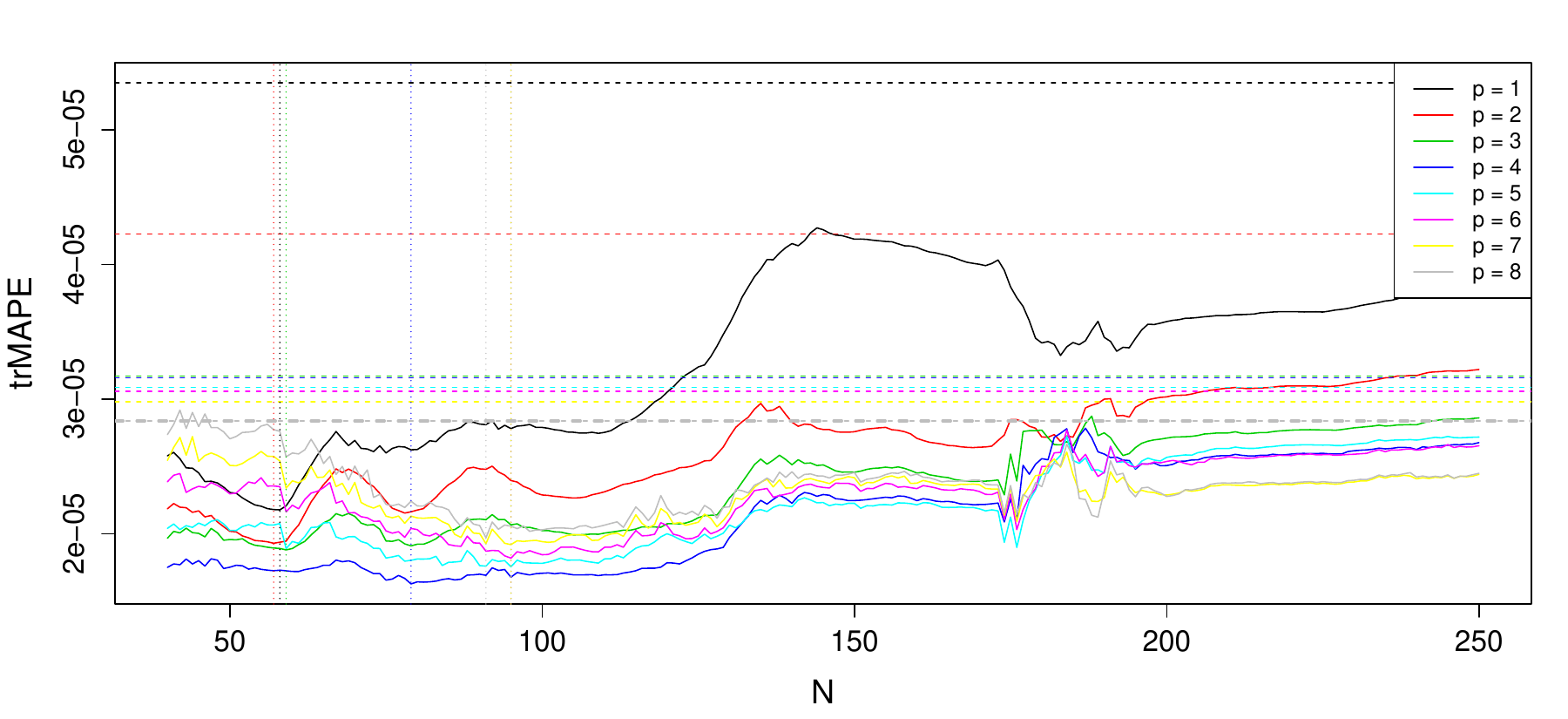}
     \includegraphics[width=\linewidth]{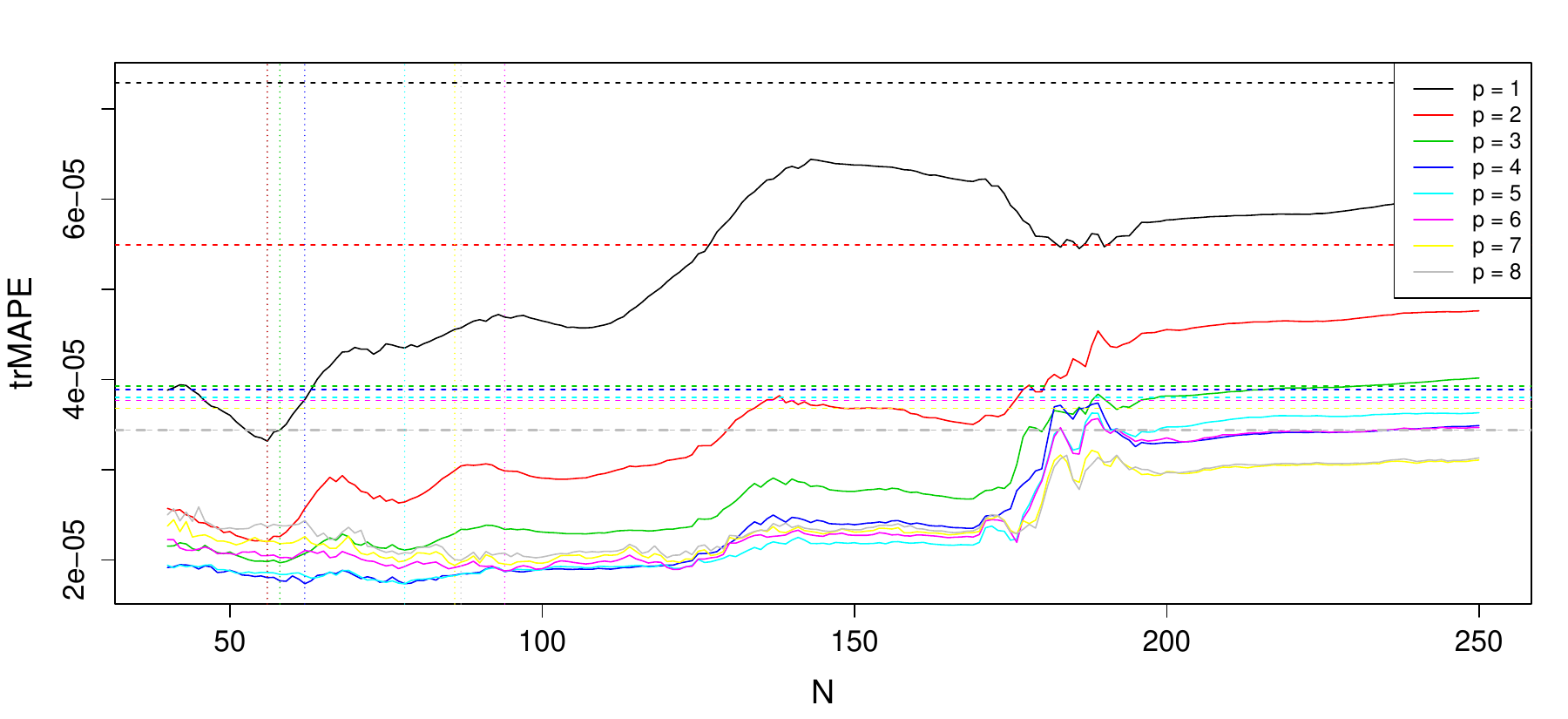}
		 \includegraphics[width=\linewidth]{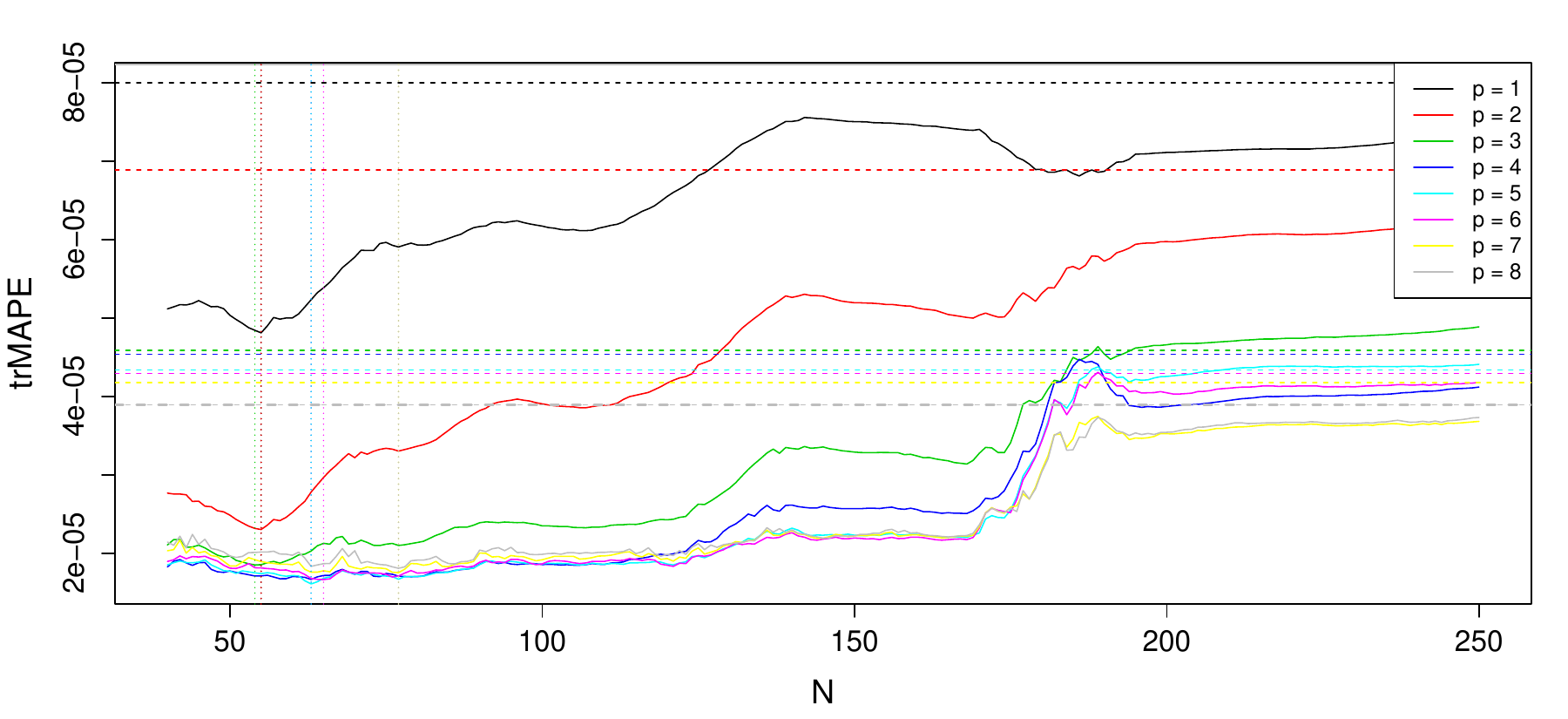}
    \end{center}
    \caption{Empirical trimmed mean absolute prediction errors (trMAPE) computed on the first validation set (predicting the observations 548 to 567) of the squared and centred FTSE returns. Top, middle and bottom panel show the trMAPE for the 1, 2 and 3-step ahead predictions, respectively. The colours indicate which $p$ was used. The colour code is described in the plot's legend. The solid lines correspond to the trMAPE for different $N$ when the locally stationary approach is used. The dashed lines show the trMAPE when the stationary approach is used. The horizontal grey line indicates the trMAPE for the trivial forecasts ($f^{\lstat}_{t,h;0,N}$ and $f^{\stat}_{t,h;0}$). The trMAPE in this case is $8.2 \times 10^{-5}$.
    \label{ftse_data_mspe}}
\end{figure}

First, we consider the trMAPEs of forecasting the 20 observations from the first validation set to determine the optimal $\bar p_{\stat}$, $\bar p_{\lstat}$ and $\bar N_{\lstat}$. We use a bar instead of the hat to indicate that the trMAPEs were used instead of the MSPEs. In Figure~\ref{ftse_data_mspe} we can see, for the 1-step, 2-step and 3-step ahead forecasts, that the lines depicting the trMAPEs have a characteristic shape: as $N$ increases the trMAPEs slightly decreases (for each $p$ at a different level) until it starts increasing around $N \approx 60$. After this follows another phase of slight decreasing and increasing with the new minimum higher than the minimum of the previous phase. We further observe that the overall level is typically lower than that of the trMAPEs of the stationary approach. The last such minimum in our plots is obtained when $N$ is around 170--180.

The observations 1 through to 373 were recorded from 2~January 2015 to 23~June 2016 (the day of the EU referendum) and observations 374 through to 607 were recorded from 24~June 2016 to 26~May 2017. This implies that the final 234 observations were recorded after the EU referendum, meaning that there are 175 observations between the EU referendum and the observations to be forecast in the first step. Thus, the last minimum of the lines, when $N$ is roughly about 170, corresponds to the time of the referendum. The sudden increase of the trMAPE indicates the change in bias of the Yule Walker estimator due to non-stationarity when pre-referendum data is starting to be used for the estimation of the prediction coefficients. Another important observation is that also the post-referendum part of the diagram ($40 \leq N \leq 175$) shows signs of non-stationarity.
More specifically, each phase of up-movement indicate that the variance is reduced less than the squared bias increases. The increase from the first (and global) minimum at around $N \approx 60$ onwards corresponds to taking data from the end of November 2016 and earlier into account and might correspond to changes due to effects of the election of the US president. The minimum trMAPE for forecasting the data from the end of the estimation set are summarised in Table~\ref{ftse_numbers}. We observe that for $h=1,2,3,4$ the optimum segment size is roughly 60 such that no observations prior to November 2016 are used for estimation. For $h=5$ the optimum segment size is 41 and thus even smaller. This implies that no observation prior to the presidential election in the US are used for estimation of the forecasting coefficients.

\begin{table} \scriptsize
\begin{center}
\begin{tabular}{|c||c|c|c|c|c|}
\hline
 $h$ & $\bar p_{\stat}$ & ${\rm trMAPE}_{T,1}^{\stat}(h)$ & $\bar p_{\lstat}$ & $\bar N_{\lstat}$ & ${\rm trMAPE}_{T,1}^{\lstat}(h)$ \\
\hline \hline
1 &   8 & 2.838118e-05  &   4  &  79 & 1.628012e-05 \\
2 &   8 & 3.440985e-05  &   5  &  78 & 1.736197e-05 \\
3 &   8 & 3.892572e-05  &   5  &  63 & 1.610966e-05 \\
4 &   8 & 4.786001e-05  &   6  &  65 & 1.731724e-05 \\
5 &   8 & 5.161963e-05  &   5  &  53 & 2.209272e-05 \\
\hline
\end{tabular}
\hspace*{-0.25cm}\begin{tabular}{|c||c|c|c||c|c|c|}
\hline
 $h$ & ${\rm trMAPE}_{T,2}^{\stat}(h)$ & ${\rm trMAPE}_{T,2}^{\lstat}(h)$ & $\bar R_{T,2}(h)$ & ${\rm trMAPE}_{T,3}^{\stat}(h)$ & ${\rm trMAPE}_{T,3}^{\lstat}(h)$ & $\bar R_{T,3}(h)$\\
\hline \hline
1 & 2.838118e-05 & 1.628012e-05 & 1.743 & 3.33206e-05  & 2.395498e-05 & 1.391 \\
2 & 3.440985e-05 & 1.736197e-05 & 1.982 & 3.817851e-05 & 2.505271e-05 & 1.524 \\
3 & 3.892572e-05 & 1.610966e-05 & 2.416 & 4.369565e-05 & 2.750278e-05 & 1.589 \\
4 & 4.786001e-05 & 1.731724e-05 & 2.764 & 4.974355e-05 & 2.81844e-05  & 1.765 \\
5 & 5.161963e-05 & 2.209272e-05 & 2.336 & 5.390384e-05 & 4.560496e-05 & 1.182 \\
%
%
%
\hline
\end{tabular}
\caption{\textit{Minimum empirical trimmed mean absolute prediction errors (trMAPE) for $h$-step ahead prediction, $h=1,2,3,4,5$, of the squared and centred FTSE~100 data. Analysis performed with $m := 20$ and $p_{\max} = 8$. Top table shows values computed on the first validation set. Bottom table shows values computed on the second validation set and on the test set.}}  \label{ftse_numbers}
\end{center}
\end{table}

Using these predictors to forecast the 20 observations from the second validation set we see, in Table~\ref{ftse_numbers}, that the trMAPE of the stationary approach are typically (2.1 to 3.7 times) larger than the trMAPE of the locally stationary approach. We thus choose to work with the locally stationary approach. In Table~\ref{ftse_numbers} we denote the ratios of the trMAPE of the stationary approach over the trMAPE of the locally stationary approach by $\bar R_{T,j}(h)$, where the bar indicates that the trMAPE and not the MSPE is used.
Forecasting the 20 observations from the test set we see that the trMAPEs of the stationary approach are again larger than those of the locally stationary approach, but not quite as much as on the second validation set. Still, following our procedure, we chose the better performing approach (the locally stationary one).

For this example, we further conducted a sensitivity analysis, by varying the parameters $m$ and $p_{\max}$. Selected results, in which we see the results are mostly stable when changing the parameters are shown in Appendix~\ref{a:extraEmpExpl} \citep{KleyEtAl2019}.

\section{Analysis of the localised Yule-Walker estimator under general conditions and local stationarity}
\label{sec:propertiesYW}

In this section we discuss the probabilistic properties of the localised Yule-Walker estimator $\hat a_{N,T}^{(p)}(t)$ defined in~\eqref{eqn:YW}. We believe the results to be of independent interest and therefore present them in this separate section. They are also key results for the proofs of the result in Section~\ref{theorysectionar1}. Our results will hold under Assumptions~\ref{a:loc_stat}--\ref{a:MomentCond} (cf. Section~\ref{theorysectionar1}). The assumptions are not restrictive and, in particular, the concentration result in this section will hold for a broad class of locally stationary processes and, in particular, does not require that the data come from a tvAR($p$) model. Further, we allow for any $1+p \leq N \leq T$ and, in particular, allow for a diverging model order $p$, as $T \rightarrow \infty$. We do not, as do for example~\cite{dahlgir1998}, require that $N = o(T)$.

The main result of this section (Theorem~\ref{thm:propertiesYW}) provides a non-asymptotic bound for the Euclidean distance of $\hat a_{N,T}^{(p)}(t)$ to the following population quantity:
\begin{equation}
\begin{split}
	\label{eqn:YW_bar}
		\bar a_{N,T}^{(p)}(t) & := \big( \E \hat\Gamma_{N,T}^{(p)}(t) \big)^{-1} \big( \E \hat\gamma_{N,T}^{(p)}(t) \big)
			= \big( \bar a_{1,N,T}^{(p)}(t), \ldots, \bar a_{p,N,T}^{(p)}(t) \big)'.
		\end{split}
\end{equation}

The Yule-Walker estimator is widely used in practice and $\hat a_{N,T}^{(p)}(t)$ and its properties have been studied in detail under various conditions. \cite{BercuEtAl1997,BercuEtAl2000} and \cite{Bercu2001} derive large deviation principles for Gaussian AR processes when the model order is 1. A simple exponential inequality, also for model order 1, is given in Section 5.2 of \cite{BercuTouati2008}. \cite{YuSi2009} prove a large deviation principle for general, but fixed, model order. \cite{Jirak2012,Jirak2014} derives simultaneous confidence bands. The cited results all require that the underlying process is stationary. \cite{dahlgir1998} analyse the bias and variance of the localised Yule-Walker estimator in the framework of local stationarity. They do not, however, provide an exponential inequality, and, as far as we are aware, no result as the one we provide below is available at present.
The exponential inequality in Theorem~\ref{thm:propertiesYW}, which we now state, is explicit in terms of all parameters and constants. We make use of the explicitness to derive Corollary~\ref{kor:propertiesYW}, by which the localised Yule-Walker estimator is strongly, uniformly consistent, even when the model order is diverging as the sample size grows.
\begin{satz}\label{thm:propertiesYW}
Let $(X_{t,T})_{t \in \IZ, T \in \IN^*}$ satisfy Assumptions~\ref{a:loc_stat}--\ref{a:MomentCond} and $\E X_{t,T} = 0$.
Then, for every $T \geq 2 C_1 p^2$, $N \geq 1 + p \geq 2$ and $\varepsilon > 0$, we have:
	\begin{equation*}
	\begin{split}
		& \IP \big( \| \hat a_{N,T}^{(p)}(t) - \bar a_{N,T}^{(p)}(t) \| > \varepsilon \big) \\
		& \quad \leq 3 p \exp\Bigg( - \frac{\Big( \frac{m_f}{4 p} \min\Big\{1, \varepsilon \frac{1}{8 C_0}\Big\} \Big)^2}{2\Big(C_{1,1} \frac{p}{N-p} + \Big( \frac{m_f}{4 p} \min\Big\{1, \varepsilon \frac{1}{8 C_0}\Big\} \Big)^{(3+4d)/(2+2d)} \Big(C_{2,1} \frac{p}{N-p} \Big)^{1/(2+2d)}\Big)} \Bigg) \\
	\end{split}
	\end{equation*}
		\begin{equation*}
	\begin{split}
		& \quad \leq \begin{cases}
		3 p \exp\Bigg( - \frac{m_f^2}{ 32 C_{1,1} \frac{p^3}{N-p} + m_f^{(3+4d)/(2+2d)} \Big(32 C_{2,1} \frac{p^2}{N-p} \Big)^{1/(2+2d)}\Big)} \Bigg)
		& \varepsilon \geq 1/(8 C_0) \\
		3 p \exp\Bigg( - \varepsilon^2 \frac{m_f^2}{2^{12} C_{1,1}} \Big( C_0^2 \frac{p^3}{N-p} \Big)^{-1} \Bigg)
			& \varepsilon \leq \min\{U_{p,N}, \frac{1}{8 C_0}\}  \\
		3 p \exp\Bigg( - \varepsilon^{1/(2+2d)} \Big(\frac{m_f}{2^{9+4d} C_{2,1}}\Big)^{1/(2+2d)} \Big( C_0 \frac{p^2}{N-p} \Big)^{-1/(2+2d)} \Bigg) 
			& \frac{1}{8 C_0} > \varepsilon \geq \min\{U_{p,N}, \frac{1}{8 C_0}\}
			\end{cases}
	\end{split}
	\end{equation*}
	where $\hat a_{N,T}^{(p)}(t)$ is defined in~\eqref{eqn:YW}, $\bar a_{N,T}^{(p)}(t)$ is defined in~\eqref{eqn:YW_bar},
	\begin{equation*}
		U_{p,N} := \frac{32 C_0}{m_f} \Big(\frac{C_{1,1}^{2+2d}}{C_{2,1}}\Big)^{1/(3+4d)} \Big( \frac{p^{(4+6d)/(3+4d)}}{(N-p)^{(1+2d)/(3+4d)}} \Big),
		\end{equation*}
	and $C_0$, $C_1$ and $C_{1,1}$, $C_{2,1}$ are defined in~\eqref{eqn:def_C0} and~\eqref{lem:exp_ineq_gamma_k2:def:C1C2}, respectively.
\end{satz}
The proof of Theorem~\ref{thm:propertiesYW} is deferred to Section~\ref{app:main} of the appendix.

Theorem~\ref{thm:propertiesYW} is a key ingredient to the proof of Lemma~\ref{lem:main} which is essential to the proof of the performance-guarantee-result (Theorem~\ref{thm:main3}) of our procedure. Further, it implies

\begin{kor}\label{kor:propertiesYW}
	Let $(X_{t,T})_{t \in \IZ, T \in \IN^*}$ satisfy Assumptions~\ref{a:loc_stat}--\ref{a:MomentCond}, $\E X_{t,T} = 0$ and let $P = P_{T}$ and $N = N_T$ be sequences of integers that satisfy
		$2 \leq 1+P \leq N \leq T$. Assume that $P = o(N^{(1+2d)/(4+6d)})$ and $N \rightarrow \infty$, as $T \rightarrow \infty$. Further, assume that there exists a sequence $R_T$ with $0 \leq R_T \rightarrow \infty$ and
		$R_T \log(T) = o\big( (N/P)^{1/(3+4d)} \big)$, as $T \rightarrow \infty$, where $d$ is the constant from Assumption~\ref{a:MomentCond}. Then, we have
	\[\sup_{p=1,\ldots,P} \sup_{t = N, \ldots, T} \| \hat a_{N,T}^{(p)}(t) - \bar a_{N,T}^{(p)}(t) \|
		= O\Bigg( P^{3/2} \Big( \frac{ \log(T)}{N} \Big)^{1/2} \Bigg), \text{ almost surely, as $T \rightarrow \infty$.}\]
\end{kor}

\begin{rem}\label{rem:propertiesYW}
	For any stationary AR($p$) model we have that $\bar a_{N,T}^{(p)}(u)$ corresponds to the vector of coefficients. This can be seen from Lemma~\ref{lem:rel_a_bar} and the fact that $C_1 = 0$ if the model is stationary. Thus, choosing $N_T = T$ and $P_T = p$, our result yields the same rate as Theorem~1 in \cite{lai1982}, by which the (least squares) estimator is strongly consistent with rate $(\log(T)/T)^{1/2}$. An early consistency result for the Yule-Walker estimate with diverging model order is Theorem~6 in \cite{Hong-ZhiEtAl1982}. Under the assumption that $P = O(\log(T)^a)$, $a > 1$ or $P = C \log T$, $C$ sufficiently large, they prove that the rate of convergence is $O\big( (\log \log T / T)^{1/2} \big)$.
\end{rem}

\section{Conclusion}\label{sec:conclusion}

In this paper, we have presented a method to choose between different forecasting procedures, based on the empirical mean squared prediction errors the procedures achieve. Using the empirical rather than the asymptotic mean squared prediction error, our procedure automatically takes into account that different models should be preferred depending on the amount of data available, which is an important difference to the Focused Information Criterion by \cite{clahjo2003}. Working in the general framework of locally stationary time series we choose from two classes of forecasts that were motivated by approximating the serial dependence of the time series by time-varying or traditional autoregressive models. The procedure implicitly balances the modelling bias (which is lower if the model is more complex) and the variance of estimation (which increases for more complex models). Our two step procedure automatically chooses the number of forecasting coefficients to be used and the segment size from which the forecasting coefficients are estimated.

In a comprehensive simulation study we have illustrated that it is often advisable to use a forecasting procedure derived from a simpler model when not a vast amount of data is available. In particular, in the tvAR models of our simulations, if the variation over time is not very pronounced and when the tangent processes are not close to being unit root it is advisable to work with the simpler stationary model, even when the data are non-stationary.

As an important side result of our rigorous theoretical analysis of the method, we have shown that the localised Yule-Walker estimator is strongly, uniformly consistent under local stationarity.

\bibliography{references}

\clearpage
\appendix

\section*{Appendix}

In Section~\ref{app:main} we provide proofs of the results in the main text. In Section~\ref{app:C:thm:main:3}, we provide a proof for Theorem~\ref{thm:main3}, the performance guarantee of our model selection procedure. The proof relies on properties of the empirical mean squared prediction errors for fixed model order and segment (Lemmas~\ref{lem:main}--\ref{lem:rel_MSPE}) which we state in Section~\ref{app:C:threeLemmas}. Theorem~\ref{thm:propertiesYW} which is about concentration properties of the localised Yule-Walker estimate under local stationarity, is proved in Section~\ref{app:C:thm:propertiesYW}. Corollary~\ref{kor:propertiesYW} which is about the strong consistency of the localised Yule-Walker estimate is proved in Section~\ref{app:C:kor:propertiesYW}. Lemmas~\ref{thm:main3:kor1} and \ref{thm:main3:kor2}, which fascilitate our discussion of the special case of our procedure are proved in Section~\ref{app:proofSimplVer}.

In Sections~\ref{sec:Lemmas_a}--\ref{sec:Cov} we provide technical results about the properties of quantities related to the second order moments. In Section~\ref{sec:Lemmas_a} we state results about the vector $\bar a_{N,T}^{(p)}(u)$, defined in~\eqref{def:a_Delta}, around which the localised Yule-Walker estimator concentrates. We also discuss how it is related to the mean square minimising 1-step ahead forecasting coefficients. In Section~\ref{app:der_thMSPE} we discuss properties of $\bar v_{N,T}^{(p,h)}(u)$, the $h$-step ahead version of $\bar a_{N,T}^{(p)}(u)$. Further, we establish properties of $g_{\Delta}^{(p,h)}(u)$ and ${\rm MSPE}_{\Delta_1, \Delta_2}^{(p,h)}(u)$ from the definition of $q(\delta)$ that is important for Assumptions~\ref{a:m_and_N} and~\ref{a:T}.
In Section~\ref{sec:Cov:mom}, we provide approximation results for expectations of Toepliz matrices of empirical localised autocovariances $\hat \gamma_{k;N,T}(t)$, defined in~\eqref{def:acf} and in Section~\ref{sec:Cov:concentration} we establish concentration results. In Section~\ref{app:TRes} we state a number of technical lemmas that we use in the proofs of our results. We state these results in a separate section, because we believe that they are useful for proving similar results in the future.

Sections~\ref{app:proofAuxLem}--\ref{a:extraSim} that are only available in the extended, arXiv'ed version of the manuscript, cf. \cite{KleyEtAl2019}, contain supplementary material. In Section~\ref{app:proofAuxLem} we provide the proofs of Lemmas~\ref{lem:main}--\ref{lem:rel_MSPE}. In Section~\ref{sec:res_sausta91} we cite two results from~\cite{sausta91} which we use for our proof in Section~\ref{sec:Cov:proof}.
In Sections~\ref{a:extraEmpExpl}--\ref{a:extraSim} we provide additional material for our simulation and empirical study.

\section{Proofs of Theorems~\ref{thm:main3} and~\ref{thm:propertiesYW}, of Corollary~\ref{kor:propertiesYW}, and of Lemmas~\ref{thm:main3:kor1} and \ref{thm:main3:kor2}}\label{app:main}

\subsection{Outlook}

\begin{figure}[t]
	\begin{center}
		\begin{tikzpicture}[
this_sec/.style = {fill=blue!30},
other_sec/.style = {fill=gray!30},
]


\node [other_sec] at (-2,0) {Corollary~\ref{kor:exp_gamma_k}};
\node (F2iii) at (-2,-0.5) {(iii)};
\node [other_sec] at (1,0) {Lemma~\ref{lem:BoundEV_M0}};
\node (F4ib) at (0.5,-0.5) {(i-b)};
\node (F4iic) at (1.5,-0.5) {(ii-c)};
\node [other_sec] (F5) at (4,0) {Lemma~\ref{lem:exp_ineq_gamma_k2}};
\node [other_sec] (G6) at (7,0) {Lemma~\ref{lem:ineqQuotientMV}};

\node [this_sec] (61) at (2.5,-2) {Theorem~\ref{thm:propertiesYW}};
\node [this_sec] (62) at (5.5,-2) {Corollary~\ref{kor:propertiesYW}};

\node [this_sec] (31) at (-0.5,-2) {Theorem~\ref{thm:main3}};
\node [other_sec] (C1) at (2.5,-4) {Lemma~\ref{lem:main}};
\node [other_sec] (C2) at (-0.5,-4) {Lemma~\ref{lem:nicerboundP}};
\node [other_sec] (C3) at (-3.5,-4) {Lemma~\ref{lem:rel_MSPE}};

\path[->] (62) edge (61);
\path[->] (61) edge (F2iii);
\path[->] (61) edge (F4ib);
\path[->] (61) edge (F4iic);
\path[->] (61) edge (F5);
\path[->] (61) edge (G6);
\path[->] (C1) edge (61);
\path[->] (31) edge (C1);
\path[->] (31) edge (C2);
\path[->] (31) edge (C3);

\end{tikzpicture}
	\end{center}
	\caption{Map of the results proved in Section~\ref{app:main}.}
	\label{fig:map_app_C}
\end{figure}
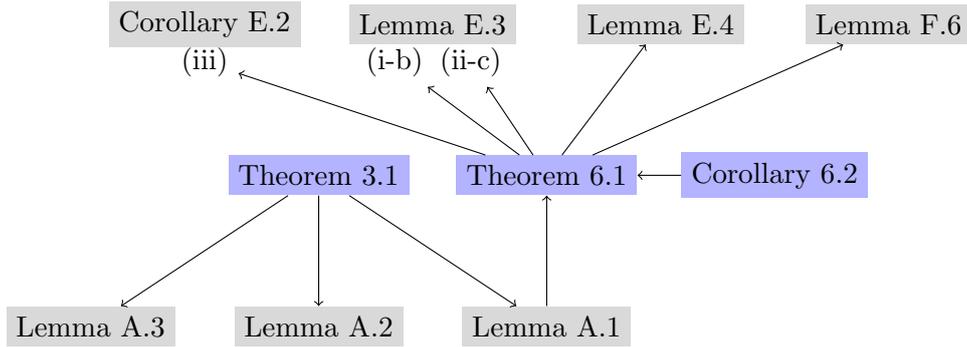

In this section we provide the proofs of the results from Sections~\ref{theorysectionar1} and~\ref{sec:propertiesYW}. In Section~\ref{app:C:threeLemmas} we state and discuss three auxiliary results (Lemmas~\ref{lem:main}--\ref{lem:rel_MSPE}) which facilitate the proof of our main result (Theorem~\ref{thm:main3}). The auxiliary results are about the empirical mean squared prediction error. Their proofs are deferred to Section~\ref{app:proofAuxLem} \citep{KleyEtAl2019}. The proof of Theorem~\ref{thm:main3}, by which our model selection procedure chooses models consistently with high probability, is then stated in Section~\ref{app:C:thm:main:3}. Because the proof of Lemma~\ref{lem:main} heavily relies on our result about the Yule-Walkers estimators (Theorem~\ref{thm:propertiesYW}), our proof of Theorem~\ref{thm:main3}, implicitly, also depends on it. The proof of Theorem~\ref{thm:propertiesYW} and its corollary (Corollary~\ref{kor:propertiesYW}), by which the localised Yule-Walker estimator is uniformly, strongly consistent, are stated in Sections~\ref{app:C:thm:propertiesYW} and~\ref{app:C:kor:propertiesYW}, respectively. For the proof of Theorem~\ref{thm:propertiesYW} we employ some of our results about the localised empirical autocovariance estimate from Section~\ref{sec:Cov} and a technical result from Section~\ref{app:TRes}. For the readers convenience, we include Figure~\ref{fig:map_app_C} in which the dependence of the various results is illustrated graphically.

\subsection{Three technical lemmas for the proof}\label{app:C:threeLemmas}
We now introduce two quantities that combine constants from the assumptions. Stating the results in terms of these constants will help to better interpret the bounds and significantly shorten otherwise complicated expressions. To this end, we define
\begin{equation}\label{eqn:def_C0}
		C_0 := (2\pi)^{1/2} M_f / m_f, \quad \text{and} \quad
		C_1 := (2\pi M'_f + C) m_f^{-1}.
\end{equation}

The constant $C_0$ can be interpreted in terms of the strength of serial correlation. Note that $C_0$ will be smaller if there is little variation (uniform in local time) of the spectral density with respect to frequency. In particular, it will be minimal if the spectral density is constant. This would corresponds to the case of white noise. The constant $C_1$ can be interpreted as divergence from stationarity. In particular, note the meaning of the two summands of the first factor. The constant $M'_f$ corresponds to the rapidity of changes in stationarity and will vanish in case of stationarity. The constant $C$ corresponds to the quality of locally approximating the correlation structure with a stationary processes correlation structure. It, also, vanishes if the underlying process is stationary.

The aim of the auxiliary results is to approximate
general mean squared prediction errors of the form
\begin{equation}
\label{eqn:MSPE}
	{\rm MSPE}_{s,m,N,T}^{(p,h)} := \frac{1}{m} \sum_{t=s+1}^{s+m} \big( X_{t+h,T} - \sum_{i=1}^p \hat v^{(p,h)}_{i;N,T}(t) X_{t-i+1,T} \big)^2,
\end{equation}
with $\hat v^{(p,h)}_{i;N,T}(t)$ defined in~\eqref{def:vhat} and $\hat v^{(p,h)}_{i;0,T}(t) := \hat v^{(p,h)}_{i;t,T}(t)$.

The first auxiliary result (Lemma~\ref{lem:main}) entails that the quantity defined in~\eqref{eqn:MSPE} is, with high probability, close to
\begin{equation}
\label{eqn:MSPE_bar}
	\overline{{\rm MSPE}}_{s,m,N,T}^{(p,h)} := \frac{1}{m} \sum_{t=s+1}^{s+m} \E \big( X_{t+h,T} - \sum_{i=1}^p \bar v^{(p,h)}_{i,N,T}(t) X_{t-i+1,T} \big)^2,
\end{equation}
with

\begin{minipage}[t]{.45\linewidth}
\begin{equation*}
\begin{split}
	& \big( \bar v_{N,T}^{(p,h)}(t) \big)' \\
				&  := \big(
				\bar v_{1;N,T}^{(p,h)}(t), \bar v_{2;N,T}^{(p,h)}(t), \cdots, \bar v_{p;N,T}^{(p,h)}(t) \big)	\\
				& := e'_1 \big( \bar A_{N,T}^{(p)}(t) \big)^h \\
				& := e'_1 \big( e_1 \big( \bar a_{N,T}^{(p)}(t) \big)' + H \big)^h,
\end{split}
\end{equation*}
\end{minipage}
\begin{minipage}[t]{.54\linewidth}
\begin{equation}\label{eqn:def_e_H}
e_1 := \begin{pmatrix} 1 \\ 0 \\ \vdots \\ 0 \end{pmatrix},
\ H := 
				\begin{pmatrix}
							0			 & 0      & \cdots & 0  & 0 \\
							1			 & 0      & \cdots & 0  & 0 \\
							0      & 1      & \cdots & 0  & 0 \\
							\vdots & \ddots & \cdots & 0  & 0 \\
							0      & 0      & \cdots & 1  & 0
						\end{pmatrix}.
\end{equation}
\end{minipage}
where $\bar a_{N,T}^{(p)}(t)$ is defined in~\eqref{eqn:YW_bar}, $e_1$ denotes the first canonical unity vector of dimension $p$ and $H$ denotes a $p \times p$ Jordan block with all eigenvalues equal to zero. The second auxiliary result (Lemma~\ref{lem:nicerboundP}) provides a simplified probability bound for the result in Lemma~\ref{lem:main} that can be applied in an especially relevant case.

By our third auxiliary result (Lemma~\ref{lem:rel_MSPE}) we have that $\overline{{\rm MSPE}}_{s,m,N,T}^{(p,h)}$ in turn can be approximated by ${\rm MSPE}_{N/T, m/T}^{(p,h)}(s/T)$, where ${\rm MSPE}_{\Delta_1, \Delta_2}^{(p,h)}(u)$ is the quantity defined in~\eqref{eqn:thMSPE}, with continuous time indices $\Delta_1$ and $\Delta_2$. Note that this quantity also appears in $q(\delta)$ defined in~\eqref{def:f} which is a relevant component of Assumptions~\ref{a:m_and_N} and~\ref{a:T}.

Some comparison of $\overline{{\rm MSPE}}_{s,m,N,T}^{(p,h)}$, defined in~\eqref{eqn:MSPE_bar}, and ${\rm MSPE}_{N/T, m/T}^{(p,h)}(s/T)$, as defined in~\eqref{eqn:thMSPE} are in order:
Note that $\overline{{\rm MSPE}}_{s,m,N,T}^{(p,h)}$ is defined as the expectation of a modified version of ${\rm MSPE}_{s,m,N,T}^{(p,h)}$, the modification being that $\hat v_{N,T}^{(p,h)}(t)$ is exchanged by $\bar v_{N,T}^{(p,h)}(t)$. As before, we will denote $\bar v_{0,T}^{(p,h)}(t) := \bar v_{t,T}^{(p,h)}(t)$.


We have that $g^{(p,h)}_{N/T}(t/T)$ approximates $\E[ ( X_{t+h,T}- f^{{\lstat}}_{t,h;p,N})^2]$, with $f^{{\lstat}}_{t,h;p,N}$ defined in~\eqref{eqn:lsEstim}. Therefore, the expectation of the empirical mean squared prediction error~\eqref{eqn:MSPE} we are considering is naturally an average of these quantities:
\[ \E [{\rm MSPE}_{s,m,N,T}^{(p,h)}] = \frac{1}{m}\sum_{t=s+1}^{s+m} \E[ ( X_{t+h,T}- f^{{\lstat}}_{t,h;p,N})^2]. \]

We now state the results that the quantities defined in~\eqref{eqn:MSPE} and~\eqref{eqn:MSPE_bar} are close, with high probability.
\begin{lem}\label{lem:main}
Let $(X_{t,T})_{t \in \IZ, T \in \IN^*}$ satisfy Assumptions~\ref{a:loc_stat}--\ref{a:MomentCond} and $\E X_{t,T} = 0$.
	Then, for every $m, h \in \IN^*$, $p \in \IN$, $N \geq 6 C_0 p^{2}$, $\varepsilon > 0$ and $T \geq 10 C_1 p^2$, with ${\rm MSPE}_{s,m,N,T}^{(p,h)}$ defined in~\eqref{eqn:MSPE} and $\overline{{\rm MSPE}}_{s,m,N,T}^{(p,h)}$ defined in~\eqref{eqn:MSPE_bar}, we have
	that
\begin{equation*}
		\IP \Big( \Big| {\rm MSPE}_{s,m,N,T}^{(p,h)} - \overline{{\rm MSPE}}_{s,m,N,T}^{(p,h)} \Big| > \varepsilon \Big) \leq P^{(p,h)}_{m,N}(\varepsilon)
\end{equation*}
and
\begin{equation*}
		\IP \Big( \Big| {\rm MSPE}_{s,m,0,T}^{(p,h)} - \overline{{\rm MSPE}}_{s,m,0,T}^{(p,h)} \Big| > \varepsilon \Big) \leq P^{(p,h)}_{m,s}(\varepsilon)
\end{equation*}
with
\begin{equation*}
	\begin{split}
& P^{(p,h)}_{m,N}(\varepsilon) :=
(1 + 4p + 2p^2) \\
		& \cdot \exp\Bigg( - \frac{\frac{\varepsilon^2}{(p+1)^4} }{8\Big(\big(2 C_0 + 1\big)^{4 h} \frac{C_{1,2} (h+p-1)}{m} + (\frac{\varepsilon}{2 (p+1)^2})^{(3+8d)/(2+4d)} \Big( \big(2 C_0 + 1\big)^{2 h} \frac{C_{2,2} (h+p-1)}{m} \Big)^{1/(2+4d)}\Big)} \Bigg) \\
		& \qquad + 6 m p^2 (p+1) 
		\exp\Bigg( - \frac{\eta^2}{2\Big(C_{1,1} \frac{p}{N-p} + \eta^{(3+4d)/(2+2d)} \Big(C_{2,1} \frac{p}{N-p} \Big)^{1/(2+2d)}\Big)} \Bigg),
	\end{split}
\end{equation*}
where
\begin{align*}
	\eta & := \frac{m_f}{4 p} \min\Big\{1, \bar\mu / (8 C_0) \Big\}, &
	\bar\mu & := 2^{1-h} \frac{\mu}{\mu + h (2 C_0)^{h-1}}, \\
	\mu & := \frac{\bar \varepsilon}{2 \Big( \big( 2 C_0 + 1 \Big)^{2h} + \bar \varepsilon \Big)^{1/2}}, &
	\bar \varepsilon & := \frac{\varepsilon/(p+1)^2}{2 \big( (6 \pi M_f c^2 24^d)^2 + \varepsilon^2/(p+1)^4 \big)^{1/4}},
	\end{align*}
and the constants $C_0$, $C_1$, and $C_{1,1}$, $C_{1,2}$, $C_{2,1}$, $C_{2,2}$, and $m_f$, $M_f$, and $c$, $d$ are defined in~\eqref{eqn:def_C0}, \eqref{lem:exp_ineq_gamma_k2:def:C1C2}, and Assumptions~\ref{a:sd_bounded} and~\ref{a:MomentCond}, respectively.
\end{lem}

In a typical application the bound $P^{(p,h)}_{m,N}(\varepsilon)$ will be small. More precisely, the following, more accessible bound for $P^{(p,h)}_{m,N}(\varepsilon)$, proved in Section~\ref{app:proofAuxLem} \citep{KleyEtAl2019}, will be useful
\begin{lem}\label{lem:nicerboundP}
There exist constants $D_1, D_2, D_3 > 0$ and $K_0 > 1$, defined in the proof, such that for any
\begin{equation}\label{cond:eps}
	\max\Big\{ \Big(\frac{h+p}{m}\Big)^{\frac{1+4d}{3+8d}} K_0^{h} p^2, \Big(\frac{p}{N-p}\Big)^{\frac{1+2d}{3+4d}} K_0^{h} p^3 h \Big\}
	 < \varepsilon \leq \min\{6\pi M_f c^2 24^d, 1\} (p+1)^2,
\end{equation}
we have
\begin{equation*}
	P^{(p,h)}_{m,N}(\varepsilon)
	\leq D_1 \Bigg[ p^2 \exp\Bigg(-D_2 \Big( \frac{m}{h+p} \Big)^{1/(3+8d)} \Bigg)
	+ m p^3 \exp\Bigg(-D_3 \Big( \frac{N-p}{p} \Big)^{1/(3+4d)} \Bigg) \Bigg].
\end{equation*}
\end{lem}

Note that we are interested in the scenario where $\varepsilon > 0$ may be small. Therefore, if we allow that $p$ and $h$ may be large, we have to require $m$ and $N$ to be of a minimum size.

We now state the result that the quantities defined in~\eqref{eqn:MSPE_bar} and~\eqref{eqn:thMSPE} are close. The quality of the approximation depends on the parameters $T$, $p$ and $h$, but is uniform with respect to $s$, $m$ and $N$:

\begin{lem} \label{lem:rel_MSPE}
	Let $(X_{t,T})_{t \in \IZ, T \in \IN^*}$ satisfy Assumptions~\ref{a:loc_stat}--\ref{a:MomentCond} and $\E X_{t,T} = 0$.
	Then, for every $m, h \in \IN^*$, $p \in \IN$, $T \geq 6 h 2^{h} C_1 p^2$, and
	$N \geq 4 h 2^h C_0 p^{2}$, with $\overline{{\rm MSPE}}_{s,m,N,T}^{(p,h)}$ defined in~\eqref{eqn:MSPE_bar} and ${\rm MSPE}_{\Delta_1, \Delta_2}^{(p,h)}(u)$ defined in~\eqref{eqn:thMSPE}, we have
	\[ \Big| \overline{{\rm MSPE}}_{s,m,N,T}^{(p,h)} - {\rm MSPE}_{N/T, m/T}^{(p,h)}(s/T)\Big|
		\leq 8 h 2^h \big( C_0 \big)^{2 h+1} \Big[ 6 (2\pi M'_f + C) \frac{p^2}{T} + \frac{p^{2}}{N} \Big]
	\]
	and
	\[ \Big| \overline{{\rm MSPE}}_{s,m,0,T}^{(p,h)} - {\rm MSPE}_{s/T, m/T}^{(p,h)}(s/T)\Big|
		\leq 8 h 2^h \big( C_0 \big)^{2 h+1} \Big[ 6 (2\pi M'_f + C) \frac{p^2}{T} + \frac{p^{2}}{N} \Big].
	\]
\end{lem}

The proofs of the three lemmas are long and technical. We therefore defer them to Section~\ref{app:proofAuxLem} \citep{KleyEtAl2019}.

A few comments about Lemma~\ref{lem:rel_MSPE} are in order. Note that the approximation error is zero in case of a stationary time series, as then $2\pi M'_f + C = 0$. Note further, that the approximation will be better, if $h$ and $p$ are small compared to $T$. More precisely, if $h (2 C_0^2)^h p^2 = o(T)$, then the difference will vanish asymptotically. In particular, if $h = O(1)$, then it would suffice to assume that $p = o(T^{1/2})$, for the approximation error to vanish asymptotically.

\subsection{Proof of Theorem~\ref{thm:main3}}\label{app:C:thm:main:3}

The constants $D_1$, $D_2$ and $D_3$ are defined as
\begin{equation}\label{def:D}
\begin{split}
	D_1 & := 12,\\
	D_2 & := \big( 2^8 \max\{C_{1,2}, C_{2,2}^{1/(2+4d)}\}\big)^{-1}, \text{ and}\\
	D_3 & := K_1^2 / \big( 2^{12} \max\{C_{1,1}, (K_1^{3+4d} C_{2,1})^{1/(2+2d)}\} \big),
\end{split}
\end{equation}
where $K_1 := m_f / (32 \min\{(6 \pi M_f c^2 24^d)^{1/2}, 1\})$ and
\begin{equation}\label{lem:exp_ineq_gamma_k2:def:C1C2}
	\begin{aligned}
		C_{1,\alpha} & := 12 \cdot 2^{10 \alpha d + 7} \alpha^{4 \alpha d} \big(\max\{c^2, 3\pi M_f, 1\} \big)^{2\alpha} {\rm e} \Big( 1 + \frac{1}{\log \rho}\Big)\Big(1 + K^{1/2} \Big), \\
		C_{2,\alpha} & := 12 \cdot 2^{4 \alpha d + 3} \alpha^{2 \alpha d} \big(\max\{c^2, 3\pi M_f, 1\} \big)^{\alpha} {\rm e} \Big( 1 + \frac{1}{\log \rho}\Big),
		\end{aligned}
\end{equation}
with $\alpha \in \{1,2\}$. In the definitions, we have $K$ and $\rho$ the constants from Assumption~\ref{a:mixing}, $M_f$ and $m_f$ the constants from Assumption~\ref{a:sd_bounded}, and $c$ and $d$ the constants from Assumption~\ref{a:MomentCond}.

To compact notation, we denote $s_2 := T-h$, ${\rm MSPE}_{s_{i},m,N,T}^{(p_1,h)}$ by $X_i$ and ${\rm MSPE}_{s_{i},m,0,T}^{(p_2,h)}$ by $Y_i$. Further, denote
${\rm MSPE}_{N/T,m/T}^{(p_1,h)}(\frac{s_{i}}{T})$ and ${\rm MSPE}_{s_1/T,m/T}^{(p_2,h)}(\frac{s_{i}}{T})$
by $\bar Y_i$ and $\bar X_i$, respectively. 
Further, we abbreviate $A := Y_1 - X_1 (1 + \delta)$ and $B := Y_1-Y_2 + (X_2 - X_1) (1 + \delta)$.

First note that Assumptions~\ref{a:m_and_N} and~\ref{a:T} imply that
\begin{equation*}
		T \geq \max\big\{ 10 C_1 (\max\mathcal{P})^2, 
		                    6 h 2^{h} C_1 (\max\mathcal{P})^2 
												\big\}, \quad
		\min\mathcal{N} \geq 4 h 2^h C_0 (\max\mathcal{P})^2 
\end{equation*}
Therefore, the conditions of Lemmas~\ref{lem:main} and~\ref{lem:rel_MSPE} are satisfied.
Further, note that since
\[\min\mathcal{N} \geq 8 h 2^h \big( C_0 \big)^{2 h+1} \big( \max\mathcal{P} \big)^2 \Big[ 6 (2\pi M'_f + C) + 1 \Big] \Big( 20 (1+\delta) / q(\delta) \Big)\]
and because $N \leq T$ for all $N \in \mathcal{N}$, we have that the bound from Lemma~\ref{lem:rel_MSPE} can again be bounded
\begin{equation}\label{proof:thm:main3:1}
	8 h 2^h \big( C_0 \big)^{2 h+1} \big( \max\mathcal{P} \big)^2 \Big[ 6 (2\pi M'_f + C) \frac{1}{T} + \frac{1}{N} \Big] \leq \frac{q(\delta)}{20 (1+\delta)} =: \varepsilon
\end{equation}

Finally, note that by Assumption~\ref{a:T}, we have
\[T \geq 4 m \big( 2 h + 1 \big) \big( C_0 \big)^{2 h+1} M'_f \frac{20 (1+\delta)}{q(\delta)} \]
which implies that (a quantity related to the bound from Lemma~\ref{lem:norm_g_der}(iv)) can be bounded
\begin{equation}\label{proof:thm:main3:2}
		 4 \big( 2 h + 1 \big) \big( C_0 \big)^{2 h+1} M'_f
		\Big| \frac{s_1 - s_2}{T} \Big| \leq \varepsilon.
\end{equation}

Now, for the proof of the Theorem, note that
\begin{align}
	& \IP\Big( (\hat R_{T,2}(h) \geq 1 + \delta \text{ and }
		\hat R_{T,3}(h) \geq 1 + \delta)
		\text{ or }
		(\hat R_{T,2}(h) < 1 + \delta \text{ and }
		\hat R_{T,3}(h) < 1 + \delta) \Big) \nonumber \\
	& \qquad \geq 1 - \sum_{p_1, p_2 \in \mathcal{P}} \sum_{N \in \mathcal{N}} \Big( \IP\big( |A| \leq q(\delta)/2 \big) + \IP\big( |B| > q(\delta)/2 \big) \Big), \label{eqn:main3:forApp}
\end{align}
which we prove in Section~\ref{sec:eq_forApp} \citep{KleyEtAl2019}.

We now bound the part of the right hand side of~\eqref{eqn:main3:forApp} that involves the quantity $A$. Using the fact that
\begin{equation*}
\begin{split}
	| \bar Y_1 - \bar X_1 (1 + \delta)|
	& = | Y_1 + \bar Y_1 - Y_1 - X_1 (1 + \delta) + (X_1 - \bar  X_1) (1 + \delta)| \\
	& \leq | Y_1 - X_1 (1 + \delta) | + | \bar Y_1 - Y_1| + |X_1 - \bar  X_1| (1 + \delta),
\end{split}
\end{equation*}
we have the first inequality of
\begin{align}
	& \IP\Big( |A| \leq q(\delta)/2 \Big) 
	= \IP\Big( |Y_1 - X_1 (1 + \delta)| \leq q(\delta)/2 \Big) \nonumber \\
	& \leq \IP\Big( |Y_1 - \bar Y_1 | + |X_1 - \bar X_1 | (1+\delta) \geq | \bar Y_1 - \bar X_1 (1 + \delta)| - q(\delta)/2 \Big) \nonumber \\
	& \leq \IP\Big( |Y_1 - \bar Y_1 | \geq \frac{1}{2} (| \bar Y_1 - \bar X_1 (1 + \delta)| - q(\delta)/2) \Big) \nonumber \\
	& \qquad + \IP\Big( |X_1 - \bar X_1 | \geq \frac{1}{2 (1 + \delta)} (| \bar Y_1 - \bar X_1 (1 + \delta)| - q(\delta)/2) \Big) \nonumber \\
	& \leq \IP\Big( |Y_1 - \bar Y_1 | > q(\delta)/10) \Big)
		+ \IP\Big( |X_1 - \bar X_1 | > \frac{q(\delta)}{10 (1 + \delta)} ) \Big) \label{prf:main3:A} \\
	& \leq P^{(p_2,h)}_{m,T-m}\Big( \frac{q(\delta)}{20} \Big)
		+ P^{(p_1,h)}_{m,N}\Big( \frac{q(\delta)}{20 (1+\delta)} \Big)
	\leq 2 P^{(p_{\max},h)}_{m,N_{\min}}\Big( \frac{q(\delta)}{20 (1+\delta)} \Big), \label{prf:main3:B}
\end{align}
where $p_{\max} := \max\mathcal{P}$ and $N_{\min} := \min\mathcal{N}$.
For the inequality in~\eqref{prf:main3:A} we have used the definition of $q(\delta)$ and $1/4 > 1/10$. For the first inequality in~\eqref{prf:main3:B} we have used Lemmas~\ref{lem:main} and~\ref{lem:rel_MSPE} and~\eqref{proof:thm:main3:1} to obtain
\begin{equation}\label{eqn:main_N}
		\IP \Big( \Big| {\rm MSPE}_{s,m,N,T}^{(p,h)} - {\rm MSPE}_{N/T,m/T}^{(p,h)}(s/T) \Big| > 2 \varepsilon \Big) \leq P^{(p,h)}_{m,N}(\varepsilon)
\end{equation}
and
\begin{equation}\label{eqn:main_0}
		\IP \Big( \Big| {\rm MSPE}_{s,m,0,T}^{(p,h)} - {\rm MSPE}_{s/T,m/T}^{(p,h)}(s/T) \Big| > 2 \varepsilon \Big) \leq P^{(p,h)}_{m,s}(\varepsilon).
\end{equation}
For the second inequality in~\eqref{prf:main3:B} we have used that
\[p_1 \leq p_2 \Rightarrow P^{(p_1,h)}_{m,N}(\varepsilon) \leq P^{(p_2,h)}_{m,N}(\varepsilon), \quad
N_1 \leq N_2 \Rightarrow P^{(p,h)}_{m,N_1}(\varepsilon) \geq P^{(p,h)}_{m,N_2}(\varepsilon),\]
and $\varepsilon_1 \leq \varepsilon_2 \Rightarrow P^{(p,h)}_{m,N}(\varepsilon_1) \geq P^{(p,h)}_{m,N}(\varepsilon_2)$.

We now bound the part of the right hand side of~\eqref{eqn:main3:forApp} that involves the quantity $B$.
We have
\begin{align}
	& \IP\Big( |B| > q(\delta)/2 \Big)
	= \IP\Big( |Y_1-Y_2 + (X_2 - X_1) (1 + \delta)| > q(\delta)/2 \Big) \nonumber \\ & \leq \IP\Big( | Y_1-Y_2 | > q(\delta)/4\Big)
		+ \IP\Big( | X_2 - X_1 | > \frac{q(\delta)}{4(1 + \delta)} \Big) \nonumber \\
	& \leq 2 P^{(p_2,h)}_{m,T-m}\Big( \frac{q(\delta)}{20} \Big)
		+ 2 P^{(p_1,h)}_{m,N}\Big( \frac{q(\delta)}{20 (1+\delta)} \Big)
	\leq 4 P^{(p_{\max},h)}_{m,N_{\min}}\Big( \frac{q(\delta)}{20 (1+\delta)} \Big). \label{prf:main3:bnd_B}
\end{align}

Note that we have
\begin{equation*} 
	\begin{split}
		& \IP \Big( \Big| {\rm MSPE}_{s_1,m,N,T}^{(p,h)} - {\rm MSPE}_{s_2,m,N,T}^{(p,h)} > 5 \varepsilon \Big) \\
		& \leq \IP \Big( \Big| {\rm MSPE}_{s_1,m,N,T}^{(p,h)} - {\rm MSPE}_{N/T,m/T}^{(p,h)}(s_1/T) \Big| > 2 \varepsilon \Big) \\
		& \quad + \IP \Big( \Big| {\rm MSPE}_{s_2,m,N,T}^{(p,h)} - {\rm MSPE}_{N/T,m/T}^{(p,h)}(s_2/T) \Big| > 2 \varepsilon \Big) \\
		& \quad + I\Big\{ \Big| {\rm MSPE}_{N/T,m/T}^{(p,h)}(s_1/T) - {\rm MSPE}_{N/T,m/T}^{(p,h)}(s_2/T) \Big| > \varepsilon \Big) \Big\},
	\end{split}
\end{equation*}
where the first two terms can be bound by an application of~\eqref{eqn:main_N} and the indicator function vanishes for all $T$ satisfying the condition of the Theorem, because
\begin{equation*}
	\begin{split}
		& \Big| {\rm MSPE}_{N/T,m/T}^{(p,h)}(s_1/T) - {\rm MSPE}_{N/T,m/T}^{(p,h)}(s_2/T) \Big|
		\leq 4 \big( 2 h + 1 \big) \big( C_0 \big)^{2 h+1} M'_f
		\Big| \frac{s_1 - s_2}{T} \Big|, \\
		\end{split}
\end{equation*}
where Lemma~\ref{lem:norm_g_der}(iv) was employed to obtain~\eqref{proof:thm:main3:2} for the last inequality.

Thus, combining~\eqref{eqn:main3:forApp}, \eqref{prf:main3:B} and~\eqref{prf:main3:bnd_B}, we have shown that
\begin{equation*}
\begin{split}
	& \IP\Big( (\hat R_{T,2}(h) \geq 1 + \delta \text{ and }
		\hat R_{T,3}(h) \geq 1 + \delta)
		\text{ or }
		(\hat R_{T,2}(h) < 1 + \delta \text{ and }
		\hat R_{T,3}(h) < 1 + \delta) \Big) \nonumber \\
	& \qquad \geq 1 - 6 |\mathcal{P}|^2 |\mathcal{N}| P^{(p_{\max},h)}_{m,N_{\min}}\Big( \frac{q(\delta)}{20 (1+\delta)} \Big).
	\end{split}
\end{equation*}

An application of Lemma~\ref{lem:nicerboundP} finishes the proof of the theorem.

\begin{rem}\label{rem:asymptEmpMSPE}
Equations~\eqref{eqn:main_N}--\eqref{eqn:main_0}, which are immediate consequences of Lemmas~\ref{lem:main} and~\ref{lem:rel_MSPE}, can be used to derive the almost sure convergence of
\[\Big| {\rm MSPE}_{s,m,N,T}^{(p,h)} - {\rm MSPE}_{N/T,m/T}^{(p,h)}(s/T) \Big|
\text{ and }
\Big| {\rm MSPE}_{s,m,0,T}^{(p,h)} - {\rm MSPE}_{s/T,m/T}^{(p,h)}(s/T) \Big|,
\]under appropriate conditions, using a classical Borel-Cantelli argument.

This asymptotic view of ${\rm MSPE}_{s,m,N,T}^{(p,h)}$ and ${\rm MSPE}_{s,m,0,T}^{(p,h)}$, in particular, implies that we may interpret ${\rm MSPE}_{\Delta_1, \Delta_2}^{(p,h)}(u)$ as an approximation of the expectation of the empirical MSPE for an $h$-step ahead linear forecast of order $p$, where observations up to (local) time $u$ have been made. The $\Delta_1$ and $\Delta_2$ are (localised) length which are related to the segment length of observations used for the estimation of the forecasting coefficients and the segment from which the observations $X_{t+h,T}$ that are being forecasted are taken, respectively.
\end{rem}

We now proceed with the proofs of the results from Section~\ref{sec:propertiesYW}.

\subsection{Proof of Theorem~\ref{thm:propertiesYW}}\label{app:C:thm:propertiesYW}

Let	$M := \hat\Gamma_{N,T}^{(p)}(t)$, $M_0 := \E M$, $v := \hat\gamma_{N,T}^{(p)}(t)$, and $v_0 := \E v$.
By Lemma~\ref{lem:BoundEV_M0}(ii-c) we deduce that $M_0$ is invertible for $T \geq 2 p^2 C_1$, because it is positive definite with smallest eigenvalue larger or equal to $m_f/2$.
An application of Lemma~\ref{lem:ineqQuotientMV}, with the spectral norm as the matrix norm and the Euclidean norm as the vector norm yields
\begin{equation*}
\begin{split}
	& \IP \big( \| \hat a_{N,T}^{(p)}(t) - \bar a_{N,T}^{(p)}(t) \| > \varepsilon \big)
	= \IP \big( \| M^{-1} v - M_0^{-1} v_0 \| > \varepsilon \big) \\
	& \leq \IP \Big( \| M - M_0 \| > \frac{1}{2 \|M_0^{-1}\|} \Big) + \IP \Big( \|v-v_0\| > \frac{\varepsilon}{4} \frac{1}{\|M_0^{-1}\|} \Big) \\
				& \qquad + \IP \Big( \| M - M_0 \|  > \frac{\varepsilon}{4} \frac{1}{(\|M_0^{-1}\|)^2 \, \| v_0 \|} \Big) I\{ \|v_0\| \neq 0\} \\
	& \leq \IP \Big( \max_{k=0,\ldots,p-1} |\hat\gamma_{k;N,T}(t) - \E \hat\gamma_{k;N,T}(t)| > \frac{1}{4 p} m_f \Big) \\
	& \qquad + \IP \Big( \max_{k=1,\ldots,p} |\hat\gamma_{k;N,T}(t) - \E \hat\gamma_{k;N,T}(t)| > \frac{\varepsilon}{8 p^{1/2}} m_f \Big) \\
				& \qquad + \IP \Big( \max_{k=0,\ldots,p-1} |\hat\gamma_{k;N,T}(t) - \E \hat\gamma_{k;N,T}(t)|  > \frac{\varepsilon}{32 (2\pi)^{1/2} M_f p } m_f^2 \Big) \\
	& \leq 3 p \max_{k=0,\ldots,p} \IP \Big( |\hat\gamma_{k;N,T}(t) - \E \hat\gamma_{k;N,T}(t)| > \frac{m_f}{4 p} \min\Big\{1, \frac{\varepsilon p^{1/2}}{2}, \frac{\varepsilon}{8 C_0} \Big\} \Big), \\
	& = 3 p \max_{k=0,\ldots,p} \IP \Big( |\hat\gamma_{k;N,T}(t) - \E \hat\gamma_{k;N,T}(t)| > \frac{m_f}{4 p} \min\Big\{1, \frac{\varepsilon}{8 C_0} \Big\} \Big), \\
\end{split}
\end{equation*}
where we have use Lemma~\ref{lem:BoundEV_M0}(ii-c) again to bound $1/\|M_0^{-1}\|$. In the last step we employed that $\frac{p^{1/2}}{2} \geq \frac{1}{8 C_0}$. Further, we have used that $M - M_0$ satisfies
\begin{multline*} \| M - M_0 \|_1 = \| M - M_0 \|_{\infty}
= \max_{1 \leq \ell \leq p} \sum_{h=1}^p |\hat\gamma_{h-\ell;N,T}(t) - \E \hat\gamma_{h-\ell;N,T}(t) | \\ 
\leq p \max_{k=0,\ldots,p-1} |\hat\gamma_{k;N,T}(t) - \E \hat\gamma_{k;N,T}(t)|.
\end{multline*}
Thus, by H\"{o}lder's inequality
\[\|M - M_0\| \leq \Big( \|M - M_0\|_1 \|M - M_0\|_{\infty} \Big)^{1/2} \leq p \max_{k=0,\ldots,p-1} |\hat\gamma_{k;N,T}(t) - \E \hat\gamma_{k;N,T}(t)|. \]

For the Euclidean norm we have used
\[\|v-v_0\| \leq p^{1/2} \|v-v_0\|_{\infty} = p^{1/2} \max_{k=1,\ldots,p} |\hat\gamma_{k;N,T}(t) - \E \hat\gamma_{k;N,T}(t)|.\]

Finally, 
by Corollary~\ref{kor:exp_gamma_k}(iii) and Lemma~\ref{lem:BoundEV_M0}(i-b), we have
\begin{equation*}
	\begin{split}
	\|v_0\| & = \| \E \hat \gamma_{N,T}^{(p)}(t/T) \|
\leq \| f_{p,N} \circ \gamma_{N/T}^{(p)}(t/T) \| + \| \E \hat \gamma_{N,T}^{(p)}(t/T) - f_{p,N} \circ \gamma_{N/T}^{(p)}(t/T) \| \\
& \leq (2\pi)^{1/2} M_f + 2 T^{-1} p^{3/2} C_1 m_f \leq 2 (2\pi)^{1/2} M_f,
	\end{split}
\end{equation*}
where the second inequality holds for $T \geq 2 \frac{p^{3/2} C_1 m_f}{(2\pi)^{1/2} M_f} = 2 C_1 p^{3/2} / C_0$, which is the case, as $T \geq 2 C_1 p^2$ is assumed. Here we also have used that $\| f_{p,N} \circ x \| \leq \|x\|$, as all entries of $f_{p,N}$ are between 0 and 1.
Applying Lemma~\ref{lem:exp_ineq_gamma_k2} yields the assertion, because
\begin{equation*}
\begin{split}
	& \IP \big( \| \hat a_{N,T}^{(p)}(t) - a_{N,T}^{(p)}(t) \| > \varepsilon \big) \\
& \leq 3 p \max_{h=0,\ldots,p} \exp\Bigg( - \frac{\eta^2}{2\Big(C_{1,1}\frac{ h_*}{N-|h|} + \eta^{(3+4d)/(2+2d)} \Big( C_{2,1}\frac{ h_*}{N-|h|} \Big)^{1/(2+2d)}\Big)} \Bigg) \\
& = 3 p \exp\Bigg( - \frac{\eta^2}{2\Big(C_{1,1} \frac{p}{N-p} + \eta^{(3+4d)/(2+2d)} \Big(C_{2,1} \frac{p}{N-p} \Big)^{1/(2+2d)}\Big)} \Bigg) \\
\end{split}
\end{equation*}
where $\eta := \frac{m_f}{4 p} \min\Big\{1, \varepsilon \frac{1}{8 C_0}\Big\}$, and the third line follows from the fact, for any two integers $N$ and $p$ with $N \geq 1+p \geq 2$ we have that $(\frac{h_*}{N-|h|})_{h=0,1,\ldots,p}$ is an increasing sequence. This is easy to see: $\frac{1}{N-0} \leq \frac{1}{N-1} \leq \ldots \leq \frac{p-1}{N-p+1} \leq \frac{p}{N-p}$.
Note that $T \geq 2 p C_1 \geq C/(\pi m_f)$, such that this condition of Lemma~\ref{lem:exp_ineq_gamma_k2} is met.\hfill$\square$

\subsection{Proof of Corollary~\ref{kor:propertiesYW}}\label{app:C:kor:propertiesYW}

Note the fact that, if $R_n \geq 0$ is a sequence with $R_n \rightarrow \infty$, as $n \rightarrow \infty$, then 
\[b_n = O(1) \Leftrightarrow b_n = o(r_n), \forall \, 0 \leq r_n \leq R_n, \text{ with } r_n \rightarrow \infty, \text{ as } n \rightarrow \infty.\]
Thus, employing the Borel-Cantelli lemma, it suffices to show that, for any given $\varepsilon > 0$
and sequence $0 \leq r_T \leq R_T^{1/2}$ with $r_T \rightarrow \infty$,
we have
\[\sum_{T=1}^{\infty} \IP \Big( \sup_{p=1,\ldots,P} \sup_{t = N, \ldots, T} \| \hat a_{N,T}^{(p)}(t) - \bar a_{N,T}^{(p)}(t) \| > \varepsilon P^{3/2} \Big( \frac{ \log(T)}{N} \Big)^{1/2} r_T \Big) < \infty.\]

This follows, since we have
	\begin{equation*}
	\begin{split}
		& \IP \Big( \sup_{p=1,\ldots,P} \sup_{t = N, \ldots, T} \| \hat a_{N,T}^{(p)}(t) - a_{N,T}^{(p)}(t) \| > \varepsilon P^{3/2} \Big( \frac{ \log(T)}{N} \Big)^{1/2} r_T \Big) \\
		& \leq P \cdot T \cdot \sup_{p=1,\ldots,P} \sup_{t = N, \ldots, T} \IP \Big(   \| \hat a_{N,T}^{(p)}(t) - a_{N,T}^{(p)}(t) \| > \varepsilon P^{3/2} \Big( \frac{ \log(T)}{N} \Big)^{1/2} r_T \Big) \\
		& \leq P \cdot T \cdot \sup_{p=1,\ldots,P} \sup_{t = N, \ldots, T} \IP \Big(   \| \hat a_{N,T}^{(p)}(t) - a_{N,T}^{(p)}(t) \| > \varepsilon p^{3/2} \tilde C^{1/2} \Big( \frac{ \log(T)}{N-p} \Big)^{1/2} r_T \Big) \\
		& \leq P \cdot T \cdot \sup_{p=1,\ldots,P} \sup_{t = N, \ldots, T} 3 p \exp\Bigg( - \varepsilon^{2} \frac{p^{3} \log(T)}{N-p} \tilde C r_T^2 \frac{m_f^2}{2^{12} C_{1,1}} \Big( C_0^2 \frac{p^3}{N-p} \Big)^{-1} \Bigg) \Bigg) \\
		& = 3 T^3 \exp\Bigg( - \varepsilon^{2} \log(T) \tilde C r_T^2 \frac{m_f^2}{2^{12} C_{1,1}} C_0^{-2} \Bigg) \Bigg) \leq 3 T^{-2},
	\end{split}
	\end{equation*}
for $T$ large enough. In the second inequality we have used the fact that, due to $P = o(N)$, there exists a $\tilde C > 0$ such that $1/N \geq \tilde C/(N-P)$, for $T$ large enough.
Note that we have $P = o(T^{1/2})$, from $N \leq T$, $P = o(N^{(1+2d)/(4+6d)})$ and $d \geq 1/2$, such that, in the third inequality, Theorem~\ref{thm:propertiesYW} can be applied, where we have also used the fact that, under the assumptions made
\begin{equation*}
	p^{3/2} \Big( \frac{ \log(T)}{N-p} \Big)^{1/2} R_T^{1/2}
	= o\Big( \frac{P^{(4+6d)(3+4d)}}{N^{(1+2d)/(3+4d)}} \Big),
\end{equation*}
implying that, for $T$ large enough, we have
\[\varepsilon p^{3/2} \Big( \frac{ \log(T)}{N-p} \Big)^{1/2} r_T \leq \min\{U_{p,N}, 1/(8 C_0)\} = U_{p,N}.\]
This completes the proof.\hfill$\square$

\subsection{Proofs of Lemmas~\ref{thm:main3:kor1} and \ref{thm:main3:kor2}}\label{app:proofSimplVer}

For the proof of Lemma~\ref{thm:main3:kor1} it suffices to show that
\[q(\delta) := \min_{N \in \mathcal{N}} \Big| {\rm MSPE}_{s_1/T,m/T}^{(1,1)}(\frac{s_1}{T})
- (1+\delta) \cdot {\rm MSPE}_{N/T,m/T}^{(1,1)}(\frac{s_1}{T}) \Big| \geq \delta \pi m_f \big( 1 - \rho^2 \big).\]
Likewise, to show Lemma~\ref{thm:main3:kor2}, we bound $q(\delta)$ with $\pi m_f D_{\inf}^2 / 2$ on the right hand side.

Denoting
\[\gamma_{k}(u,\Delta) := \int_0^1 \gamma_k(u+\Delta (x-1)) {\rm d}x
= \Delta^{-1} \int_{u - \Delta}^u \gamma_k(v) {\rm d}v \]
we have, by definition~\eqref{eqn:thMSPE}, that
\begin{equation*}
\begin{split}
	& {\rm MSPE}_{\Delta_1,\Delta_2}^{(1,1)}(u)
	= \int_0^1 g^{(1,1)}_{\Delta_1}\big( u + \Delta_2 (1-x) \big) {\rm d}x \\
	& = \Delta_2^{-1} \int_u^{u+\Delta_2} \Big( \gamma_0( w )
			-  2 \frac{\gamma_1\big( w; \Delta_1 \big)}{\gamma_0\big( w; \Delta_1 \big)} \gamma_1( w )
			+ \Big(\frac{\gamma_1\big( w; \Delta_1 \big)}{\gamma_0\big( w; \Delta_1 \big)}\Big)^2 \gamma_0( w )				 \Big) {\rm d}w.
\end{split}
\end{equation*}

To find the lower bound we want, it therefore suffices to proof lower bounds, for every $w \in [s_1/T, (s_1+m)/T]$, of the following difference
\begin{equation}\label{eqn:lem:xy:D}
\begin{split}
	& \Big(\Big( \gamma_0( w ) -  2 \frac{\gamma_1\big( w; s_1/T \big)}{\gamma_0\big( w; s_1/T \big)} \gamma_1( w )
			+ \Big(\frac{\gamma_1\big( w; s_1/T \big)}{\gamma_0\big( w; s_1/T \big)}\Big)^2 \gamma_0( w ) \Big) \\
			& - (1+\delta) \Big( \gamma_0( w ) -  2 \frac{\gamma_1\big( w; N/T \big)}{\gamma_0\big( w; N/T \big)} \gamma_1( w )
			+ \Big(\frac{\gamma_1\big( w; N/T \big)}{\gamma_0\big( w; N/T \big)}\Big)^2 \gamma_0( w ) \Big) \Big).
\end{split}
\end{equation}
For Lemma~\ref{thm:main3:kor2} we will bound $-1 \times \eqref{eqn:lem:xy:D}$.
For notational convenience we omit the $w$'s and denote
\[E := \frac{\gamma_1(w, N/T)}{\gamma_0(w, N/T)}, \text{ and }
  F := \frac{\gamma_1(w, s_1/T)}{\gamma_0(w, s_1/T)}.\]

By elementary considerations it can be shown that
\begin{equation}
	\eqref{eqn:lem:xy:D}
	= \gamma_0 \Bigg( \Big( F - \frac{\gamma_1}{\gamma_0}\Big)^2 - \Big( \frac{\gamma_1}{\gamma_0} - E \Big)^2 - \delta \Big( 1 - \Big( \frac{\gamma_1}{\gamma_0} \Big)^2 + \Big( \frac{\gamma_1}{\gamma_0} - E  \Big)^2 \Big) \Bigg).\label{eqn:lem:xy:D:ineq}
\end{equation}

By~\eqref{def:D_inf}, we have $|F - \frac{\gamma_1}{\gamma_0}| \geq D_{\inf}$ and by~\eqref{def:D_sup}, we have $|F - \frac{\gamma_1}{\gamma_0}| \leq D_{\sup}$. Further, we have that $|\frac{\gamma_1}{\gamma_0} - E| \leq M'_f N / (m_f T)$, uniformly with respect to $\omega$, which can be seen as follows: first, note that
\begin{equation*}
	\begin{split}
	\Big| \gamma_{k}(w, N/T) - \gamma_{k}(w, 0) \Big| 
	& \leq \int_0^1 \Big| \gamma_k(w) - \gamma_k(w-\frac{N}{T} (1-x)) \Big| {\rm d}x \\
	& \leq 2 \pi M'_f \int_0^1 \frac{N}{T} (1-x) {\rm d}x  =  \pi M'_f \frac{N}{T} \\
	\end{split}
\end{equation*}

Further, note that we have $\frac{x}{y} - \frac{x_0}{y_0} = \frac{1}{y_0} \big( \frac{x}{y} (y_0 - y) + (x - x_0) \big)$
and thus
\begin{equation}\label{eqn:lem:xy:E}
\begin{split}
	\Big| \frac{\gamma_1}{\gamma_0} - E \Big|
		& \leq \frac{1}{\gamma_0(w; N/T)} \Big( \frac{|\gamma_1|}{\gamma_0} + 1 \Big) \pi M'_f \frac{N}{T}
		\leq \frac{M'_f}{m_f} \frac{N}{T},
\end{split}
\end{equation}
where we have used that
$2 \pi m_f \leq \gamma_0(w; \Delta) := \int_0^1 \gamma_0(w+\Delta (x-1)) {\rm d}x$
and $|\gamma_1|/\gamma_0 \leq 1$.
Employing~\eqref{eqn:lem:xy:D:ineq}, we have now brought the tools together to prove Lemma~\ref{thm:main3:kor1}:
\begin{equation*}
\begin{split}
	-1 \times \eqref{eqn:lem:xy:D}
	& = \gamma_0 \Bigg( \delta \Big( 1 - \Big( \frac{\gamma_1}{\gamma_0} \Big)^2 + \Big( \frac{\gamma_1}{\gamma_0} - E  \Big)^2 \Big) - \Big( F - \frac{\gamma_1}{\gamma_0}\Big)^2 + \Big( \frac{\gamma_1}{\gamma_0} - E \Big)^2 \Bigg) \\
	& \geq 2 \pi m_f \Big( 1 - \rho^2 \Big) \Big( \delta/2 + \delta/2  - D_{\sup}^2 / \big( 1 - \rho^2 \big) \Big)  \Big) \Big)  \geq \pi m_f \delta  \big( 1 - \rho^2 \big). \\
\end{split}
\end{equation*}
For the first inequality we have used the fact that $(\gamma_1/\gamma_0 - E)^2 \geq 0$ and the definitions of $\rho$ and $D_{\sup}$. For the second inequality we have used the condition imposed on $\delta$.

Finally, employing~\eqref{eqn:lem:xy:D:ineq} again, we prove Lemma~\ref{thm:main3:kor2}:
\begin{equation*}
\begin{split}
	\eqref{eqn:lem:xy:D}
	& = \gamma_0 \Bigg( \Big( F - \frac{\gamma_1}{\gamma_0}\Big)^2 - \Big( \frac{\gamma_1}{\gamma_0} - E \Big)^2 - \delta \Big( 1 - \Big( \frac{\gamma_1}{\gamma_0} \Big)^2 + \Big( \frac{\gamma_1}{\gamma_0} - E  \Big)^2 \Big) \Bigg) \\
	& \geq 2 \pi m_f \Big( \Big( F - \frac{\gamma_1}{\gamma_0}\Big)^2 - \Big( \frac{\gamma_1}{\gamma_0} - E \Big)^2 - 2 \delta \Big) \\
	& \geq 2 \pi m_f \Big( D_{\inf}^2 - \Big( \frac{M'_f}{m_f} \frac{N}{T} \Big)^2 - 2 \delta \Big)
	\geq 2 \pi m_f \Big( D_{\inf}^2 / 2 - 2 \delta \Big) \geq \pi m_f D_{\inf}^2 / 2, \\
\end{split}
\end{equation*}
where in the first inequality we have used
\[\Big( \frac{\gamma_1}{\gamma_0} - E  \Big)^2 \leq \Big( \frac{M'_f}{m_f} \frac{N}{T} \Big)^2 \leq 1,\]
as we have $D_{\inf} \leq 2$ and thus $\max\mathcal{N} \leq (m_f / M'_f) T$ follows from condition~\eqref{thm:main3:kor2:assN}. For the second inequality we have used the definition of $D_{\inf}$ and again condition~\eqref{thm:main3:kor2:assN}, by which we have $D_{\inf}^2 / 2 \geq \big( M'_f N / (m_f T) \big)^2$. Finally, for the third inequality we have used that by assumption in the Corollary $2 \delta \leq D_{\inf}^2 / 4$.

The first bound, $q(\delta) \geq \delta \pi m_f ( 1 - \rho^2)$, implies that if
\begin{equation*}
	m > 2 \Big(\frac{\pi m_f ( 1 - \rho^2)}{20 K_0} \frac{\delta}{1+\delta} \Big)^{\frac{3+8d}{1+4d}}
	\quad \text{and} \quad
	\min\mathcal{N} > \Big(\frac{\pi m_f ( 1 - \rho^2)}{20 K_0} \frac{\delta}{1+\delta}\Big)^{\frac{3+4d}{1+2d}} + 1
\end{equation*}
and
\begin{equation*}
	\min\mathcal{N} \geq 16 \big( C_0 \big)^{3} \max\Big\{\frac{20 (1+\delta)}{\delta \pi m_f ( 1 - \rho^2)}, 1\Big\} \Big[ 6 (2\pi M'_f + C) + 1 \Big] 
\end{equation*}
then Assumption~\ref{a:m_and_N} holds, and if
\begin{equation*}
	\begin{split}
		T \geq \max\Big\{ 12 C_1, 12 m \big( C_0 \big)^{3} M'_f \frac{20}{ \pi m_f ( 1 - \rho^2)} \frac{1+\delta}{\delta}\Big\},
	\end{split}
	\end{equation*}
then Assumption~\ref{a:T} holds. Hence, we have proven Lemma~\ref{thm:main3:kor1} where the constants can be chosen as
\[K_1 := 2 \Big(\frac{\pi m_f ( 1 - \rho^2)}{20 K_0} \Big)^{\frac{3+8d}{1+4d}},\]
\begin{multline*}
	K_2 := \max\Big\{\Big(\frac{\pi m_f ( 1 - \rho^2)}{20 K_0} \Big)^{\frac{3+4d}{1+2d}} + 1, \\
		16 \big( C_0 \big)^{3} \max\Big\{\frac{20 (1+(1-\rho^2)/(2 D_{\sup}^2))}{\pi m_f ( 1 - \rho^2)}, 1\Big\} \Big[ 6 (2\pi M'_f + C) + 1 \Big] \Big\}
\end{multline*}
and
\[K_3 := \max\Big\{ 12 C_1, 12 \big( C_0 \big)^{3} M'_f \frac{20}{ \pi m_f ( 1 - \rho^2)} \Big(1+\frac{1-\rho^2}{2 D_{\sup}^2}\Big)\Big\}.\]
The second bound, $q(\delta) \geq \pi D_{\inf}^2 m_f / 2$, implies that if
\begin{equation*}
	m > 2 \Big(\frac{\pi (M'_f)^2}{20 K_0 m_f (1+\delta)} \Big( \frac{\max\mathcal{N}}{T}\Big)^2\Big)^{\frac{3+8d}{1+4d}}
\end{equation*}
and
\begin{equation*}
\begin{split}
	\min\mathcal{N} & > \max\Big\{\Big(\frac{\pi (M'_f)^2}{20 K_0 m_f (1+\delta)} \Big( \frac{											\max\mathcal{N}}{T}\Big)^2\Big)^{\frac{3+4d}{1+2d}} + 1, \\
		& \qquad\qquad 16 \big( C_0 \big)^{3} \max\Big\{\frac{20 (1+\delta)}{\pi D_{\inf}^2 m_f / 2}, 1\Big\} \Big[ 6 (2\pi M'_f + C) + 1 \Big]  \Big\}
\end{split}
\end{equation*}
then Assumption~\ref{a:m_and_N} holds, and if
\begin{equation*}
		T \geq \max\Big\{ 12 C_1, 12 m \big( C_0 \big)^{3} M'_f \frac{20 (1+\delta)}{\pi D_{\inf}^2 m_f / 2} \Big\}.
	\end{equation*}
then Assumption~\ref{a:T} holds.
Hence, we have proven Lemma~\ref{thm:main3:kor2} where the constants can be chosen as
\[K_4 := 2 \Big(\frac{\pi (M'_f)^2}{20 K_0 m_f} \Big)^{\frac{3+8d}{1+4d}},\]
\begin{equation*}
\begin{split}
	K_5 & := \max\Big\{\Big(\frac{\pi (M'_f)^2}{20 K_0 m_f} \Big)^{\frac{3+4d}{1+2d}} + 1,
		16 \big( C_0 \big)^{3} \max\Big\{\frac{20 (1+\frac{1}{8} D_{\inf}^2))}{\pi D_{\inf}^2 m_f / 2}, 1\Big\} \Big[ 6 (2\pi M'_f + C) + 1 \Big] \Big\}
\end{split}
\end{equation*}
and
\begin{equation*}
	K_6 := 12 \max\Big\{ C_1, \big( C_0 \big)^{3} M'_f \frac{20 (1+\frac{1}{8} D_{\inf}^2)}{\pi D_{\inf}^2 m_f / 2} \Big\}.
	\end{equation*}
This finishes the proof of Lemmas~\ref{thm:main3:kor1} and~\ref{thm:main3:kor2}.
\hfill$\square$

\section[Lemmas regarding a]{Lemmas regarding $a$}
\label{sec:Lemmas_a}

\subsection{Outlook}

In this section we state and discuss results relating quantities that are encountered in connection with the localised Yule-Walker estimator.
In Section~\ref{sec:Lemmas_a_state} we state and discuss three lemmas. In Lemma~\ref{lem:rel_a_1} we make precise that $a_{0}^{(p)}(t/T)$ approximates the \emph{time-varying $1$-step linear prediction coefficients} which, for $p \in \IN^*$ and $t=1,\ldots,T$, are defined as
\[ \tilde a_{T}^{(p)}(t) := \arg\min_{a = (a_1, \ldots, a_p)' \in \IR^p}
		\E \Big[\Big( X_{t,T} - \sum_{j=1}^p a_j X_{t-j,T} \Big)^2\Big] = \big( \tilde\Gamma_{T}^{(p)}(t)\big)^{-1} \tilde\gamma^{(p)}_{T}(t),\]
		where
		\begin{equation}\label{eqn:tildeGammaT}
			\begin{split}
			\tilde\gamma_{T}^{(p)}(t) & := (\Cov\left(X_{t,T}, X_{t-1,T}\right), \ldots, \Cov\left(X_{t,T}, X_{t-p,T}\right))', \\
			\tilde\Gamma_{T}^{(p)}(t) & := ( \Cov\left(X_{t-i,T}, X_{t-j,T}\right);\,i,j=1,\ldots,d).
			\end{split}
		\end{equation}

In Lemma~\ref{lem:rel_a_bar} we make precise that $\bar a_{N,T}^{(p)}(t)$, defined in~\eqref{eqn:YW_bar}, is related to $a_{\Delta}^{(p)}(u)$, defined in~\eqref{def:a_Delta}, in the sense that $a_{0}^{(p)}(t/T)$ and $a_{N/T}^{(p)}(t/T)$ approximate $\bar a_{N,T}^{(p)}(t)$. In Lemma~\ref{lem:norm_a} a bound for the norm of $a_{\Delta}^{(p)}(u)$ is provided, which is independent of $p$, $\Delta$ and $u$.

Proofs of the results in Section~\ref{sec:Lemmas_a_state} are provided in Section~\ref{sec:Lemmas_a_proofs} \citep{KleyEtAl2019}. 
The proofs rely on results about expectations of localised autocovariance estimates from Section~\ref{sec:Cov} and an approximation result for inverses of matrices (Lemma~\ref{lem:stoerungslemma2}). For the readers convenience, we include Figure~\ref{fig:map_app_F} in which the dependence of the various results is illustrated graphically.

\begin{figure}[t]
	\begin{center}
		\begin{tikzpicture}[
this_sec/.style = {fill=blue!30},
other_sec/.style = {fill=gray!30},
]



\node [this_sec] (D1) at (9,10) {Lemma~\ref{lem:rel_a_1}};
\node [this_sec] (D2) at (6,10) {Lemma~\ref{lem:rel_a_bar}};
\node [this_sec] (D3) at (3,10) {Lemma~\ref{lem:norm_a}};
\node (D2i) at (5.5,9.5) {(i)};
\node (D2ii) at (6.5,9.5) {(ii)};

\node [other_sec] (F2) at (3,8) {Corollary~\ref{kor:exp_gamma_k}};
\node (F2ii) at (3,7.5) {(i)};
\node (F2iii) at (3,7) {(iii)};
\node [other_sec] (F4) at (3,6) {Lemma~\ref{lem:BoundEV_M0}};
\node (F4i) at (3,5.5) {(i)};

\node [other_sec] (G2) at (8,8) {Lemma~\ref{lem:stoerungslemma2}};

\path[->] (D2i) edge [bend left=50] (F2ii);
\path[->] (D2i) edge [bend left=60] (F4i);
\path[->] (D3) edge [bend right=90] (F4i);
\path[->] (D1) edge [bend left=90] (F4i);
\path[->] (D2i) edge (G2);
\path[->] (D2ii) edge (G2);
\path[->] (D1) edge (G2);
\path[->] (D2ii) edge [bend left=30](F2iii);
\path[->] (D2ii) edge [bend left=60](F4i);

\end{tikzpicture}
	\end{center}
	\vspace*{-2cm}\caption{Map of the lemmas in Section~\ref{sec:Lemmas_a}.}
	\label{fig:map_app_F}
\end{figure}

\subsection{Statement of the lemmas}\label{sec:Lemmas_a_state}
		
The following two lemmas discuss approximation properties of $\bar a_{N,T}^{(p)}(t)$ and $\tilde a_T^{(p)}(t)$:

\begin{lem}\label{lem:rel_a_1}
	Let $(X_{t,T})_{t \in \IZ, T \in \IN^*}$ satisfy Assumptions~\ref{a:loc_stat}, \ref{a:sd_bounded}, \ref{a:sd_der_bounded}, and $\E X_{t,T} = 0$. Define $C_0$ and $C_1$ as in~\eqref{eqn:def_C0}. Then, if $T \geq 2 p^2 C_1$, we have
	\[\| \tilde a_{T}^{(p)}(t) - a_0^{(p)}(t/T) \|
	\leq \frac{1}{T}  \Big( 5 C_0 C_1 \, p^2\Big). \]
\end{lem}

\cite{rousan16} prove a similar bound (Lemma~3):
\[\| \tilde a_{T}^{(p)}(t) - a_0^{(p)}(t/T) \| \leq \frac{D_1}{T}, \quad D_1 := \frac{C p^{1/2} (p 2^p + 1)}{\pi m_f},\]
for $T \geq T_0 := \frac{C p^{3/2}}{\pi m_f}$. Note that (for larger $p$) their constant $D_1$ can be substantially larger than the constant in Lemma~\ref{lem:rel_a_1}, which is largely due to a different representations of $\tilde a_{T}^{(p)}(t) - a_0^{(p)}(t/T)$ in their proof.

It is worth mentioning that in case of a stationary process, where $C_1 = 0$, Lemma~\ref{lem:rel_a_1} implies that $\tilde a_{T}^{(p)}(t)$ and $a_0^{(p)}(t/T)$ coincide.

\begin{lem}\label{lem:rel_a_bar}
	Let $(X_{t,T})_{t \in \IZ, T \in \IN^*}$ satisfy Assumptions~\ref{a:loc_stat}, \ref{a:sd_bounded}, \ref{a:sd_der_bounded}, and $\E X_{t,T} = 0$. Define $C_0$ and $C_1$ as in~\eqref{eqn:def_C0}. Then, if
	
	(i) $T \geq 8 p N C_1$ and $N \geq 4 p^{2} \frac{M_f}{m_f}$, then
	$\| \bar a_{N,T}^{(p)}(t) - a_0^{(p)}(t/T) \|
	\leq \big(9 C_0 C_1 \big)\frac{p \, N}{T} + \big( 3 C_0^2 \big)\frac{p^{2}}{N}$.
	
	(ii) $T \geq 4 p^2 C_1$ and $N \geq 4 p^{2} \frac{M_f}{m_f}$, then
	$\| \bar a_{N,T}^{(p)}(t) - a_{N/T}^{(p)}(t/T) \|
	\leq \big( 5 C_0 C_1 \big) \frac{p^2}{T} + \big( 3 C_0^2 \big)\frac{p^{2}}{N}$.
	
\end{lem}

Note that, if $p^2 = o(T)$, as $T \rightarrow \infty$, then we have, by Lemma~\ref{lem:rel_a_1}, that $\tilde a_T^{(p)}(t)$ and $a_0^{(p)}(t/T)$ are asymptotically equivalent in the sense that the Euclidean norm of the difference vanishes asymptotically. For $N p = o(T)$ and $p = o(N^{1/2})$ we have, by Lemma~\ref{lem:rel_a_bar}(i), that $\bar a_{N,T}^{(p)}(t)$ and $a_0^{(p)}(t/T)$ are asymptotically equivalent, too. Therefore, since $0 \leq p^2 \leq N p$, we have: if $N p = o(T)$ and $p = o(N^{1/2})$, then $\bar a_{N,T}^{(p)}(t)$ and $\tilde a_T^{(p)}(t)$ are asymptotically equivalent.
Note further, that in the case of a tvAR($p$) model, the quantity $\tilde a_T^{(p)}(t)$ coincides with the vector of coefficients $(a_1(t/T), \ldots, a_p(t/T))$, as is evident from the Yule-Walker equations.

It is worth mentioning that in case of a stationary process, where $C_1 = 0$, the bounds in Lemmas~\ref{lem:rel_a_1} and~\ref{lem:rel_a_bar} are independent of $T$.

We will also need the following result that bounds the norm of $a_{\Delta}^{(p)}(u)$:
\begin{lem}\label{lem:norm_a}
	Let $(X_{t,T})_{t \in \IZ, T \in \IN^*}$ satisfy Assumptions~\ref{a:loc_stat}, \ref{a:sd_bounded}, \ref{a:sd_der_bounded}, and $\E X_{t,T} = 0$. Then, for $u \in \IR$, $p \in \IN^*$ and $\Delta \geq 0$, we have
	\[ \| a_{\Delta}^{(p)}(u) \| \leq (2\pi)^{1/2} M_f / m_f =: C_0.\]
\end{lem}

By Lemma~2 in \cite{rousan16} we have $\| a^{(p)}_0(u) \| \leq 2^p$. Their proof adapts arguments from Lemma~4.2 in \cite{dahlgir1998} where $\| \hat a^{(p)}_0(u) \| \leq 2^p$ almost surely is proven. 
We choose to work with the bound from Lemma~\ref{lem:norm_a}, because it has the advantage that it does not depend on $p$. Further, note that neither of the bounds is sharp, as by Cauchy-Schwarz inequality we clearly have $\| a^{(1)}_0(u) \| \leq 1$.

In Lemmas~\ref{lem:rel_v_bar}(i) and(ii) we show similar bounds for the approximation of $\bar v_{N,T}^{(p,h)}(t)$ with $v^{(p,h)}_{0}(t/T)$ or $v^{(p,h)}_{N/T}(t/T)$.

\section[Lemmas regarding v, g and MSPE]{Lemmas regarding $v$, $g$ and MSPE}\label{app:der_thMSPE}

\subsection{Outlook}
In this section we state and discuss results relating quantities that are encountered in connection with the $h$-step ahead forecasting coefficients and the empirical mean squared prediction errors. In particular, this are the quantities $v_{\Delta}^{(p,h)}(u)$, $g_{\Delta}^{(p,h)}(u)$ and ${\rm MSPE}_{\Delta_1, \Delta_2}^{(p,h)}(u)$. 
In Section~\ref{app:der_thMSPE:state} we state and discuss four lemmas. 
In Lemma~\ref{lem:rel_v_bar} we make precise that $\bar v_{N,T}^{(p,h)}(t)$ can be approximated by $v^{(p,h)}_{0}(t/T)$ or $v^{(p,h)}_{N/T}(t/T)$, where $\bar v^{(p,h)}_{\Delta}(u)$ was defined in~\eqref{eqn:def_e_H}. In Lemma~\ref{lem:norm_v_der} we state bounds for norms of $v_{\Delta}^{(p,h)}(u)$ and its derivatives with with respect to $u$ or $\Delta$. In Lemma~\ref{lem:norm_a_bar}, we state bounds for norms of $\bar a_{N,T}^{(p,h)}(t)$. In Lemma~\ref{lem:norm_g_der}(i)--(iii) we state bounds for $g_{\Delta}^{(p,h)}(u)$ and its derivates with respect to $u$ or $\Delta$. In Lemma~\ref{lem:norm_g_der}(iv)--(vi) we state bounds for the derivatives of ${\rm MSPE}_{\Delta_1, \Delta_2}^{(p,h)}(u)$ with respect to $u$, $\Delta_1$ or $\Delta_2$.

Proofs of the results in Section~\ref{app:der_thMSPE:state} are provided in Section~\ref{app:der_thMSPE:proof} \citep{KleyEtAl2019}. The proofs rely on some analogous bounds for the quantities related to the Yule-Walker estimator (Section~\ref{sec:Lemmas_a}), on results on expectations of localised autocovariance estimates (Section~\ref{sec:Cov}) and an approximation result for powers of matrices (Lemma~\ref{lem:expansionMatrixPower}). For the readers convenience, we include Figure~\ref{fig:map_app_G} in which the dependence of the various results is illustrated graphically.

\begin{figure}[t]
	\begin{center}
		\begin{tikzpicture}[
this_sec/.style = {fill=blue!30},
other_sec/.style = {fill=gray!30},
]


\node [this_sec] at (-2,0) {Lemma~\ref{lem:rel_v_bar}};
\node (E1i) at (-2,-0.5) {(i)};
\node (E1ii) at (-2,-1) {(ii)};
\node [this_sec] at (1,0) {Lemma~\ref{lem:norm_v_der}};
\node (E2i) at (1,-0.5) {(i)};
\node (E2ii) at (1,-1) {(ii)};
\node (E2iii) at (1,-1.5) {(iii)};
\node [this_sec] (F3) at (4,0) {Lemma~\ref{lem:norm_a_bar}};
\node (E3i) at (4,-0.5) {(i)};
\node (E3ii) at (4,-1) {(ii)};
\node [this_sec] at (7,0) {Lemma~\ref{lem:norm_g_der}};
\node (E4i) at (7,-0.5) {(i)};
\node (E4ii) at (7,-1) {(ii)};
\node (E4iii) at (7,-1.5) {(iii)};
\node (E4iv) at (7,-2.5) {(iv)};
\node (E4v) at (7,-2) {(v)};
\node (E4vi) at (7,-3) {(vi)};

\node [other_sec] (D3) at (-0.5,2) {Lemma~\ref{lem:norm_a}};
\node [other_sec] (D2) at (2.5,2) {Lemma~\ref{lem:rel_a_bar}};
\node (D2i) at (2,1.5) {(i)};
\node (D2ii) at (3,1.5) {(ii)};
\node [other_sec] (D1) at (5.5,2) {Lemma~\ref{lem:rel_a_1}};

\node [other_sec] (G3) at (-3.5,-3) {Lemma~\ref{lem:expansionMatrixPower}};

\node [other_sec] (F4) at (1,-3) {Lemma~\ref{lem:BoundEV_M0}};
\node (F4i) at (1,-3.5) {(i)};

\path[->] (E1i) edge [bend right = 60] (G3);
\path[->] (E1ii) edge [bend right = 30] (G3);
\draw[->] (E1i) .. controls (0,-1) and (-1, 1) .. (D2i);
\draw[->] (E1i) .. controls (0,-.5) and (-1,0) .. (D3);
\draw[->] (E2i) .. controls (-0.5,-0.5) .. (D3);
\draw[->] (E1ii) .. controls (0,-1) and (-0.7,0.6) .. (0.7,0.7) .. controls (2,1) .. (D2ii);
\draw[->] (E3i) .. controls (2.5,-0.5) .. (D2ii);
\draw[->] (E3i) .. controls (2,-1) and (2.5,1) .. (D3);
\draw[->] (E2ii) .. controls (-1.5,-1.5) and (-1,-4).. (F4i);
\draw[->] (E2iii) .. controls (-1,-1.5) and (-1,-3.5).. (F4i);
\path[->] (E3ii) edge [bend right = 5] (E2i);
\draw[->] (E4i) .. controls (4,-0.8) and (5,-2) .. (E2i);
\draw[->] (E4i) .. controls (4,-1.5) and (5,-2) .. (E2ii);
\draw[->] (E4i) .. controls (4,-2) and (3,-3.5) .. (F4i);
\draw[->] (E4ii) .. controls (4,-1.2) and (4,-2) .. (E2i);
\draw[->] (E4ii) .. controls (4,-1.5) and (5,-2.5) .. (E2ii);
\draw[->] (E4ii) .. controls (4,-2) and (3,-4) .. (F4i);
\draw[->] (E4iii) .. controls (4,-1.5) and (5,-2.5) .. (E2iii);
\draw[->] (E3ii) .. controls (3,-1) and (2,-3) .. (E1ii);
\path[->] (E4iv) edge [out=0, in=0] (E4ii);
\path[->] (E4v) edge [out=0, in=0] (E4iii);
\path[->] (E4vi) edge [out=0, in=0] (E4i);


\end{tikzpicture}
	\end{center}
	\caption{Map of the lemmas in Section~\ref{app:der_thMSPE}.}
	\label{fig:map_app_G}
\end{figure}


\subsection{Statement of the lemmas}\label{app:der_thMSPE:state}
The following lemma is constructed analogously to Lemma~\ref{lem:rel_a_bar}, but for the $h$-step ahead coefficients.
\begin{lem}\label{lem:rel_v_bar}
	Let $(X_{t,T})_{t \in \IZ, T \in \IN^*}$ satisfy Assumptions~\ref{a:loc_stat}, \ref{a:sd_bounded}, \ref{a:sd_der_bounded}, and $\E X_{t,T} = 0$. Define $C_0$ and $C_1$ as in~\eqref{eqn:def_C0}. Then, we have, for $\bar v_{N,T}^{(p,h)}(t)$ defined in~\eqref{eqn:def_e_H},
	
	(i) if $T \geq 18 C_1 p N$ and $N \geq 6 p^{2} C_0$, with $v_0^{(p,h)}(t/T)$ defined in~\eqref{def:v_Delta}, that
	\begin{equation*}
	\begin{split}
	\| \bar v_{N,T}^{(p,h)}(t) - v_0^{(p,h)}(t/T) \|
	&  \leq h \big( 2 C_0 \big)^h \Big(5 C_1\frac{p N}{T} + 2 \frac{p^{2}}{N} C_0 \Big).
		\end{split}
\end{equation*}
	
	(ii) if $T \geq 10 C_1 p^2$ and $N \geq 6 p^{2} C_0$, with $v_{N/T}^{(p,h)}(t/T)$ defined in~\eqref{def:v_Delta}, that
		\begin{equation*}
	\begin{split}
	\| \bar v_{N,T}^{(p,h)}(t) - v_{N/T}^{(p,h)}(t/T) \|
	&  \leq h \big( 2 C_0 \big)^h \Big(3 C_1\frac{p^2}{T} + 2 \frac{p^{2}}{N} C_0 \Big).
		\end{split}
\end{equation*}
	
	
\end{lem}

Further, for the norms of $u \mapsto v_{\Delta}^{(p,h)}(u)$ and $\Delta \mapsto v_{\Delta}^{(p,h)}(u)$ and the derivatives, we have the following bounds:

\begin{lem}\label{lem:norm_v_der}
	Let $(X_{t,T})_{t \in \IZ, T \in \IN^*}$ satisfy Assumptions~\ref{a:loc_stat}, \ref{a:sd_bounded}, \ref{a:sd_der_bounded}, and $\E X_{t,T} = 0$. $C_0$ as in~\eqref{eqn:def_C0} and $m_f$, $M_f$, $M'_f$ from the assumptions. Then, with $v_{\Delta}^{(p,h)}(u)$ defined in~\eqref{def:v_Delta}, we have
	
	(i) $\| v_{\Delta}^{(p,h)}(u) \| \leq \big( C_0 \big)^h$,
	
	(ii) $v_{\Delta}^{(p,h)}(\cdot)$ is continuously differentiable, with
	\[
		\Big\| \frac{\partial}{\partial u} v_{\Delta}^{(p,h)}(u) \Big\|
		\leq h \big( C_0 \big)^h M'_f (m_f^{-1} + M_f^{-1} ),
	\]
	
	(iii) $\Delta \mapsto v_{\Delta}^{(p,h)}(u)$, $\Delta > 0$, is continuously differentiable, with
	\[
		\Big\| \frac{\partial}{\partial \Delta} v_{\Delta}^{(p,h)}(u) \Big\|
		\leq 2 h \big( C_0 \big)^h M'_f (m_f^{-1} + M_f^{-1} ) / \Delta.
	\]
\end{lem}

Lemma~\ref{lem:norm_v_der} also holds for $h=1$. Part~(i) thus extends Lemma~\ref{lem:norm_a}.

Finally, we use Lemmas~\ref{lem:rel_a_1}, \ref{lem:rel_a_bar} and~\ref{lem:norm_a} to bound the norm of $\bar a_{N,T}^{(p)}(t)$ and $\bar v_{N,T}^{(p)}(t,h)$.
\begin{lem}\label{lem:norm_a_bar}
	Let $(X_{t,T})_{t \in \IZ, T \in \IN^*}$ satisfy Assumptions~\ref{a:loc_stat}, \ref{a:sd_bounded}, \ref{a:sd_der_bounded}, and $\E X_{t,T} = 0$. Define $C_0$ and $C_1$ as in~\eqref{eqn:def_C0}. Then,
	
	(i) for	$T \geq 10 C_1 p^2$ and $N \geq 6 C_0 p^{2}$ we have, for $\bar a_{N,T}^{(p,h)}(t)$ defined in~\eqref{eqn:YW_bar},
	\[\| \bar a_{N,T}^{(p)}(t) \|
	\leq 2 C_0,
	\text{ and }
	\| \bar v_{N,T}^{(p,h)}(t) \|_{\infty}
	\leq \Big( 2 C_0 + 1 \Big)^h.\]
	
	Further,
	(ii) for $T \geq 6 h 2^h C_1 p^2$ and $N \geq 4 h 2^h C_0 p^{2}$
	we have, for $\bar v_{N,T}^{(p,h)}(t)$ defined in~\eqref{eqn:def_e_H},
	\[\| \bar v_{N,T}^{(p,h)}(t) \|	\leq 2 \big( C_0 \big)^h.\]
\end{lem}

Note that Lemma~\ref{lem:norm_a_bar}(i) implies that we have
$\| \bar v_{N,T}^{(p,h)}(t) \| \leq p^{1/2} \Big( 2 C_0 + 1 \Big)^h$. The bound in Lemma~\ref{lem:norm_a_bar}(ii) does not depend on $p$, but require larger $T$ and $N$.

An important observation is that, as a function of $u$, ${\rm MSPE}_{N/T, n/T}^{(p,h)}(u)$ is differentiable with bounded derivative

\begin{lem} \label{lem:norm_g_der}
	Let $(X_{t,T})_{t \in \IZ, T \in \IN^*}$ satisfy Assumptions~\ref{a:loc_stat}, \ref{a:sd_bounded}, \ref{a:sd_der_bounded}, and $\E X_{t,T} = 0$. Define $C_0$ as in~\eqref{eqn:def_C0} and the other constants from the assumptions. Then, the function
	$g^{(p,h)}_{\Delta}$, defined in~\eqref{eqn:defg}, is continuously differentiable and the derivatives are bounded. More precisely, We have
	
	(i)
	\[ \big| g^{(p,h)}_{\Delta}(u) \big|
		\leq 4 M_f \big( C_0 \big)^{2 h},
	\]
	
	(ii)
	\[ \Big| \frac{\partial}{\partial u} g^{(p,h)}_{\Delta}(u) \Big|
		\leq 4 \big( 2 h + 1 \big) \big( C_0 \big)^{2 h+1} M'_f,
	\]
	
	(iii)
	\[ \Big| \frac{\partial}{\partial \Delta} g^{(p,h)}_{\Delta}(u) \Big|
		\leq 8 \big( 2 h + 1 \big) \big( C_0 \big)^{2 h+1} M'_f / \Delta.
	\]
	
	In particular, this implies, or ${\rm MSPE}_{\Delta_1, \Delta_2}^{(p,h)}(u)$ defined in~\eqref{eqn:thMSPE}, that
	
	(iv)
	\[ \Big| \frac{\partial}{\partial u} {\rm MSPE}_{\Delta_1, \Delta_2}^{(p,h)}(u)\Big|
		\leq 4 \big( 2 h + 1 \big) \big( C_0 \big)^{2 h+1} M'_f.
	\]
	
	(v)
	\[ \Big| \frac{\partial}{\partial \Delta_1} {\rm MSPE}_{\Delta_1, \Delta_2}^{(p,h)}(u)\Big|
		\leq 8 \big( 2 h + 1 \big) \big( C_0 \big)^{2 h+1} M'_f / \Delta_1.
	\]
	
	(vi)
	\[ \Big| \frac{\partial}{\partial \Delta_2} {\rm MSPE}_{\Delta_1, \Delta_2}^{(p,h)}(u)\Big|
		\leq 8 M_f \big( C_0 \big)^{2 h} / \Delta_2.
	\]
	
\end{lem}

\section{Properties of empirical localised autocovariances}\label{sec:Cov}

\subsection{Outlook}\label{sec:Cov:outlook}

In this section we establish properties of the empirical localised autocovariances under local stationarity. In Section~\ref{sec:Cov:mom} we state three lemmas about the estimators' moments and in Section~\ref{sec:Cov:concentration} we state two lemmas about the concentration of the estimators. More precisely, in Lemma~\ref{lem:exp_gamma_k} we approximate the expectation of the empirical autocovariance and state bounds for the approximation error. In Corollary~\ref{kor:exp_gamma_k} we employ the approximation results from Lemma~\ref{lem:exp_gamma_k} to approximate matrices of such expectations and bound the approximation error (in spectral norm). In Lemma~\ref{lem:BoundEV_M0} we establish lower and upper bounds for the eigenvalues of $\Gamma_{\Delta}^{(p)}(u)$ and $\E \hat\Gamma_{N,T}^{(p)}(t)$. In Lemma~\ref{lem:exp_ineq_gamma_k2} we establish a concentration result for the localised empirical autocovariance. Lemma~\ref{lem:exp_ineq_gamma_k2} follows as a special case from Lemma~\ref{lem:expIneqGenSum} where a concentration result for generalised sums under local stationarity is established.

Proofs of the results are proved in Section~\ref{sec:Cov:proof} \citep{KleyEtAl2019}. The proofs rely on technical results to bound the matrix norm of perturbed inverse matrices and approximation of sums by integrals (cf. Section~\ref{app:TRes}) as well as on general concentration results from \cite{sausta91} which we cite in Section~\ref{sec:res_sausta91} \citep{KleyEtAl2019}.
For the readers convenience, we include Figure~\ref{fig:map_app_H} in which the dependence of the various results is illustrated graphically.

\begin{figure}[t]
	\begin{center}
		\begin{tikzpicture}[
this_sec/.style = {fill=blue!30},
other_sec/.style = {fill=gray!30},
]


\node [this_sec] at (-2,0) {Lemma~\ref{lem:exp_gamma_k}};
\node (F1i) at (-2,-0.5) {(i)};
\node (F1ii) at (-2,-1) {(ii)};
\node (F1iii) at (-2,-1.5) {(iii)};
\node [this_sec] at (1,0) {Corollary~\ref{kor:exp_gamma_k}};
\node (F2i) at (1,-0.5) {(i)};
\node (F2ii) at (1,-1) {(ii)};
\node (F2iii) at (1,-1.5) {(iii)};
\node [this_sec] (F3) at (7,0) {Lemma~\ref{lem:FpnEV}};
\node [this_sec] at (4,0) {Lemma~\ref{lem:BoundEV_M0}};
\node (F4ia) at (4,-0.5) {(i-a)};
\node (F4ib) at (4,-1) {(i-b)};
\node (F4iia) at (4,-1.5) {(ii-a)};
\node (F4iib) at (4,-2) {(ii-b)};
\node (F4iic) at (4,-2.5) {(ii-c)};

\node [this_sec] (F5) at (1,-5) {Lemma~\ref{lem:exp_ineq_gamma_k2}};
\node [this_sec] (F6) at (4,-5) {Lemma~\ref{lem:expIneqGenSum}};

\node [other_sec] (G8) at (-0.5,1) {Lemma~\ref{lem:sumApprox}};
\node [other_sec] (G1) at (6,-3) {Lemma~\ref{lem:stoerungslemma}};

\node [other_sec] (H1) at (2.5,-6) {Lemma~\ref{lem:sausta91:Lemma24}};
\node [other_sec] (H2) at (5.5,-6) {Lemma~\ref{lem:sausta91:Theorem417}};

\path[->] (F1i) edge [out = 0, in = -90] (G8);
\path[->] (F1ii) edge [out=180, in=200] (F1i);
\path[->] (F1iii) edge [out=180, in=160] (F1i);
\path[->] (F2iii) edge (F1ii);
\path[->] (F2ii) edge (F1iii);
\path[->] (F4iia) edge [out=0, in=0] (F4ia);
\path[->] (F4iia) edge [out=-20, in=-90] (F3);
\path[->] (F4iib) edge [out = 0, in = 90] (G1);
\path[->] (F4iia) edge [bend right=10] (F2iii);
\path[->] (F4iib) edge [bend left = 10] (F2iii);
\path[->] (F4iic) edge [out=180, in=200] (F4iib);
\path[->] (F5) edge (F6);
\path[->] (F6) edge (H1);
\path[->] (F6) edge (H2);

\end{tikzpicture}
	\end{center}
	\caption{Map of the lemmas in Section~\ref{sec:Cov}.}
	\label{fig:map_app_H}
\end{figure}

\subsection{Approximations for moments}\label{sec:Cov:mom}

\begin{lem}
	\label{lem:exp_gamma_k}
	Let $(X_{t,T})_{t \in \IZ, T \in \IN^*}$ satisfy Assumptions~\ref{a:loc_stat} and \ref{a:sd_der_bounded}, and $\E X_{t,T} = 0$. Then, with $\hat \gamma_{k;N,T}(t)$ defined in~\eqref{def:acf}, $f(u, \lambda)$ and $C$ from Assumption~\ref{a:loc_stat}, and $M'_f$ from Assumption~ \ref{a:sd_der_bounded}, we have:
	(i)
	\begin{equation*}
		\Big| \E \hat \gamma_{k;N,T}(t) - \frac{N-|k|}{N} \int_{-\pi}^{\pi} \Bigg[ \int_0^1 f\big(\frac{t-N+|k|}{T} + \frac{N-|k|}{T} u, \lambda\big) {\rm d}u \Bigg] {\rm e}^{{\rm i} |k| \lambda} {\rm d}\lambda \Big|
		\leq \frac{2\pi M'_f + C}{T}
	\end{equation*}
	and (ii)
	\begin{equation*}
		\Big| \E \hat \gamma_{k;N,T}(t) - \frac{N-|k|}{N} \int_{-\pi}^{\pi} \Bigg[ \int_0^1 f\big(\frac{t-N}{T} + \frac{N}{T} u, \lambda\big) {\rm d}u \Bigg] {\rm e}^{{\rm i} |k| \lambda} {\rm d}\lambda
		\Big| \leq \frac{2\pi (|k|+1) M'_f + C}{T}
	\end{equation*}
	and (iii)
	\begin{equation*}
		\Big| \E \hat \gamma_{k;N,T}(t) - \frac{N-|k|}{N} \gamma_{k}(t/T)\Big| \\
		\leq \frac{2 \pi M'_f (N-|k|+1) + C}{T}.
		\end{equation*}
\end{lem}

\begin{kor}\label{kor:exp_gamma_k}
	Under the conditions of Lemma~\ref{lem:exp_gamma_k}, with $\tilde \Gamma_{T}^{(p)}(t)$ and $\tilde \gamma_{T}^{(p)}(t)$ defined in~\eqref{eqn:tildeGammaT}, $\Gamma_{\Delta}^{(p)}(u)$ and $\gamma_{\Delta}^{(p)}(u)$ defined in~\eqref{eqn:gammaDelta_u}, $\hat\Gamma_{N,T}^{(p)}(t)$ and $\hat\gamma_{N,T}^{(p)}(t)$ defined in~\eqref{eqn:hatGamma_t}, and $F_{p,n}$ and $f_{p,n}$ defined for any $n=1,2,\ldots$ and $p = 1, \ldots, n$ as
	\[F_{p,N} := \big( 1 - |j-k|/N \big)_{j,k=1,\ldots,p}, \quad \text{and} \quad
	  f_{p,N} := \big( 1-1/N, \ldots, 1-p/N \big)',\]
	we have: (i)
	\begin{align*}
		\| \tilde \Gamma_{T}^{(p)}(t) - \Gamma_0^{(p)}(t/T) \| & \leq \frac{p^{2}}{T} (2 \pi M'_f + C) &
		\| \tilde \gamma_{T}^{(p)}(t) - \gamma_0^{(p)}(t/T) \| & \leq \frac{p^{1/2}}{T}  C
	\end{align*}
	and (ii)
	\begin{align*}
		\| \E \hat\Gamma_{N,T}^{(p)}(t) - F_{p,N} \circ \Gamma_0^{(p)}(t/T) \| & \leq \frac{p}{T} (2 \pi M'_f (N+1) + C)  \\
		\| \E \hat\gamma_{N,T}^{(p)}(t) - f_{p,N} \circ \gamma_0^{(p)}(t/T) \| & \leq \frac{p^{1/2}}{T}  (2 \pi M'_f N + C) 
	\end{align*}
	and (iii)
	\begin{align*}
		\| \E \hat\Gamma_{N,T}^{(p)}(t) - F_{p,N} \circ \Gamma_{N/T}^{(p)}(t/T) \| & \leq \frac{p^{2}}{T} (2 \pi M'_f + C) \\
		\| \E \hat\gamma_{N,T}^{(p)}(t) - f_{p,N} \circ \gamma_{N/T}^{(p)}(t/T) \| & \leq 2 \frac{p^{3/2}}{T} \Big( 2 \pi M'_f + C \Big)
	\end{align*}
	\end{kor}

\begin{lem}\label{lem:BoundEV_M0}
		Let $(X_{t,T})_{t \in \IZ, T \in \IN^*}$ satisfy Assumptions~\ref{a:loc_stat}, \ref{a:sd_bounded}, and \ref{a:sd_der_bounded}, and $\E X_{t,T} = 0$. Then, we have:
	
	(i-a) the matrices $\Gamma^{(p)}(u)$ and $\Gamma_{\Delta}^{(p)}(u)$ are positive definite, hence invertible, for $u \in \IR$ and $\Delta \geq 0$, with their eigenvalues between $m_f$ and $M_f$. In other words, the norms of the matrices and their inverses are uniformly bounded:
	\[m_f \leq 1/\|\Gamma_{\Delta}^{(p)}(u)^{-1}\| \leq \|\Gamma_{\Delta}^{(p)}(u)\| \leq M_f.\]
	(i-b) the norms of the respective vectors are uniformly bounded:
	\[\|\gamma_{\Delta}^{(p)}(u)\| \leq (2\pi)^{1/2} M_f.\]
	
	(ii-a) The largest eigenvalue of $\E \hat\Gamma_{N,T}^{(p)}(t)$ satisfies the following bound:
	\[\| \E \hat\Gamma_{N,T}^{(p)}(t) \| \leq M_f + \frac{p^2}{T} (2 \pi M'_f + C).\]
	(ii-b) if $T > m_f^{-1} p^2 (2 \pi M'_f + C)$, then the matrix $\E \hat\Gamma_{N,T}^{(p)}(t)$ is positive definite, and the smallest eigenvalue satisfies the following bound:
	\[m_f - \frac{p^2}{T} (2 \pi M'_f + C) \leq \frac{1}{\| \big(\E \hat\Gamma_{N,T}^{(p)}(t)\big)^{-1}\|}\]
	(ii-c) in particular, if $T \geq 2 m_f^{-1} p^2 (2 \pi M'_f + C)$ we thus have
	\[\frac{1}{2} m_f \leq \frac{1}{\| \big(\E \hat\Gamma_{N,T}^{(p)}(t)\big)^{-1}\|} \leq \| \E \hat\Gamma_{N,T}^{(p)}(t) \| \leq \frac{3}{2} M_f.\]
\end{lem}

\subsection{Exponential inequalities for empirical covariances}\label{sec:Cov:concentration}

We now state an exponential inequalities for the empirical covariances:
\begin{lem}
	\label{lem:exp_ineq_gamma_k2}
	Let $(X_{t,T})_{t \in \IZ, T \in \IN^*}$ satisfy Assumptions~\ref{a:loc_stat}, \ref{a:mixing}, \ref{a:sd_bounded}, and \ref{a:MomentCond} and $\E X_{t,T} = 0$. Then, for $T \geq C/(\pi m_f)$, $n \in \IN^*$, $h \in \IN$ and $\varepsilon > 0$, we have
	\begin{equation*}
	\begin{split}
	& \IP\Big(\Big| \hat\gamma_{h;N,T}(t) - \E \hat\gamma_{h;N,T}(t) \Big| \geq \varepsilon\Big)
	 \leq 
	 \exp\Bigg( - \frac{\varepsilon^2}{2\Big(\frac{C_{1,1} h_*}{N-|h|} + \varepsilon^{(3+4d)/(2+2d)} \Big( \frac{C_{2,1} h_*}{N-|h|} \Big)^{1/(2+2d)}\Big)} \Bigg) \\
	& \quad \leq \begin{cases}
		\exp\Big( -\cfrac{\varepsilon^2}{4} \cfrac{N-|h|}{C_{1,1} h_*} \Big) & \varepsilon \leq \Big(\cfrac{h_*}{N-|h|} \Big)^{(1+2d)/(3+4d)} \Big(\cfrac{C_{1,1}^{2+2d}}{C_{2,1}}\Big)^{1/(3+4d)} \\ 
		\exp\Big( -\cfrac{1}{4} \Big( \varepsilon \cfrac{N-|h|}{C_{2,1} h_*} \Big)^{1/(2+2d)} \Big) & \varepsilon \geq \Big(\cfrac{h_*}{N-|h|} \Big)^{(1+2d)/(3+4d)} \Big(\cfrac{C_{1,1}^{2+2d}}{C_{2,1}}\Big)^{1/(3+4d)},
	\end{cases}
	\end{split}
	\end{equation*}
where $h_* := |h| + I\{h=0\}$, $\hat\gamma_{h;N,T}(t)$ is defined in~\eqref{def:acf}, and the constants $C_{1,1}$ and $C_{2,1}$ are defined in~\eqref{lem:exp_ineq_gamma_k2:def:C1C2}.
\end{lem}
Note that the right hand side does not depend on $t$.

\begin{lem}\label{lem:expIneqGenSum}
Let $(X_{t,T})_{t \in \IZ, T \in \IN^*}$ satisfy Assumptions~\ref{a:loc_stat}, \ref{a:mixing}, \ref{a:sd_bounded}, and \ref{a:MomentCond} and $\E X_{t,T} = 0$. Let~$a_t$, $t=b, \ldots, b+n-1$ be a bounded sequence of numbers; i.\,e., $|a_t| \leq A$. Then, for $\alpha \in \IN^*$, $T \geq C/(\pi m_f)$, $n \in \IN^*$, $b \in \IZ$, $h \in \IN$ and $\varepsilon > 0$, we have
\begin{equation*}
	\begin{split}
	& \IP \Big( \Big| n^{-1} \sum_{t=b}^{b+n-1} a_t (X_{t,T}^{\alpha} X_{t+h,T}^{\alpha} - \E (X_{t,T}^{\alpha} X_{t+h,T}^{\alpha})) \Big| > \varepsilon \Big) \\
	& \quad \leq 
	 \exp\Bigg( - \frac{\varepsilon^2}{2\Big(\frac{C_{1,\alpha} A^2 h_*}{n} + \varepsilon^{(3+4 \alpha d)/(2+2\alpha d)} \Big( \frac{C_{2,\alpha} A h_*}{n} \Big)^{1/(2+2\alpha d)}\Big)} \Bigg) \\
	& \quad \leq \begin{cases}
		\exp\Big( -\cfrac{(\varepsilon/A)^2}{4} \cfrac{n}{C_{1,\alpha} h_*} \Big) & \varepsilon \leq A \Big( \cfrac{C_{1,\alpha}^{2+2 \alpha d}}{C_{2,\alpha}} \Big)^{1/(3+4\alpha d)} \Big(\cfrac{h_*}{n}\Big)^{(1+2 \alpha d)/(3+4\alpha d)} \\ 
		\exp\Big( -\cfrac{1}{4} \Big( \cfrac{\varepsilon}{A} \cfrac{n}{C_{2,\alpha} h_*} \Big)^{1/(2+2\alpha d)} \Big) & \varepsilon \geq A \Big(\cfrac{C_{1,\alpha}^{2+2 \alpha d}}{C_{2,\alpha}} \Big)^{1/(3+4\alpha d)} \Big(\cfrac{h_*}{n}\Big)^{(1+2 \alpha d)/(3+4\alpha d)},
	\end{cases}
	\end{split}
	\end{equation*}
	where $h_* := |h| + I\{h=0\}$ and the constants $C_{1,\alpha}$ and $C_{2,\alpha}$ are defined in~\eqref{lem:exp_ineq_gamma_k2:def:C1C2} in the proof [depending only on $\alpha$, $d$, $C$, $M_f$, $\rho$, and $K$]. 
\end{lem}
Note the important fact that the bounds in the inequality do not depend on~$b$.

\section{Technical Results}\label{app:TRes}

In the previous sections we used the following general results, which are not restricted to locally stationary processes. In some of these technical lemmas we denote by $\| \cdot \|_M$ or $\| \cdot \|_v$ an arbitrary matrix or vector norm, respectively. Special properties we require include \emph{submultiplicativity} of a matrix norm, and \emph{compatibility} of a matrix norm with a vector norm. A matrix norm which satisfies $\| A B \|_M \leq \|A \|_M \|B \|_M$ for all square matrices ($m=n$), is said to be submultiplicative. A matrix norm $\| \cdot \|_M$ and vector norm $\| \cdot \|_v$ are said to be compatible if $\| A x \|_v \leq \| A \|_M \| x \|_v$ for all square matrices $A$ and vectors $x$ (of sizes that allow for the matrix product).

\begin{lem}
	\label{lem:stoerungslemma}
	Let $A \in \IR^{p \times p}$ be an invertible matrix and $\Delta \in \IR^{p \times p}$ be a matrix with $\| A^{-1} \|_M \cdot \| \Delta \|_M < 1$ for a submultiplicative matrix norm $\| \cdot \|_M$. Then, the matrix $A + \Delta$ is invertible and we have
	\[ \| (A + \Delta)^{-1} \|_M \leq \frac{\| A^{-1} \|_M}{1 - \| A^{-1} \|_M \cdot \| \Delta \|_M}\]
\end{lem}

An important corollary to the above lemma is the following:
\begin{lem}
	\label{lem:stoerungslemma2}
	Let $A \in \IR^{p \times p}$ be an invertible matrix and $\Delta \in \IR^{p \times p}$ be a matrix with $\| A^{-1} \|_M \cdot \| \Delta \|_M \leq c < 1$ for a submultiplicative matrix norm $\| \cdot \|_M$. Then, the matrix $A + \Delta$ is invertible and we have
	\[ \| (A + \Delta)^{-1} - A^{-1} \|_M \leq \| \Delta \|_M \frac{\| A^{-1} \|_M^2}{1 - c}.\]
\end{lem}

\begin{lem}\label{lem:expansionMatrixPower}
	Let $A$ and $A_0$ be two square matrices and $\| \cdot \|_M$ be a submultiplicative matrix norm. Then, for any $h \in \IN$,
\begin{equation*}
	\|A^h - A_0^h\|_M
	\leq h \| A-A_0 \|_M \big( \|A-A_0\|_M + \|A_0\|_M \big)^{h-1}.
\end{equation*}
\end{lem}

\begin{lem}\label{lem:ineq_productOfTwo}
	Let $u$ and $v$ be two real-valued random variables. Further, let $u_0$ and $v_0$ be two real numbers. Then, for all $\varepsilon > 0$
\begin{equation*}
	\IP( | u v - u_0 v_0 | > \varepsilon)
		\leq \IP \Big( |u - u_0| > \frac{1}{2} \frac{\varepsilon}{(|v_0|^2 + \varepsilon)^{1/2}} \Big)
		+ \IP \Big( |v - v_0| > \frac{1}{2} \frac{\varepsilon}{(|u_0|^2 + \varepsilon)^{1/2}} \Big).
\end{equation*}
\end{lem}

For the proof in the main part we need the following lemma:
\begin{lem}\label{lem:disentangle}
	Let $X_t$ and $\hat a_t$, $t=1,\ldots,n$, be two sequences of random variables, and $\alpha_t$, $t=1,\ldots,n$ be a sequence of numbers. Assume that there exists a constant $m_2^2 > 0$ such that $\max_{t=1,\ldots,n} \E X_t^2 \leq m_2^2 < \infty$. Then, for any $\varepsilon > 0$, we have
	\begin{equation*}
		\begin{split}
		& \IP \Big( \Big| \sum_{t=1}^n \big( \hat a_t X_t - \alpha_t \E(X_t) \big) \Big| > n \varepsilon \Big) \\
		& \leq \IP \Big( \sup_{t=1,\ldots,n} | \hat a_t - \alpha_t |
		> \frac{\varepsilon}{2 \big( (2 m_2^2)^2 + \varepsilon^2 \big)^{1/4}}  \Big) \\
		& \qquad + \IP \Big( \Big| \sum_{t=1}^n ( X_t^2 - \E X_t^2) \Big| > n \varepsilon / 2 \Big)
			+ \IP \Big( \Big| \sum_{t=1}^n \alpha_t (X_t - \E X_t) \Big| > n \varepsilon / 2 \Big).
	\end{split}
	\end{equation*}
\end{lem}

We will further use the following lemmas:
\begin{lem}\label{lem:ineqQuotientMV}
	Let $M \in \IR^{p \times p}$ be a random $p \times p$ matrix with existing expectation $M_0 := \E M$, which is assumed to be invertible. Further, let $v$ be a $\IR^p$-valued random vector $v$ with existing expectation $\E v := v_0$. Then, for every submultiplicative matrix norm $\| \cdot \|_M$ that is compatible with the (vector) norm $\|\cdot \|_v$, we have: for every $\varepsilon > 0$
	\begin{equation*}
		\begin{split}
			& \IP \Big( \Big\| M^{-1} v - M_0^{-1} v_0 \Big\|_v > \varepsilon \Big) \\
				& \leq \IP \Big( \| M - M_0 \|_M > \frac{1}{2 \|M_0^{-1}\|_M} \Big) + \IP \Big( \|v-v_0\|_v > \frac{\varepsilon}{4} \frac{1}{\|M_0^{-1}\|_M} \Big) \\
				& \qquad + \IP \Big( \| M - M_0 \|_M  > \frac{\varepsilon}{4} \frac{1}{(\|M_0^{-1}\|_M)^2 \, \| v_0 \|_v} \Big) I\{ \|v_0\|_v \neq 0\}.
		\end{split}
	\end{equation*}
\end{lem}

\begin{lem}\label{lem:PowerYW}
Let $x = (x_1, \ldots, x_p)$ be a random vector and $x_0 = (x_{0,1}, \ldots, x_{0,p})$ be a deterministic vector. Define two $p \times p$ matrices
\[
A := \begin{pmatrix}
							x_1 & x_2 & \cdots & x_{p-1} & x_p \\
							1						 & 0            & \cdots & 0              & 0 \\
							0            & 1            & \cdots & 0              & 0 \\
							\vdots       & \ddots       & \cdots & 0              & 0 \\
							0            & 0            & \cdots & 1              & 0
						\end{pmatrix}
\text{ and }
A_0 := \begin{pmatrix}
							x_{0,1} & x_{0,2} & \cdots & x_{0,p-1} & x_{0,p} \\
							1						 & 0            & \cdots & 0              & 0 \\
							0            & 1            & \cdots & 0              & 0 \\
							\vdots       & \ddots       & \cdots & 0              & 0 \\
							0            & 0            & \cdots & 1              & 0
						\end{pmatrix}
\]
For any $h=1,2,\ldots$ define $v := ( 1, 0, \ldots, 0 ) A^h$ and $v_0 := ( 1, 0, \ldots, 0 ) A_0^h$. Then, for every $\varepsilon > 0$,
\begin{equation*}
	\IP( \| v - v_0\| > \varepsilon) \leq
	\IP\Big( \| x - x_0\| > 2^{1-h} \cfrac{\varepsilon}{ \varepsilon + h (\max\{\|x_0\|, 1\})^{h-1}} \Big).
\end{equation*}

\end{lem}

The following lemma ensures that $b$-sub-Gaussian processes satisfy Assumption~\ref{a:MomentCond}.
\begin{lem}\label{lem:MomentCond}
	Let	$(X_{t,T})_{t\in\IZ, T\in\IN^{*}}$ satisfy Assumptions~\ref{a:loc_stat} and~\ref{a:sd_bounded}. Assume that $T \geq \frac{C}{\pi m_f}$ and that the standardized variables $X_{t,T} / \sigma_{t,T}$ are $b$-sub-Gaussian ($b > 0$); i.\,e., 
\[\E\Big( \exp\big(\xi X_{t,T} / \sigma_{t,T} \big) \Big) \leq \exp\Big(\frac{b^2 |\xi|^2}{2}\Big), \quad \xi \in \IR.\]
Then, the process $(X_{t,T})_{t\in\IZ, T\in\IN^{*}}$ satisfies Assumption~\ref{a:MomentCond} with $c := 6 \pi b M_f$ and $d := 1/2$.
\end{lem}

\begin{lem}\label{lem:sumApprox}
Let $f: [0,1] \rightarrow \IR$ be continuous and differentiable on $(0,1)$. Then, for every $A, B = 0, \ldots, T$, $T \in \IN_*$, $A < B$, we have
	\begin{equation*}
	\begin{split}
		\Big| \frac{1}{B-A} \sum_{\ell=A+1}^B f(\ell/T) - \int_0^1 f \Big(\frac{A}{T} + \frac{B-A}{T}u \Big) {\rm d}u \Big| \leq \frac{1}{T} \sup_{A/T < u < B/T} | f'(u) |.
	\end{split}
	\end{equation*}
\end{lem}

\clearpage
\section*{Supplementary Material}

\section{Proofs of Lemmas~\ref{lem:main}--\ref{lem:rel_MSPE}}\label{app:proofAuxLem}

\subsection{Outlook}

In this section we provide the proofs of the results from Section~\ref{app:main}.

The proof of Lemma~\ref{lem:main} is given in Section~\ref{app:proofAuxLem:1}. For the proof we first expand the quantity ${\rm MSPE}_{s,n,N,T}^{(p,h)}$ using simple algebra to then decompose the probability from the assertion. To this end we employ (a) elementary considerations, (b) Lemma~\ref{lem:disentangle} to disentangle the combinations of observations and forecasting coefficients, and (c) Lemmas~\ref{lem:ineq_productOfTwo} and~\ref{lem:norm_a_bar}(i) to treat products of different components of the forecasting coefficients. The decomposition yields two groups of terms. Each group is a sum of probabilites with respect to concentration. The first group, \eqref{eqn:D:a}, concerns generalised sums of observations. The second group, \eqref{eqn:D:b}, concerns the $h$-step ahead forecasting coefficients. We bound the terms in~\eqref{eqn:D:a} by applying Lemma~\ref{lem:expIneqGenSum} and we bound the terms in~\eqref{eqn:D:b} by first employing Lemma~\ref{lem:PowerYW} and treating the resulting term by employing Theorem~\ref{thm:propertiesYW}. Lemma~\ref{lem:nicerboundP}, which provides a more convenient bound for the probability $P^{(p,h)}_{m,N}(\varepsilon)$ in Lemma~\ref{lem:main}, is proved in Section~\ref{app:proofAuxLem:2} by tedious, but elementary considerations.

The proof of Lemma~\ref{lem:rel_MSPE} is given in Section~\ref{app:proofAuxLem:3}. For the proof, we decompose $\overline{{\rm MSPE}}_{s,n,N,T}^{(p,h)}$ into three sums and ${\rm MSPE}_{N/T, n/T}^{(p,h)}(s/T)$ into three integrals. To bound the respective approximation errors we employ bounds on derivatives derived in Section~\ref{app:der_thMSPE}, and bounds on eigenvalues of Toeplitz matrices of moments of empirical autocovariances from Lemma~\ref{lem:BoundEV_M0}.

For the readers convenience, we include Figure~\ref{fig:map_app_E} in which the dependence of the various results is illustrated graphically.
\begin{figure}[b]
	\begin{center}
		\begin{tikzpicture}[
this_sec/.style = {fill=blue!30},
other_sec/.style = {fill=gray!30},
]



\node [other_sec] (E3) at (3,12) {Lemma~\ref{lem:norm_a_bar}};
\node (E3i) at (2.5,11.5) {(i)};
\node (E3ii) at (3.5,11.5) {(ii)};
\node [other_sec] (E4) at (6,12) {Lemma~\ref{lem:norm_g_der}};
\node (E4ii) at (6,11.5) {(ii)};
\node [other_sec] (E1) at (9,12) {Lemma~\ref{lem:rel_v_bar}};
\node (E1ii) at (9,11.5) {(ii)};
\node [other_sec] (E2) at (12,12) {Lemma~\ref{lem:norm_v_der}};
\node (E2i) at (12,11.5) {(ii)};

\node [other_sec] (G4) at (0,8) {Lemma~\ref{lem:ineq_productOfTwo}};
\node [other_sec] (G5) at (3,8) {Lemma~\ref{lem:disentangle}};
\node [other_sec] (G7) at (6,8) {Lemma~\ref{lem:PowerYW}};
\node [other_sec] (G8) at (9,8) {Lemma~\ref{lem:sumApprox}};

\node [other_sec] (61) at (0,10) {Theorem~\ref{thm:propertiesYW}};

\node [this_sec] (A2) at (3,10) {Lemma~\ref{lem:main}};
\node [this_sec] (A3) at (12,10) {Lemma~\ref{lem:rel_MSPE}};

\node [other_sec] (F4) at (9,10) {Lemma~\ref{lem:BoundEV_M0}};
\node (F4ia) at (9,9.5) {(i-a)};
\node (F4ib) at (9,9) {(i-b)};
\node [other_sec] (F6) at (6,10) {Lemma~\ref{lem:expIneqGenSum}};

\path[->] (A2) edge (E3i);
\path[->] (A2) edge (G4);
\path[->] (A2) edge (G5);
\path[->] (A2) edge (G7);
\path[->] (A2) edge (F6);
\path[->] (A2) edge (61);

\path[->] (A3) edge (E3ii);
\path[->] (A3) edge (E4ii);
\path[->] (A3) edge (E1ii);
\path[->] (A3) edge (E2i);
\path[->] (A3) edge (F4ia);
\path[->] (A3) edge (F4ib);
\path[->] (A3) edge (G8);

\end{tikzpicture}
	\end{center}
	\caption{Map of the results proved in Section~\ref{app:proofAuxLem}.}
	\label{fig:map_app_E}
\end{figure}

\subsection{Proof of Lemma~\ref{lem:main}}\label{app:proofAuxLem:1}
We prove the first equation in full detail and comment on how the proof needs to be adapted for the second inequality in the end. First, we consider the case when $p \geq 1$. In the end of the proof we will comment on the (easier) case $p=0$. Denote $\hat v_{N,T}^{(p,h)}(t) =: (\hat v_{1,t}, \ldots, \hat v_{p,t})'$ and note that, omitting the second index of $X_{t,T}$ for the sake of brevity, we have
\begin{equation*}
\begin{split}
		& {\rm MSPE}_{s,n,N,T}^{(p,h)} \\
		& = n^{-1} \sum_{t=s+1}^{s+n} \Big( X_{t+h}^2
		- 2 X_{t+h} \sum_{i=1}^p \hat v_{i,t} X_{t-i+1}
		+ \sum_{i_1, i_2=1}^p \hat v_{i_1,t} \hat v_{i_2,t} X_{t-i_1+1} X_{t-i_2+1} \Big)
\end{split}
\end{equation*}
and an analogous equation for $\overline{{\rm MSPE}}_{s,n,N,T}^{(p,h)}$, with $\bar v_{N,T}^{(p,h)}(t) =: (\bar v_{1,t}, \ldots, \bar v_{p,t})'$.

With this notation we have
\begin{equation*}
\begin{split}
		& \IP \Big( \Big| {\rm MSPE}_{s,n,N,T}^{(p,h)} - \overline{{\rm MSPE}}_{s,n,N,T}^{(p,h)} \Big| > \varepsilon \Big) \\
		& = \IP \Big( \Big| n^{-1} \sum_{t=s+1}^{s+n} \big( X_{t+h}^2 - \E(X_{t+h}^2) \big) \\
		& \quad - 2 \sum_{i=1}^p \Big( n^{-1} \sum_{t=s+1}^{s+n} \big( \hat v_{i,t} X_{t+h} X_{t-i+1}  - \bar v_{i,t} \E( X_{t+h} X_{t-i+1} ) \big) \Big) \\
		& \quad + \sum_{i_1, i_2=1}^p n^{-1} \sum_{t=s+1}^{s+n} \big( \hat v_{i_1,t} \hat v_{i_2,t} X_{t-i_1+1} X_{t-i_2+1} - \bar v_{i_1,t} \bar v_{i_2,t} \E( X_{t-i_1+1} X_{t-i_2+1} ) \big) \Big| > \varepsilon \Big),
\end{split}
\end{equation*}

Therefore, we have
\begin{align}
		& \IP \Big( \Big| {\rm MSPE}_{s,n,N,T}^{(p,h)} - \overline{{\rm MSPE}}_{s,n,N,T}^{(p,h)} \Big| > \varepsilon \Big) \nonumber \\
		& \leq \IP \Big( \Big| n^{-1} \sum_{t=s+1}^{s+n} \big( X_{t+h}^2 - \E(X_{t+h}^2) \big) \Big| > \varepsilon/(p+1)^2 \Big) \label{eqn:A} \\
		& \quad + 2 \sum_{i=1}^p \IP \Big( \Big| n^{-1} \sum_{t=s+1}^{s+n} \big( \hat v_{i,t} X_{t+h} X_{t-i+1}  - \bar v_{i,t} \E( X_{t+h} X_{t-i+1} ) \big) \Big| > \varepsilon/(p+1)^2 \Big) \label{eqn:B} \\
		& \quad + \sum_{i_1, i_2=1}^p \IP \Big( \Big| n^{-1} \sum_{t=s+1}^{s+n} \big( \hat v_{i1,t} \hat v_{i2,t} X_{t-i_1+1} X_{t-i_2+1} \nonumber \\
		& \qquad\qquad\qquad - \bar v_{i1,t} \bar v_{i2,t} \E( X_{t-i_1+1} X_{t-i_2+1} ) \big) \Big| > \varepsilon/(p+1)^2 \Big). \label{eqn:C}
	\end{align}

\pagebreak
We now use Lemma~\ref{lem:disentangle} to disentangle~\eqref{eqn:B} and~\eqref{eqn:C}. More precisely, for the $i$th addend ($i = 1, \ldots, p$) in~\eqref{eqn:B} we have
\begin{align}
		& \IP \Big( \Big| n^{-1} \sum_{t=s+1}^{s+n} \big( \hat v_{i,t} X_{t+h} X_{t-i+1}  - \bar v_{i,t} \E( X_{t+h} X_{t-i+1} ) \big) \Big| > \varepsilon/(p+1)^2 \Big) \nonumber \\
		& \leq n \sup_{t=s+1,\ldots,s+n} \IP \Big( | \hat v_{i,t} - \bar v_{i,t} |
			> \frac{\varepsilon/(p+1)^2}{2 \big( (6 \pi M_f c^2 24^d)^2 + \varepsilon^2/(p+1)^4 \big)^{1/4}} \Big) \label{eqn:B:a} \\
		& \qquad + \IP \Big( \Big| \sum_{t=s+1}^{s+n} ( X_{t+h}^2 X_{t-i+1}^2 - \E (X_{t+h}^2 X_{t-i+1}^2)) \Big| > \frac{n \varepsilon}{2 (p+1)^2} \Big) \label{eqn:B:b} \\
		& \qquad + \IP \Big( \Big| \sum_{t=s+1}^{s+n} (\bar v_{i,t}) (X_{t+h} X_{t-i+1} - \E (X_{t+h} X_{t-i+1})) \Big| > \frac{n \varepsilon}{2 (p+1)^2} \Big), \label{eqn:B:c}
	\end{align}
where we have used the Cauchy-Schwarz inequality, Assumption~\ref{a:MomentCond} and~\eqref{eqn:cons:sd_bounded}, which holds as $T \geq \frac{C}{\pi m_f}$, to (uniformly) bound the second moments
\begin{equation*}
	\begin{split}
		& \E (X_{t+h}^2 X_{t-i+1}^2) \leq \big( \E (X_{t+h}^4) \E (X_{t-i+1}^4) \big)^{1/2} 
		\leq \big( c^4 (4!)^{2d} \sigma_{t+h,T}^2 \sigma_{t-i+1,T}^2) \big)^{1/2}
		\leq 3 \pi M_f c^2 24^d
	\end{split}
\end{equation*}
Note that we have used $m_2^2 := 3 \pi M_f c^2 24^d$ in the application of Lemma~\ref{lem:disentangle}. For the $(i_1, i_2)$th addend ($i_1, i_2 = 1, \ldots, p$) in~\eqref{eqn:C} we analogously have
\begin{align}
		& \IP \Big( \Big| n^{-1} \sum_{t=s+1}^{s+n} \big( \hat v_{i1,t} \hat v_{i2,t} X_{t-i_1+1} X_{t-i_2+1}  - \bar v_{i1,t} \bar v_{i2,t} \E( X_{t-i_1+1} X_{t-i_2+1} ) \big) \Big| > \varepsilon/(p+1)^2 \Big) \nonumber \\
		& \leq n \sup_{t=s+1,\ldots,s+n} \IP \Big( | \hat v_{i1,t} \hat v_{i2,t} - \bar v_{i1,t} \bar v_{i2,t} | > \frac{\varepsilon/(p+1)^2}{2 \big( (6 \pi M_f c^2 24^d)^2 + \varepsilon^2/(p+1)^4 \big)^{1/4}}  \Big) \label{eqn:C:a}\\
		& \qquad + \IP \Big( \Big| \sum_{t=s+1}^{s+n} ( X_{t-i_1+1}^2 X_{t-i_2+1}^2 - \E ( X_{t-i_1+1}^2 X_{t-i_2+1}^2 ) \Big| > \frac{n \varepsilon}{2 (p+1)^2} \Big) \label{eqn:C:b} \\
		& \qquad + \IP \Big( \Big| \sum_{t=s+1}^{s+n} \big(\bar v_{i1,t} \bar v_{i2,t} \big) (X_{t-i_1+1} X_{t-i_2+1} - \E ( X_{t-i_1+1} X_{t-i_2+1} )) \Big| > \frac{n \varepsilon}{2 (p+1)^2} \Big) \label{eqn:C:c}
	\end{align}
where we bounded $\sup_{t=s,\ldots,s+n-1} \E (X_{t-i_1+1}^2 X_{t-i_2+1}^2)$ by $m_2^2$ as before.

We now use Lemma~\ref{lem:ineq_productOfTwo} to bound~\eqref{eqn:C:a} by a sum of two quantities resembling those from~\eqref{eqn:B:a}. Applying Lemmas~\ref{lem:ineq_productOfTwo} and~\ref{lem:norm_a_bar}(i), where for the second we have required that $T \geq 10 C_1 p^2$ and $N \geq 6 C_0 p^{2}$, yields that~\eqref{eqn:C:a} can be bounded by
\begin{equation}\label{eqn:C:a2}
	2 n \sup_{i \in \{i_1, i_2 \}}\sup_{t=s+1,\ldots,s+n} \IP \Big( | \hat v_{i,t} - \bar v_{i,t} | > \mu \Big).
\end{equation}

In conclusion, we have shown
\begin{align}
		& \IP \Big( \Big| {\rm MSPE}_{s,n,N,T}^{(p,h)} - \overline{{\rm MSPE}}_{s,n,N,T}^{(p,h)} \Big| > \varepsilon \Big) \nonumber \\
		& \leq \eqref{eqn:A} 
		+ 2 p \sup_{i=1,\ldots,p} \Big( \eqref{eqn:B:b}_i + \eqref{eqn:B:c}_i) \Big) 
		+ p^2 \sup_{i_1, i_2 =1,\ldots,p} \Big( \eqref{eqn:C:b}_{i_1, i_2} + \eqref{eqn:C:c}_{i_1, i_2} \Big) \label{eqn:D:a} \\
		& \quad + 2 p \sup_{i=1,\ldots,p} \eqref{eqn:B:a}_i + p^2 \sup_{i_1, i_2 =1,\ldots,p} \eqref{eqn:C:a2}_{i_1, i_2}, \label{eqn:D:b}
\end{align}
where the subscript at the equation numbers indicate that the respective expressions depend on $i$ or $i_1, i_2$. We now discuss how to bound~\eqref{eqn:D:a} and~\eqref{eqn:D:b}.

\begin{table}
	\centering
		\begin{tabular}{cccccc}
											& $b$	                & $a_t$                               & $\alpha$  & $h$		      & $\varepsilon$                              \\ \hline
			\eqref{eqn:A}		& $s+h+1$               & $1$                                 & $2$       & $0$         & $\varepsilon/(p+1)^2$                      \\
			\eqref{eqn:B:b} & $s-i+2$             & $1$                                 & $2$       & $h+i-1$     & $\frac{1}{2}\varepsilon/(p+1)^2$           \\
			\eqref{eqn:B:c} & $s-i+2$             & $\bar v_{i,t}$              & $1$       & $h+i-1$     & $\frac{1}{2} \varepsilon/(p+1)^2$ \\
			\eqref{eqn:C:b} & $s-\max\{i_1,i_2\}+2$ & $1$                                 & $2$       & $|i_1-i_2|$ & $\frac{1}{2}\varepsilon/(p+1)^2$         \\
			\eqref{eqn:C:c} & $s-\max\{i_1,i_2\}+2$ & $\bar v_{i_1,t} \bar v_{i_2,t}$ & $1$       & $|i_1-i_2|$ & $\frac{1}{2} \varepsilon/(p+1)^2$ \\
		\end{tabular}
		\caption{Parameters for the application of Lemma~\ref{lem:expIneqGenSum}.}
		\label{tab:par}
\end{table}

To bound~\eqref{eqn:D:a} we note that~\eqref{eqn:A}, \eqref{eqn:B:b}, \eqref{eqn:B:c}, \eqref{eqn:C:b}, and \eqref{eqn:C:c} are all of the type to which Lemma~\ref{lem:expIneqGenSum} can be applied. We do so with the quantities $b$, $a_t$, $\alpha$, $h$, $\varepsilon$ of Lemma~\ref{lem:expIneqGenSum} chosen as in Table~\ref{tab:par}.

Note that for the discussion of \eqref{eqn:B:c} and \eqref{eqn:C:c}, we have, by Lemma~\ref{lem:norm_a_bar}(i) that
\[|a_{t,\eqref{eqn:C:c}}| \leq \big( 2 C_0 + 1 \big)^{h} =: A_{\eqref{eqn:B:c}},
\text{ and }
|a_{t,\eqref{eqn:C:c}}| \leq \big( 2 C_0 + 1 \big)^{2h} =: A_{\eqref{eqn:C:c}},\]
for \eqref{eqn:B:c} and \eqref{eqn:C:c}, respectively. (We wrote $a_{t,(\#)}$ for the sequence of equation $(\#)$.)

Note that (i) the bound from Lemma~\ref{lem:expIneqGenSum} is independent of $b$, that (ii) increasing $A$, $\alpha$ or $h_*$ will increase the bound and (iii) decreasing $\varepsilon$ will also increase the bound.

Thus, \eqref{eqn:D:a} can be bounded by
\begin{equation*}
	\begin{split}
		& (1+4p + 2p^2) \\
		& \cdot \exp\Bigg( - \frac{\frac{\varepsilon^2}{(p+1)^4} }{8\Big(\big(2 C_0 + 1\big)^{4 h} \frac{C_{1,2} (h+p-1)}{n} + (\frac{\varepsilon}{2 (p+1)^2})^{(3+8d)/(2+4d)} \Big( \big(2 C_0 + 1\big)^{2 h} \frac{C_{2,2} (h+p-1)}{n} \Big)^{1/(2+4d)}\Big)} \Bigg).
	\end{split}
\end{equation*}

It remains to bound~\eqref{eqn:D:b}. Note that
\begin{equation}
\eqref{eqn:D:b} \leq (2 p n + 2 p^2 n) \sup_{i=1,\ldots,p} \sup_{t=s+1,\ldots,s+n} \IP \Big( | \hat v_{i,t} - \bar v_{i,t} | > \mu \Big),\label{eqn:D:b2}
\end{equation}
where we have used the fact that $\mu \leq \bar \varepsilon$, because of $\big( 2 C_0 + 1 \big)^{2h} \geq 1$.

Thus, to bound~\eqref{eqn:D:b} we can employ Lemma~\ref{lem:PowerYW}, which yields that
\begin{equation}\label{eqn:bound_vit}
	\begin{split}
		& \IP \Big( | \hat v_{i,t} - \bar v_{i,t} | > \mu \Big)
		\leq \IP \Big( \| \hat v_{N,T}^{(p,h)}(t) - \bar v_{N,T}^{(p,h)}(t) | > \mu \Big)
		\\
		& \leq
			\IP\Big( \| \hat a_{N,T}^{(p)}(t) - \bar a_{N,T}^{(p)}(t)\|
			> 2^{1-h} \frac{\mu}{\mu + h (2 C_0)^{h-1}} \Big)
	\end{split}
\end{equation}
where we have used that $\max\{\| \bar a_{N,T}^{(p)}(t) \|, 1\} \leq \max\{2 C_0, 1\} = 2 C_0$ by Lemma~\ref{lem:norm_a_bar}(i) and the fact that $C_0 \geq (2\pi)^{1/2}$.

Thus, employing~\eqref{eqn:D:b2}, \eqref{eqn:bound_vit} and Theorem~\ref{thm:propertiesYW} (note: $T \geq 10 C_1 p^2 \geq 2 C_1 p^2$), we have that~\eqref{eqn:D:b} is bounded by
\begin{align}
		& 2 n p (p+1) \IP\Big( \| \hat a_{i,N,T}^{(p)}(t) - \bar a_{i,N,T}^{(p)}(t)\|
			> 2^{1-h} \frac{\mu}{\mu + h (2 C_0)^{h-1}} \Big) \nonumber \\
		& \leq 6 n p^2 (p+1) 
		\exp\Bigg( - \frac{\eta^2}{2\Big(C_{1,2} \frac{p}{N-p} + \eta^{(3+4d)/(2+2d)} \Big(C_{2,2} \frac{p}{N-p} \Big)^{1/(2+2d)}\Big)} \Bigg). \label{eqn:D:b3}
	\end{align}
	This completes the proof of the first inequality of Lemma~\ref{lem:main}.
	For the second inequality denote $\hat v_{t,T}^{(p,h)}(t) =: (\hat v_{1,t}, \ldots, \hat v_{p,t})'$ and apply the same stream of arguments, but set $N = t$ and instead of~\eqref{eqn:D:b3} use
	\begin{align*}
		& 2 n p (p+1) \IP\Big( \| \hat a_{i,t,T}^{(p)}(t) - \bar a_{i,t,T}^{(p)}(t)\|
			> 2^{1-h} \frac{\mu}{\mu + h (2 C_0)^{h-1}} \Big) \nonumber \\
		& \leq 6 n p^2 (p+1) 
		\exp\Bigg( - \frac{\eta^2}{2\Big(C_{1,2} \frac{p}{t-p} + \eta^{(3+4d)/(2+2d)} \Big(C_{2,2} \frac{p}{t-p} \Big)^{1/(2+2d)}\Big)} \Bigg) \\
		& \leq 6 n p^2 (p+1) 
		\exp\Bigg( - \frac{\eta^2}{2\Big(C_{1,2} \frac{p}{s-p} + \eta^{(3+4d)/(2+2d)} \Big(C_{2,2} \frac{p}{s-p} \Big)^{1/(2+2d)}\Big)} \Bigg),
	\end{align*}
	as $s \leq t$. The second inequality follows.
	
	Lastly, for the case where $p=0$, note that it suffices to bound~\eqref{eqn:A}.

\subsection{Proof of Lemma~\ref{lem:nicerboundP}}\label{app:proofAuxLem:2}
We begin by proving that the more general inequality
\begin{equation}\label{eqn:nicer_bound}
\begin{split}
	P^{(p,h)}_{m,N}(\varepsilon) & \leq D_1 \Bigg[ p^2 \exp\Bigg(-D_2 \Big( \varepsilon\frac{m}{K_0^h p^2 (h+p)} \Big)^{1/(2+4d)} \Bigg) I\{ \varepsilon > K_0^h \Big(\frac{h+p}{m} \Big)^{\frac{1+4d}{3+8d}} p^2 \} \\
	& \qquad\qquad + p^2 \exp\Bigg(-D_2 \Big(\varepsilon \frac{m^{1/2}}{K_0^h p^2 (h+p)^{1/2}}\Big)^2 \Bigg) I\{ \varepsilon \leq K_0^h \Big(\frac{h+p}{m} \Big)^{\frac{1+4d}{3+8d}} p^2 \} \\
	& \qquad\qquad + m p^3 \exp\Bigg(-D_3 \Big( \varepsilon\frac{N-p}{K_0^h p^4 h} \Big)^{1/(2+2d)} \Bigg) I\{\varepsilon > \Big( \frac{p}{N-p} \Big)^{\frac{1+2d}{3+4d}} p^3 h K_0^h\} \\
	& \qquad\qquad + m p^3 \exp\Bigg(-D_3 \Big(\varepsilon \frac{(N-p)^{1/2}}{p^{7/2} h K_0^h} \Big)^2 \Bigg) I\{\varepsilon \leq \Big( \frac{p}{N-p} \Big)^{\frac{1+2d}{3+4d}} p^3 h K_0^h \} 	\Bigg]
\end{split}
\end{equation}
holds with the constants defined in~\eqref{def:D}:
\begin{equation*}
\begin{split}
	D_1 & := 12,\\
	D_2 & := \big( 2^8 \max\{C_{1,2}, C_{2,2}^{1/(2+4d)}\}\big)^{-1}, \text{ and}\\
	D_3 & := K_1^2 / \big( 2^{12} \max\{C_{1,1}, (K_1^{3+4d} C_{2,1})^{1/(2+2d)}\} \big),
\end{split}
\end{equation*}
where $C_{1,1}, C_{1,2}, C_{2,1}, C_{2,2}$ are defined in~\eqref{lem:exp_ineq_gamma_k2:def:C1C2}.
In the end of this section we will prove that this implies the assertion of Lemma~\ref{lem:nicerboundP}. Now, for the proof of~\eqref{eqn:nicer_bound} we will derive bounds for $\eta$ in terms of $\varepsilon$.

Let $K_0 := 4 C_0 (2 C_0+1)$ and $K_1 := m_f / (32 \min\{(6 \pi M_f c^2 24^d)^{1/2}, 1\})$. Then, we have
\begin{equation}\label{bnd1_eta}
	\eta \leq K_1 \frac{1}{p^3 h K_0^h} \ \varepsilon.
\end{equation}
Further, if $\varepsilon \leq \min\{ 6 \pi M_f c^2 24^d, 1 \} \cdot (p+1)^2$, then we have
\begin{equation}\label{bnd2_eta}
	\eta \geq \frac{K_1}{32} \frac{1}{p^3 h K_0^h} \ \varepsilon.
\end{equation}

For the proof of~\eqref{bnd1_eta} and~\eqref{bnd2_eta}, denote $C_2 := 1/(2 \min\{(6\pi M_f c^2 24^d)^{1/2}, 1\})$. Then, we have the following upper bounds:
\begin{equation}
\label{bnd1}
\eta \leq \frac{m_f}{4 p} \frac{\bar\mu}{8 C_0}, \quad
\bar\mu \leq \frac{\mu}{h (4 C_0)^{h-1}}, \quad
\mu \leq \frac{\bar\varepsilon}{2 (2 C_0 + 1)^h},
\end{equation}
and
\begin{equation}
\label{bnd2}
\bar\varepsilon \leq C_2 \min\Big\{\frac{\varepsilon}{(p+1)^2}, \frac{\varepsilon^{1/2}}{p+1} \Big\} \leq C_2 \frac{\varepsilon}{(p+1)^2} \leq C_2 \frac{\varepsilon}{p^2}
\end{equation}

We therefore obtain~\eqref{bnd1_eta}, because
\begin{equation*}
\begin{split}
	\eta & \leq \frac{m_f}{4 p} \frac{1}{8 C_0} \frac{1}{h (4 C_0)^{h-1}} \frac{1}{2 (2 C_0 + 1)^h} C_2 \frac{\varepsilon}{p^2} \leq \frac{m_f}{16} C_2 \frac{1}{h (4 C_0)^h (2 C_0 + 1)^h} \frac{\varepsilon}{p^3} \\
\end{split}
\end{equation*}

For~\eqref{bnd2_eta} note that from the assumption that $\varepsilon \leq \min\{ 6 \pi M_f c^2 24^d, 1 \} \cdot (p+1)^2$, we have
\begin{align}
	\varepsilon & \leq (6 \pi M_f c^2 24^d) (p+1)^2 \label{eqn:eps1} \\
	\varepsilon & \leq \min\{ (6 \pi M_f c^2 24^d)^{1/2}, 1 \} (p+1)^2 = \frac{1}{2 C_2} (p+1)^2 \label{eqn:eps2}
\end{align}
where~\eqref{eqn:eps1} follows as $\min\{a, b\} \leq a$ and~\eqref{eqn:eps2} follows since $\min\{K, 1\} \leq \min\{K^{1/2}, 1\}$.

Now note that from~\eqref{eqn:eps1}, we have
\begin{equation*}
\begin{split}
	\bar\varepsilon & = \frac{\varepsilon/(p+1)^2}{2 \big( (6 \pi M_f c^2 24^d)^2 + \varepsilon^2/(p+1)^4 \big)^{1/4}} \geq \frac{\varepsilon/(p+1)^2}{2 \big( 2 (6 \pi M_f c^2 24^d)^2 \big)^{1/4}} \geq \frac{C_2}{8} \frac{\varepsilon}{p^2}.
\end{split}
\end{equation*}

Further, note that from~\eqref{bnd2}, \eqref{eqn:eps2} and the fact that $2 C_0 + 1 \geq 1$, we have
\begin{equation}\label{bnd3}
\bar\varepsilon \leq C_2 \frac{\varepsilon}{(p+1)^2} \leq \frac{1}{2} \leq (2 C_0 + 1)^{2h}
\end{equation}
which implies that
\begin{equation*}
\begin{split}
	\mu = \frac{\bar \varepsilon}{2 \Big( \big( 2 C_0 + 1 \Big)^{2h} + \bar \varepsilon \Big)^{1/2}} \geq \frac{\bar \varepsilon}{2 \Big( 2 \big( 2 C_0 + 1 \Big)^{2h}\Big)^{1/2}} \geq \frac{\bar \varepsilon}{4 \big( 2 C_0 + 1 \big)^{h}}
	\geq \frac{C_2}{32} \frac{\varepsilon}{p^2 \big( 2 C_0 + 1 \big)^{h}} ,
\end{split}
\end{equation*}

Further, from~\eqref{bnd1}, \eqref{bnd3}, $(2C_0+1)^h \geq 1$ and $h (2C_0)^{h-1} \geq 1$, we have
\begin{equation}
\label{bnd4}
\mu \leq \frac{\bar\varepsilon}{2 (2 C_0 + 1)^h} \leq \frac{1}{4 (2 C_0 + 1)^h} \leq 1 \leq h (2 C_0)^{h-1}
\end{equation}
which implies that
\begin{equation*}
\begin{split}
	\bar\mu & = 2^{1-h} \frac{\mu}{\mu + h (2 C_0)^{h-1}}
		\geq \frac{\mu}{2 h (4 C_0)^{h-1}}
		\geq \frac{C_2}{64} \frac{1}{h p^2}
		 \frac{\varepsilon}{ (4 C_0)^{h-1} \big( 2 C_0 + 1 \big)^{h}}
\end{split}
\end{equation*}

Finally, from~\eqref{bnd1}, \eqref{bnd4}, $h (4 C_0)^{h-1} \geq 1$ and $8 C_0 \geq 1$ we have
\[\bar\mu \leq \frac{\mu}{h (4 C_0)^{h-1}} \leq \frac{1}{h (4 C_0)^{h-1}} \leq 8 C_0,\]
therefore
\begin{equation*}
\begin{split}
	\eta & = \frac{m_f}{4 p} \min\Big\{1, \bar\mu / (8 C_0) \Big\}
		= \frac{m_f}{4 p} \frac{\bar\mu}{8 C_0}
		\geq \frac{m_f}{4 p} \frac{1}{8 C_0}
		\frac{C_2}{64} \frac{1}{h p^2}
		 \frac{\varepsilon}{ (4 C_0)^{h-1} \big( 2 C_0 + 1 \big)^{h}} \\
		& = \frac{C_2 m_f}{512} \frac{1}{h p^3}
		 \frac{\varepsilon}{ (4 C_0)^h \big( 2 C_0 + 1 \big)^{h}}.
\end{split}
\end{equation*}
This finished the proof of~\eqref{bnd2_eta}.

Now we are going to bound the two exponentials in the expression for $P^{(p,h)}_{m,N}(\varepsilon)$. For the first one recall that $1/D_2 = 2^8 \max\{C_{1,2}, C_{2,2}^{1/(2+4d)}\}$, such that we have
\begin{equation*}
	\begin{split}
		& \exp\Bigg( - \frac{\frac{\varepsilon^2}{(p+1)^4} }{8\Big(\big(2 C_0 + 1\big)^{4 h} \frac{C_{1,2} (h+p-1)}{m} + (\frac{\varepsilon}{2 (p+1)^2})^{(3+8d)/(2+4d)} \Big( \big(2 C_0 + 1\big)^{2 h} \frac{C_{2,2} (h+p-1)}{m} \Big)^{1/(2+4d)}\Big)} \Bigg) \\
		& \leq \exp\Bigg( - \frac{\frac{\varepsilon^2}{(p+1)^4} }{8\Big(K_0^{2h} \frac{C_{1,2} (h+p)}{m} + (\frac{\varepsilon}{2 (p+1)^2})^{(3+8d)/(2+4d)} \Big( K_0^h \frac{C_{2,2} (h+p)}{m} \Big)^{1/(2+4d)}\Big)} \Bigg) \\
		& \leq \exp\Bigg( - \frac{\frac{\varepsilon^2}{p^4} }{2^7 \max\{C_{1,2}, C_{2,2}^{1/(2+4d)}\} \Big(K_0^{2h} \frac{(h+p)}{m} + (\frac{\varepsilon}{p^2})^{(3+8d)/(2+4d)} \Big( K_0^h \frac{ (h+p)}{m} \Big)^{1/(2+4d)}\Big)} \Bigg) \\
		& \leq \begin{cases}
				\exp\Big( - \varepsilon^2 \frac{m}{K_0^{2h} p^4 (h+p)} D_2 \Big) & \varepsilon \leq K_0^h \Big(\frac{h+p}{m} \Big)^{(1+4d)/(3+8d)} p^2 \\
				\exp\Big( - \varepsilon^{1/(2+4d)} \Big(\frac{m}{K_0^h p^2 (h+p)} \Big)^{1/(2+4d)} D_2 \Big) & \varepsilon > K_0^h \Big(\frac{h+p}{m} \Big)^{(1+4d)/(3+8d)} p^2
			\end{cases} \\
& \leq \begin{cases}
				\exp\Big( - \Big(\varepsilon \frac{m^{1/2}}{K_0^h p^2 (h+p)^{1/2}}\Big)^2 D_2 \Big) & \varepsilon \leq K_0^h \Big(\frac{h+p}{m} \Big)^{(1+4d)/(3+8d)} p^2 \\
				\exp\Big( - \Big(\varepsilon \frac{m}{K_0^h p^2 (h+p)} \Big)^{1/(2+4d)} D_2 \Big) & \varepsilon > K_0^h \Big(\frac{h+p}{m} \Big)^{(1+4d)/(3+8d)} p^2
			\end{cases}
	\end{split}
\end{equation*}

where we have used $(2 C_0+1)^2 \leq (2 C_0+2 C_0) (2 C_0+1) = 4 C_0 (2 C_0+1) = K_0$ and the two cases hold, because
\begin{equation*}
\begin{split}
	& (\frac{\varepsilon}{p^2})^{(3+8d)/(2+4d)} \Big( K_0^h \frac{h+p}{m} \Big)^{1/(2+4d)} \leq K_0^{2h} \frac{h+p}{m} \\
	& \Leftrightarrow (\frac{\varepsilon}{p^2})^{3+8d} K_0^h \frac{h+p}{m} \leq K_0^{2h(2+4d)} \Big(\frac{h+p}{m}\Big)^{2+4d} \\
	& \Leftrightarrow (\frac{\varepsilon}{p^2})^{3+8d} \leq K_0^{h(3+8d)} \Big(\frac{h+p}{m}\Big)^{1+4d} \\
	& \Leftrightarrow \varepsilon \leq K_0^h \Big(\frac{h+p}{m}\Big)^{(1+4d)/(3+8d)} p^2\\
\end{split}
\end{equation*}

For the second one, letting $D_3 := K_1^2 / \big( 2^{12} \max\{C_{1,1}, (K_1^{3+4d} C_{2,1})^{1/(2+2d)}\} \big)$ and by employing~\eqref{bnd1_eta} and~\eqref{bnd2_eta}, we obtain
\begin{equation*}
	\begin{split}
		&	\exp\Bigg( - \frac{\eta^2}{2\Big(C_{1,1} \frac{p}{N-p} + \eta^{(3+4d)/(2+2d)} \Big(C_{2,1} \frac{p}{N-p} \Big)^{1/(2+2d)}\Big)} \Bigg) \\
		&	\leq \exp\Bigg( - K_1^2 \frac{\big(\varepsilon / (p^3 h K_0^h) \big)^2}{2^{11} \max\{C_{1,1}, (K_1^{3+4d} C_{2,1})^{1/(2+2d)}\}\Big( \frac{p}{N-p} + \big( \varepsilon / (p^3 h K_0^h) \big)^{(3+4d)/(2+2d)} \Big( \frac{p}{N-p} \Big)^{1/(2+2d)}\Big)} \Bigg) \\
		&	= \exp\Bigg( - 2 D_3 \frac{\big(\varepsilon / (p^3 h K_0^h) \big)^2}{\Big( \frac{p}{N-p} + \big( \varepsilon / (p^3 h K_0^h) \big)^{(3+4d)/(2+2d)} \Big( \frac{p}{N-p} \Big)^{1/(2+2d)}\Big)} \Bigg) \\
		& \leq \begin{cases}
			\exp\Bigg( - D_3 \frac{\big(\varepsilon / (p^3 h K_0^h) \big)^2}{\Big( \frac{p}{N-p}\Big)} \Bigg) & \varepsilon \leq \Big( \frac{p}{N-p} \Big)^{(1+2d)/(3+4d)} p^3 h K_0^h \\
			\exp\Bigg( - D_3 \frac{\big(\varepsilon / (p^3 h K_0^h) \big)^2}{\Big(\big( \varepsilon / (p^3 h K_0^h) \big)^{(3+4d)/(2+2d)} \Big( \frac{p}{N-p} \Big)^{1/(2+2d)}\Big)} \Bigg) & \Big( \frac{p}{N-p} \Big)^{(1+2d)/(3+4d)} p^3 h K_0^h \\
			& \quad < \varepsilon \leq \min\{ 6 \pi M_f c^2 24^d, 1 \} \cdot (p+1)^2 
		\end{cases} \\
		& \leq \begin{cases}
			\exp\Bigg( - D_3 \Big(\varepsilon \frac{(N-p)^{1/2}}{p^{7/2} h K_0^h} \Big)^2 \Bigg) & \varepsilon \leq \Big( \frac{p}{N-p} \Big)^{(1+2d)/(3+4d)} p^3 h K_0^h \\
			\exp\Bigg( - D_3 \Big(\varepsilon \frac{N-p}{p^4 h K_0^h} \Big)^{1/(2+2d)} \Bigg) & \Big( \frac{p}{N-p} \Big)^{(1+2d)/(3+4d)} p^3 h K_0^h \\
			& \quad < \varepsilon \leq \min\{ 6 \pi M_f c^2 24^d, 1 \} \cdot (p+1)^2 
		\end{cases}
	\end{split}
\end{equation*}

where the two cases in the last bound hold, because
\begin{equation*}
\begin{split}
	& \big( \varepsilon / (p^3 h K_0^h) \big)^{(3+4d)/(2+2d)} \Big( \frac{p}{N-p} \Big)^{1/(2+2d)} \leq \frac{p}{N-p} \\
	& \Leftrightarrow \big( \varepsilon / (p^3 h K_0^h) \big)^{3+4d}  \leq \Big( \frac{p}{N-p} \Big)^{1+2d} \\
	& \Leftrightarrow \varepsilon \leq \Big( \frac{p}{N-p} \Big)^{(1+2d)/(3+4d)} p^3 h K_0^h \\
\end{split}
\end{equation*}

This finishes the proof of~\eqref{eqn:nicer_bound}, from which we will now deduce that if~\eqref{cond:eps} holds, we have
\begin{equation}\label{eqn:nicer_bound2}
\begin{split}
	& P^{(p,h)}_{m,N}(\varepsilon) \\
	& \leq D_1 \Bigg[ p^2 \exp\Bigg(-D_2 \Big( \frac{m}{h+p} \Big)^{1/(3+8d)} \Bigg)
	+ m p^3 \exp\Bigg(-D_3 \Big( \frac{N-p}{p} \Big)^{1/(3+4d)} \Bigg) \Bigg].
\end{split}
\end{equation}
Note that if~\eqref{cond:eps} holds, then only the first and third line of~\eqref{eqn:nicer_bound} are relevant. Also, since
\[\varepsilon > \Big(\frac{h+p}{m}\Big)^{\frac{1+4d}{3+8d}} K_0^{h} p^2,\]
we have
\begin{equation*}
\begin{split}
	& \exp\Bigg(-\Big( \varepsilon\frac{m}{K_0^h p^2 (h+p)} \Big)^{1/(2+4d)} \Bigg)
	\leq \exp\Bigg(-\Big( \Big(\frac{h+p}{m}\Big)^{\frac{1+4d}{3+8d}} \frac{m}{(h+p)} \Big)^{1/(2+4d)} \Bigg) \\
	& = \exp\Bigg(-\Big( \Big(\frac{m}{h+p}\Big)^{\frac{3+8d}{3+8d} - \frac{1+4d}{3+8d}} \Big)^{1/(2+4d)} \Bigg)
	= \exp\Bigg(-\Big(\frac{m}{h+p}\Big)^{\frac{1}{3+8d}} \Bigg).
\end{split}
\end{equation*}

Further, since
\[\varepsilon > \Big(\frac{p}{N-p}\Big)^{\frac{1+2d}{3+4d}} K_0^{h} p^3 h,\]
we have
\begin{equation*}
\begin{split}
	& \exp\Bigg(-\Big( \varepsilon\frac{N-p}{K_0^h p^4 h} \Big)^{1/(2+2d)} \Bigg)
	\leq \exp\Bigg(-\Big( \Big(\frac{p}{N-p}\Big)^{\frac{1+2d}{3+4d}} \frac{N-p}{ p } \Big)^{1/(2+2d)} \Bigg) \\
	& = \exp\Bigg(-\Big( \Big(\frac{N-p}{p}\Big)^{\frac{3+4d}{3+4d}-\frac{1+2d}{3+4d}}  \Big)^{1/(2+2d)} \Bigg)
	= \exp\Bigg(-\Big(\frac{N-p}{p}\Big)^{\frac{1}{3+4d}} \Bigg).
\end{split}
\end{equation*}

This completes the proof of~\eqref{eqn:nicer_bound2} and hence of the assertion of the lemma.

\subsection{Proof of Lemma~\ref{lem:rel_MSPE}}\label{app:proofAuxLem:3}
Note the two representations
\begin{equation*}
	\begin{split}
	\overline{{\rm MSPE}}_{s,n,N,T}^{(p,h)}
	& = \frac{1}{n} \sum_{t=s+1}^{s+n} \E \big( X_{t+h,T} - \sum_{i=1}^p \bar v^{(p,h)}_{i;N,T}(t) X_{t-i+1,T} \big)^2 \\
	& = \frac{1}{n} \sum_{t=s+1}^{s+n}
		\Big ( g^{(1)}_{h,T}(t) + g^{(2)}_{p,h,N,T}(t) + g^{(3)}_{p,h,N,T}(t) \Big) \\
	\end{split}
\end{equation*}
where
\begin{equation*}
	\begin{split}
		g^{(1)}_{h,T}(t) & := \E \big( X_{t+h,T}^2 \big) \\
		g^{(2)}_{p,h,N,T}(t) & := - 2 \sum_{i=1}^p \bar v^{(p,h)}_{i;N,T}(t) \E \big(X_{t-i+1,T} X_{t+h,T} \big) \\
		g^{(3)}_{p,h,N,T}(t) & := \sum_{i_1=1}^p \sum_{i_2=1}^p \bar v^{(p,h)}_{i_1;N,T}(t) \bar v^{(p,h)}_{i_2;N,T}(t) \E \big( X_{t-i_1+1,T} X_{t-i_2+1,T} \big),
			\end{split}
\end{equation*}
and
\begin{equation*}
	\begin{split}
	& {\rm MSPE}_{N/T,n/T}^{(p,h)}(s/T)
	= \int_0^1 g_{N/T}^{(p,h)} \Big( \frac{s}{T} + \frac{n}{T} x \Big) {\rm d}x \\
	& \qquad = \frac{1}{n} \sum_{t=s+1}^{s+n}
		\Big( g^{(1)}_{h}(t/T) + g^{(2)}_{p,h;N/T}(t/T) + g^{(3)}_{p,h;N/T}(t/T) \Big) + R_{p,h,N,T}(t)\\
	\end{split}
\end{equation*}
where $g^{(p,h)}_{\Delta}(u) := g^{(1)}_{h}(u) + g^{(2)}_{p,h;\Delta}(u) + g^{(3)}_{p,h;\Delta}(u)$, and
\begin{equation*}
	\begin{split}
		g^{(1)}_{h}(u) & := \gamma_0(u) \\
		g^{(2)}_{p,h;\Delta}(u) & := - 2 \sum_{i=1}^p v^{(p,h)}_{i;\Delta}(u) \gamma_{i+h-1}(u) \\
			& = -2 \big( v^{(p,h)}_{\Delta}(u) \big)' \gamma_0^{(p,h)}(u) \\
		g^{(3)}_{p,h;\Delta}(u) & := \sum_{i_1=1}^p \sum_{i_2=1}^p v^{(p,h)}_{i_1;\Delta}(u) v^{(p,h)}_{i_2;\Delta}(u) \gamma_{i_2-i_1}(u) \\
			& = \big( v^{(p,h)}_{\Delta}(u) \big)' \Gamma^{(p)}(u) v^{(p,h)}_{\Delta}(u).
			\end{split}
\end{equation*}

Note that we have
\begin{equation}\label{lem:rel_MSPE:prf:A}
	\begin{split}
			& \overline{{\rm MSPE}}_{s,n,N,T}^{(p,h)} - {\rm MSPE}_{N/T, n/T}^{(p,h)}(s/T) \\
			& = R_{p,h,N,T}(t) + \frac{1}{n} \sum_{t=s+1}^{s+n} 
				\Big ( r^{(1)}_{h,T}(t) + r^{(2a)}_{p,h,N,T}(t) + r^{(3a)}_{p,h,N,T}(t) \Big), \\
	\end{split}
\end{equation}
and
\begin{equation*}
	\begin{split}
			& \overline{{\rm MSPE}}_{s,n,0,T}^{(p,h)} - {\rm MSPE}_{s/T, n/T}^{(p,h)}(s/T) \\
			& = R_{p,h,N,T}(t) + \frac{1}{n} \sum_{t=s+1}^{s+n} 
				\Big ( r^{(1)}_{h,T}(t) + r^{(2)}_{p,h,t,T}(t) + r^{(3)}_{p,h,t,T}(t) \Big), \\
	\end{split}
\end{equation*}
with
\begin{equation*}
	\begin{split}
			r^{(1)}_{h,T}(t) & := g^{(1)}_{h,T}(t) - g^{(1)}_{h}(t/T),\\
			r^{(2)}_{p,h,N,T}(t) & := g^{(2)}_{p,h,N,T}(t) - g^{(2)}_{p,h;N/T}(t/T), \\
			r^{(3)}_{p,h,N,T}(t) & := g^{(3)}_{p,h,N,T}(t) - g^{(3)}_{p,h;N/T}(t/T), \\
	\end{split}
\end{equation*}

Because, by Lemma~\ref{lem:sumApprox} and Lemma~\ref{lem:norm_g_der}(ii) we have
\[| R_{p,h,N,T}(t)| \leq \frac{1}{T} \sup_{\frac{s}{T} < u < \frac{s+n}{T}} | \frac{\partial}{\partial u} g_{N/T}^{(p,h)}(u) | \leq \frac{4}{T} \big( 2 h + 1 \big) \big( C_0 \big)^{2 h+1} M'_f \]
it now suffices to prove bounds (uniform with respect to $t$) for
$r^{(1)}_{h,T}(t)$,
$r^{(2)}_{p,h,N,T}(t)$, and
$r^{(3)}_{p,h,N,T}(t)$. For $r^{(1)}_{h,T}(t)$ we have
\begin{equation*}
	\begin{split}
		g^{(1)}_{h,T}(t) & := \E \big( X_{t+h,T}^2 \big) = \tilde\gamma_{0,T}(t+h) \\
			&  = \gamma_0(\frac{t}{T}) + \gamma_0(\frac{t+h}{T}) - \gamma_0(\frac{t}{T}) + \tilde\gamma_{0,T}(t+h) - \gamma_0(\frac{t+h}{T}) \\
			& =: g^{(1)}_{h}(t/T) + r^{(1)}_{h,T}(t),
			\end{split}
\end{equation*}
where $|r^{(1)}_{h,T}(t)| \leq (2 \pi M'_f h + C)/T$, by~\eqref{eqn:sd_der_bounded} and~\eqref{eqn:def:locally_stationary}.

For $r^{(2)}_{p,h,N,T}(t)$ we have
\begin{equation*}
	\begin{split}
		& |r^{(2)}_{p,h,N,T}(t)| 
		:= \Big| - 2 \sum_{i=1}^p \Big(
			\bar v^{(p,h)}_{i;N,T}(t) \E \big(X_{t-i+1,T} X_{t+h,T} \big) \\
			& \qquad\qquad - v^{(p,h)}_{i;N/T}(t/T) \gamma_{i+h-1}(\frac{t}{T}) \Big) \Big|\\
		& = 2 \Big| \sum_{i=1}^p \Big(
			\big( \bar v^{(p,h)}_{i;N,T}(t) - v^{(p,h)}_{i;N/T}(t/T) \big) \gamma_{i+h-1}(\frac{t}{T}) \Big) \\
			& \qquad + \sum_{i=1}^p \Big(
			 \bar v^{(p,h)}_{i;N,T}(t) \Big( \E \big(X_{t-i+1,T} X_{t+h,T} \big) - \gamma_{i+h-1}(\frac{t+h}{T}) \\
			& \qquad\qquad\qquad\qquad\qquad + \gamma_{i+h-1}(\frac{t+h}{T}) -  \gamma_{i+h-1}(\frac{t}{T}) \Big) \Big) \Big| \\
			& \leq 2 \Big(  \| \bar v^{(p,h)}_{N,T}(t) - v^{(p,h)}_{N/T}(t/T) \| \cdot \| \gamma_0^{(p,h)}(\frac{t}{T}) \| + p \| \bar v^{(p,h)}_{N,T}(t) \|_{\infty} (C + 2\pi h M'_f) / T \Big) \\
			& \leq 2 \Big(  3 h \big( 2 C_0 \big)^h C_1 \frac{p^2}{T} \cdot (2\pi)^{1/2} M_f 
			+ 2 h \big( 2 C_0 \big)^h \frac{p^{2}}{N} C_0
			+ p 2 \big( C_0 \big)^h \cdot (C + 2\pi  M'_f) h / T \Big) \\
			& \leq 2 \big( C_0 \big)^h \Big[ \frac{1}{T} \Big(  h 2^h C_1 \Big( 3 p^2 \Big) (2\pi)^{1/2} M_f + 2 h p (C + 2\pi M'_f) \Big) + h 2^{h+1} \frac{p^{2}}{N} C_0 \Big] \\
			& \leq 2 \big( C_0 \big)^h \Big[ \frac{1}{T} \Big(  h 2^h \big( 3 p^2 \big) C_0 (2\pi M'_f + C) + 2 h p^2 C_0 (2\pi M'_f+C) \Big) + h 2^{h+1} \frac{p^{2}}{N} C_0 \Big] \\
			& \leq 2 \big( C_0 \big)^{h+1} \Big[ (  3 h 2^h + 2 h ) \frac{1}{T} p^2 (2\pi M'_f + C)  + h 2^{h+1} \frac{p^{2}}{N} \Big] \\
			& \leq  h \big( 2 C_0 \big)^{h+1} \Big[ 5 (2\pi M'_f + C) \frac{p^2}{T} + \frac{p^{2}}{N} \Big]    \\
	\end{split}
\end{equation*}
where the second inequality holds by employing Lemmas~\ref{lem:rel_v_bar}(ii) and~\ref{lem:norm_a_bar}(ii), both of which can be applied, since $T \geq 6 h 2^{h} C_1 p^2 \geq 10 C_1 p^2$ and $N \geq 4 h 2^h C_0 p^{2} \geq 6 p^{2} C_0$, due to the assumptions and $h \geq 1$.
Also note that $\| \gamma_0^{(p,h)}(\frac{t+h}{T}) \| \leq \| \gamma_0^{(p+h)}(\frac{t+h}{T}) \|$ to which we apply Lemma~\ref{lem:BoundEV_M0}(i-b). We also employed the bound on $\gamma'_{i+h-1}(u)$ from below~\eqref{eqn:gammaPrime} and used the mean value theorem.

For $r^{(3)}_{p,h,N,T}(t)$, employing $a b c - a_0 b_0 c_0 = b_0 c (a-a_0) + a c (b-b_0) + a_0 b_0 (c-c_0)$, we have
\begin{equation*}
	\begin{split}
		& |r^{(3)}_{p,h,N,T}(t)| 
		:= \Big| \sum_{i_1=1}^p \sum_{i_2=1}^p \Big( 
		\bar v^{(p,h)}_{i_2;N,T}(t) \gamma_{i_2-i_1}(\frac{t}{T})
			\Big( v^{(p,h)}_{i_1;N/T}(t/T) - \bar v^{(p,h)}_{i_1;N,T}(t)\Big) \\
		& \qquad\qquad + v^{(p,h)}_{i_1;N/T}(t/T) \gamma_{i_2-i_1}(\frac{t}{T})
			\Big( v^{(p,h)}_{i_2;N/T}(t) - \bar v^{(p,h)}_{i_2;N,T}(t/T) \Big) \\
		& \qquad\qquad + \bar v^{(p,h)}_{i_1;N,T}(t) \bar v^{(p,h)}_{i_2;N,T}(t)
			\Big( \gamma_{i_2-i_1}(\frac{t}{T}) - \E \big( X_{t-i_1+1,T} X_{t-i_2+1,T}\big) \Big) \Big) \Big| \\
		& \leq \big| (\bar v^{(p,h)}_{N,T}(t))' \Gamma^{(p)}_0(t/T) \big( v^{(p,h)}_{N/T}(t/T) - \bar v^{(p,h)}_{N,T}(t) \big) \\
		& \qquad + (v^{(p,h)}_{N/T}(t/T))' \Gamma^{(p)}_0(t/T) \big( v^{(p,h)}_{N/T}(t/T) - \bar v^{(p,h)}_{N,T}(t) \big) \big| \\
		& \qquad + \| \bar v^{(p,h)}_{N,T}(t) \|_1^2 (2\pi M'_f (p-1) + C ) / T \\
		& \leq \big( \|\bar v^{(p,h)}_{N,T}(t)\| + \| v^{(p,h)}_{N/T}(t/T) \| \big) \| \Gamma^{(p)}_0(t/T) \| \cdot \| v^{(p,h)}_{N/T}(t/T) - \bar v^{(p,h)}_{N,T}(t) \| \\
		& \qquad + p \| \bar v^{(p,h)}_{N,T}(t) \|^2 (2\pi M'_f (p-1) + C ) / T
	\end{split}
\end{equation*}
\begin{equation*}
	\begin{split}
		& \leq \big( 2 \big( C_0 \big)^h + \big( C_0 \big)^h \big) M_f \Big( 3 h \big( 2 C_0 \big)^h C_1 \frac{p^2}{T} + 2 h \big( 2 C_0 \big)^h C_0 \frac{p^{2}}{N} \Big) \\
		& \qquad + p \Big( 2 \big( C_0 \big)^h \Big)^2 (2\pi M'_f (p-1) + C ) / T \\
		& = T^{-1} \big( C_0 \big)^{2 h} \Big( 9 h 2^h M_f C_1 p^2
							+ 4 p (2\pi M'_f (p-1) + C ) \Big) + 6 \big( C_0 \big)^{2 h} 2^h h C_0 \frac{p^{2}}{N}\\
		& \leq T^{-1} \big( C_0 \big)^{2 h} \Big( 9 h 2^h M_f C_1 p^2
							+ 4 p^2 (2\pi M'_f + C ) \Big) + 6 \big( C_0 \big)^{2 h} 2^h h C_0 \frac{p^{2}}{N} \\
		& = T^{-1} \big( C_0 \big)^{2 h} p^2 (2\pi M'_f + C )  \Big( 9 (2 \pi)^{-1/2} h 2^h C_0	+ 4 \Big) + 6 \big( C_0 \big)^{2 h} 2^h h C_0 \frac{p^{2}}{N} \\
		& \leq 8 (2\pi M'_f + C) h 2^h \big( C_0 \big)^{2 h+1} T^{-1} p^2 + 6 \big( C_0 \big)^{2 h} 2^h h C_0 \frac{p^{2}}{N} \\
	\end{split}
\end{equation*}
where the last inequality is due to Lemma~\ref{lem:norm_a_bar}(ii), \ref{lem:norm_v_der}(i), \ref{lem:BoundEV_M0}(i-a), \ref{lem:rel_v_bar}(ii), and~\ref{lem:norm_a_bar}(ii) to the respective terms. Note that Lemma~\ref{lem:norm_a_bar}(ii) can be applied as we have assumed that $T \geq 3 h 2^{h} C_1 p^2$ and $N \geq 4 h 2^h C_0 p^{2}$, which implies the condition for Lemma~\ref{lem:rel_v_bar}(ii), as $h \geq 1$.

For the first inequality, we have used the fact that
\begin{equation*}
	\begin{split}
		& | \gamma_{i_2-i_1}(\frac{t}{T}) - \E \big( X_{t-i_1+1,T} X_{t-i_2+1,T}\big) | \\
		& \leq | \gamma_{i_2-i_1}(\frac{t}{T}) - \gamma_{i_2-i_1}(\frac{t-i_1 + 1}{T}) | 
		+ | \gamma_{i_2-i_1}(\frac{t-i_1+1}{T}) - \E \big( X_{t-i_1+1,T} X_{t-i_2+1,T}\big) | \\
	\end{split}
\end{equation*}

For the very last inequality we used $9/(2\pi)^{1/2}+4 \leq 8$.

Substituting these four results into~\eqref{lem:rel_MSPE:prf:A}, we have
\begin{equation*}
	\begin{split}
			& \Big| \overline{{\rm MSPE}}_{s,n,N,T}^{(p,h)} - {\rm MSPE}_{N/T, n/T}^{(p,h)}(s/T) \Big| \\
			& \leq \frac{4}{T} \big( 2 h + 1 \big) \big( C_0 \big)^{2 h+1} M'_f
				+ (2 \pi M'_f h + C)/T \\
			& \qquad
				+ 5 (2\pi M'_f + C) h \big( 2 C_0 \big)^{h+1} \frac{1}{T} p^2 + 8  (2\pi M'_f + C) h 2^h \big( C_0 \big)^{2 h+1} T^{-1} p^2 \\
			& \qquad
				 + h \big( 2 C_0 \big)^{h+1} \frac{p^2}{N} + 6 \big( C_0 \big)^{2 h} 2^h h C_0 \frac{p^{2}}{N} \\
			& \leq 48 (2\pi M'_f + C) h 2^h \big( C_0 \big)^{2 h+1} \frac{p^2}{T}
				+ 8 h 2^h \big( C_0 \big)^{2 h+1} \frac{p^{2}}{N} \\
			& = 8 h 2^h \big( C_0 \big)^{2 h+1} \Big[ 6 (2\pi M'_f + C) \frac{p^2}{T} + \frac{p^{2}}{N} \Big],
	\end{split}
\end{equation*}
which completes this proof.

\clearpage
\section{Proofs ommitted in the Appendix}\label{sec:techDetails}

\subsection[Proof of (43)]{Proof of~\eqref{eqn:main3:forApp}}\label{sec:eq_forApp}
With the notation from the proof of Theorem~\ref{thm:main3}, we have
\begin{equation*}
\begin{split}
	& \IP\Big( (\hat R_{T,2}(h) \geq 1 + \delta \text{ and }
		\hat R_{T,3}(h) \geq 1 + \delta) \\
	& \qquad \text{ or }
		(\hat R_{T,2}(h) < 1 + \delta \text{ and }
		\hat R_{T,3}(h) < 1 + \delta) \Big) \\
	& = 1 - \IP\Big( \bigcup_{p_1, p_2 \in \mathcal{P}} \bigcup_{N \in \mathcal{N}} \{ (\hat R_{T,2}(h) < 1 + \delta \text{ or }
		\hat R_{T,3}(h) < 1 + \delta) \\
	& \qquad \text{ and }
		(\hat R_{T,2}(h) \geq 1 + \delta \text{ or }
		\hat R_{T,3}(h) \geq 1 + \delta), \ \hat p_{\lstat} = p_1, \ \hat N_{\lstat} = N, \ \hat p_{\stat} = p_2 \}\Big) \\
	& \geq 1 - \sum_{p_1, p_2 \in \mathcal{P}} \sum_{N \in \mathcal{N}} \IP\Big( \big( Y_1 < X_1 (1 + \delta) \text{ or } Y_2 < X_2 (1 + \delta) \big) \\
	& \qquad\qquad\qquad \text{ and }
		\big( Y_1 \geq X_1 (1 + \delta) \text{ or } Y_2 \geq X_2 (1 + \delta) \big) \Big) \\
	& = 1 - \sum_{p_1, p_2 \in \mathcal{P}} \sum_{N \in \mathcal{N}} \IP\Big( \big( Y_1 < X_1 (1 + \delta) \text{ or } Y_1 < X_1 (1 + \delta) + Y_1-Y_2 + (X_2 - X_1) (1 + \delta) \big) \\
	& \qquad\qquad\qquad \text{ and }
		\big( Y_1 \geq X_1 (1 + \delta) \text{ or } Y_1 \geq X_1 (1 + \delta) + Y_1-Y_2 + (X_2 - X_1) (1 + \delta) \big) \Big) \\
	& = 1 - \sum_{p_1, p_2 \in \mathcal{P}} \sum_{N \in \mathcal{N}} \IP\Big( \big( A < 0 \text{ or } A < B \big)
	\text{ and }
		\big( A \geq 0 \text{ or } A \geq B \big) \Big) \\
	& = 1 - \sum_{p_1, p_2 \in \mathcal{P}} \sum_{N \in \mathcal{N}}\IP\Big( A < \max\{0, B\}, \ A \geq \min\{0, B\} \Big) \\
	& = 1 - \sum_{p_1, p_2 \in \mathcal{P}} \sum_{N \in \mathcal{N}} \Big( \IP\Big( A < \max\{0, B\}, \ A \geq \min\{0, B\}, \ |B| > \varepsilon \Big) \\
	& \qquad + \IP\Big( A < \max\{0, B\}, \ A \geq \min\{0, B\}, \ -\varepsilon \leq B < 0 \Big) \\
	& \qquad + \IP\Big( A < \max\{0, B\}, \ A \geq \min\{0, B\}, \ 0 \leq B \leq \varepsilon \Big) \Big) \\
	& = 1 - \sum_{p_1, p_2 \in \mathcal{P}} \sum_{N \in \mathcal{N}} \Big( \IP\Big( A < \max\{0, B\}, \ A \geq \min\{0, B\}, \ |B| > \varepsilon \Big) \\
	& \qquad + \IP\Big( A < 0, \ A \geq B, \ -\varepsilon \leq B < 0 \Big) \\
	& \qquad + \IP\Big( A < B, \ A \geq 0, \ 0 \leq B \leq \varepsilon \Big) \Big) \\
	& \geq 1 - \sum_{p_1, p_2 \in \mathcal{P}} \sum_{N \in \mathcal{N}} \Big( \IP\Big( |B| > \varepsilon \Big) + \IP\Big( -\varepsilon \leq A < 0 \Big) + \IP\Big( 0 \leq A < \varepsilon \Big) \Big)\\
	& \geq 1 - \sum_{p_1, p_2 \in \mathcal{P}} \sum_{N \in \mathcal{N}} \Big( \IP\Big( |B| > \varepsilon \Big) + \IP\Big( |A| \leq \varepsilon \Big) \Big).  \\
\end{split}
\end{equation*}
Note that in~\eqref{eqn:main3:forApp} we use $\varepsilon = q(\delta)/2$.\hfill$\square$

\clearpage
\subsection{Proof of Lemmas~\ref{lem:rel_a_1}--\ref{lem:norm_a}}\label{sec:Lemmas_a_proofs}
This section concludes with the proof of the lemmas above.

\textbf{Proofs of Lemma~\ref{lem:norm_a}.}
Note that
\begin{equation*}
	\begin{split}
		\| a_{\Delta}^{(p)}(u) \|
		& = \| \Gamma_{\Delta}^{(p)}(u)^{-1} \gamma_{\Delta}^{(p)}(u) \|
		\leq \| \Gamma_{\Delta}^{(p)}(u)^{-1}\| \| \gamma_{\Delta}^{(p)}(u) \|.
	\end{split}
\end{equation*}
The assertion follows from Lemma~\ref{lem:BoundEV_M0}(i).
\hfill$\square$

\textbf{Outline for the proof of Lemma~\ref{lem:rel_a_1} and Lemma~\ref{lem:rel_a_bar}.}

Note that the norm to be bounded is, in all cases, of a difference of quantities of the form $a_1 := \Gamma_1^{-1} \gamma_1$ and $a_2 := \Gamma_2^{-1} \gamma_2$, where the norms of the components of $a_2$ and norms of the differences of the components can be bounded using results from Section~\ref{sec:Cov:mom}. We denote these bounds by $\|\Gamma_2^{-1}\| \leq K$, $\|\gamma_2\| \leq \kappa$, $\|\Gamma_1 - \Gamma_2\| \leq D$, and $\|\gamma_1 - \gamma_2\| \leq \delta$.

For the bound of the norm of interest note that
\begin{equation*}
	\begin{split}
		a_1 - a_2 & = \Gamma_1^{-1} \gamma_1 - \Gamma_2^{-1} \gamma_2
		= (\Gamma_1^{-1} - \Gamma_2^{-1}) \gamma_1 + \Gamma_2^{-1} (\gamma_1 - \gamma_2) 
	\end{split}
\end{equation*}
Thus, we have
		$ \|a_1 - a_2\| 
		\leq \|\Gamma_1^{-1} - \Gamma_2^{-1}\| \cdot \|\gamma_1\| + \|\Gamma_2^{-1}\| \cdot \|\gamma_1 - \gamma_2\|$.
For each of the inequalities to be proven, we will provide condition~1 (bounds for $T$ and $N$) such that
\begin{equation}\label{lem:rel_a:prf:cond1}
\| \Gamma_2 - \Gamma_1 \| \cdot \|\Gamma_2^{-1}\| \leq D \cdot  K \leq 1/2.
\end{equation}
Under condition~1, we have, by Lemma~\ref{lem:stoerungslemma2}, that
\[ \|\Gamma_1^{-1} - \Gamma_2^{-1}\| = \|(\Gamma_2 + \Gamma_1 - \Gamma_2) ^{-1} - \Gamma_2^{-1}\|
\leq 2 \| \Gamma_1 - \Gamma_2 \| \cdot \| \Gamma_2^{-1} \|^2, \]
which yields
\begin{equation}\label{lem:rel_a:prf:bound:pre}
\|a_1 - a_2\| \leq 
		 2 \| \Gamma_1 - \Gamma_2 \| \cdot \| \Gamma_2^{-1} \|^2 \cdot \|\gamma_1\| + \|\Gamma_2^{-1}\| \cdot \|\gamma_1 - \gamma_2\|.
\end{equation}

Finally, we will provide condition~2 (bounds for $T$ and $N$) such that
\begin{equation}\label{lem:rel_a:prf:cond2}
	\kappa - \delta \geq 0.
\end{equation}
Under condition~2, we have
$ \|\gamma_1\| \leq \kappa + \delta \leq 2 \kappa$, by the triangle inequality.
It is implied that, if condition~1 and~2 hold, the norm $\|a_1 - a_2\|$ does not exceed
\begin{equation}\label{lem:rel_a:prf:bound}
4 \ K^2 \ \kappa \ D + K \ \delta.
\end{equation}
For the proof of the individual cases it thus suffices to refer to Section~\ref{sec:Cov:mom} where appropriate bounds (i.\,e., $K$, $\kappa$, $D$, and $\delta$) can be found, from which the combined bounds for $T$ and $N$, and~\eqref{lem:rel_a:prf:bound} are obtained.

\textbf{Proof of Lemma~\ref{lem:rel_a_1}.}

We proceed as outlined and choose $\tilde a_{T}^{(p)}(t)$ to be $a_1$ and $a_0^{(p)}(t/T)$ as $a_2$. Lemma~\ref{lem:BoundEV_M0}(i) and Corollary~\ref{kor:exp_gamma_k}(i) provide the needed bounds. We have
\[
\| \tilde \Gamma_{T}^{(p)}(t) - \Gamma_0^{(p)}(t/T) \| \cdot \| \Gamma_0^{(p)}(u)^{-1} \| \leq T^{-1} p^2 (2\pi M'_f + C) m_f^{-1}
\]
Thus, \eqref{lem:rel_a:prf:cond1} holds if $T \geq T_{0,1} := 2 p^2 (2\pi M'_f + C) m_f^{-1}$. Further, we have
\[\| \tilde \gamma_{T}^{(p)}(t) - \gamma_0^{(p)}(t/T) \| \leq T^{-1} p^{1/2} C, \quad \| \gamma_0^{(p)}(t/T) \| \leq (2\pi)^{1/2} M_f,\]
such that~\eqref{lem:rel_a:prf:cond2} holds if $T \geq T_{0,2} := \frac{C}{M_f} (\frac{p}{2\pi})^{1/2}$.
Consequently, by \eqref{lem:rel_a:prf:bound}, we have the following bound
\begin{equation*}
	\begin{split}
		\| \tilde a_{T}^{(p)}(t) - a_0^{(p)}(t/T) \|
		& \leq 4 T^{-1} p^2 (2\pi M'_f + C) m_f^{-2} (2 \pi)^{1/2} M_f
						+ m_f^{-1} T^{-1} p^{1/2} C \\
		& = \frac{1}{T}  \Big( 4 C_0 C_1 \, p^2	+ p^{1/2} C/m_f \Big),
	\end{split}
\end{equation*}
which holds, if $T \geq  \max\{2 p^2 C_1 , p^{1/2} C_0^{-1} C/m_f \}$.
The assertion follows, as we have $p^{1/2} C / m_f \leq C_1 p^2$ and $C_0^{-1} < 1$.\hfill$\square$

\textbf{Proof of Lemma~\ref{lem:rel_a_bar}.}
We proceed as outlined. For (i) we choose $\bar a_{N,T}^{(p)}(t)$ to be $a_1$ and $a_0^{(p)}(t/T)$ as $a_2$. For the needed bounds note that we have
\begin{equation*}
\begin{split}
	& \| \E \hat \Gamma_{N,T}^{(p)}(t) - \Gamma_0^{(p)}(t/T) \| \cdot \| \Gamma_0^{(p)}(u)^{-1} \| \\
	& \leq \big( \| \E \hat \Gamma_{N,T}^{(p)}(t) - F_{p,N} \circ \Gamma_0^{(p)}(t/T) \| + \| F_{p,N} \circ \Gamma_0^{(p)}(t/T) - \Gamma_0^{(p)}(t/T) \| \big) \cdot \| \Gamma_0^{(p)}(u)^{-1} \| \\
	& \leq \Big( \frac{p}{T} (2\pi M'_f (N+1)+C) + \frac{p^{2}}{N} M_f \Big) m_f^{-1}.
	\end{split}
\end{equation*}
The first inequality follows from the triangle inequality. In the second inequality we have employed Lemma~\ref{lem:BoundEV_M0}(i) and Corollary~\ref{kor:exp_gamma_k}(ii), and the fact that
\begin{equation*}
\begin{split}
	& \| F_{p,N} \circ \Gamma_0^{(p)}(t/T) - \Gamma_0^{(p)}(t/T) \|
	= \| \big( |i - j| / N \big)_{i,j=1,\ldots,p} \circ \Gamma_0^{(p)}(t/T) \| \\
	& \leq N^{-1} \| \big( |i - j| \big)_{i,j=1,\ldots,p}\| \cdot \| \Gamma_0^{(p)}(t/T) \|
	\leq N^{-1} \| \big( |i - j| \big)_{i,j=1,\ldots,p}\|_F \, M_f
	\leq \frac{p^{2}}{N} M_f,
	\end{split}
\end{equation*}	
where the first inequality is a consequence of the submultiplicativity of the spectral norm of the Hadamard product; cf. (1.1) in \cite{JohnsonNylen1990}, and $\sum_{i,k=1}^p |i-j|^2 = p^2 (p^2-1) / 6$.

Thus, \eqref{lem:rel_a:prf:cond1} holds if $T \geq 4 p (2\pi (N+1) M'_f + C) m_f^{-1}$ and $N \geq 4 p^{2} M_f / m_f$. 
Further, note that we have $\| \gamma_0^{(p)}(t/T) \| \leq (2\pi)^{1/2} M_f$ and 
\begin{equation*}
\begin{split}
	& \| \E \hat \gamma^{(p)}_{N,T}(t) - \gamma_0^{(p)}(t/T) \|
	\leq \| \E \hat \gamma^{(p)}_{N,T}(t) - f_{p,N} \circ \gamma_0^{(p)}(t/T) \|
		+ \| f_{p,N} \circ \gamma_0^{(p)}(t/T) - \gamma_0^{(p)}(t/T)\| \\
	& \leq \frac{p^{1/2}}{T} (2\pi M'_f N + C) + \frac{p}{N} (2\pi)^{1/2} M_f,
	\end{split}
\end{equation*}
where we have used the fact that
\begin{equation*}
	\| f_{p,N} \circ \gamma_0^{(p)}(t/T) - \gamma_0^{(p)}(t/T)\|
	= \Big( \sum_{k=1}^p \big(\frac{k}{N} \gamma_k(t/T)\big)^2 \Big)^{1/2} \leq \frac{p}{N} \|\gamma_0^{(p)}(t/T)\| \leq \frac{p}{N} (2\pi)^{1/2} M_f.
\end{equation*}
In consequence, we see that~\eqref{lem:rel_a:prf:cond2} holds if $T \geq 2 \frac{2\pi M'_f N + C}{M_f} (\frac{p}{2\pi})^{1/2}$ and $N \geq 2 p$.
Consequently, by \eqref{lem:rel_a:prf:bound}, we have the following bound

\begin{equation*}
	\begin{split}
		\| \bar a_{N,T}^{(p)}(t) - a^{(p)}(t/T) \|
		& \leq 4 m_f^{-2} (2\pi)^{1/2} M_f \Big( \frac{p}{T} (2\pi M'_f (N+1)+C) + \frac{p^{2}}{N} M_f \Big) \\
		& \quad + m_f^{-1} \Big(\frac{p^{1/2}}{T} (2\pi M'_f N + C) + \frac{p}{N} (2\pi)^{1/2} M_f\Big) \\
		& = 4 T^{-1} p (2 \pi M'_f (N+1) + C) m_f^{-2} (2 \pi)^{1/2} M_f \\
		& \quad + m_f^{-1} T^{-1} p^{1/2} (2 \pi M'_f N + C) \\
		& \quad + 4 m_f^{-2} (2\pi)^{1/2} M_f \frac{p^{2}}{N} M_f + \frac{p}{N} (2\pi)^{1/2} M_f m_f^{-1}\\
		& \leq 4 T^{-1} p 2 N (2 \pi M'_f + C) m_f^{-2} (2 \pi)^{1/2} M_f \\
		& \quad + m_f^{-1} T^{-1} p^{1/2} N (2 \pi M'_f + C)
					+ 2 C_0^2 \frac{p^{2}}{N} + C_0 \frac{p}{N} \\
		& = C_1 \frac{N}{T} \Big( 8 C_0 p + p^{1/2} \Big) + 2 C_0^2 \frac{p^{2}}{N} + C_0 \frac{p}{N} \\
	\end{split}
\end{equation*}
for $T \geq 2 \max\{2 p (2 \pi M'_f (N+1) + C) m_f^{-1}, \frac{2\pi M'_f N + C}{M_f} (\frac{p}{2\pi})^{1/2}\}$ and $N \geq 4 p^{2} M_f / m_f$. We choose to state the result to hold for $T \geq 2 N C_1 \max\{4 p, p^{1/2} C_0^{-1}\}$, which is more restrictive, but allows for the more compact expression, as $C_0 > 1$.

For the proof of (ii) we choose $\bar a_{N,T}^{(p)}(t)$ to be $a_1$ and $a_{N/T}^{(p)}(t/T)$ as $a_2$. Arguments as in the proof of (i), together with Lemma~\ref{lem:BoundEV_M0}(i) and Corollary~\ref{kor:exp_gamma_k}(iii) provide the needed bounds.
We have
\[
\| \E \hat \Gamma_{N,T}^{(p)}(t) - \Gamma_{N/T}^{(p)}(t/T) \| \cdot \| \Gamma_{N/T}^{(p)}(t/T)^{-1} \| \leq \Big( T^{-1} p^2 (2 \pi M'_f + C) + \frac{p^{2}}{N} M_f \Big) m_f^{-1}
\]
Thus, \eqref{lem:rel_a:prf:cond1} holds if $T \geq 4 p^2 (2 \pi M'_f + C) m_f^{-1}$ and $N \geq 4 p^{2} M_f / m_f$. 
Further, we have
\[\| \E \hat \gamma^{(p)}_{N,T}(t) - \gamma_{N/T}^{(p)}(t/T) \| \leq 2 \frac{p^{3/2}}{T} (2\pi M'_f + C) + \frac{p}{N} (2\pi)^{1/2} M_f, \quad \| \gamma_{N/T}^{(p)}(t/T) \| \leq (2\pi)^{1/2} M_f,\]
such that~\eqref{lem:rel_a:prf:cond2} holds if $T \geq \frac{4 p^{3/2} (2\pi M'_f + C)}{(2\pi)^{1/2} M_f}$ and $N \geq 2 p$.
Consequently, by \eqref{lem:rel_a:prf:bound}, we have the following bound
\begin{equation*}
	\begin{split}
		\| \bar a_{N,T}^{(p)}(t) - a_{N/T}^{(p)}(t/T) \|
		& \leq 4 m_f^{-2} (2\pi)^{1/2} M_f \Big( T^{-1} p^2 (2 \pi M'_f + C) + \frac{p^{2}}{N} M_f \Big) \\
		& \quad + m_f^{-1} \Big( 2 \frac{p^{3/2}}{T} (2\pi M'_f + C) + \frac{p}{N} (2\pi)^{1/2} M_f \Big) \\
		& \leq 4 m_f^{-2} (2\pi)^{1/2} M_f T^{-1} p^2 (2 \pi M'_f + C) \\
		& \qquad				+ 2 m_f^{-1} T^{-1} p^{3/2} (2\pi M'_f + C) + 2 C_0^2 \frac{p^{2}}{N} + C_0 \frac{p}{N} \\
		& = C_1 \frac{1}{T} \Big( 4 C_0 p^2  + 2 p^{3/2} \Big) + 2 C_0^2 \frac{p^{2}}{N} + C_0 \frac{p}{N},
	\end{split}
\end{equation*}
for $T \geq C_1 \max\{ 4 p^2, 4 p^{3/2} C_0^{-1}\}$ and $N \geq 4 p^{2} M_f / m_f$.
Thus, part (ii) of the lemma follows from $C_0^{-1} < 1$ and $p^{3/2} < C_0 p^2 / 2$.\hfill$\square$

\subsection{Proofs of Lemmas~\ref{lem:rel_v_bar}--\ref{lem:norm_g_der}}\label{app:der_thMSPE:proof}
\textbf{Proof of Lemma~\ref{lem:rel_v_bar}.} By the definitions of $\bar A_{N,T}^{(p)}(t)$ and $A^{(p)}_{\Delta}(t/T)$, which we now denote by $A$ and $A_0$, respectively.
We have
\begin{equation*}
	\begin{split}
		& \| \bar v_{N,T}^{(p,h)}(t) - v_0^{(p,h)}(t/T) \|
		= \big\| (1, 0, \ldots, 0 ) \big( A^h - A_0^h \big) \big\| \\
		& \leq \Big\| A^h - A_0^h \big\|
		\leq h \| A-A_0 \| \big( \|A-A_0\| + \|A_0\| \big)^{h-1} \\
		& \leq h \| \bar a_{N,T}^{(p)}(t) - a_0^{(p)}(t/T) \| \big( \| \bar a_{N,T}^{(p)}(t) - a_0^{(p)}(t/T) \| + \| a_0^{(p)}(t/T) \| \big)^{h-1} \\
		& \leq h ( 9 \frac{N}{T} C_0 C_1 p + 3 \frac{p^{2}}{N} C_0^2) \big( 9 \frac{N}{T} C_0 C_1 p + 3 \frac{p^{2}}{N} C_0^2 + C_0 \big)^{h-1} \\
		& \leq h \Big( 9 \frac{N}{T} C_0 C_1 p + 3 \frac{p^{2}}{N} C_0^2 \Big) \big( 2 C_0 \big)^{h-1}
			\leq h \Big( 5 C_1 p \frac{N}{T} + 2 \frac{p^{2}}{N}C_0 \Big) \big( 2 C_0 \big)^h, \\
	\end{split}
\end{equation*}
where the first inequality is due to the submultiplicativity of the spectral norm, the second inequality is due to Lemma~\ref{lem:expansionMatrixPower}, the third inequality is due to the fact that $\|A-A_0\| \leq \|A-A_0\|_F$ and $\| A_0 \| = \| e_1 a^{(p)}(t/T) \|$, because adding/subtracting the Jordan block doesn't change the eigenvalues (see also the proof of Lemma~\ref{lem:norm_v_der}).

The fourth inequality is due to Lemma~\ref{lem:norm_a} and Lemma~\ref{lem:rel_a_bar}(i), which can be applied since $18 C_1 p N \geq 8 p N C_1$ and $6 p^{2} C_0 \geq 4 p^{2} M_f / m_f$. The fifth inequality holds since
$T \geq 18 C_1 p N$ implies $9 \frac{N}{T} C_0 C_1 p \leq C_0/2$ and $N \geq 6 p^{2} C_0$ implies $3 \frac{p^{2}}{N} C_0^2 \leq C_0/2$.

For (ii) we have, by a similar argument as above that
\begin{equation*}
	\begin{split}
		& \| \bar v_{N,T}^{(p,h)}(t) - v_{N/T}^{(p,h)}(t/T) \|
		\leq h \big(C_1 \frac{1}{T} ( 5 C_0 p^2 )+ 3 \frac{p^{2}}{N} C_0^2 \Big) \big( C_1 \frac{1}{T} ( 5 C_0 p^2 ) + 3 \frac{p^{2}}{N} C_0^2 + C_0 \big)^{h-1} \\
		& \leq h \Big( C_1 \frac{1}{T} ( 5 C_0 p^2 ) + 3 \frac{p^{2}}{N} C_0^2 \Big) \big( 2 C_0 \big)^{h-1}
			\leq h \Big( 3 C_1 \frac{p^2}{T} + 2 \frac{p^{2}}{N} C_0 \Big) \big( 2 C_0 \big)^{h} \\
	\end{split}
\end{equation*}
where the second inequality holds due to
$T \geq 10 C_1 p^2$ implies $C_1 \frac{1}{T} (5 C_0 p^2) \leq C_0/2$ and $N \geq 6 p^{2} C_0$ implies $3 \frac{p^{2}}{N} C_0^2 \leq C_0/2$.
Note that $T \geq 4 C_1 p^2$ and that $N \geq 4 p^{2} M_f / m_f$, which are the conditions required to apply Lemma~\ref{lem:rel_a_bar}(ii).\hfill$\square$


\textbf{Proof of Lemma~\ref{lem:norm_v_der}.}
Let $e_1$ and $H$ be defined as in~\eqref{eqn:def_e_H}.
We derive a compact expression for ${\rm d}v$ where $v(u)$ is short for
\[ \big( v^{(p,h)}_{\Delta}(u) \big)' := e'_1 \big( e_1 (a(u))' + H \big)^h =: e'_1 X(u)^h,\]
with $a(u)$ being short for
$a^{(p)}_{\Delta}(u) := \Gamma_{\Delta}^{(p)}(u)^{-1} \gamma_{\Delta}^{(p)}(u)$.
We will further abbreviate $\Gamma(u) := \Gamma_{\Delta}^{(p)}(u)$ and $\gamma(u) := \gamma_{\Delta}^{(p)}(u)$.

By Weyl's inequality we have that the eigenvalues of $\mu_1 \leq \ldots \leq \mu_p$ of $e_1 a(u)' + H$ fulfil
$\lambda_1 + \rho_i \leq \mu_i \leq \rho_i + \lambda_p$,
where $\rho_1 \leq \ldots \leq \rho_p$ denote the eigenvalues of $e_1 a(u)'$ and $\lambda_1 \leq \ldots \leq \lambda_p$ denote the eigenvalues of $H$.
Note further, that the eigenvalues of $H$ are $\lambda_1 = \ldots = \lambda_p = 0$ (it is a Jordan block). Therefore: $\mu_i = \rho_i$. In particular: $\| e_1 a(u)' + H \| = \| e_1 a(u)' \| \leq \| a(u) \| \leq C_0$, where the last inequality follows from Lemma~\ref{lem:norm_a}. This observation, obviously, implies that
$\| v \| = \| e'_1 X^h\| \leq \|X\|^h \leq \|a\|^h$,
which yields (i).

For the proof of (ii) recall that the notation reflects that $v$, $a$, $\Gamma$, and $\gamma$ are functions of the variable $u$ (for a fixed $\Delta \geq 0$).
By elementary rule~(15), p.\,148 and Theorem 3, p.\,151 from~\cite{magneu88}, we have, that
\begin{equation}\label{eqn:a_der}
	{\rm d} a = - \Gamma^{-1} ({\rm d} \Gamma) \Gamma^{-1} \gamma + \Gamma^{-1} {\rm d} \gamma
\end{equation}

Iterating elementary rule (15), p.\,148, we get that, for every square matrix function $X$ and $h = 1, 2, \ldots$ that
\[{\rm d} (X^h) = \sum_{k=1}^h X^{k-1} ({\rm d} X) X^{h-k}.\]

Obviously we have ${\rm d} X = e_1 ({\rm d} a)'$, which yields
\begin{equation}\label{eqn:D_v}
	\begin{split}
		{\rm d} v' & = e'_1 \sum_{k=1}^h X^{k-1} ({\rm d} X) X^{h-k} \\
			& = e'_1 \sum_{k=1}^h \big( e_1 a' + H \big)^{k-1} (e_1 ({\rm d} a)') \big( e_1 a' + H \big)^{h-k}
	\end{split}
\end{equation}

For $\| {\rm d} a \|$ note that, by~\eqref{eqn:a_der} and employing Lemma~\ref{lem:BoundEV_M0}(i), we have
\begin{equation*}
	\begin{split}
		\| {\rm d} a \| & \leq \|\Gamma^{-1}\|^2 \| {\rm d} \Gamma \| \| \gamma \| + \| \Gamma^{-1}\| \| {\rm d} \gamma \| \\
			& \leq m_f^{-2} M'_f ( (2\pi)^{1/2} M_f ) + m_f^{-1} ( (2\pi)^{1/2} M'_f ) \\
			& = (2\pi)^{1/2} m_f^{-1} M'_f \big( M_f / m_f + 1 \big).
	\end{split}
\end{equation*}

Here we have used that, by the assumed continuity of $\frac{\partial}{\partial u} f(u,\lambda)$, we have
\[
	({\rm d} \Gamma_{\Delta}^{(p)})_{i,j}(u) = \int_{-\pi}^{\pi} {\rm e}^{{\rm i} (i-j) \lambda} \Big( \int_0^1 \frac{\partial}{\partial u} f\big(u+\Delta(x-1),\lambda\big) {\rm d} x \Big) {\rm d}\lambda,
\]
such that an application of Lemma~4.1 from~\cite{gra09} (cf. the proof of Lemma~\ref{lem:BoundEV_M0}(i)) yields that $\| {\rm d} \Gamma \| \leq M'_f$. Further, since
\[
	({\rm d} \gamma_{\Delta}^{(p)})_{j}(u) = \int_{-\pi}^{\pi} {\rm e}^{{\rm i} j \lambda} \Big( \int_0^1 \frac{\partial}{\partial u} f\big(u+\Delta(x-1),\lambda\big) {\rm d} x \Big) {\rm d}\lambda,
\]
we have by Bessel's inequality (see, again, the proof of Lemma~\ref{lem:BoundEV_M0}(i)) that $\| {\rm d} \gamma\| \leq (2\pi)^{1/2} M'_f$. In conclusion, [cf. \eqref{eqn:D_v}] we have shown
\begin{equation*}
	\begin{split}
		\| {\rm d} v \| & \leq h \|a \|^{h-1} \|{\rm d}a\|
			\leq h \big( C_0 \big)^{h-1} (2\pi)^{1/2} m_f^{-1} M'_f \big( M_f / m_f + 1 \big). \\
	\end{split}
\end{equation*}

For the proof of (iii) consider $\Delta$ to be the argument of the functions (instead of $u$). Let this be reflected by changing the notation to
$v(\Delta) := v^{(p,h)}_{\Delta}(u)$,
$a(\Delta) := a^{(p)}_{\Delta}(u)$,
$\Gamma(\Delta) := \Gamma^{(p)}_{\Delta}(u)$, and
$\gamma(\Delta) := \gamma^{(p)}_{\Delta}(u)$ (for a fixed $u \in \IR$).
Note that all previous arguments remain the same, but for the derivation of ${\rm d} \Gamma$ and ${\rm d}\gamma$. To this end, denote $f_{\lambda}(u) := f(\lambda, u)$ and note that we have
\begin{equation*}
\begin{split}
	\frac{\partial}{\partial \Delta} (\gamma_{\Delta}^{(p)}(u))_{j}
		& = \int_{-\pi}^{\pi} {\rm e}^{{\rm i} j \lambda} \Big( \int_0^1 \frac{\partial}{\partial \Delta} f\big(u+\Delta(x-1),\lambda\big) {\rm d} x \Big) {\rm d}\lambda \\
		& = \int_{-\pi}^{\pi} {\rm e}^{{\rm i} j \lambda} \Big( \int_0^1 f'_{\lambda}(u+\Delta(x-1)) x \, {\rm d} x \Big) {\rm d}\lambda \\
		& = \int_{-\pi}^{\pi} {\rm e}^{{\rm i} j \lambda} \Big( \Big[ \frac{f_{\lambda}(u+\Delta(x-1))}{\Delta} x \Big]_0^1 - \int_0^1 \frac{f_{\lambda}(u+\Delta(x-1))}{\Delta} \, {\rm d} x \Big) {\rm d}\lambda \\
		& = \int_{-\pi}^{\pi} {\rm e}^{{\rm i} j \lambda} \Big( \frac{1}{\Delta} \int_0^1 \big( f_{\lambda}(u) - f_{\lambda}(u+\Delta(x-1)) \big) \, {\rm d} x \Big) {\rm d}\lambda. \\
\end{split}
\end{equation*}

Thus,
\[
\Big\| \frac{\partial}{\partial \Delta} \, \gamma_{\Delta}^{(p)}(u) \Big\| \leq 2 (2\pi)^{1/2} M_f / \Delta, \text{ and }
\Big\| \frac{\partial}{\partial \Delta} \, \Gamma_{\Delta}^{(p)}(u) \Big\| \leq 2 M_f / \Delta.
\]
\hfill$\square$

\textbf{Proof of Lemma~\ref{lem:norm_a_bar}.}
This follows, because
\begin{equation*}
	\begin{split}
		\| \bar a_{N,T}^{(p)}(t) \|
		& = \| a^{(p)}_{N/T}(t/T) + \bar a_{N,T}^{(p)}(t) - a^{(p)}_{N/T}(t/T) \| \\
		& \leq \| a^{(p)}_{N/T}(t/T) \| + \| \bar a_{N,T}^{(p)}(t) - a^{(p)}_{N/T}(t/T) \| \\
		& \leq C_0 + \big( 5 C_0 C_1 \big) \frac{p^2}{T} + \big( 3 C_0^2 \big)\frac{p^{2}}{N}
		\leq 2 C_0,
	\end{split}
\end{equation*}
The second inequality follows by application of Lemmas~\ref{lem:norm_a} and~\ref{lem:rel_a_bar}(ii). Note that Lemma~\ref{lem:rel_a_bar}(ii) can be applied since
$T \geq 10 C_1 p^2 \geq 4 C_1 p^2$ and $N \geq 6 C_0 p^{2} \geq 4 p^{2} M_f / m_f$.
The third inequality follows, since $T \geq 10 C_1 p^2$ and $N \geq 6 C_0 p^{2} $, as assumed.

For the second inequality note that, similar to the definition of $\hat v_{N,T}^{(p,h)}(t)$ in Step 3 of Section~\ref{sec:prec} and~\eqref{def:vhat}, we have the following recursive relationship
\begin{equation*}
\begin{split}
	\big( \bar v^{(p,h)}_{N,T}(t) \big)' & = e'_1 \big( e_1 (\bar a^{(p)}_{N,T}(t))' + H \big)^h
				= \big( \bar v^{(p,h-1)}_{N,T}(t) \big)' \big( e_1 (\bar a^{(p)}_{N,T}(t))' + H \big) \\
			& = \big( \bar v^{(p,h-1)}_{N,T}(t) \big)' e_1 (\bar a^{(p)}_{N,T}(t))' + \big( \bar v^{(p,h-1)}_{N,T}(t) \big)' H. \\
\end{split}
\end{equation*}
In other terms, we have
\[\bar v_{i,N,T}^{(p,1)}(t) = \bar a_{i,N,T}^{(p)}(t)\]
and, for every h = 2, 3, \ldots, we have
\begin{equation*}
	\begin{split}
		\bar v_{i,N,T}^{(p)}(t,h)
			= \bar a_{i,N,T}^{(p)}(t) \bar v_{1,N,T}^{(p,h-1)}(t)
				+ \bar v_{i+1,N,T}^{(p,h-1)}(t) I\{i \leq p-1\},
	\end{split}
\end{equation*}
which implies that
\[\|\bar v_{N,T}^{(p,h)}(t)\|_{\infty} \leq (\| \bar a_{N,T}^{(p)}(t) \|_{\infty} + 1) \| \bar v_{N,T}^{(p,h-1)}(t) \|_{\infty} \leq \big(\| \bar a_{N,T}^{(p)}(t) \|_{\infty} + 1 \big)^{h-1} \| \bar a_{N,T}^{(p)}(t) \|_{\infty}. \]

Employing the first part of this lemma we thus have
\begin{equation*}
	\begin{split}\|\bar v_{N,T}^{(p,h)}(t)\|_{\infty}
	& \leq \big( 2 C_0 + 1 \big)^{h-1} 2 C_0 \leq \big( 2 C_0 + 1 \big)^h. \\
	\end{split}
\end{equation*}

For the proof of (ii), as in the proof of the bound for $\| \bar a_{N,T}^{(p)}(t) \|$ in (i), we have
\begin{equation*}
	\begin{split}
		\| \bar v_{N,T}^{(p,h)}(t) \|
		& \leq \| v^{(p,h)}_{N/T}(t/T) \| + \| \bar v_{N,T}^{(p,h)}(t) - v^{(p,h)}_{N/T}(t/T) \| \\
		& \leq \big( C_0 \big)^h
					  + h \big( 2 C_0 \big)^h \Big(3 C_1\frac{p^2}{T} + 2 \frac{p^{2}}{N} C_0 \Big)
			\leq 2 \big( C_0 \big)^h,
	\end{split}
\end{equation*}
where, for the second inequality, we have employed Lemma~\ref{lem:norm_v_der}(i) and~\ref{lem:rel_v_bar}(ii), which can be applied since $T \geq 6 h 2^h C_1 p^2 \geq 10 C_1 p^2$ and $N \geq 4 h 2^h p^{2} C_0 \geq 6 p^{2} C_0$ by assumption. The third inequality holds, since
\begin{equation*}
	\begin{split}
T & \geq h \big( 2 C_0 \big)^h \Big( 3 C_1\frac{p^2}{T} \Big) / \big( (C_0)^h / 2 \big) 
						= 6 h 2^h C_1 p^2, \text{ and}\\
		N & \geq 2 h \big( 2 C_0 \big)^h p^{2} C_0 / \big( (C_0)^h / 2 \big)
			 = 4 h 2^h p^{2} C_0,
	\end{split}
\end{equation*}
which, for all $h=1,2,\ldots$, implies the condition required to apply Lemma~\ref{lem:rel_v_bar}(ii).
\hfill$\square$

\textbf{Proof of Lemma~\ref{lem:norm_g_der}.}
For the sake of brevity, we omit sub and superscripts that are not needed for this particular argument. More precisely, in this section we denote $g^{(p,h)}_{\Delta}(u)$ by $g$, $v_{\Delta}^{(p,h)}(u)$ by $v$, $\gamma_0(u)$ by $\gamma_0$, $\gamma_0^{(p,h)}(u)$ by $\gamma$, and $\Gamma_0^{(p)}(u)$ by $\Gamma$, respectively.

For (i) note that we have, by~\eqref{equation:local_covariance_function} together with Assumption~\ref{a:sd_bounded}, by Lemma~\ref{lem:norm_v_der}(i)+(ii), and by Lemma~\ref{lem:BoundEV_M0}(i),
\begin{equation*}
	\begin{split}
		\| g \| & \leq \| \gamma_0 \| +  2 \| v' \| \cdot \| \gamma \| + \| v' \| \cdot \| \Gamma \| \cdot \| v \| \\
			& \leq 2\pi M_f +  2 (2\pi)^{1/2} M_f \big( C_0 \big)^h  + \big( C_0 \big)^{2h} M_f  \\
			& = \big((2\pi)^{1/2} + \big( C_0 \big)^h \big)^2 M_f
			\leq 4 M_f \big( C_0 \big)^{2 h} \\
	\end{split}
\end{equation*}

We now compute the derivatives of $g^{(p,h)}_{\Delta}(u)$.
We have:
\begin{equation*}
	\begin{split}
		{\rm d} g & = {\rm d} \big( \gamma_0 - 2 v' \gamma + v' \Gamma v \big) \\
			& = {\rm d} \gamma_0 - 2 \big( ({\rm d} v') \gamma + v' {\rm d} \gamma \big) \\
			& \quad + ( {\rm d} v') \Gamma v + v' ( {\rm d} \Gamma) v + v' \Gamma ( {\rm d} v ) \\
			& = {\rm d} \gamma_0 + 2 ( {\rm d} v') ( \Gamma v - \gamma) + v' \big( ( {\rm d} \Gamma) v - 2 {\rm d} \gamma \big).\\
	\end{split}
\end{equation*}
Note that, for $\Delta = 0$ and $h=1$, we have $v=a$ and by the definition of $a$, it follows that $ \Gamma v - \gamma = 0$, such that the second term of the derivative will only show if $h \geq 2$.

In consequence, for (ii) we bound ${\rm d} g_{p,h}$ as follows:
\begin{equation*}
	\begin{split}
		& |{\rm d} g_{p,h}| \\
			& = \| {\rm d} \gamma_0 + 2 ( {\rm d} v') ( \Gamma v - \gamma) + v' \big( ( {\rm d} \Gamma) v - 2 {\rm d} \gamma \big) \| \\
			& \leq \| {\rm d} \gamma_0 \| + 2 \| {\rm d} v' \| (\| \Gamma \| \| v \|  + \|\gamma\|) + \| v' \| \big( \|{\rm d} \Gamma \| \| v \| + 2 \| {\rm d} \gamma \| \big) \\
			& \leq 2 \pi M'_f + 2  h \big( C_0 \big)^h M'_f (m_f^{-1} + M_f^{-1} ) \Big( M_f \big( C_0 \big)^h + (2 \pi)^{1/2} M_f \Big) \\
			& \quad + \big( C_0 \big)^h \Big( M'_f \big( C_0 \big)^h + 2 (2 \pi)^{1/2} M'_f \Big) \\
			& = M'_f \Big[ 2 \pi + 2  h \big( C_0 \big)^h (M_f / m_f + 1 ) \Big( \big( C_0 \big)^h + (2 \pi)^{1/2} \Big)
			+ \big( C_0 \big)^h \Big( \big( C_0 \big)^h + 2 (2 \pi)^{1/2} \Big) \Big] \\
			& \leq M'_f \Big[ \big( C_0 \big)^{2 h} + 2  h \big( C_0 \big)^h (M_f / m_f + 1 ) \Big( \big( C_0 \big)^h + \big(C_0\big)^h \Big)
			+ \big( C_0 \big)^h \Big( \big( C_0 \big)^h + 2 \big(C_0\big)^h \Big) \Big] \\
			& = M'_f \Big[ \big( C_0 \big)^{2 h} + \big( C_0 \big)^{2 h} \Big( 4  h  (M_f / m_f + 1 ) + 3 \Big) \Big]
			\leq 4 M'_f \big( C_0 \big)^{2 h} \Big( h  (M_f / m_f + 1 ) + 1 \Big) \\
			& \leq 4 M'_f \big( C_0 \big)^{2 h} \Big( h  (C_0 + C_0 ) + C_0 \Big)
			\leq 4 M'_f \big( C_0 \big)^{2 h+1} \big( 2 h + 1 \big).
	\end{split}
\end{equation*}
In the above, for the second inequality, we have used~\eqref{eqn:gammaPrime} together with Assumption~\ref{a:sd_der_bounded}, Lemma~\ref{lem:norm_v_der}(i)+(ii), and Lemma~\ref{lem:BoundEV_M0}(i).

For the proof of (iii) use the partial derivatives of $\gamma$ and $\Gamma$ with respect to $\Delta$ that were derived in the proof of Lemma~\ref{lem:norm_v_der}(iii). For the proof of (iii) use Lemma~\ref{lem:norm_v_der}(iii) instead of Lemma~\ref{lem:norm_v_der}(ii) and use the bounds for $\| {\rm d} \Gamma \|$ and $\| {\rm d} \gamma \|$ that were derived in the end of the proof of Lemma~\ref{lem:norm_v_der}(iii).
For the proof of (iv) note that the partial derivative of $g$ is continuous by (i) and employ Leibniz's integral rule. (v) follows analogously from (ii). For the proof of (vi) we note that, by the argument in the proof of Lemma~\ref{lem:norm_v_der}(iii), we have
\begin{equation*}
\begin{split}
	& \frac{\partial}{\partial \Delta_2} {\rm MSPE}_{\Delta_1, \Delta_2}^{(p,h)}(u)
		= \frac{1}{\Delta_2} \int_0^1 \big( g_{\Delta_1}^{(p,h)}(u) - g_{\Delta_1}^{(p,h)}(u+\Delta_2(x-1)) \big) \, {\rm d} x.
\end{split}
\end{equation*}
(vi) then follows from (i) of this lemma. \hfill$\square$

\subsection{Proofs of the Lemmas~\ref{lem:exp_gamma_k}--\ref{lem:expIneqGenSum}}\label{sec:Cov:proof}

\textbf{Proof of Lemma~\ref{lem:exp_gamma_k}.} We have
\begin{equation*}
\begin{split}
	\E \hat \gamma_{k;N,T}(t) &  = \frac{1}{N} \sum_{\ell=t-N+|k|+1}^{t} \E (X_{\ell-|k|,T} X_{\ell,T}) \\
	& = \frac{N-|k|}{N} \frac{1}{N-|k|} \sum_{\ell=t-N+|k|+1}^{t} \big( \int_{-\pi}^{\pi} f(\ell/T, \lambda) {\rm e}^{{\rm i} |k| \lambda} {\rm d}\lambda + R^{(1)}_{\ell,|k|,T} \big), \\
	& = \frac{N-|k|}{N} \int_{-\pi}^{\pi} \Bigg[ \frac{1}{N-|k|} \sum_{\ell=t-N+|k|+1}^{t} f(\ell/T, \lambda) \Bigg] {\rm e}^{{\rm i} |k| \lambda} {\rm d}\lambda + R^{(2)}_{t,N,|k|,T}, \\
	& = \frac{N-|k|}{N} \int_{-\pi}^{\pi} \Bigg[ \int_0^1 f\big(\frac{t-N+|k|}{T} + \frac{N-|k|}{T} u, \lambda\big) {\rm d}u \\
	& \qquad\qquad + R^{(3)}_{t-N+|k|,t,T}(\lambda) \Bigg] {\rm e}^{{\rm i} |k| \lambda} {\rm d}\lambda + R^{(2)}_{t,N,|k|,T} \\
	& = \frac{N-|k|}{N} \int_{-\pi}^{\pi} \Bigg[ \int_0^1 f\big(\frac{t-N+|k|}{T} + \frac{N-|k|}{T} u, \lambda\big) {\rm d}u \Bigg] {\rm e}^{{\rm i} |k| \lambda} {\rm d}\lambda + R^{(4)}_{t,N,|k|,T},
\end{split}
\end{equation*}
where we have
\begin{equation*}
	\begin{split}
		|R^{(1)}_{\ell,|k|,T}| & \leq C/T, \quad \text{by~\eqref{eqn:time_varying_covariance_function}, \eqref{equation:local_covariance_function} and~\eqref{eqn:def:locally_stationary},}\\
		|R^{(2)}_{t,N,|k|,T}| & = \frac{1}{N} \Big| \sum_{\ell=t-N+|k|+1}^{t} R^{(1)}_{\ell,|k|,T} \Big| \leq \frac{N-|k|}{N} C/T \leq C/T \\
		\big| R^{(3)}_{t-N+|k|,t,T}(\lambda) \big| & \leq \sup_{(t-N+|k|)/T < u < t/T} \Big| \frac{\partial}{\partial u} f(u,\lambda) \Big| \frac{1}{T} \leq M'_f / T,
	\end{split}
\end{equation*}
by Lemma~\ref{lem:sumApprox}, where the differentiability (with respect to $u$) is guaranteed by Assumption~\ref{a:sd_der_bounded}, and
\begin{equation*}
	\begin{split}
		|R^{(4)}_{t,N,|k|,T}| & = \Big| \frac{N-|k|}{N} \int_{-\pi}^{\pi} R^{(3)}_{t-N+|k|,t,T}(\lambda) {\rm e}^{{\rm i} |k| \lambda}  {\rm d}\lambda +  R^{(2)}_{t,N,|k|,T} \Big| \leq \frac{2\pi M'_f + C}{T}
	\end{split}
\end{equation*}

For the second equation in the lemma, note that (for fixed $t, N, k, u$) we have
	\begin{equation}\label{lem:exp_gamma_k:eqn:A}
	\begin{split}
		& f\big(\frac{t-N+|k|}{T} + \frac{N-|k|}{T} u, \lambda\big) \\
		& = f\big(\frac{t-N}{T} + \frac{N}{T} u, \lambda\big) \\
		& \quad + f\big(\frac{t-N+|k|}{T} + \frac{N-|k|}{T} u, \lambda\big) - f\big(\frac{t-N}{T} + \frac{N}{T} u, \lambda\big) \\
		& = f\big(\frac{t-N}{T} + \frac{N}{T} u, \lambda\big) + \frac{\partial}{\partial x} f\big(x, \lambda\big)\Big|_{x=\xi} \Big( \frac{|k|}{T} (1-u)\Big), \\
	\end{split}
	\end{equation}
where $\xi \in \Big(\frac{t-N}{T} + \frac{N}{T} u, \frac{t-N+|k|}{T} + \frac{N-|k|}{T} u\Big)$.

Thus, by Assumption~\ref{a:sd_der_bounded} and because $u \in [0,1]$, we have
	\begin{equation*}
	\begin{split}
		\Big| f\big(\frac{t-N+|k|}{T} + \frac{N-|k|}{T} u, \lambda\big) - f\big(\frac{t-N}{T} + \frac{N}{T} u, \lambda\big) \Big|
		& \leq M'_f \Big( \frac{|k|}{T}\Big).\\
	\end{split}
	\end{equation*}
	
Thus,
	\begin{equation*}
	\begin{split}
		& \Big| \int_{-\pi}^{\pi} \Bigg[ \int_0^1 \Big( f\big(\frac{t-N+|k|}{T} + \frac{N-|k|}{T} u, \lambda\big) - f\big(\frac{t-N}{T} + \frac{N}{T} u, \lambda\big) \Big) {\rm d}u \Bigg] {\rm e}^{{\rm i} |k| \lambda} {\rm d}\lambda \Big| \\
		& \leq  \int_{-\pi}^{\pi} \int_0^1 \Bigg| f\big(\frac{t-N+|k|}{T} + \frac{N-|k|}{T} u, \lambda\big) - f\big(\frac{t-N}{T} + \frac{N}{T} u, \lambda\big) \Bigg| {\rm d}u {\rm d}\lambda \leq 2 \pi M'_f \Big( \frac{|k|}{T}\Big).
	\end{split}
	\end{equation*}

For the proof of the third part of the Lemma, note that
(for fixed $t, N, k, u$) we have, in analogy to~\eqref{lem:exp_gamma_k:eqn:A}
	\begin{equation*}
	\begin{split}
		f\big(\frac{t-N+|k|}{T} + \frac{N-|k|}{T} u, \lambda\big)
		& = f\big(\frac{t}{T}, \lambda\big) + \frac{\partial}{\partial x} f\big(x, \lambda\big)\Big|_{x=\xi} \Big( \frac{|k|-N}{T} (1-u)\Big), \\
	\end{split}
	\end{equation*}
where $\xi \in \Big(\frac{t-N+|k|}{T} + \frac{N-|k|}{T} u, \frac{t}{T}\Big)$.

Thus, by Assumption~\ref{a:sd_der_bounded} and because $u \in [0,1]$, we have
	\begin{equation*}
	\begin{split}
		\Big| f\big(\frac{t-N+|k|}{T} + \frac{N-|k|}{T} u, \lambda\big) - f\big(\frac{t}{T}, \lambda\big) \Big|
		& \leq M'_f \Big( \frac{N-|k|}{T}\Big).\\
	\end{split}
	\end{equation*}

Thus,
	\begin{equation*}
	\begin{split}
		& \Big| \int_{-\pi}^{\pi} \Bigg[ \int_0^1 \Big( f\big(\frac{t-N+|k|}{T} + \frac{N-|k|}{T} u, \lambda\big) - f\big(\frac{t}{T}, \lambda\big) \Big) {\rm d}u \Bigg] {\rm e}^{{\rm i} |k| \lambda} {\rm d}\lambda \Big| \\
		& \leq  \int_{-\pi}^{\pi} \int_0^1 \Bigg| f\big(\frac{t-N+|k|}{T} + \frac{N-|k|}{T} u, \lambda\big) - f\big(\frac{t}{T}, \lambda\big) \Bigg| {\rm d}u {\rm d}\lambda \leq 2 \pi M'_f \Big( \frac{N-|k|}{T}\Big),
	\end{split}
	\end{equation*}
which completes the proof of the lemma.\hfill$\square$

\textbf{Proof of Corollary~\ref{kor:exp_gamma_k}.}

For (i), note that, letting
\begin{equation*}
	\begin{split}
		\tilde\Gamma_{T}^{(p)}(t) & - \Gamma_0^{(p)}(t/T) =: \big( \delta_{j,k}(t,T) \big)_{j,k=1,\ldots,p}\\
		\tilde\gamma^{(p)}_{T}(t) & - \gamma_0^{(p)}(t/T) =: \big( \eta_{j}(t,T) \big)_{j=1,\ldots,p}
	\end{split}
\end{equation*}
we have
\begin{equation*}
	\begin{split}
		|\delta_{j,k}(t,T)|
			& = \Big| \tilde\gamma_{k-j,T}(t-j) - \gamma_{j-k}\big(\frac{t-j}{T}\big)
					+ \gamma_{j-k}\big(\frac{t-j}{T}\big) - \gamma_{j-k}\big(\frac{t}{T}\big) \Big| \\
			& \leq C/T + \Big| \int\limits_{-\pi}^{\pi}\exp\left({\rm i}(j-k)\lambda\right) \Big[ f(\frac{t-j}{T},\lambda) - f(\frac{t}{T},\lambda) \Big] {\rm d}\lambda \Big| \\
			& \leq (C + 2\pi M'_f j)/T
		\end{split}
\end{equation*}
where the first part of the first inequality is due to~\eqref{eqn:def:locally_stationary} and the second part is due to~\eqref{eqn:sd_der_bounded} and the mean value theorem.

We bound the norms of $\Delta := \big( \delta_{j,k}(t,T) \big)_{j,k=1,\ldots,p}$
\begin{equation*}
	\begin{split}
		\| \Delta \|_1 & = \max_{1 \leq k \leq p} \sum_{j=1}^p | \delta_{j,k} |
					\leq \max_{1 \leq k \leq p} \sum_{j=1}^p \frac{C + 2\pi M'_f j}{T}
					= T^{-1} p ( C + \pi M'_f (p+1) ) \\
		\| \Delta \|_{\infty} & = \max_{1 \leq j \leq p} \sum_{k=1}^p | \delta_{j,k} |
					\leq \max_{1 \leq j \leq p} \sum_{k=1}^p \frac{C + 2\pi M'_f j}{T}
					= p \frac{C + 2\pi M'_f p}{T} \\
		\end{split}
\end{equation*}

Then, we have, by H\"{o}lder's inequality
\begin{equation*}
	\begin{split}
		\| \Delta \|_2 & \leq \Big( \| \Delta \|_1 \| \Delta \|_{\infty} \Big)^{1/2} \\
			& \leq T^{-1} p \Big( ( C + \pi M'_f (p+1) ) (C + 2\pi M'_f p) \Big)^{1/2} \\
			& \leq T^{-1} p (C + 2\pi M'_f p) \leq T^{-1} p^{2} (C + 2\pi M'_f)
		\end{split}
\end{equation*}

Further, note that
\begin{equation*}
	\begin{split}
		|\eta_{j}(t,T)|
			& = \Big| \tilde\gamma_{j,T}(t) - \gamma_j\big(\frac{t}{T}\big) \Big| \leq C/T
		\end{split}
\end{equation*}
and we get the second part of (i) by bounding the Euclidean norm by $p^{1/2}$ times the maximum norm.

For (ii) note that, by Lemma~\ref{lem:exp_gamma_k}(iii), we have
\begin{equation*}
	\begin{split}
		\E \hat \Gamma_{N,T}^{(p)}(t) - F_{p,N} \circ \Gamma_0^{(p)}(t/T) & =: \big( \delta_{j-k}(t,T) \big)_{j,k=1,\ldots,p}  =: \Delta\\
		\E \hat \gamma^{(p)}_{N,T}(t) - f_{p,N} \circ \gamma_0^{(p)}(t/T) & =: \big( \delta_{j}(t,T) \big)_{j=1,\ldots,p}  =: \delta
	\end{split}
\end{equation*}
with
\begin{equation*}
	\begin{split}
		|\delta_{k}(t,T)| & \leq \Big|\frac{2 \pi M'_f (N-|k|+1) + C}{T}\Big|. 
		\end{split}
\end{equation*}

We bound the spectral norm of $\Delta$ by its Frobenius norm:
\begin{equation*}
	\begin{split}
		\| \Delta \|_2^2 & \leq \| \Delta \|_F^2 
			= \sum_{i=1}^p \sum_{j=1}^p |\delta_{i-j}|^2
			\leq \Big| \frac{2 \pi M'_f (N+1) + C}{T}\Big|^2 p^2 \\
		\end{split}
\end{equation*}

For the second part of (ii) note, that
\begin{equation*} 
	\begin{split}
		& \Big\| \E \hat \gamma^{(p)}_{N,T}(t) - f_{p,N} \circ \gamma_0^{(p)}(t/T) \Big\|
		= \Big( \sum_{j=1}^p |\delta_{j}(t,T)|^2 \Big)^{1/2} \\
		& \leq \Big( \sum_{j=1}^p \Big|\frac{2 \pi M'_f N + C}{T}\Big|^2 \Big)^{1/2}
		\leq p^{1/2} \Big|\frac{2 \pi M'_f N + C}{T}\Big| \\
	\end{split}
\end{equation*}

For the proof of (iii) note that, letting
\begin{equation*}
	\begin{split}
		\E \Gamma_{N,T}^{(p)}(t) & - F_{p,N} \circ \Gamma_{N/T}^{(p)}(t/T) =: \big( \delta_{j-k}(t,T) \big)_{j,k=1,\ldots,p}\\
		\E \gamma^{(p)}_{N,T}(t) & - f_{p,N} \circ \gamma_{N/T}^{(p)}(t/T) =: \big( \eta_{j}(t,T) \big)_{j=1,\ldots,p}
	\end{split}
\end{equation*}
we have
\begin{equation*}
\begin{split}
	\|\Delta\|^2 &
	\leq \|\Delta\|_F^2
	= \sum_{i=1}^p \sum_{j=1}^p |\delta_{i-j}|^2
	\leq \sum_{i=1}^p \sum_{j=1}^p \Big( \frac{2\pi (|i-j|+1) M'_f}{T} + \frac{C}{T}\Big)^2 \\
	& \leq \Big(\frac{2\pi M'_f + C}{T}\Big)^2 \Big( \sum_{i=1}^p \sum_{j=1}^p |i-j|^2 + 2 \sum_{i=1}^p \sum_{j=1}^p |i-j| + \sum_{i=1}^p \sum_{j=1}^p 1 \Big) \\
	& = \Big(\frac{2\pi M'_f + C}{T}\Big)^2 \frac{p}{6} \Big( (p^2 - 1) (p+4) + 6p \Big)
	\leq \Big(\frac{2\pi M'_f + C}{T}\Big)^2 p^4.
\end{split}
\end{equation*}

For the other bound note
\begin{equation*}
\begin{split}
	& \Big( \sum_{j=1}^{p} (\frac{2\pi(j+1)M'_f+C}{T})^2 \Big)^{1/2}
	\leq \frac{1}{T} \Big( 2\pi M'_f \big( \sum_{j=1}^{p} (j+1)^2 \big)^{1/2} + p^{1/2} C \Big) \\
	& \leq \frac{1}{T} \Big( 2\pi M'_f \big( p (p+1)^2 \big)^{1/2} + p^{1/2} C \Big)
	\leq \frac{1}{T} \Big( 4\pi M'_f p^{3/2} + p^{1/2} C \Big). \\
\end{split}
\end{equation*}

\hfill$\square$

For the proof of Lemma~\ref{lem:BoundEV_M0} we will use the following lemma which is derived using a result from \cite{BapatSunder1985}.

\begin{lem}\label{lem:FpnEV}
	Let $M$ be a positive definite $p \times p$ matrix. For $n=1,2,\ldots$ and $p = 1, \ldots, n$ define
	\[F_{p,n} := \big( 1 - |j-k|/n \big)_{j,k=1,\ldots,p}, \quad \text{and} \quad
	  f_{p,n} := \big( 1-1/n, \ldots, 1-p/n \big)'.\]
	Then, $F_{p,n}$ is positive semidefinite and $F_{p,n} \circ M$ is positive definite. Further, the eigenvalues of $F_{p,n} \circ M$ are bounded by the smallest eigenvalue of $M$ from below and by the largest eigenvalue of $M$ from above; i.\,e.,
	\[0 < 1/\|M^{-1}\| \leq 1/\|(F_{p,n} \circ M)^{-1}\| \leq \|F_{p,n} \circ M\| \leq \|M\|.\]
\end{lem}

\textbf{Proof of Lemma~\ref{lem:FpnEV}.}

First, $F_{p,n}$ is positive semidefinite, as we have $F_{p,n} = A_{p,n} A'_{p,n}$, where
\[A_{p,n} := n^{-1/2} \begin{pmatrix}
0      & \cdots & 0      & 0      & 1      & \cdots & 1      & 1      & \cdots & 1      & 1 \\
0      & \cdots & 0      & 1      & 1      & \cdots & 1      & 1      & \cdots & 1      & 0 \\
\vdots & \ddots & \vdots & \vdots & \vdots & \ddots & \vdots & \vdots & \ddots & \vdots & \vdots \\
1      & \cdots & 1      & 1      & 1      & \cdots & 1      & 0      & \cdots & 0      & 0
\end{pmatrix}\]
is a $p \times (n+p)$ matrix.
Then, we have $x' F_{p,n} x = \| A'_{p,n} x \|^2 \geq 0$.

We now employ Theorem~3(i) from \cite{BapatSunder1985}, which states the following: if $A$ and $B$ are two positive semidefinite $p \times p$ matrices and $k=1,\ldots,p$, then
\begin{equation}\label{eqn:lem:FpnEV:a}
	\sum_{i=1}^k \lambda_i(A \circ B) \leq \sum_{i=1}^k \lambda_i(A) B_{ii}, \quad k = 1, \ldots, p.
\end{equation}
If the diagonal entries of $B$ all equal 1, this means that the spectrum of $A \circ B$ is (weakly) majorised by that of $A$. (It follows that it's majorised and not just weakly majorised, because ${\rm trace}(A \circ B) = {\rm trace}(A)$ if $B_{ii} = 1$.)

Further, if $\lambda_i(A \circ B) \geq 0$ and $\lambda_i(A) B_{ii} \geq 0$, for all $i=1,\ldots,p$, and if $\sum_{i=1}^p \lambda_i(A \circ B) = \sum_{i=1}^p \lambda_i(A) B_{ii}$, then
\begin{equation}\label{eqn:lem:FpnEV:b}
	\sum_{i=k}^p \lambda_i(A \circ B) \geq \sum_{i=k}^p \lambda_i(A) B_{ii}, \quad k = 1, \ldots, p.
\end{equation}

We employ these inequalities with $A = M$ and $B = F_{p,n}$. Note that both are positive semidefinite ($M$ by assumption and $F_{p,n}$ as we already showed) and that we further have $B_{ii} = 1$. From \eqref{eqn:lem:FpnEV:b}, with $k = 1$, we obtain $\|F_{p,n} \circ M\| \leq \|M\|$ and from \eqref{eqn:lem:FpnEV:b}, with $k = p$, we see $1/\|M^{-1}\| \leq 1/\|(F_{p,n} \circ M)^{-1}\|$.\hfill$\square$

\textbf{Proof of Lemma~\ref{lem:BoundEV_M0}.} We use the notation from~\cite{gra09}. For every function $f$ defined on $[0,2\pi]$, with a Fourier series that has absolutely summable Fourier coefficients and for every $n \in \IN^*$, he defines, in~(4.8), the $p \times p$ Toeplitz matrix
\[T_p(f) := \Bigg[ \int_{-\pi}^{\pi} f(\lambda) {\rm e}^{-{\rm i}(k-j)\lambda} {\rm d}\lambda \ : \ k,j = 0,1,\ldots,p-1 \Bigg].\]

We will apply Lemma~4.1 from~\cite{gra09} according to which, if $T_p(f)$ is Hermitian (which is the case if $f$ is real-valued), its eigenvalues lie between the essential infimum and the essential supremum of $f$.

Letting $g_{\Delta}(u,\lambda) := \int_0^1 f\big(u + \Delta(x-1), \lambda\big) {\rm d}x$, we have $\Gamma_{\Delta}^{(p)}(u) = T_p\big( g_{\Delta}(u,\cdot) \big)$. The proof of (i-a) is thus completed, as $g_{\Delta}(u,\IR) \subset [m_f, M_f]$, by Assumption~\ref{a:sd_bounded}.
(i-b) follows from Bessel's inequality, as we have
\begin{equation*}
	\begin{split}
	\|\gamma_{\Delta}^{(p)}(u)\|^2
	 & = \sum_{\ell=1}^p \Big| \int\limits_{-\pi}^{\pi} {\rm e}^{{\rm i} \ell \lambda} g_{\Delta}(u,\lambda){\rm d}\lambda \Big|^2
	\leq \int\limits_{-\pi}^{\pi} \big| g_{\Delta}(u,\lambda)\big|^2 {\rm d}\lambda 
	\leq 2\pi M_f^2. \\
	\end{split}
\end{equation*}

For the proof of (ii-a) apply the triangle inequality:
\[ \|\E \hat\Gamma_{N,T}^{(p)}(t) \| \leq \| F_{p,N} \circ \Gamma_{N/T}^{(p)}(t/T) \| + \|\E \hat\Gamma_{N,T}^{(p)}(t) - F_{p,N} \circ \Gamma_{N/T}^{(p)}(t/T) \| \]
and note that the first term on the right hand side can be bounded by
\[\| F_{p,N} \circ \Gamma_{N/T}^{(p)}(t/T) \| \leq \| \Gamma_{N/T}^{(p)}(t/T) \| \leq M_f,\]
where we have employed Lemma~\ref{lem:FpnEV} and (i-a) of this lemma. The second part can be bounded by employing Corollary~\ref{kor:exp_gamma_k}(iii).

For the proof of (ii-b) we will employ Lemma~\ref{lem:stoerungslemma}. Following the notation there, we denote $A := F_{p,n} \circ \Gamma_{N/T}^{(p)}(t/T)$ and $\Delta := \E \Gamma_{N,T}^{(p)}(t) - F_{p,N} \circ \Gamma_{N/T}^{(p)}(t/T)$. By Corollary~\ref{kor:exp_gamma_k}(iii) we have $\| \Delta \| \leq \frac{p^2}{T} (C + 2 \pi M'_f)$. Thus, by assuming $T > m_f^{-1} p^2 (C + 2 \pi M'_f)$ we have
\[\| A^{-1}\| \| \Delta\| \leq m_f^{-1} \frac{p^2}{T} (C + 2 \pi M'_f) < 1. \]
Under this condition Lemma~\ref{lem:stoerungslemma} asserts that $\E \Gamma_{N,T}^{(p)}(t)$ is invertible with
\[\frac{1}{\|\E \Gamma_{N,T}^{(p)}(t)^{-1}\|} \geq \frac{1 - \| A^{-1}\| \cdot \|\Delta\|}{\|A^{-1}\|} = \frac{1}{\|A^{-1}\|} - \|\Delta\| \geq m_f - \frac{p^2}{T} (2 \pi M'_f + C).\]

Finally, because $2 m_f^{-1} p^2 (2 \pi M'_f + C) \geq 2 M_f^{-1} p^2 (2 \pi M'_f + C)$ we see that (ii-c), the ``in particular'', holds.\hfill$\square$

\textbf{Proof of Lemma~\ref{lem:exp_ineq_gamma_k2}.} First, note that $\frac{N}{N-|h|} \geq 1$ and therefore,
\begin{equation}
\label{lem:exp_ineq_gamma_k2:eqn1}\IP\Big(\Big| \hat\gamma_{h;N,T}(t) - \E \hat\gamma_{h;N,T}(t) \Big| \geq \varepsilon\Big)
\leq \IP\Big(\frac{N}{N-|h|} \Big| \hat\gamma_{h;N,T}(t) - \E \hat\gamma_{h;N,T}(t) \Big| \geq \varepsilon\Big).
\end{equation}
We next state and prove Lemma~\ref{lem:expIneqGenSum}, which entails a bound for the right hand side of~\eqref{lem:exp_ineq_gamma_k2:eqn1} as a special case (let $a_t := 1$, $\alpha = 1$, $b := t-N+1$, and $n := N-|h|$).\hfill$\square$

\textbf{Proof of Lemma~\ref{lem:expIneqGenSum}.}

We proceed along the lines of the proof of Lemma~2 in \cite{golzee01} and apply the two results from~\cite{sausta91} which are cite in Section~\ref{sec:res_sausta91}. As \cite{golzee01} pointed out, several other exponential-type inequalities are available in the literature; e.\,g., \cite{Doukh1994} and \cite{bos98}. These result could be used to derive the type of bound we need, but under our conditions the results from \cite{sausta91} seem most suitable here.

The assertion of Lemma~\ref{lem:expIneqGenSum} follows from Lemma~\ref{lem:sausta91:Lemma24} if we show that~\eqref{lem:sausta91:Lemma24:cond} holds with $\gamma := 1 + 2 \alpha d$,
\begin{equation*}
	\begin{split}
		\xi & := \sum_{t=1}^{n} n^{-1} a_t \big(X_{t+b-1,T}^\alpha X_{t+b-1+h,T}^\alpha - \E (X_{t+b-1,T}^\alpha X_{t+b-1+h}^\alpha)\big),
	\end{split}
	\end{equation*}
	and
\begin{equation*}
	\begin{aligned}
		H &:= \frac{C_{1,\alpha} A^2 h_*}{n}, \quad
		\bar\Delta & := \frac{n}{C_{2,\alpha} A h_*},
		\end{aligned}
\end{equation*}
where $h_* := |h| + I\{h=0\}$, and
\begin{equation*}
	\begin{aligned}
		C_{1,\alpha} & := 12 \cdot 2^{10 \alpha d + 7} \alpha^{4 \alpha d} \big(\max\{c^2, 3\pi M_f, 1\} \big)^{2\alpha} {\rm e} \Big( 1 + \frac{1}{\log \rho}\Big)\Big(1 + K^{1/2} \Big), \\
		C_{2,\alpha} & := 12 \cdot 2^{4 \alpha d + 3} \alpha^{2 \alpha d} \big(\max\{c^2, 3\pi M_f, 1\} \big)^{\alpha} {\rm e} \Big( 1 + \frac{1}{\log \rho}\Big),
		\end{aligned}
\end{equation*}
as defined in~\eqref{lem:exp_ineq_gamma_k2:def:C1C2}
We show this by employing Lemma~\ref{lem:sausta91:Theorem417} with
\[Y_{\ell,T} := n^{-1} a_t \big(X_{\ell+b-1,T}^\alpha X_{\ell+b-1+h,T}^\alpha - \E (X_{\ell+b-1,T}^\alpha X_{\ell+b-1+h,T}^\alpha)\big),\]
$\ell = 1, \ldots, n$. The first step in the application of Lemma~\ref{lem:sausta91:Theorem417} is to show that~\eqref{lem:sausta91:Theorem417:cond} holds with appropriately chosen $H_1$ and $\gamma_1$. The second step will be to show that the bound we get from Lemma~\ref{lem:sausta91:Theorem417} can again be bounded by the right hand side in~\eqref{lem:sausta91:Lemma24:cond} with $H$, $\bar\Delta$, and $\gamma$ as defined before.

For the first step, observe that for any $\ell = 1, \ldots, n$ and $k = 2, 3, \ldots$ we have
\begin{align}
		& n^k \E|Y_{\ell,T}|^k \nonumber \\
		& = |a_t|^k \E \big| X_{\ell+b-1,T}^\alpha X_{\ell+b-1+h,T}^\alpha - \E (X_{\ell+b-1,T}^\alpha X_{\ell+b-1+h,T}^\alpha) \big|^k \nonumber \\
			& \leq A^k 2^{k-1} \E \Big( |X_{\ell+b-1,T}^\alpha X_{\ell+b-1+h,T}^\alpha|^k + \big| (\E X_{\ell+b-1,T}^{2\alpha}) (\E X_{\ell+b-1+h,T}^{2\alpha} ) \big|^{k/2} \Big) \nonumber \\
			& \leq A^k 2^{k} (k!)^{2 \alpha d}  \big( (2 \alpha)^{2 \alpha d} (\max\{c^2, 3 \pi M_f, 1\})^{\alpha} \big)^k. \label{lem:expIneqGenSum:eqn:A}
	\end{align}
For the first part of~\eqref{lem:expIneqGenSum:eqn:A} we have used that
\begin{equation*}
	\begin{split}
		& \E | X_{\ell+b-1,T}^\alpha X_{\ell+b-1+h,T}^\alpha |^k = \E (|X_{\ell+b-1,T}|^{\alpha k} |X_{\ell+b-1+h,T}|^{\alpha k}) \\
		& \leq (\E |X_{\ell+b-1,T}|^{2 \alpha k})^{1/2} ( \E|X_{\ell+b-1+h,T}|^{2 \alpha k})^{1/2} \\
		& \leq \big( c^{2 \alpha k - 2} \big( (2 \alpha k)! \big)^d \sigma_{\ell+b-1,T}^2 c^{2 \alpha k - 2} \big( (2 \alpha k)! \big)^d \sigma_{\ell+b-1+h,T}^{2} \big)^{1/2} \\
		& \leq  c^{2 \alpha k - 2} (2\alpha)^{2 \alpha d k}   ( k! )^{2 \alpha d}  (3\pi M_f)
		\leq  (k!)^{2 \alpha d} \big( (2 \alpha)^{2 \alpha d} (\max\{c^2, 3 \pi M_f, 1\})^{\alpha} \big)^{k}.
	\end{split}
\end{equation*}
In the above we have used Cauchy-Schwarz inequality, Assumption~\ref{a:MomentCond}, Assumptions~\ref{a:loc_stat}, \ref{a:sd_bounded}, which, if $T \geq C/(\pi m_f)$, imply~\eqref{eqn:cons:sd_bounded}, and the elementary inequality $(L k)! \leq L^{L k} (k!)^L$, which holds for $L, k \in \IN^*$. For the second part of~\eqref{lem:expIneqGenSum:eqn:A} we have used that
\begin{equation*}
	\begin{split}
		& \big| (\E X_{\ell+b-1}^{2\alpha}) (\E X_{t-s-1+h,T}^{2\alpha} ) \big|^{k/2}
			\leq | c^{2\alpha-2} ((2\alpha)!)^d \sigma_{\ell+b-1,T}^{2} c^{2\alpha-2} ((2\alpha)!)^d \sigma_{\ell+b-1+h,T}^{2} |^{k/2} \\
			& \leq \big( c^{2\alpha-2} (2 \alpha)^{2 \alpha d} (3 \pi  M_f ) \big)^k
			\leq (k!)^{2 \alpha d} \big( (2 \alpha)^{2 \alpha d} (\max\{c^2, 3 \pi M_f, 1\})^{\alpha} \big)^k,
	\end{split}
\end{equation*}
where we have employed arguments as before, $(2\alpha)! \leq (2\alpha)^{2\alpha}$, and $1 \leq (k!)^{2 \alpha d}$.

Thus we have shown that~\eqref{lem:sausta91:Theorem417:cond} holds, with
$H_1 := 2 A (2 \alpha)^{2 \alpha d} (\max\{c^2, 3 \pi M_f, 1\})^{\alpha}\} / n$ and $\gamma_1 := 2 \alpha d - 1$.

For the second step we turn our attention to $\Delta_n(\nu)$ defined in Lemma~\ref{lem:sausta91:Theorem417}.
Note that, for $s < t$ and $h \geq 0$, we have
$\sigma(Y_{\ell,T} : 1 \leq \ell \leq s) \subseteq \sigma(X_{\ell,T} : \ell \leq b+s-1+h)$ and
$\sigma(Y_{\ell,T} : \ell \geq t) \subseteq \sigma(X_{\ell,T} : \ell \geq t+b-1)$, which implies that
\begin{equation*}
\begin{split}
	\alpha^Y(s,t) & := \sup_{A \in \sigma(Y_{\ell,T} : 1 \leq \ell \leq s)} \sup_{B \in \sigma(Y_{\ell,T} : \ell \geq t)} \big| \IP(A \cap B) - \IP(A) \IP(B) \big| \\
	& \leq  \sup_{A \in \sigma(X_{\ell,T} : \ell \leq s+b-1+h)} \sup_{B \in \sigma(X_{\ell,T} : \ell \geq t+b-1)}
	\big| \IP(A \cap B) - \IP(A) \IP(B) \big| \leq \alpha(t-s-|h|),
	\end{split}
\end{equation*}
where $\alpha$ is the function from Assumption~\ref{a:mixing}. The case where $h < 0$ is derived analogously.
Further, note that the trivial inequality $\alpha(t-s-|h|) \leq 1$ holds for $t \leq s+|h|$.
Thus, adopting the argument from step~3 in the proof of Lemma~2 in~\cite{golzee01}, we have
\begin{equation*}
	\begin{split}
		\sum_{t=s}^n \alpha^Y(s,t)^{1/\nu}
		& \leq \sum_{t=s}^n \alpha(t-s-|h|)^{1/\nu}
		\leq \sum_{t=s}^{s+|h|} 1 + \sum_{t=s+|h|+1}^{n} \alpha(t-s-|h|)^{1/\nu} \\
		& \leq |h| + \sum_{t=1}^{n-s-|h|} K^{1/\nu} \rho^{-t/\nu}
		\leq |h| + K^{1/\nu} \sum_{t=0}^{\infty} \rho^{-t/\nu} \\
		& = |h| + K^{1/\nu} \frac{\rho^{1/\nu}}{\rho^{1/\nu}-1}
		\leq |h| + K^{1/\nu} \Big( 1 + \frac{\nu}{\log \rho}\Big), \\
	\end{split}
\end{equation*}
where the last inequality follows from $\exp(x) - 1 \geq x$, for $x \geq 0$.

Thus, for $\delta = 1$ and $k=2,3,\ldots$, we have
\begin{equation*}
	\begin{split}
		& \Delta_n((1+1/\delta)(k-1))^{k-1}
		\leq \Big(|h| + K^{\frac{1}{2(k-1)}} \Big( 1 + \frac{2(k-1)}{\log \rho}\Big) \Big)^{k-1} \\
		& \leq 2^{k-2} \Big(|h|^{k-1} + K^{1/2} \Big( 1 + \frac{2(k-1)}{\log \rho}\Big)^{k-1} \Big) \\
		& \leq 2^{k-2} \Big(h_*^{k-1} + K^{1/2} \Big( 1 + \frac{2}{\log \rho}\Big)^{k-1}(k-1)^{k-1} \Big) \\
		& \leq 2^{k-2} \Big(h_*^{k-1} + K^{1/2} \Big( 1 + \frac{2}{\log \rho}\Big)^{k-1} (k-1)! {\rm e}^{k-1} \Big) \\
		& \leq 2^{k-2} \Big(h_*^{k-1} \Big( 1 + \frac{2}{\log \rho}\Big)^{k-1} (k-1)! {\rm e}^{k-1} \Big) \Big(1 + K^{1/2} \Big) \\
		& \leq \Big( 4 \ {\rm e} \ h_* \Big( 1 + \frac{1}{\log \rho}\Big) \Big)^{k-1} k!  \Big(1 + K^{1/2} \Big). \\
	\end{split}
\end{equation*}
where we have used the fact that $p^p \leq p! \, {\rm e}^{p}$, for $p \in \IN$.

Recall that $\gamma := 1 + 2 \alpha d$ and $\gamma_1 := 2 \alpha d - 1$. Thus, by Lemma~\ref{lem:sausta91:Theorem417}, we have that
\begin{equation*}
	\begin{split}
		& \Big| \cum_k \Big( \sum_{t=1}^n Y_t \Big) \Big| \\
		& \leq 2(k!)^{2+2 \alpha d-1} 12^{k-1} \big( H_1 \big)^k 2^{2 \alpha d k}
		\Big( 4 \ {\rm e} \ h_* \Big( 1 + \frac{1}{\log \rho}\Big) \Big)^{k-1} k!  \Big(1 + K^{1/2} \Big) n \\
		& = 12^{k-1} 2^{1+2 \alpha d k + 2+2 \alpha d + 2(k-1)} (k!/2)^{2+2 \alpha d} \Big( n H_1 \Big)^k	\Big( {\rm e} \ h_* \Big( 1 + \frac{1}{\log \rho}\Big) / n \Big)^{k-1} \Big(1 + K^{1/2} \Big) \\
		& = 12 \cdot 12^{k-2} 2^{2 (\alpha d + 1) k + 2 \alpha d + 1} (k!/2)^{2+2 \alpha d} \Big( n H_1 \Big)^k	\Big( {\rm e} \ h_* \Big( 1 + \frac{1}{\log \rho}\Big) / n \Big)^{k-1} \Big(1 + K^{1/2} \Big) \\
		& = \Big(\frac{k!}{2}\Big)^{1+\gamma} \frac{12 \cdot 2^{4 (\alpha d + 1) + 2 \alpha d + 1} \Big( n H_1 \Big)^2{\rm e} \ h_* \Big( 1 + \frac{1}{\log \rho}\Big)\Big(1 + K^{1/2} \Big) / n }{\Big( n 12^{-1} 2^{-2 (\alpha d + 1)} \big(n H_1\big)^{-1} \Big({\rm e} \ h_* \Big( 1 + \frac{1}{\log \rho}\Big) \Big)^{-1} \Big)^{k-2} }
		\leq \Big(\frac{k!}{2}\Big)^{1+\gamma} \frac{H}{\bar\Delta^{k-2}},
	\end{split}
\end{equation*}
where we have used the fact that
$n H_1 \leq 2^{2 \alpha d+1} A \alpha^{2 \alpha d} (\max\{c^2, 3 \pi M_f,1 \})^{\alpha}$.
Applying Lemma~\ref{lem:sausta91:Lemma24} yields the assertion.\hfill$\square$

\subsection{Proofs of the Lemmas~\ref{lem:stoerungslemma}--\ref{lem:sumApprox}}\label{sec:TRes:proof}

\textbf{Proof of Lemma~\ref{lem:stoerungslemma}.} Let $X := I - A^{-1} (A + \Delta) = - A^{-1} \Delta$. Then, because $\| X \|_M \leq \| A^{-1} \|_M \| \Delta \|_M < 1$, Lemma~2.1 in~\cite{dem97} can be applied, which asserts that $I - X = A^{-1} (A + \Delta)$ is invertible and that
\[\| (I - X)^{-1} \|_M \leq \frac{1}{1-\|X\|_M}.\]
We deduce that $A + \Delta$ is invertible, because we have
\[\det (A + \Delta) = \det(A) \det(I - X) \neq 0,\]
where we have used that $A$ is invertible by assumption and $I - X$ is invertible by Lemma~2.1 in~\cite{dem97}. Finally, the bound on the matrix norm of the inverse follows from the assumed submultiplicativity
\begin{equation*}
\begin{split}
	\| (A + \Delta)^{-1}\|_M
	& = \| (I - X)^{-1} A^{-1}\|_M
	\leq \| (I - X)^{-1} \|_M \cdot \| A^{-1} \|_M \\
	& \leq \frac{\| A^{-1} \|_M}{1-\|A^{-1} \Delta\|_M}
	\leq \frac{\| A^{-1} \|_M}{1-\|A^{-1} \|_M \|\Delta\|_M}.
\end{split}
\end{equation*}
This finishes the proof of the lemma.\hfill$\square$

\textbf{Proof of Lemma~\ref{lem:stoerungslemma2}.} Note that Lemma~\ref{lem:stoerungslemma} can be applied, which yields that $A + \Delta$ is invertible as an immediate consequence. Further, note that
\[(A + \Delta)^{-1} - A^{-1} = (A + \Delta)^{-1}\big( A - (A + \Delta) \big) A^{-1}.\]
Employing the submultiplicativity of the norm and the inequality from Lemma~\ref{lem:stoerungslemma} yields
\[ \| (A + \Delta)^{-1} - A^{-1} \|_M \leq \frac{\| A^{-1} \|_M}{1 - \| A^{-1} \|_M \cdot \| \Delta \|_M} \cdot \| -\Delta \|_M \cdot \| A^{-1} \|_M.\]
The  assertion then follows, because we have $\| -\Delta \|_M = \| \Delta \|_M$ and, by assumption, $\| A^{-1} \|_M \cdot \| \Delta \|_M \leq c$ holds.\hfill$\square$

\textbf{Proof of Lemma~\ref{lem:expansionMatrixPower}.} The statement for $h \in \{0,1\}$ is obvious. For $h = 2,3, \ldots$ note that by the binomial theorem we have
\begin{equation*}
\begin{split}
	 A^h = \big( (A - A_0) + A_0 \big)^h & = \sum_{\substack{(i_1, \ldots, i_h) \\ \in \{0,1\}^h}} \prod_{\ell = 1}^h (A-A_0)^{1-i_{\ell}} A_0^{i_{\ell}}.
\end{split}
\end{equation*}
This obviously implies
\begin{equation*}
\begin{split}
	\|A^h - A_0^h\|_M & = \Big\| \sum_{\substack{(i_1, \ldots, i_h) \\ \in \{0,1\}^h \setminus \{(1,\ldots,1)\}}} \prod_{\ell = 1}^h (A-A_0)^{1-i_{\ell}} A_0^{i_{\ell}} \Big\|_M \\
		& \leq \sum_{j=0}^{h-1} \binom{h}{j+1} \|A-A_0\|_M^{j+1} \|A_0\|_M^{h-j-1} \\
		& = \sum_{j=0}^{h-1} \frac{h}{j+1} \frac{(h-1)!}{j!(h-1-j)!} \|A-A_0\|_M^{j+1} \|A_0\|_M^{h-j-1} \\
		& \leq h \| A-A_0 \|_M \sum_{j=0}^{h-1} \binom{h-1}{j} \|A-A_0\|_M^{j} \|A_0\|_M^{h-1-j} \\
		& = h \| A-A_0 \|_M \big( \|A-A_0\|_M + \|A_0\|_M \big)^{h-1}. \\
\end{split}
\end{equation*}

\textbf{Proof of Lemma~\ref{lem:ineq_productOfTwo}.}
Note that
\begin{equation*}
\begin{split}
	& \IP( | u v - u_0 v_0 | > \varepsilon) 
	= \IP( | (u - u_0) (v - v_0) + u_0 (v - v_0) + (u - u_0) v_0 | > \varepsilon) \\
	& \leq \IP( |u - u_0| |v - v_0| + |u_0| |v - v_0| + |u - u_0| |v_0| > \varepsilon) \\
	& \leq \IP( \frac{1}{2} |u - u_0|^2 + \frac{1}{2} |v - v_0|^2
			+ |u_0| |v - v_0| + |v_0| |u - u_0| > \varepsilon) \\
	& \leq \IP( |u - u_0|^2 + 2 |v_0| |u - u_0| > \varepsilon)
			+ \IP( |v - v_0|^2 + 2 |u_0| |v - v_0| > \varepsilon), \\
\end{split}
\end{equation*}
where the first inequality is due to the triangle inequality and the second is due to Young's inequality. We now bound the first term and note that the second term can be handled analogously. Note that
\begin{equation}\label{lem:ineq_productOfTwo:arg}
\begin{split}
	\IP( |u - u_0|^2 + 2 |v_0| |u - u_0| > \varepsilon)
	& = \IP( |u - u_0| > \sqrt{\varepsilon + |v_0|^2} - |v_0|) \\
\end{split}
\end{equation}
The assertion then follows, because
\begin{equation}\label{lem:ineq_productOfTwo:arg2}
	(|v_0|^2+\varepsilon)^{1/2} - (|v_0|^2)^{1/2} = \frac{1}{2} \xi^{-1/2} \varepsilon \geq \frac{1}{2} (|v_0|^2 + \varepsilon)^{-1/2} \varepsilon.
\end{equation}
The equality in the above holds due to the mean value theorem, for some $\xi \in [|v_0|^2, |v_0|^2 + \varepsilon]$. The inequality is due to $\xi \leq |v_0|^2 + \varepsilon$.\hfill$\square$

\textbf{Proof of Lemma~\ref{lem:disentangle}}. Note that, denoting $\mu_t := \E X_t$, we have
\begin{equation*}
	\begin{split}
		& \IP \Big( \Big| \sum_{t=1}^n \big( \hat a_t X_t - \alpha_t \E(X_t) \big) \Big| > n \varepsilon \Big)
		= \IP \Big( \Big| \sum_{t=1}^n ( \hat a_t - \alpha_t) X_t
			+ \sum_{t=1}^n \alpha_t (X_t - \mu_t) \Big| > n \varepsilon \Big) \\
		& \leq \IP \Big( \Big| \sum_{t=1}^n ( \hat a_t - \alpha_t) X_t \Big|
			+ \Big| \sum_{t=1}^n \alpha_t (X_t - \mu_t) \Big| > n \varepsilon \Big) \\
		& \leq \IP \Big( \Big| \sum_{t=1}^n ( \hat a_t - \alpha_t)^2\Big|^{1/2} \ \Big|\sum_{t=1}^n X_t^2 \Big|^{1/2}
			+ \Big| \sum_{t=1}^n \alpha_t (X_t - \mu_t) \Big| > n \varepsilon \Big) \\
		& \leq \IP \Big( \Big| \sum_{t=1}^n ( \hat a_t - \alpha_t)^2 \Big|
			\Big| \sum_{t=1}^n ( X_t^2 - \E X_t^2 + \E X_t^2) \Big| > \big(\frac{n \varepsilon}{2}\big)^2 \Big)
			+ \IP \Big( \Big| \sum_{t=1}^n \alpha_t (X_t - \mu_t) \Big| > \frac{n \varepsilon}{2} \Big) \\
	\end{split}
\end{equation*}
Further, we have
\begin{equation*}
	\begin{split}
	& \IP \Big( \Big| \sum_{t=1}^n ( \hat a_t - \alpha_t)^2 \Big|
			\cdot \Big| \sum_{t=1}^n ( X_t^2 - \E X_t^2 + \E X_t^2) \Big| > (n \varepsilon/2)^2 \Big) \\
		& \leq \IP \Big( \Big| \sum_{t=1}^n ( \hat a_t - \alpha_t)^2 \Big|
			\cdot \Big| \sum_{t=1}^n ( X_t^2 - \E X_t^2) \Big| + \Big| \sum_{t=1}^n ( \hat a_t - \alpha_t)^2 \Big| \cdot \Big| \E  \sum_{t=1}^n X_t^2 \Big| > (n \varepsilon/2)^2 \Big) \\	
		& \leq \IP \Big( \frac{1}{2} \Big| \sum_{t=1}^n ( \hat a_t - \alpha_t)^2 \Big|^2
			+ \frac{1}{2} \Big| \sum_{t=1}^n ( X_t^2 - \E X_t^2) \Big|^2 + \Big| \sum_{t=1}^n ( \hat a_t - \alpha_t)^2 \Big| n m_2^2 > (n \varepsilon/2)^2 \Big) \\
		& \leq \IP \Big( \Big| \sum_{t=1}^n ( \hat a_t - \alpha_t)^2 \Big|^2
			 + 2 n m_2^2 \Big| \sum_{t=1}^n ( \hat a_t - \alpha_t)^2 \Big| > \frac{(n \varepsilon)^2}{4} \Big)
		+ \IP \Big( \Big| \sum_{t=1}^n ( X_t^2 - \E X_t^2) \Big| > \frac{n \varepsilon}{2} \Big) \\
		& \leq \IP \Big( \Big| \sum_{t=1}^n ( \hat a_t - \alpha_t)^2 \Big|
									> \sqrt{ (n \varepsilon/2)^2 + (n m_2^2)^2 } - n m_2^2 \Big)
		+ \IP \Big( \Big| \sum_{t=1}^n ( X_t^2 - \E X_t^2) \Big| > \frac{n \varepsilon}{2} \Big)
	\end{split}
\end{equation*}
where for the second inequality we employed Young's inequality and in the fourth inequality we have used the argument in~\eqref{lem:ineq_productOfTwo:arg} from the proof of Lemma~\ref{lem:ineq_productOfTwo}. The assertion then follows, because using~\eqref{lem:ineq_productOfTwo:arg2} from the proof of Lemma~\ref{lem:ineq_productOfTwo}, we have
\begin{equation*}
		\sqrt{ (n \varepsilon/2)^2 + (n m_2^2)^2 } - n m_2^2
		\geq \frac{1}{2} (n \varepsilon/2)^2 \big( (n m_2^2)^2 + (n \varepsilon/2)^2 \big)^{-1/2}
		= \frac{n}{2} 2 (\varepsilon/2)^2 \big( (2 m_2^2)^2 + \varepsilon^2 \big)^{-1/2},
\end{equation*}
which finishes the proof.\hfill$\square$

\textbf{Proof of Lemma~\ref{lem:ineqQuotientMV}.}
	For every $\delta \in (0, 1 / \| M_0^{-1} \|_M )$, we have
	\begin{equation*}
		\begin{split}
		& \IP \Big( \Big\| M^{-1} v - M_0^{-1} v_0  \Big\|_v > \varepsilon \Big) \\
		& = \IP \Big( \Big\| M^{-1} v - M_0^{-1} v_0 \Big\|_v > \varepsilon \ \Big| \ \| M - M_0 \|_M > \delta \Big) \IP ( \| M - M_0 \|_M > \delta) \\
		& \quad + \IP \Big( \Big\| M^{-1} v - M_0^{-1} v_0  \Big\|_v > \varepsilon, \ \| M - M_0 \|_M \leq \delta \Big) \\
		& \stackrel{(*)}{\leq} \IP ( \| M - M_0 \|_M > \delta) \\
		& \quad + \IP \Big( \|v-v_0\|_v > \frac{\varepsilon}{2} \frac{1-\| M_0^{-1} \|_M \, \delta}{\|M_0^{-1}\|_M} \Big) \\
		& \quad + \IP \Big( \| M - M_0 \|_M \| v_0 \|_v > \frac{\varepsilon}{2} \frac{1-\| M_0^{-1} \|_M \, \delta}{(\|M_0^{-1}\|_M)^2} \Big), \\
	\end{split}
	\end{equation*}
where for $(*)$ we have used the fact that any conditional probability is $\leq 1$ for the first part and the following argument for the second part:
	\begin{equation*}
		\begin{split}
		& \IP \Big( \Big\| M^{-1} v - M_0^{-1} v_0  \Big\|_v > \varepsilon, \ \| M - M_0 \|_M \leq \delta \Big) \\
		& = \IP \Big( \Big\| M^{-1} v - M^{-1} v_0 + M^{-1} v_0 - M_0^{-1} v_0  \Big\|_v > \varepsilon, \ \| M - M_0 \|_M \leq \delta \Big) \\
		& \leq \IP \Big( \|M^{-1}\|_M \|v-v_0\|_v + \| M^{-1} \|_M \, \| M_0 - M \|_M \, \| M_0^{-1} \|_M \| v_0 \|_v > \varepsilon, \ \| M - M_0 \|_M \leq \delta \Big) \\
		& \leq \IP \Big( \frac{\|M_0^{-1}\|_M}{1-\| M_0^{-1} \|_M \, \delta} \|v-v_0\|_v + \frac{(\|M_0^{-1}\|_M)^2}{1-\| M_0^{-1} \|_M \, \delta} \| M - M_0 \|_M \| v_0 \|_v > \varepsilon, \ \| M - M_0 \|_M \leq \delta \Big) \\
		& \leq \IP \Big( \|v-v_0\|_v > \frac{\varepsilon}{2} \frac{1-\| M_0^{-1} \|_M \, \delta}{\|M_0^{-1}\|_M}, \ \| M - M_0 \|_M \leq \delta \Big) \\
		& \quad + \IP \Big( \| M - M_0 \|_M \| v_0 \|_v > \frac{\varepsilon}{2} \frac{1-\| M_0^{-1} \|_M \, \delta}{(\|M_0^{-1}\|_M)^2}, \ \| M - M_0 \|_M \leq \delta \Big) \\
	\end{split}
	\end{equation*}

In the above inequalities we have used Lemma~\ref{lem:stoerungslemma}, by which we have:
\[\| M^{-1} \|_M
		= \| (M_0 + M - M_0)^{-1} \|_M
		\leq \frac{\|M_0^{-1}\|_M}{1-\| M_0^{-1} \|_M \, \| M - M_0 \|_M}
		\leq \frac{\|M_0^{-1}\|_M}{1-\| M_0^{-1} \|_M \, \delta} \]
The assertion then follows by choosing $\delta = \frac{1}{2 \|M_0^{-1}\|_M}$.\hfill$\square$

\textbf{Proof of Lemma~\ref{lem:PowerYW}.}
Note that, for $h=1$ we have $v = x$ and $v_0 = x_0$. Thus, the assertion holds due to $\varepsilon/(\varepsilon + 1) \leq \varepsilon$. For $h = 2, 3, \ldots$, we have
\[\| v - v_0\| = \| (1, 0, \ldots, 0) (A^h - A_0^h) \| \leq \| A^h - A_0^h \|,\]

Thus, with this and by Lemma~\ref{lem:expansionMatrixPower}, we have
\begin{equation*}
\begin{split}
	& \IP( \| v - v_0\| > \varepsilon) \leq \IP( \|A^h - A_0^h\| > \varepsilon) \\
	& \leq \IP( \| A-A_0 \| \big( \|A-A_0\| + \|A_0\| \big)^{h-1} > \varepsilon/h) \\
	& \leq \IP( \| x-x_0 \| \big( \|x-x_0\| + \|x_0\| \big)^{h-1} > \varepsilon/h)
	=: (*),
	\end{split}
\end{equation*}
where the second inequality uses $\| A-A_0 \| \leq \| A-A_0 \|_F = \|x-x_0\|$ and the fact that $\|A_0\| = \| e_1 x'_0 \| \leq \|x_0\|$. For $\|A_0\| = \| e_1 x'_0 \|$ note that subtracting the Jordan block doesn't change the eigenvalues (see also the argument involving Weyl's inequalty regarding the norm of $X$ in the proof of Lemma~\ref{lem:norm_v_der}). For the rest of the derivation denote $y := \| x - x_0 \|^{1/(h-1)}$. Then, we have
\begin{equation}\label{eqn:polyn}
	(*) = \IP( y^{h} + \|x_0\| y - (\varepsilon/h)^{1/(h-1)} > 0)
	\leq \IP\Big( y > y_{\min} \Big),
\end{equation}
where $y_{\min} := \min\{ |y_j| \, : \, j=1,\ldots,h\}$ is the minimum of the absolute values of the (complex) roots $y_1, \ldots, y_h$ of the polynomial function $p(y) := y^{h} + \|x_0\| y - (\varepsilon/h)^{1/(h-1)}$. The inequality in~\eqref{eqn:polyn} holds due to the following argument: Due to Descartes' rule of sign, the polynomial $p(y)$ has exactly one positive root, say $y_0$. Since $p(0) < 0$ we have by the intermediate value theorem that if $p(y) > 0$ then $y > y_0 \geq y_{\min}$, which implies~\eqref{eqn:polyn}.

Then, since by a standard argument involving Rouch\'{e}'s theorem, we have
\[y_{\min} \geq \frac{(\varepsilon/h)^{1/(h-1)}}{(\varepsilon/h)^{1/(h-1)} + \max\{\|x_0\|, 1\}},\]
we have
\begin{equation*}
\begin{split}
	\IP( \| v - v_0\| > \varepsilon)
	&	\leq \IP\Big( \| x - x_0\|^{1/(h-1)} > \frac{(\varepsilon/h)^{1/(h-1)}}{(\varepsilon/h)^{1/(h-1)} + \max\{\|x_0\|, 1\}} \Big) \\
	&	= \IP\Big( \| x - x_0\| > \frac{\varepsilon/h}{\big( (\varepsilon/h)^{1/(h-1)} + \max\{\|x_0\|, 1\} \big)^{h-1}} \Big) \\
	&	= \IP\Big( \| x - x_0\| > \frac{\varepsilon}{\big( \varepsilon^{1/(h-1)} + h^{1/(h-1)} \max\{\|x_0\|, 1\} \big)^{h-1}} \Big) \\
	&	\leq \IP\Big( \| x - x_0\| > \frac{\varepsilon}{ \max\{2^{h-2},1\} \big( \varepsilon + h \max\{\|x_0\|, 1\}^{h-1} \big)} \Big) \\
	\end{split}
\end{equation*}

Thus, we have shown
\begin{equation*}
	\IP( \| v - v_0\| > \varepsilon) \leq
	\begin{cases}
	\IP\Big( \| x - x_0\| > \varepsilon \Big)& \text{if $h=1$,} \\
	\IP\Big( \| x - x_0\| > 2^{2-h} \cfrac{\varepsilon}{ \varepsilon + h \max\{\|x_0\|, 1\}^{h-1}} \Big) & \text{if $h = 2, 3, \ldots$}.
	\end{cases}
\end{equation*}

We report a bound that is larger, to have a compact expression that is valid for all $h$.\hfill$\square$

\textbf{Proof of Lemma~\ref{lem:MomentCond}.} It is well know that (see, for example, Exercise 2.3.25(i) in \cite{str11}), if $X$ is a $b$-sub-Gaussian random variable, then
\[\E |X|^p \leq 2^{p/2+1} \Gamma(p/2+1) b^p,\]
for every $p \in (0,\infty)$. Thus, we have
\begin{equation*}
	\begin{split}
		\E |X_{t,T}|^k
			&	= \sigma_{t,T}^{k-2} \E \Big|\frac{X_{t,T}}{\sigma_{t,T}}\Big|^k \sigma_{t,T}^2
			\leq (3 \pi M_f)^{(k-2)/2} 2^{k/2+1} \Gamma\big( k / 2 + 1  \big) b^k \sigma_{t,T}^2 \\
			& \leq (3 \pi M_f)^{k} 2^{k} (k!)^{1/2} b^k \sigma_{t,T}^2 = (6 \pi b M_f)^{k} (k!)^{1/2} \sigma_{t,T}^2, \\
	\end{split}
\end{equation*}
where we have used~\eqref{eqn:cons:sd_bounded} for the first inequality and for the second inequality assumed that without loss of generality $3\pi M_f \geq 1$ and we have used that, for $p \in \IN^*$, we have
\[
	\Gamma(p/2 + 1) = 
	\begin{cases}
		(p/2)! & \text{ if $p$ is even}, \\
		\sqrt{\pi} \prod_{j=0}^{\lfloor p/2 \rfloor} \frac{p-2 j}{2} & \text{ if $p$ is odd},
	\end{cases}
\]
which implies $\Gamma(p/2 + 1)^2 \leq p!$, for $p \geq 3$, which finishes the proof, as $\E X_{t,T} = 0$.\hfill$\square$

\textbf{Proof of Lemma~\ref{lem:sumApprox}.}

By the mean value theorem and the assumed conditions we have that for every $\ell \in \IZ$ there exists $\xi_{\ell} \in [(\ell-1)/T, \ell/T]$ and $\zeta_{\ell} \in (\xi_{\ell}, \ell/T)$ such that
\begin{equation*}
\begin{split}
	\int_{(\ell-1)/T}^{\ell/T} f(u) {\rm d}u
	& = f(\xi_{\ell}) (\frac{\ell}{T} - \frac{\ell-1}{T})
	= f(\ell/T) \frac{1}{T} + (f(\xi_\ell) - f(\ell/T)) \frac{1}{T} \\
	& = f(\ell/T) \frac{1}{T} + f'(\zeta_\ell) (\xi_\ell - \ell/T) \frac{1}{T}
\end{split}
\end{equation*}

Thus, for $A, B = 0, \ldots, T$ we have
\begin{equation}\label{lem:sumApprox:eqn:A}
\begin{split}
	\int_{A/T}^{B/T} f(u) {\rm d}u
	= \sum_{\ell=A+1}^B
	\int_{(\ell-1)/T}^{\ell/T} f(u) {\rm d}u
	= \frac{1}{T} \sum_{\ell=A+1}^B f(\ell/T) + R_T(A,B)
\end{split}
\end{equation}
where
\begin{equation}\label{lem:sumApprox:eqn:B}
\begin{split}
	|R_T(A,B)| & = \Big| \sum_{\ell=A+1}^B f'(\zeta_\ell) (\xi_\ell - \ell/T) \frac{1}{T} \Big| 
	\leq \frac{B-A}{T^2} \sup_{A/T < u < B/T} | f'(u) | . \\
\end{split}
\end{equation}
In the above we have used that $|\ell/T-\xi_{\ell}| \leq 1/T$, if $\xi_{\ell} \in [(\ell-1)/T, \ell/T]$.

Integration by substitution yields that
\begin{equation}\label{lem:sumApprox:eqn:C}
\begin{split}
	\int_{A/T}^{B/T} f(u) {\rm d}u
	= \frac{B-A}{T} \int_0^1 f\big( \frac{A}{T} + \frac{B-A}{T}u \big) {\rm d}u.
\end{split}
\end{equation}

Finally, substitute~\eqref{lem:sumApprox:eqn:C} into~\eqref{lem:sumApprox:eqn:A} and multiply by $T/(B-A)$, then the assertion follows from~\eqref{lem:sumApprox:eqn:B}.\hfill$\square$

\clearpage
\section[Results from Saulis and Statulevicus (1991)]{Results from~\cite{sausta91}}\label{sec:res_sausta91}

In this section we cite two results of \cite{sausta91} that are needed for the proofs in the previous sections.
\begin{lem}[Lemma~2.4, \cite{sausta91}]
	\label{lem:sausta91:Lemma24}
	For an arbitrary ran\-dom variable $\xi$ with $\E \xi = 0$, let there exist $\gamma \geq 0$, $H > 0$ and $\bar\Delta > 0$ such that
	\begin{equation}
		\label{lem:sausta91:Lemma24:cond}
		|\cum_k(\xi)| \leq \Big(\frac{k!}{2}\Big)^{1+\gamma} \frac{H}{\bar\Delta^{k-2}}, \quad k=2,3,\ldots.
	\end{equation}
	Then, for all $x \geq 0$,\footnote{Note that in~\cite{sausta91} there is a typo!}
	\begin{equation*}
	\begin{split}
	\IP(|\xi| \geq x) &
	\leq \exp\Big( - \frac{x^2}{2(H + (x/\bar\Delta^{1/(1+2\gamma)})^{(1+2\gamma)/(1+\gamma)})} \Big) \\
	& \leq \begin{cases}
		\exp\Big( -\frac{x^2}{4H} \Big) & 0 \leq x \leq (H^{1+\gamma} \bar\Delta)^{1/(1+2\gamma)} \\ 
		\exp\Big( -\frac{1}{4} (x \bar\Delta)^{1/(1+\gamma)} \Big) & x \geq (H^{1+\gamma} \bar\Delta)^{1/(1+2\gamma)}
	\end{cases}
	\end{split}
	\end{equation*}
\end{lem}

\begin{lem}[Theorem~4.17, \cite{sausta91}]
	\label{lem:sausta91:Theorem417}
	Let $X_t$, $t=1,2,\ldots$ be a random process on a probability space $(\Omega, \mathcal{F}, \IP)$ and define the strong mixing coefficients by
	\[\alpha^Y(s,t) := \sup_{A \in \sigma(X_{u} : 1 \leq u \leq s)} \sup_{B \in \sigma(X_{u} : u\geq t+k)} \big| \IP(A \cap B) - \IP(A) \IP(B) \big|.\]
	[notation taken from \cite{sausta91}, p.\,60]
	If for some $\gamma_1 \geq 0$ and $H_1 > 0$
	\begin{equation}
		\label{lem:sausta91:Theorem417:cond}
		\E |Y_t|^k \leq (k!)^{1+\gamma_1} H_1^k, \quad t=1, \ldots, n, \ k=2,3, \ldots,
	\end{equation}
	then for all $\delta > 0$
	\begin{equation*}
		\Big| \cum_k \Big( \sum_{t=1}^n Y_t \Big) \Big|
		\leq 2(k!)^{2+\gamma_1} 12^{k-1} H_1^k (1+\delta)^{(1+\gamma_1)k}
		\Delta_n((1+1/\delta)(k-1))^{k-1} n,
	\end{equation*}
	where
	\[\Delta_n(\nu) := \max \Big\{ 1, \max_{1 \leq s \leq n} \sum_{t=s}^n \alpha^Y(s,t)^{1/\nu} \Big\}.\]
\end{lem}

\clearpage

\section{Further output for the analysis in Section~\ref{data_example_brexit}}\label{a:extraEmpExpl}

\begin{table}[h] \scriptsize
\begin{center}
\begin{tabular}{|c||c|c|c|c|c|}
\hline
 $h$ & $\bar p_{\stat}$ & ${\rm trMAPE}_{T,1}^{\stat}(h)$ & $\bar p_{\lstat}$ & $\bar N_{\lstat}$ & ${\rm trMAPE}_{T,1}^{\lstat}(h)$ \\
\hline \hline
1 &   8 & 2.383838e-05  &   8  &  50 & 1.022455e-05 \\
2 &   8 & 2.903981e-05  &   8  &  49 & 1.044652e-05 \\
3 &   8 & 3.285705e-05  &   5  &  51 & 9.672596e-06 \\
4 &   8 & 4.172401e-05  &   8  &  54 & 1.036038e-05 \\
5 &   8 & 4.433286e-05  &   8  &  50 & 1.362222e-05 \\

\hline
\end{tabular}
\hspace*{-0.25cm}\begin{tabular}{|c||c|c|c||c|c|c|}
\hline
 $h$ & ${\rm trMAPE}_{T,2}^{\stat}(h)$ & ${\rm trMAPE}_{T,2}^{\lstat}(h)$ & $\bar R_{T,2}(h)$ & ${\rm trMAPE}_{T,3}^{\stat}(h)$ & ${\rm trMAPE}_{T,3}^{\lstat}(h)$ & $\bar R_{T,3}(h)$\\
\hline \hline
1 & 7.241481e-05 & 6.797737e-05 & 1.065 & 3.344593e-05 & 2.686727e-05 & 1.245 \\
2 & 6.785736e-05 & 6.202493e-05 & 1.094 & 3.703221e-05 & 2.704054e-05 & 1.370 \\
3 & 6.484127e-05 & 4.875103e-05 & 1.330 & 4.158508e-05 & 2.692850e-05 & 1.544 \\
4 & 6.798448e-05 & 5.305174e-05 & 1.281 & 4.630539e-05 & 2.780716e-05 & 1.665 \\
5 & 7.546445e-05 & 6.995483e-05 & 1.079 & 4.901377e-05 & 4.419283e-05 & 1.109 \\

\hline
\end{tabular}
\caption{\textit{Minimum empirical trimmed mean absolute prediction errors (trMAPE) for $h$-step ahead prediction, $h=1,2,3,4,5$, of the squared and centred FTSE~100 data. Analysis performed with $m := 15$ and $p_{\max} = 8$. Top table shows values computed on the first validation set. Bottom table shows values computed on the second validation set and on the test set.}}  \label{ftse_numbers_m14_p8}
\end{center}
\end{table}

\begin{table}[h] \scriptsize
\begin{center}
\begin{tabular}{|c||c|c|c|c|c|}
\hline
 $h$ & $\bar p_{\stat}$ & ${\rm trMAPE}_{T,1}^{\stat}(h)$ & $\bar p_{\lstat}$ & $\bar N_{\lstat}$ & ${\rm trMAPE}_{T,1}^{\lstat}(h)$ \\
\hline \hline
1 &    8 & 3.336108e-05 &    4 &   45 & 7.131298e-06 \\
2 &    8 & 4.070901e-05 &    5 &   44 & 7.368486e-06 \\
3 &    8 & 4.576347e-05 &    6 &   40 & 8.365509e-06 \\
4 &    8 & 5.561484e-05 &    6 &   40 & 8.420542e-06 \\
5 &    8 & 5.970097e-05 &    8 &   43 & 1.131462e-05 \\

\hline
\end{tabular}
\hspace*{-0.25cm}\begin{tabular}{|c||c|c|c||c|c|c|}
\hline
 $h$ & ${\rm trMAPE}_{T,2}^{\stat}(h)$ & ${\rm trMAPE}_{T,2}^{\lstat}(h)$ & $\bar R_{T,2}(h)$ & ${\rm trMAPE}_{T,3}^{\stat}(h)$ & ${\rm trMAPE}_{T,3}^{\lstat}(h)$ & $\bar R_{T,3}(h)$\\
\hline \hline
1 & 2.571029e-05 & 1.459625e-05 & 1.761 & 3.811591e-05 & 2.966141e-05 & 1.285 \\
2 & 2.997352e-05 & 1.283254e-05 & 2.336 & 4.169124e-05 & 2.962695e-05 & 1.407 \\
3 & 3.426832e-05 & 1.207928e-05 & 2.837 & 4.697350e-05 & 3.324401e-05 & 1.413 \\
4 & 4.366208e-05 & 1.329643e-05 & 3.284 & 5.145149e-05 & 3.032581e-05 & 1.697 \\
5 & 4.678991e-05 & 1.527034e-05 & 3.064 & 5.664217e-05 & 5.630906e-05 & 1.006 \\

\hline
\end{tabular}
\caption{\textit{Minimum empirical trimmed mean absolute prediction errors (trMAPE) for $h$-step ahead prediction, $h=1,2,3,4,5$, of the squared and centred FTSE~100 data. Analysis performed with $m := 25$ and $p_{\max} = 8$. Top table shows values computed on the first validation set. Bottom table shows values computed on the second validation set and on the test set.}}  \label{ftse_numbers_m16_p8}
\end{center}
\end{table}

\begin{table}[t] \scriptsize
\begin{center}
\begin{tabular}{|c||c|c|c|c|c|}
\hline
 $h$ & $\bar p_{\stat}$ & ${\rm trMAPE}_{T,1}^{\stat}(h)$ & $\bar p_{\lstat}$ & $\bar N_{\lstat}$ & ${\rm trMAPE}_{T,1}^{\lstat}(h)$ \\
\hline \hline
1 &    4 & 3.143868e-05 &    4 &   59 & 1.510606e-05 \\
2 &    4 & 3.871644e-05 &    4 &   62 & 1.516250e-05 \\
3 &    4 & 4.518262e-05 &    4 &   58 & 1.430857e-05 \\
4 &    4 & 5.634271e-05 &    4 &   56 & 1.499452e-05 \\
5 &    4 & 6.160170e-05 &    4 &   52 & 1.989825e-05 \\
\hline
\end{tabular}
\hspace*{-0.25cm}\begin{tabular}{|c||c|c|c||c|c|c|}
\hline
 $h$ & ${\rm trMAPE}_{T,2}^{\stat}(h)$ & ${\rm trMAPE}_{T,2}^{\lstat}(h)$ & $\bar R_{T,2}(h)$ & ${\rm trMAPE}_{T,3}^{\stat}(h)$ & ${\rm trMAPE}_{T,3}^{\lstat}(h)$ & $\bar R_{T,3}(h)$\\
\hline \hline
1 & 5.122471e-05 & 2.225700e-05 & 2.302 & 3.273401e-05 & 2.402349e-05 & 1.363 \\
2 & 4.983237e-05 & 1.824545e-05 & 2.731 & 3.820289e-05 & 2.312136e-05 & 1.652 \\
3 & 5.130980e-05 & 1.426141e-05 & 3.598 & 4.443704e-05 & 2.535910e-05 & 1.752 \\
4 & 5.863523e-05 & 1.563907e-05 & 3.749 & 5.138214e-05 & 2.594492e-05 & 1.980 \\
5 & 6.549139e-05 & 1.959563e-05 & 3.342 & 5.460377e-05 & 4.529423e-05 & 1.206 \\
\hline
\end{tabular}
\caption{\textit{Minimum empirical trimmed mean absolute prediction errors (trMAPE) for $h$-step ahead prediction, $h=1,2,3,4,5$, of the squared and centred FTSE~100 data. Analysis performed with $m := 20$ and $p_{\max} = 5$. Top table shows values computed on the first validation set. Bottom table shows values computed on the second validation set and on the test set.}}  \label{ftse_numbers_m15_p6}
\end{center}
\end{table}

\begin{table}[t] \scriptsize
\begin{center}
\begin{tabular}{|c||c|c|c|c|c|}
\hline
 $h$ & $\bar p_{\stat}$ & ${\rm trMAPE}_{T,1}^{\stat}(h)$ & $\bar p_{\lstat}$ & $\bar N_{\lstat}$ & ${\rm trMAPE}_{T,1}^{\lstat}(h)$ \\
\hline \hline
1 &    8 & 2.827214e-05 &    4 &   59 & 1.510606e-05 \\
2 &    8 & 3.418102e-05 &   15 &  122 & 1.386033e-05 \\
3 &    8 & 3.864440e-05 &    5 &   63 & 1.294023e-05 \\
4 &    8 & 4.753937e-05 &   15 &   74 & 1.370509e-05 \\
5 &    8 & 5.128567e-05 &   15 &   65 & 1.484254e-05 \\
\hline
\end{tabular}
\hspace*{-0.25cm}\begin{tabular}{|c||c|c|c||c|c|c|}
\hline
 $h$ & ${\rm trMAPE}_{T,2}^{\stat}(h)$ & ${\rm trMAPE}_{T,2}^{\lstat}(h)$ & $\bar R_{T,2}(h)$ & ${\rm trMAPE}_{T,3}^{\stat}(h)$ & ${\rm trMAPE}_{T,3}^{\lstat}(h)$ & $\bar R_{T,3}(h)$\\
\hline \hline
1 & 4.589402e-05 & 2.225700e-05 & 2.062 & 3.323796e-05 & 2.402349e-05 & 1.384 \\
2 & 4.407642e-05 & 1.556693e-05 & 2.831 & 3.807648e-05 & 2.773136e-05 & 1.373 \\
3 & 4.403332e-05 & 1.334418e-05 & 3.300 & 4.357933e-05 & 2.648961e-05 & 1.645 \\
4 & 4.998899e-05 & 1.253916e-05 & 3.987 & 4.962460e-05 & 3.145873e-05 & 1.577 \\
5 & 5.496721e-05 & 1.354236e-05 & 4.059 & 5.377520e-05 & 5.260044e-05 & 1.022 \\
\hline
\end{tabular}
\caption{\textit{Minimum empirical trimmed mean absolute prediction errors (trMAPE) for $h$-step ahead prediction, $h=1,2,3,4,5$, of the squared and centred FTSE~100 data. Analysis performed with $m := 20$ and $p_{\max} = 15$. Top table shows values computed on the first validation set. Bottom table shows values computed on the second validation set and on the test set.}}  \label{ftse_numbers_m15_p10}
\end{center}
\end{table}

\clearpage

\section{Further simulation results for Section~\ref{simulationsection1}}\label{a:extraSim}

Now we define the models.
The first two tvAR(1) models (the innovations $Z_t$ are i.\,i.\,d Gaussian white noise, in all 15 models) are defined by two periodic coefficient functions, namely the models are
\begin{align}
\nonumber
X_{t,T} & = (0.8+0.19\sin(4\pi\frac{t}{T}))X_{t-1,T}+Z_t, \\
\nonumber
X_{t,T} & = (0.3+0.19\sin(4\pi\frac{t}{T}))   X_{t-1,T}+Z_t.
\end{align}

We then look at six tvAR(1) models where the AR coefficient increases linearly, namely\\
\begin{minipage}[t]{.5\textwidth}
\vspace*{-.5cm}
\begin{align}
\label{increasing2}
X_{t,T} & = (0.5+0.19\frac{t}{T})   X_{t-1,T}+Z_t, \\
\label{increasing4}
X_{t,T} & = (0.5+0.09\frac{t}{T})   X_{t-1,T}+Z_t, \\
\label{increasing1}
X_{t,T} & = (0.8+0.19\frac{t}{T}) X_{t-1,T}+Z_t,
\end{align}
\end{minipage}%
\begin{minipage}[t]{.5\textwidth}
\vspace*{-.5cm}
\begin{align}
\label{increasing3}
X_{t,T} & = (0.9+0.09\frac{t}{T})   X_{t-1,T}+Z_t, \\
\label{increasing5}
X_{t,T} & = (0.5+0.49\frac{t}{T})   X_{t-1,T}+Z_t, \\
\label{increasing6}
X_{t,T} & = (0.5+0.4\frac{t}{T})   X_{t-1,T}+Z_t. 
\end{align}
\end{minipage}

Finally, we consider the stationary AR(1) model with coefficient $-0.6$ and two models with independent observations, one with and one without heteroscedasticity, namely\\
\begin{minipage}{.5\textwidth}
\vspace*{-.5cm}
\begin{align} 
X_{t,T} & = -0.6 X_{t-1,T}+Z_t, \label{stationaryAR}
\end{align}
\end{minipage}%
\begin{minipage}{.5\textwidth}
\vspace*{-.5cm}
\begin{align}
X_{t,T} & = Z_t, \label{indepNonHetero} \\
X_{t,T} & = \big(5-16 |\frac{t}{T} - 0.5 |^2 \big) Z_t,  \label{indepHetero}
\end{align}
\end{minipage}

We further consider the following tvAR(2) model from~\cite{dah97}
\begin{equation} \label{mdl11}
X_{t,T} = 1.8 \cos(1.5 - \cos(4 \pi \frac{t}{T})) X_{t-1,T} - 0.81 X_{t-2,T} + Z_t
\end{equation}
and one of its tangent processes (cf. the comment after Assumption~\ref{a:loc_stat})
\begin{equation} \label{mdl12}
X_{t,T} = X_{t-1,T} - 0.81 X_{t-2,T} + Z_t.
\end{equation}

Further, we consider two tvAR(1) models where the coefficient decreases linearly:
\begin{align} \label{decreasing1}
X_{t,T} & = (0.99-0.49t/T)X_{t-1,T}+Z_t \\
\label{decreasing2}
X_{t,T} & = (0.5-t/T)   X_{t-1,T}+Z_t.
\end{align}

Note the organisation of the tables and figures. For each model we have three pages with one figure, in which eight plots are displayed that showing the ratio of median performances, and four tables.

We now, summarise our findings.
For the different models we observe that for $n=100$ the best stationary predictor usually (in all models, except the strongly non-stationary~\eqref{decreasing1} and~\eqref{decreasing2}) outperforms the locally stationary approach, which can be seen from the figures where the lines are below 1 and also in the proportions from the tables. In 11 of the 15 models the 1-step ahead stationary predictor yields the better performance in more than $65\%$ of the cases. Even in the highly non-stationary models \eqref{periodic1}, \eqref{increasing5}, \eqref{decreasing1}, and~\eqref{decreasing2} where the tangent processes for the observations to be predicted are close to the unit root it is still more than $47\%$ in which the stationary approach yields the better performance. As a general conclusion we thus see that locally stationary modelling might not be advantageous for forecasting short time series. Further, for those models, when $n=1000$, we can see that the locally stationary approach will result in a significantly reduced MSPE as the sample size increases. For example, with model~\eqref{increasing5} the 1-step ahead locally stationary approach outperforms the stationary approach by about $20\%$ on average and in almost $92\%$ of the cases.

We further observe that in all of the figures the plots in the bottom row, corresponding to $n=4000, 6000, 8000$ and $10000$ are of a generally similar shape (in~\eqref{mdl11} we see this only for $n \geq 6000$). For an explanation note that, as $n$ grows larger, $m/n$ tends to zero and consequently the data in the validation and test set behave more similarly. For $n=100$, for example, the validation and test set make up $24\%$ of the data, while for $n=10000$ it is only $12.54\%$. We observe two different kinds of shapes: (a) for some of the models, the locally stationary and the stationary approach appear to be performing similarly well as $h$ increases, and (b) in other models we find that the locally stationary approach outperforms the stationary approach even more as $h$ increases. More precisely, for
models~\eqref{periodic1}--\eqref{increasing4}, \eqref{mdl11}, \eqref{decreasing1}, and \eqref{decreasing2},
we see that the locally stationary approach excels for small $h$ and both approaches are almost equally good as $h$ gets larger, while for
models~\eqref{stationaryAR}--\eqref{indepHetero}
the stationary approach is advantageous for smaller $h$ and the two approaches are of similar performance for the larger $h$'s. In model~\eqref{mdl12} the stationary procedure is always better (uniformly with respect to $h$), but as $T$ gets larger the advantage is less pronounced. Note that 
models~\eqref{stationaryAR}, \eqref{indepNonHetero}, and \eqref{mdl12}
are stationary and~\eqref{indepNonHetero} and \eqref{indepHetero} are the independent observations. It is thus not very surprising that the stationary approach outperforms the locally stationary one. Finally, in
models~\eqref{increasing1}--\eqref{increasing6},
where the tangent processes for the observations from the validation and test set are close to the unit root, we see that as $h$ gets larger the ratio of MSPEs first increases and then (slowly) decreases. These observations corroborate the rules of thumbs we have given before.

\clearpage

\begin{figure}
\centering 
\includegraphics[width=0.24\textwidth]{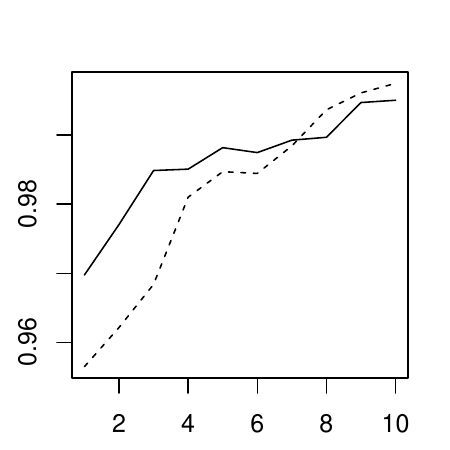}~
\includegraphics[width=0.24\textwidth]{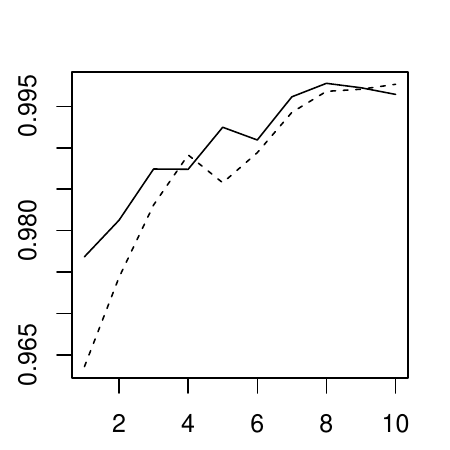}~
\includegraphics[width=0.24\textwidth]{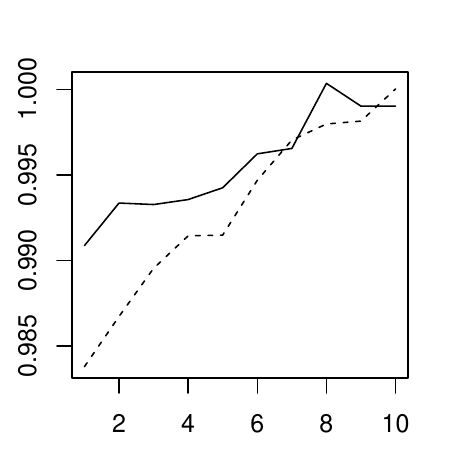}~ 
\includegraphics[width=0.24\textwidth]{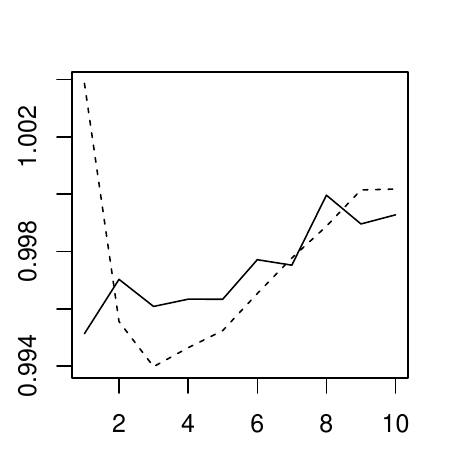} \\
\includegraphics[width=0.24\textwidth]{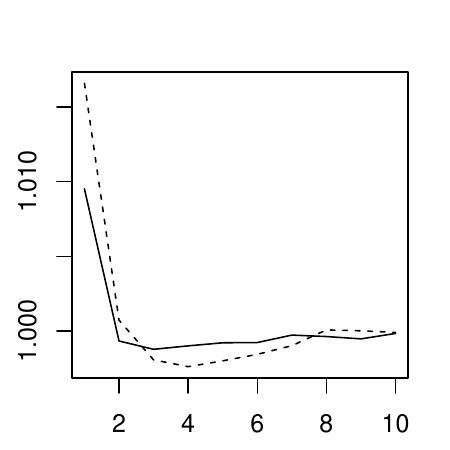}~
\includegraphics[width=0.24\textwidth]{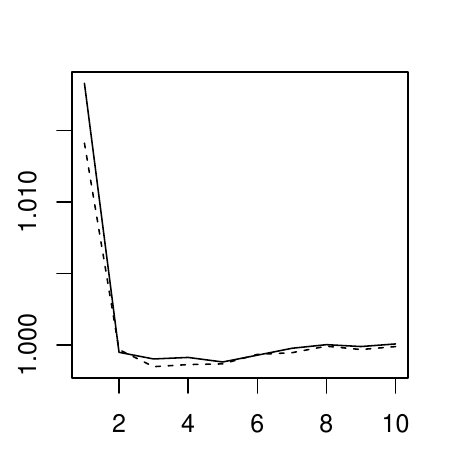}~ 
\includegraphics[width=0.24\textwidth]{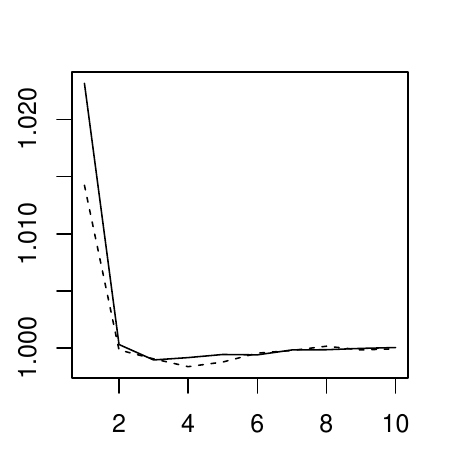}~ 
\includegraphics[width=0.24\textwidth]{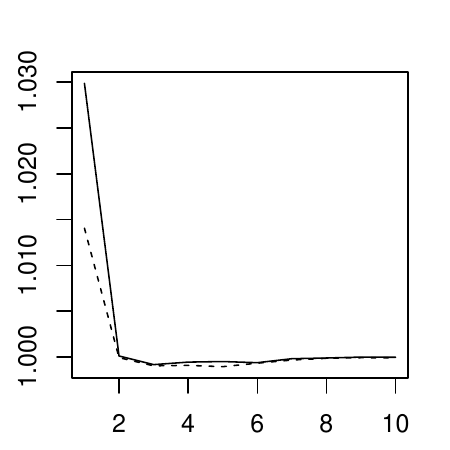}
   \caption{\it Plot of $h \mapsto R_{T,i}(h)$ for model \eqref{periodic2} and different values of $n$ (from left to right: n=100, n=200, n=500, n=1000 [first row], n=4000, n=6000, n=8000, n=10000 [second row]). Solid line: $i=3$ (test set), dashed line: $i=2$ (validation set 2). }  \label{MSPEperiodic2}
\end{figure}

\begin{table} \scriptsize
\begin{center}
\begin{tabular}{|c||c|c|c||c|c|c|}
\hline
 \multirow{2}{*}{$n$} & & \multicolumn{2}{|c||}{$R^{{\rm s}, {\rm ls}}_{T,2}(1)$}  & & \multicolumn{2}{|c|}{$R^{{\rm s}, {\rm ls}}_{T,2}(5)$} \\ 
 & & $\geq 1.01$  & $< 1.01$  & & $\geq 1.01$  & $< 1.01$ \\ 
 \hline 
 \multirow{2}{*}{$100$} & $R^{{\rm s}, {\rm ls}}_{T,3}(1) \geq 1.01$  & 0.0927 & 0.2075  & $R^{{\rm s}, {\rm ls}}_{T,3}(5) \geq 1.01$  & 0.0789 & 0.1455 \\ 
& $R^{{\rm s}, {\rm ls}}_{T,3}(1) < 1.01$  & 0.1651 & 0.5347 & $R^{{\rm s}, {\rm ls}}_{T,3}(5) < 1.01$  & 0.1384 & 0.6372 \\

\hline\hline
 \multirow{2}{*}{$n$} & & \multicolumn{2}{|c||}{$R^{{\rm s}, {\rm ls}}_{T,2}(1)$}  & & \multicolumn{2}{|c|}{$R^{{\rm s}, {\rm ls}}_{T,2}(5)$} \\ 
 & & $\geq 1.2$  & $< 1.2$  & & $\geq 1.2$  & $< 1.2$ \\ 
 \hline 
 \multirow{2}{*}{$1000$} & $R^{{\rm s}, {\rm ls}}_{T,3}(1) \geq 1.2$  & 0 & 2e-04  & $R^{{\rm s}, {\rm ls}}_{T,3}(5) \geq 1.2$  & 0 & 0 \\ 
& $R^{{\rm s}, {\rm ls}}_{T,3}(1) < 1.2$  & 5e-04 & 0.9993 & $R^{{\rm s}, {\rm ls}}_{T,3}(5) < 1.2$  & 1e-04 & 0.9999 \\

\hline\hline
 \multirow{2}{*}{$n$} & & \multicolumn{2}{|c||}{$R^{{\rm s}, {\rm ls}}_{T,2}(1)$}  & & \multicolumn{2}{|c|}{$R^{{\rm s}, {\rm ls}}_{T,2}(5)$} \\ 
 & & $\geq 1$  & $< 1$  & & $\geq 1$  & $< 1$ \\ 
 \hline 
 \multirow{2}{*}{$10000$} & $R^{{\rm s}, {\rm ls}}_{T,3}(1) \geq 1$  & 0.8504 & 0.0695  & $R^{{\rm s}, {\rm ls}}_{T,3}(5) \geq 1$  & 0.476 & 0.1447 \\ 
& $R^{{\rm s}, {\rm ls}}_{T,3}(1) < 1$  & 0.0799 & 2e-04 & $R^{{\rm s}, {\rm ls}}_{T,3}(5) < 1$  & 0.1324 & 0.2469 \\

\hline\hline
 \multirow{2}{*}{$n$} & & \multicolumn{2}{|c||}{$R^{{\rm s}, {\rm ls}}_{T,2}(1)$}  & & \multicolumn{2}{|c|}{$R^{{\rm s}, {\rm ls}}_{T,2}(5)$} \\ 
 & & $\geq 1.05$  & $< 1.05$  & & $\geq 1.05$  & $< 1.05$ \\ 
 \hline 
 \multirow{2}{*}{$10000$} & $R^{{\rm s}, {\rm ls}}_{T,3}(1) \geq 1.05$  & 0.0015 & 0.2483  & $R^{{\rm s}, {\rm ls}}_{T,3}(5) \geq 1.05$  & 0 & 0 \\ 
& $R^{{\rm s}, {\rm ls}}_{T,3}(1) < 1.05$  & 0.0208 & 0.7294 & $R^{{\rm s}, {\rm ls}}_{T,3}(5) < 1.05$  & 0 & 1 \\

\hline\hline
 \multirow{2}{*}{$n$} & & \multicolumn{2}{|c||}{$R^{{\rm s}, {\rm ls}}_{T,2}(1)$}  & & \multicolumn{2}{|c|}{$R^{{\rm s}, {\rm ls}}_{T,2}(5)$} \\ 
 & & $\geq 1.1$  & $< 1.1$  & & $\geq 1.1$  & $< 1.1$ \\ 
 \hline 
 \multirow{2}{*}{$10000$} & $R^{{\rm s}, {\rm ls}}_{T,3}(1) \geq 1.1$  & 0 & 0.0024  & $R^{{\rm s}, {\rm ls}}_{T,3}(5) \geq 1.1$  & 0 & 0 \\ 
& $R^{{\rm s}, {\rm ls}}_{T,3}(1) < 1.1$  & 0 & 0.9976 & $R^{{\rm s}, {\rm ls}}_{T,3}(5) < 1.1$  & 0 & 1 \\

\hline\hline
 \multirow{2}{*}{$n$} & & \multicolumn{2}{|c||}{$R^{{\rm s}, {\rm ls}}_{T,2}(1)$}  & & \multicolumn{2}{|c|}{$R^{{\rm s}, {\rm ls}}_{T,2}(5)$} \\ 
 & & $\geq 1.15$  & $< 1.15$  & & $\geq 1.15$  & $< 1.15$ \\ 
 \hline 
 \multirow{2}{*}{$10000$} & $R^{{\rm s}, {\rm ls}}_{T,3}(1) \geq 1.15$  & 0 & 0  & $R^{{\rm s}, {\rm ls}}_{T,3}(5) \geq 1.15$  & 0 & 0 \\ 
& $R^{{\rm s}, {\rm ls}}_{T,3}(1) < 1.15$  & 0 & 1 & $R^{{\rm s}, {\rm ls}}_{T,3}(5) < 1.15$  & 0 & 1 \\

\hline
\end{tabular}
\caption{\textit{Proportions of the individual events in~\eqref{samedecision} for the process \eqref{periodic2} and selected combinations of $n$ and $\delta$.}}  \label{MSPEanalysisperiodic2c}
\end{center}
\end{table}

\begin{table} \scriptsize
\begin{center}
	\begin{tabular}{|c||c|c|c|c|c|c|c|c|c|}\hline 
 &\multicolumn{9}{|c|}{\eqref{samedecision} holds for $h=1$} \\ 
\hline 
$\delta$ & $n$=100 & $n$=200 & $n$=500 & $n$=1000 & $n$=2000 & $n$=4000 & $n$=6000 & $n$=8000 & $n$=10000 \\ 
\hline \hline 
0  & 0.6063 & 0.6284 & 0.6129 & 0.606 & 0.6486 & 0.6964 & 0.7479 & 0.8021 & 0.8506 \\ 
0.01  & 0.6274 & 0.6627 & 0.6685 & 0.6366 & 0.6146 & 0.5723 & 0.5818 & 0.6013 & 0.6305 \\ 
0.05  & 0.7088 & 0.7881 & 0.8526 & 0.863 & 0.8584 & 0.8393 & 0.8064 & 0.7607 & 0.7309 \\ 
0.1  & 0.7935 & 0.8883 & 0.9556 & 0.975 & 0.9883 & 0.995 & 0.9959 & 0.9965 & 0.9976 \\ 
0.15  & 0.8615 & 0.942 & 0.9875 & 0.9947 & 0.9994 & 1 & 1 & 0.9999 & 1 \\ 
0.2  & 0.9033 & 0.97 & 0.9955 & 0.9993 & 1 & 1 & 1 & 1 & 1 \\ 
0.4  & 0.9774 & 0.997 & 0.9999 & 1 & 1 & 1 & 1 & 1 & 1 \\ 
0.6  & 0.9927 & 0.9998 & 1 & 1 & 1 & 1 & 1 & 1 & 1 \\ 
\hline 
\end{tabular} 
	\begin{tabular}{|c||c|c|c|c|c|c|c|c|c|}\hline 
 &\multicolumn{9}{|c|}{\eqref{samedecision} holds for $h=5$} \\ 
\hline 
$\delta$ & $n$=100 & $n$=200 & $n$=500 & $n$=1000 & $n$=2000 & $n$=4000 & $n$=6000 & $n$=8000 & $n$=10000 \\ 
\hline \hline 
0  & 0.6829 & 0.6818 & 0.682 & 0.6886 & 0.7109 & 0.7095 & 0.7166 & 0.7227 & 0.7229 \\ 
0.01  & 0.7161 & 0.7384 & 0.7957 & 0.8514 & 0.9135 & 0.9617 & 0.9747 & 0.9846 & 0.9899 \\ 
0.05  & 0.8041 & 0.8675 & 0.9483 & 0.98 & 0.9964 & 0.9998 & 1 & 1 & 1 \\ 
0.1  & 0.8771 & 0.937 & 0.9838 & 0.9971 & 0.9999 & 1 & 1 & 1 & 1 \\ 
0.15  & 0.9156 & 0.9638 & 0.9928 & 0.9992 & 1 & 1 & 1 & 1 & 1 \\ 
0.2  & 0.9398 & 0.9775 & 0.9964 & 0.9999 & 1 & 1 & 1 & 1 & 1 \\ 
0.4  & 0.9763 & 0.9936 & 0.9998 & 1 & 1 & 1 & 1 & 1 & 1 \\ 
0.6  & 0.9863 & 0.998 & 1 & 1 & 1 & 1 & 1 & 1 & 1 \\ 
\hline 
\end{tabular} 
\caption{\textit{Proportion of \eqref{samedecision}  being fulfilled for the process \eqref{periodic2} and different values of $h$, $\delta$ and $n$.}}  \label{MSPEanalysisperiodic2b}
\end{center}
\end{table}

\begin{table} \scriptsize
\begin{center}
	\begin{tabular}{|c||c|c|c|c|c|c|c|c|c|}\hline 
 &\multicolumn{9}{|p{12.5cm}|}{\centering Value of $q(\delta)$, defined in \eqref{cond:f}, for $h=1$} \\ 
\hline 
$\delta$ & $n$=100 & $n$=200 & $n$=500 & $n$=1000 & $n$=2000 & $n$=4000 & $n$=6000 & $n$=8000 & $n$=10000 \\ 
\hline \hline 
0  & 0 & 0 & 0 & 0 & 0 & 0 & 0 & 0 & 0 \\ 
0.01  & 2.8e-05 & 1e-05 & 3.1e-05 & 5.4e-05 & 0.00011 & 0.0028 & 0.0047 & 0.0054 & 0.0058 \\ 
0.05  & 0.04 & 0.026 & 0.017 & 0.012 & 0.008 & 0.0051 & 0.0036 & 0.0021 & 9e-04 \\ 
0.1  & 0.091 & 0.077 & 0.067 & 0.062 & 0.058 & 0.055 & 0.054 & 0.052 & 0.051 \\ 
0.15  & 0.14 & 0.13 & 0.12 & 0.11 & 0.11 & 0.11 & 0.1 & 0.1 & 0.1 \\ 
0.2  & 0.19 & 0.18 & 0.17 & 0.16 & 0.16 & 0.16 & 0.15 & 0.15 & 0.15 \\ 
0.4  & 0.4 & 0.38 & 0.37 & 0.36 & 0.36 & 0.36 & 0.36 & 0.35 & 0.35 \\ 
0.6  & 0.6 & 0.58 & 0.57 & 0.57 & 0.56 & 0.56 & 0.56 & 0.56 & 0.55 \\ 
\hline 
\end{tabular} 
	\begin{tabular}{|c||c|c|c|c|c|c|c|c|c|}\hline 
 &\multicolumn{9}{|p{12.5cm}|}{\centering Value of $q(\delta)$, defined in \eqref{cond:f}, for $h=5$} \\ 
\hline 
$\delta$ & $n$=100 & $n$=200 & $n$=500 & $n$=1000 & $n$=2000 & $n$=4000 & $n$=6000 & $n$=8000 & $n$=10000 \\ 
\hline \hline 
0  & 0 & 0 & 0 & 0 & 0 & 0 & 0 & 0 & 0 \\ 
0.01  & 0.01 & 0.01 & 0.01 & 0.01 & 0.01 & 0.01 & 0.01 & 0.01 & 0.01 \\ 
0.05  & 0.051 & 0.051 & 0.052 & 0.052 & 0.052 & 0.052 & 0.053 & 0.053 & 0.053 \\ 
0.1  & 0.1 & 0.1 & 0.1 & 0.1 & 0.1 & 0.1 & 0.11 & 0.11 & 0.11 \\ 
0.15  & 0.15 & 0.15 & 0.15 & 0.16 & 0.16 & 0.16 & 0.16 & 0.16 & 0.16 \\ 
0.2  & 0.2 & 0.2 & 0.21 & 0.21 & 0.21 & 0.21 & 0.21 & 0.21 & 0.21 \\ 
0.4  & 0.41 & 0.41 & 0.41 & 0.42 & 0.42 & 0.42 & 0.42 & 0.42 & 0.42 \\ 
0.6  & 0.61 & 0.61 & 0.62 & 0.62 & 0.63 & 0.63 & 0.63 & 0.63 & 0.63 \\ 
\hline 
\end{tabular} 
\caption{\textit{Values of $q(\delta)$, defined in \eqref{cond:f}, for the process \eqref{periodic2} and different values of $h$, $\delta$ and $n$.}} \label{MSPEanalysisperiodic2d}
\end{center}
\end{table}

\begin{table} \scriptsize
\begin{center}
	\begin{tabular}{|c||c|c|c|c|c|c|c|c|c|}\hline 
 &\multicolumn{9}{|c|}{\eqref{samedecision} holds for $h=1$, $i=2$} \\ 
\hline 
$\delta$ & $n$=100 & $n$=200 & $n$=500 & $n$=1000 & $n$=2000 & $n$=4000 & $n$=6000 & $n$=8000 & $n$=10000 \\ 
\hline \hline 
0  & 0.2844 & 0.2535 & 0.2795 & 0.4076 & 0.5431 & 0.5882 & 0.6001 & 0.651 & 0.7006 \\ 
0.01  & 0.2578 & 0.2148 & 0.2146 & 0.3167 & 0.4514 & 0.4971 & 0.4838 & 0.507 & 0.5306 \\ 
0.05  & 0.1706 & 0.1102 & 0.0758 & 0.094 & 0.1075 & 0.0897 & 0.0571 & 0.0385 & 0.0223 \\ 
0.1  & 0.1077 & 0.0531 & 0.0242 & 0.017 & 0.0074 & 7e-04 & 2e-04 & 1e-04 & 0 \\ 
0.15  & 0.0674 & 0.0267 & 0.0075 & 0.0035 & 4e-04 & 0 & 0 & 0 & 0 \\ 
0.2  & 0.0444 & 0.0137 & 0.0029 & 5e-04 & 0 & 0 & 0 & 0 & 0 \\ 
0.4  & 0.0101 & 0.0015 & 1e-04 & 0 & 0 & 0 & 0 & 0 & 0 \\ 
0.6  & 0.0039 & 1e-04 & 0 & 0 & 0 & 0 & 0 & 0 & 0 \\ 
\hline 
\end{tabular} 
	\begin{tabular}{|c||c|c|c|c|c|c|c|c|c|}\hline 
 &\multicolumn{9}{|c|}{\eqref{samedecision} holds for $h=1$, $i=3$} \\ 
\hline 
$\delta$ & $n$=100 & $n$=200 & $n$=500 & $n$=1000 & $n$=2000 & $n$=4000 & $n$=6000 & $n$=8000 & $n$=10000 \\ 
\hline \hline 
0  & 0.3267 & 0.3175 & 0.348 & 0.3572 & 0.3721 & 0.4678 & 0.5606 & 0.6277 & 0.6902 \\ 
0.01  & 0.3002 & 0.2729 & 0.2707 & 0.2651 & 0.2706 & 0.354 & 0.4658 & 0.5503 & 0.6307 \\ 
0.05  & 0.2106 & 0.1529 & 0.0934 & 0.0642 & 0.0445 & 0.0786 & 0.1407 & 0.2046 & 0.2498 \\ 
0.1  & 0.1398 & 0.0734 & 0.0226 & 0.0096 & 0.0043 & 0.0043 & 0.0039 & 0.0034 & 0.0024 \\ 
0.15  & 0.0869 & 0.0365 & 0.0054 & 0.0018 & 2e-04 & 0 & 0 & 1e-04 & 0 \\ 
0.2  & 0.0593 & 0.0181 & 0.0016 & 2e-04 & 0 & 0 & 0 & 0 & 0 \\ 
0.4  & 0.0131 & 0.0015 & 0 & 0 & 0 & 0 & 0 & 0 & 0 \\ 
0.6  & 0.0034 & 1e-04 & 0 & 0 & 0 & 0 & 0 & 0 & 0 \\ 
\hline 
\end{tabular} 
	\begin{tabular}{|c||c|c|c|c|c|c|c|c|c|}\hline 
 &\multicolumn{9}{|c|}{\eqref{samedecision} holds for $h=5$, $i=2$} \\ 
\hline 
$\delta$ & $n$=100 & $n$=200 & $n$=500 & $n$=1000 & $n$=2000 & $n$=4000 & $n$=6000 & $n$=8000 & $n$=10000 \\ 
\hline \hline 
0  & 0.2676 & 0.2555 & 0.2566 & 0.2562 & 0.2353 & 0.2273 & 0.2314 & 0.2236 & 0.2203 \\ 
0.01  & 0.2173 & 0.1758 & 0.1275 & 0.0825 & 0.044 & 0.0176 & 0.0107 & 0.0054 & 0.0034 \\ 
0.05  & 0.1208 & 0.0732 & 0.0254 & 0.0083 & 5e-04 & 0 & 0 & 0 & 0 \\ 
0.1  & 0.0686 & 0.0314 & 0.0068 & 0.001 & 0 & 0 & 0 & 0 & 0 \\ 
0.15  & 0.0464 & 0.0174 & 0.0024 & 1e-04 & 0 & 0 & 0 & 0 & 0 \\ 
0.2  & 0.0318 & 0.0101 & 7e-04 & 1e-04 & 0 & 0 & 0 & 0 & 0 \\ 
0.4  & 0.0117 & 0.0027 & 0 & 0 & 0 & 0 & 0 & 0 & 0 \\ 
0.6  & 0.0068 & 4e-04 & 0 & 0 & 0 & 0 & 0 & 0 & 0 \\ 
\hline 
\end{tabular} 
	\begin{tabular}{|c||c|c|c|c|c|c|c|c|c|}\hline 
 &\multicolumn{9}{|c|}{\eqref{samedecision} holds for $h=5$, $i=3$} \\ 
\hline 
$\delta$ & $n$=100 & $n$=200 & $n$=500 & $n$=1000 & $n$=2000 & $n$=4000 & $n$=6000 & $n$=8000 & $n$=10000 \\ 
\hline \hline 
0  & 0.2725 & 0.2793 & 0.2736 & 0.271 & 0.2486 & 0.2414 & 0.2334 & 0.2307 & 0.2326 \\ 
0.01  & 0.2244 & 0.1884 & 0.1358 & 0.0919 & 0.0521 & 0.0229 & 0.015 & 0.0102 & 0.0067 \\ 
0.05  & 0.1301 & 0.0811 & 0.0289 & 0.0119 & 0.0031 & 2e-04 & 0 & 0 & 0 \\ 
0.1  & 0.0773 & 0.0356 & 0.0094 & 0.0019 & 1e-04 & 0 & 0 & 0 & 0 \\ 
0.15  & 0.0488 & 0.0194 & 0.0048 & 7e-04 & 0 & 0 & 0 & 0 & 0 \\ 
0.2  & 0.0342 & 0.0124 & 0.0029 & 0 & 0 & 0 & 0 & 0 & 0 \\ 
0.4  & 0.0126 & 0.0037 & 2e-04 & 0 & 0 & 0 & 0 & 0 & 0 \\ 
0.6  & 0.0071 & 0.0016 & 0 & 0 & 0 & 0 & 0 & 0 & 0 \\ 
\hline 
\end{tabular} 
\caption{\textit{Proportion of \eqref{decisionrule}  being fulfilled for the process \eqref{periodic2} and different values of $h$, $\delta$ and $n$.}}  \label{MSPEanalysisperiodic2}
\end{center}
\end{table}

\clearpage

\begin{figure}
\centering 
\includegraphics[width=0.24\textwidth]{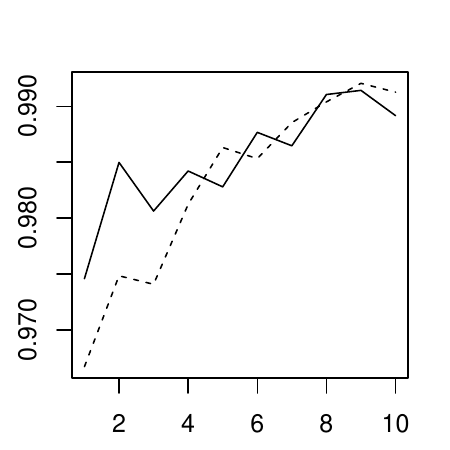}~
\includegraphics[width=0.24\textwidth]{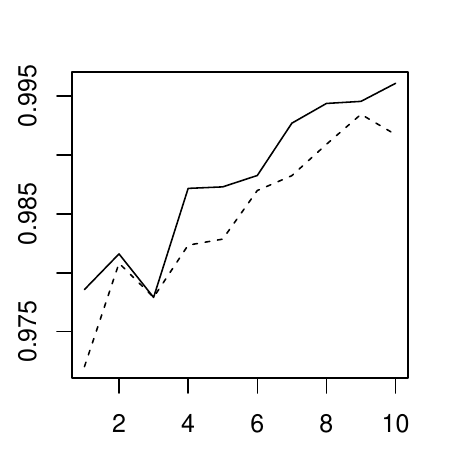}~
\includegraphics[width=0.24\textwidth]{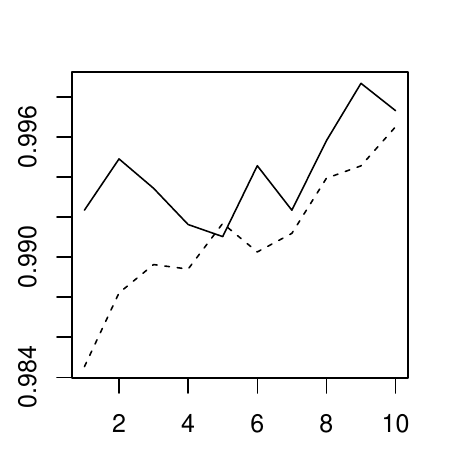}~ 
\includegraphics[width=0.24\textwidth]{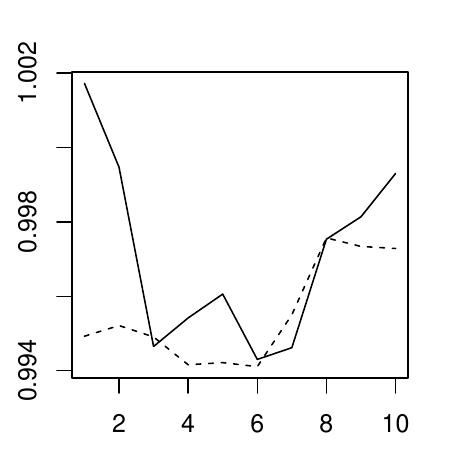}   \\
\includegraphics[width=0.24\textwidth]{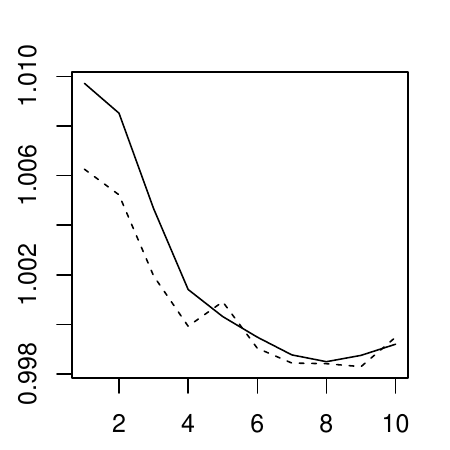}~
\includegraphics[width=0.24\textwidth]{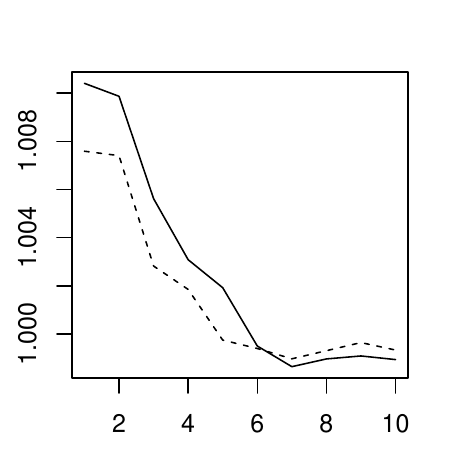}~ 
\includegraphics[width=0.24\textwidth]{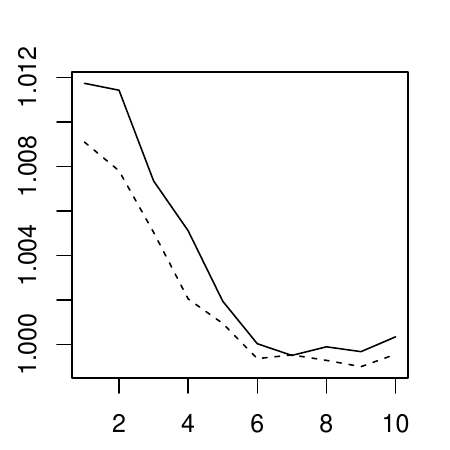}~ 
\includegraphics[width=0.24\textwidth]{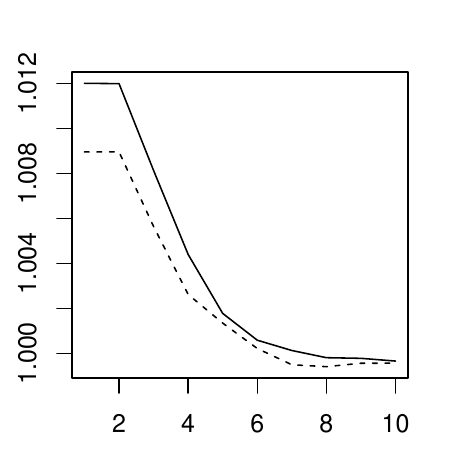}
   \caption{\it Plot of $h \mapsto R_{T,i}(h)$ for model \eqref{increasing2} and different values of $n$ (from left to right: n=100, n=200, n=500, n=1000 [first row], n=4000, n=6000, n=8000, n=10000 [second row]). Solid line: $i=3$ (test set), dashed line: $i=2$ (validation set 2). } \label{MSPEincreasing2}
\end{figure}

\begin{table} \scriptsize
\begin{center}
\begin{tabular}{|c||c|c|c||c|c|c|}
\hline
 \multirow{2}{*}{$n$} & & \multicolumn{2}{|c||}{$R^{{\rm s}, {\rm ls}}_{T,2}(1)$}  & & \multicolumn{2}{|c|}{$R^{{\rm s}, {\rm ls}}_{T,2}(5)$} \\ 
 & & $\geq 1.01$  & $< 1.01$  & & $\geq 1.01$  & $< 1.01$ \\ 
 \hline 
 \multirow{2}{*}{$100$} & $R^{{\rm s}, {\rm ls}}_{T,3}(1) \geq 1.01$  & 0.1236 & 0.1993  & $R^{{\rm s}, {\rm ls}}_{T,3}(5) \geq 1.01$  & 0.0822 & 0.1439 \\ 
& $R^{{\rm s}, {\rm ls}}_{T,3}(1) < 1.01$  & 0.1851 & 0.492 & $R^{{\rm s}, {\rm ls}}_{T,3}(5) < 1.01$  & 0.1386 & 0.6353 \\

\hline\hline
 \multirow{2}{*}{$n$} & & \multicolumn{2}{|c||}{$R^{{\rm s}, {\rm ls}}_{T,2}(1)$}  & & \multicolumn{2}{|c|}{$R^{{\rm s}, {\rm ls}}_{T,2}(5)$} \\ 
 & & $\geq 1.2$  & $< 1.2$  & & $\geq 1.2$  & $< 1.2$ \\ 
 \hline 
 \multirow{2}{*}{$1000$} & $R^{{\rm s}, {\rm ls}}_{T,3}(1) \geq 1.2$  & 3e-04 & 2e-04  & $R^{{\rm s}, {\rm ls}}_{T,3}(5) \geq 1.2$  & 0 & 0.0026 \\ 
& $R^{{\rm s}, {\rm ls}}_{T,3}(1) < 1.2$  & 0 & 0.9995 & $R^{{\rm s}, {\rm ls}}_{T,3}(5) < 1.2$  & 0.0021 & 0.9953 \\

\hline\hline
 \multirow{2}{*}{$n$} & & \multicolumn{2}{|c||}{$R^{{\rm s}, {\rm ls}}_{T,2}(1)$}  & & \multicolumn{2}{|c|}{$R^{{\rm s}, {\rm ls}}_{T,2}(5)$} \\ 
 & & $\geq 1$  & $< 1$  & & $\geq 1$  & $< 1$ \\ 
 \hline 
 \multirow{2}{*}{$10000$} & $R^{{\rm s}, {\rm ls}}_{T,3}(1) \geq 1$  & 0.7701 & 0.1313  & $R^{{\rm s}, {\rm ls}}_{T,3}(5) \geq 1$  & 0.4071 & 0.2352 \\ 
& $R^{{\rm s}, {\rm ls}}_{T,3}(1) < 1$  & 0.0825 & 0.0161 & $R^{{\rm s}, {\rm ls}}_{T,3}(5) < 1$  & 0.183 & 0.1747 \\

\hline\hline
 \multirow{2}{*}{$n$} & & \multicolumn{2}{|c||}{$R^{{\rm s}, {\rm ls}}_{T,2}(1)$}  & & \multicolumn{2}{|c|}{$R^{{\rm s}, {\rm ls}}_{T,2}(5)$} \\ 
 & & $\geq 1.01$  & $< 1.01$  & & $\geq 1.01$  & $< 1.01$ \\ 
 \hline 
 \multirow{2}{*}{$10000$} & $R^{{\rm s}, {\rm ls}}_{T,3}(1) \geq 1.01$  & 0.2386 & 0.2901  & $R^{{\rm s}, {\rm ls}}_{T,3}(5) \geq 1.01$  & 0.0322 & 0.1069 \\ 
& $R^{{\rm s}, {\rm ls}}_{T,3}(1) < 1.01$  & 0.1966 & 0.2747 & $R^{{\rm s}, {\rm ls}}_{T,3}(5) < 1.01$  & 0.0789 & 0.782 \\

\hline\hline
 \multirow{2}{*}{$n$} & & \multicolumn{2}{|c||}{$R^{{\rm s}, {\rm ls}}_{T,2}(1)$}  & & \multicolumn{2}{|c|}{$R^{{\rm s}, {\rm ls}}_{T,2}(5)$} \\ 
 & & $\geq 1.05$  & $< 1.05$  & & $\geq 1.05$  & $< 1.05$ \\ 
 \hline 
 \multirow{2}{*}{$10000$} & $R^{{\rm s}, {\rm ls}}_{T,3}(1) \geq 1.05$  & 0 & 0.0016  & $R^{{\rm s}, {\rm ls}}_{T,3}(5) \geq 1.05$  & 1e-04 & 0.0024 \\ 
& $R^{{\rm s}, {\rm ls}}_{T,3}(1) < 1.05$  & 9e-04 & 0.9975 & $R^{{\rm s}, {\rm ls}}_{T,3}(5) < 1.05$  & 0.0025 & 0.995 \\

\hline\hline
 \multirow{2}{*}{$n$} & & \multicolumn{2}{|c||}{$R^{{\rm s}, {\rm ls}}_{T,2}(1)$}  & & \multicolumn{2}{|c|}{$R^{{\rm s}, {\rm ls}}_{T,2}(5)$} \\ 
 & & $\geq 1.1$  & $< 1.1$  & & $\geq 1.1$  & $< 1.1$ \\ 
 \hline 
 \multirow{2}{*}{$10000$} & $R^{{\rm s}, {\rm ls}}_{T,3}(1) \geq 1.1$  & 0 & 0  & $R^{{\rm s}, {\rm ls}}_{T,3}(5) \geq 1.1$  & 0 & 0 \\ 
& $R^{{\rm s}, {\rm ls}}_{T,3}(1) < 1.1$  & 0 & 1 & $R^{{\rm s}, {\rm ls}}_{T,3}(5) < 1.1$  & 0 & 1 \\

\hline
\end{tabular}
\caption{\textit{Proportions of the individual events in~\eqref{samedecision} for the process \eqref{increasing2} and selected combinations of $n$ and $\delta$.}}  \label{MSPEanalysisincreasing2c}
\end{center}
\end{table}

\begin{table} \scriptsize
\begin{center}
	\begin{tabular}{|c||c|c|c|c|c|c|c|c|c|}\hline 
 &\multicolumn{9}{|c|}{\eqref{samedecision} holds for $h=1$} \\ 
\hline 
$\delta$ & $n$=100 & $n$=200 & $n$=500 & $n$=1000 & $n$=2000 & $n$=4000 & $n$=6000 & $n$=8000 & $n$=10000 \\ 
\hline \hline 
0  & 0.5924 & 0.582 & 0.556 & 0.5518 & 0.5554 & 0.6074 & 0.6722 & 0.7385 & 0.7862 \\ 
0.01  & 0.6156 & 0.621 & 0.6174 & 0.6156 & 0.5915 & 0.5529 & 0.5331 & 0.5259 & 0.5133 \\ 
0.05  & 0.7113 & 0.7678 & 0.8578 & 0.9019 & 0.9478 & 0.981 & 0.9911 & 0.9955 & 0.9975 \\ 
0.1  & 0.7958 & 0.8813 & 0.9656 & 0.9902 & 0.9988 & 0.9999 & 1 & 0.9999 & 1 \\ 
0.15  & 0.8602 & 0.9396 & 0.9908 & 0.9991 & 1 & 1 & 1 & 1 & 1 \\ 
0.2  & 0.9026 & 0.9675 & 0.9976 & 0.9998 & 1 & 1 & 1 & 1 & 1 \\ 
0.4  & 0.9673 & 0.9911 & 0.9995 & 0.9999 & 1 & 1 & 1 & 1 & 1 \\ 
0.6  & 0.9803 & 0.9922 & 0.9989 & 0.9997 & 1 & 1 & 1 & 1 & 1 \\ 
\hline 
\end{tabular} 
	\begin{tabular}{|c||c|c|c|c|c|c|c|c|c|}\hline 
 &\multicolumn{9}{|c|}{\eqref{samedecision} holds for $h=5$} \\ 
\hline 
$\delta$ & $n$=100 & $n$=200 & $n$=500 & $n$=1000 & $n$=2000 & $n$=4000 & $n$=6000 & $n$=8000 & $n$=10000 \\ 
\hline \hline 
0  & 0.702 & 0.6926 & 0.6856 & 0.6671 & 0.6412 & 0.6209 & 0.6147 & 0.5964 & 0.5818 \\ 
0.01  & 0.7175 & 0.7238 & 0.7427 & 0.7496 & 0.7618 & 0.7803 & 0.8009 & 0.787 & 0.8142 \\ 
0.05  & 0.7865 & 0.8308 & 0.8874 & 0.925 & 0.9543 & 0.9788 & 0.9873 & 0.9915 & 0.9951 \\ 
0.1  & 0.8444 & 0.9005 & 0.9508 & 0.9733 & 0.9895 & 0.9983 & 0.9991 & 0.9997 & 1 \\ 
0.15  & 0.8869 & 0.9357 & 0.9751 & 0.9892 & 0.9973 & 0.9999 & 0.9999 & 1 & 1 \\ 
0.2  & 0.9123 & 0.9541 & 0.9865 & 0.9953 & 0.9993 & 1 & 1 & 1 & 1 \\ 
0.4  & 0.9628 & 0.984 & 0.9973 & 0.9998 & 1 & 1 & 1 & 1 & 1 \\ 
0.6  & 0.9792 & 0.994 & 0.9991 & 1 & 1 & 1 & 1 & 1 & 1 \\ 
\hline 
\end{tabular} 
\caption{\textit{Proportion of \eqref{samedecision}  being fulfilled for the process \eqref{increasing2} and different values of $h$, $\delta$ and $n$.}}  \label{MSPEanalysisincreasing2b}
\end{center}
\end{table}

\begin{table} \scriptsize
\begin{center}
	\begin{tabular}{|c||c|c|c|c|c|c|c|c|c|}\hline 
 &\multicolumn{9}{|p{12.5cm}|}{\centering Value of $q(\delta)$, defined in \eqref{cond:f}, for $h=1$} \\ 
\hline 
$\delta$ & $n$=100 & $n$=200 & $n$=500 & $n$=1000 & $n$=2000 & $n$=4000 & $n$=6000 & $n$=8000 & $n$=10000 \\ 
\hline \hline 
0  & 0 & 0 & 0 & 0 & 0 & 0 & 0 & 0 & 0 \\ 
0.01  & 0.00015 & 1.8e-05 & 0.00017 & 0.00093 & 0.0015 & 0.002 & 0.0023 & 0.0024 & 0.0026 \\ 
0.05  & 0.04 & 0.039 & 0.039 & 0.038 & 0.038 & 0.037 & 0.037 & 0.037 & 0.037 \\ 
0.1  & 0.09 & 0.089 & 0.089 & 0.088 & 0.088 & 0.087 & 0.087 & 0.087 & 0.087 \\ 
0.15  & 0.14 & 0.14 & 0.14 & 0.14 & 0.14 & 0.14 & 0.14 & 0.14 & 0.14 \\ 
0.2  & 0.19 & 0.19 & 0.19 & 0.19 & 0.19 & 0.19 & 0.19 & 0.19 & 0.19 \\ 
0.4  & 0.39 & 0.39 & 0.39 & 0.39 & 0.39 & 0.39 & 0.39 & 0.39 & 0.39 \\ 
0.6  & 0.23 & 0.24 & 0.25 & 0.26 & 0.27 & 0.27 & 0.27 & 0.28 & 0.28 \\ 
\hline 
\end{tabular} 
	\begin{tabular}{|c||c|c|c|c|c|c|c|c|c|}\hline 
 &\multicolumn{9}{|p{12.5cm}|}{\centering Value of $q(\delta)$, defined in \eqref{cond:f}, for $h=5$} \\ 
\hline 
$\delta$ & $n$=100 & $n$=200 & $n$=500 & $n$=1000 & $n$=2000 & $n$=4000 & $n$=6000 & $n$=8000 & $n$=10000 \\ 
\hline \hline 
0  & 0 & 0 & 0 & 0 & 0 & 0 & 0 & 0 & 0 \\ 
0.01  & 0.011 & 0.011 & 0.011 & 0.01 & 0.0096 & 0.0093 & 0.0091 & 0.009 & 0.0089 \\ 
0.05  & 0.058 & 0.055 & 0.053 & 0.052 & 0.051 & 0.051 & 0.05 & 0.05 & 0.05 \\ 
0.1  & 0.15 & 0.15 & 0.14 & 0.14 & 0.14 & 0.14 & 0.14 & 0.14 & 0.14 \\ 
0.15  & 0.23 & 0.23 & 0.23 & 0.23 & 0.23 & 0.23 & 0.23 & 0.23 & 0.23 \\ 
0.2  & 0.32 & 0.32 & 0.33 & 0.33 & 0.33 & 0.33 & 0.33 & 0.33 & 0.33 \\ 
0.4  & 0.68 & 0.68 & 0.69 & 0.69 & 0.69 & 0.69 & 0.69 & 0.69 & 0.69 \\ 
0.6  & 1 & 1 & 1 & 1.1 & 1.1 & 1.1 & 1.1 & 1.1 & 1.1 \\ 
\hline 
\end{tabular} 
\caption{\textit{Values of $q(\delta)$, defined in \eqref{cond:f}, for the process \eqref{increasing2} and different values of $h$, $\delta$ and $n$.}} \label{MSPEanalysisincreasing2d}
\end{center}
\end{table}

\begin{table} \scriptsize
\begin{center}
	\begin{tabular}{|c||c|c|c|c|c|c|c|c|c|}\hline 
 &\multicolumn{9}{|c|}{\eqref{samedecision} holds for $h=1$, $i=2$} \\ 
\hline 
$\delta$ & $n$=100 & $n$=200 & $n$=500 & $n$=1000 & $n$=2000 & $n$=4000 & $n$=6000 & $n$=8000 & $n$=10000 \\ 
\hline \hline 
0  & 0.3466 & 0.3397 & 0.372 & 0.4311 & 0.5369 & 0.6691 & 0.748 & 0.8076 & 0.8526 \\ 
0.01  & 0.3087 & 0.2828 & 0.2716 & 0.2824 & 0.3117 & 0.3525 & 0.3829 & 0.4138 & 0.4352 \\ 
0.05  & 0.2038 & 0.1503 & 0.0755 & 0.0466 & 0.0243 & 0.0074 & 0.0026 & 0.0018 & 9e-04 \\ 
0.1  & 0.1294 & 0.0754 & 0.0169 & 0.0042 & 2e-04 & 1e-04 & 0 & 0 & 0 \\ 
0.15  & 0.0858 & 0.0416 & 0.0046 & 5e-04 & 0 & 0 & 0 & 0 & 0 \\ 
0.2  & 0.0594 & 0.0265 & 0.0024 & 3e-04 & 0 & 0 & 0 & 0 & 0 \\ 
0.4  & 0.0214 & 0.0098 & 0.002 & 2e-04 & 0 & 0 & 0 & 0 & 0 \\ 
0.6  & 0.0111 & 0.0054 & 0.0013 & 0 & 0 & 0 & 0 & 0 & 0 \\ 
\hline 
\end{tabular} 
	\begin{tabular}{|c||c|c|c|c|c|c|c|c|c|}\hline 
 &\multicolumn{9}{|c|}{\eqref{samedecision} holds for $h=1$, $i=3$} \\ 
\hline 
$\delta$ & $n$=100 & $n$=200 & $n$=500 & $n$=1000 & $n$=2000 & $n$=4000 & $n$=6000 & $n$=8000 & $n$=10000 \\ 
\hline \hline 
0  & 0.3584 & 0.3607 & 0.4132 & 0.4985 & 0.6017 & 0.7405 & 0.8162 & 0.8751 & 0.9014 \\ 
0.01  & 0.3229 & 0.3072 & 0.3086 & 0.3476 & 0.3714 & 0.4384 & 0.4798 & 0.5205 & 0.5287 \\ 
0.05  & 0.2229 & 0.1663 & 0.0967 & 0.0639 & 0.0317 & 0.012 & 0.0063 & 0.0027 & 0.0016 \\ 
0.1  & 0.1458 & 0.0819 & 0.0235 & 0.0072 & 0.001 & 0 & 0 & 1e-04 & 0 \\ 
0.15  & 0.1 & 0.0446 & 0.0088 & 0.001 & 0 & 0 & 0 & 0 & 0 \\ 
0.2  & 0.0736 & 0.0268 & 0.0038 & 5e-04 & 0 & 0 & 0 & 0 & 0 \\ 
0.4  & 0.0283 & 0.0109 & 0.0017 & 3e-04 & 0 & 0 & 0 & 0 & 0 \\ 
0.6  & 0.017 & 0.0076 & 0.001 & 3e-04 & 0 & 0 & 0 & 0 & 0 \\ 
\hline 
\end{tabular} 
	\begin{tabular}{|c||c|c|c|c|c|c|c|c|c|}\hline 
 &\multicolumn{9}{|c|}{\eqref{samedecision} holds for $h=5$, $i=2$} \\ 
\hline 
$\delta$ & $n$=100 & $n$=200 & $n$=500 & $n$=1000 & $n$=2000 & $n$=4000 & $n$=6000 & $n$=8000 & $n$=10000 \\ 
\hline \hline 
0  & 0.2545 & 0.2395 & 0.2448 & 0.2632 & 0.29 & 0.3548 & 0.3959 & 0.4407 & 0.4726 \\ 
0.01  & 0.2208 & 0.1887 & 0.1716 & 0.1654 & 0.1438 & 0.1389 & 0.1211 & 0.125 & 0.1111 \\ 
0.05  & 0.1389 & 0.0991 & 0.059 & 0.04 & 0.0211 & 0.0116 & 0.0067 & 0.0049 & 0.0026 \\ 
0.1  & 0.0861 & 0.0513 & 0.0246 & 0.0124 & 0.0045 & 8e-04 & 3e-04 & 1e-04 & 0 \\ 
0.15  & 0.0584 & 0.0311 & 0.011 & 0.0043 & 0.001 & 1e-04 & 0 & 0 & 0 \\ 
0.2  & 0.0432 & 0.0203 & 0.0062 & 0.0021 & 5e-04 & 0 & 0 & 0 & 0 \\ 
0.4  & 0.0176 & 0.0062 & 9e-04 & 2e-04 & 0 & 0 & 0 & 0 & 0 \\ 
0.6  & 0.0104 & 0.0022 & 2e-04 & 0 & 0 & 0 & 0 & 0 & 0 \\ 
\hline 
\end{tabular} 
	\begin{tabular}{|c||c|c|c|c|c|c|c|c|c|}\hline 
 &\multicolumn{9}{|c|}{\eqref{samedecision} holds for $h=5$, $i=3$} \\ 
\hline 
$\delta$ & $n$=100 & $n$=200 & $n$=500 & $n$=1000 & $n$=2000 & $n$=4000 & $n$=6000 & $n$=8000 & $n$=10000 \\ 
\hline \hline 
0  & 0.2527 & 0.2503 & 0.2542 & 0.2781 & 0.3186 & 0.3751 & 0.427 & 0.4913 & 0.5248 \\ 
0.01  & 0.2261 & 0.2075 & 0.1841 & 0.177 & 0.17 & 0.1578 & 0.1496 & 0.1596 & 0.1391 \\ 
0.05  & 0.1552 & 0.1113 & 0.0698 & 0.046 & 0.0286 & 0.0122 & 0.0084 & 0.0044 & 0.0025 \\ 
0.1  & 0.1031 & 0.064 & 0.0292 & 0.0153 & 0.0064 & 9e-04 & 6e-04 & 2e-04 & 0 \\ 
0.15  & 0.0749 & 0.0408 & 0.0153 & 0.0069 & 0.0019 & 0 & 1e-04 & 0 & 0 \\ 
0.2  & 0.0577 & 0.0288 & 0.0081 & 0.0026 & 2e-04 & 0 & 0 & 0 & 0 \\ 
0.4  & 0.0234 & 0.01 & 0.0018 & 0 & 0 & 0 & 0 & 0 & 0 \\ 
0.6  & 0.0122 & 0.0038 & 7e-04 & 0 & 0 & 0 & 0 & 0 & 0 \\ 
\hline 
\end{tabular} 
\caption{\textit{Proportion of \eqref{decisionrule}  being fulfilled for the process \eqref{increasing2} and different values of $h$, $\delta$ and $n$.}}  \label{MSPEanalysisincreasing2}
\end{center}
\end{table}

\clearpage

\begin{figure}
\centering 
\includegraphics[width=0.24\textwidth]{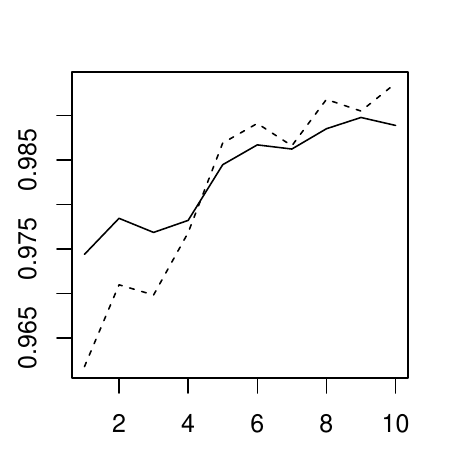}~
\includegraphics[width=0.24\textwidth]{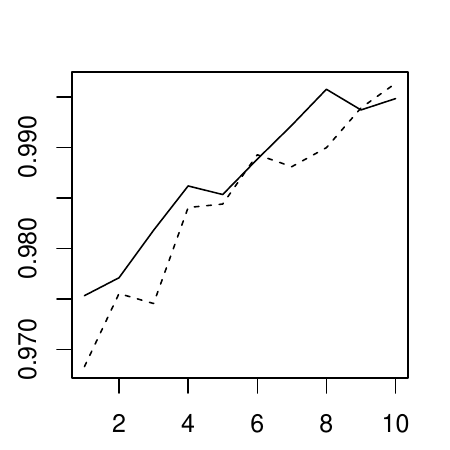}~
\includegraphics[width=0.24\textwidth]{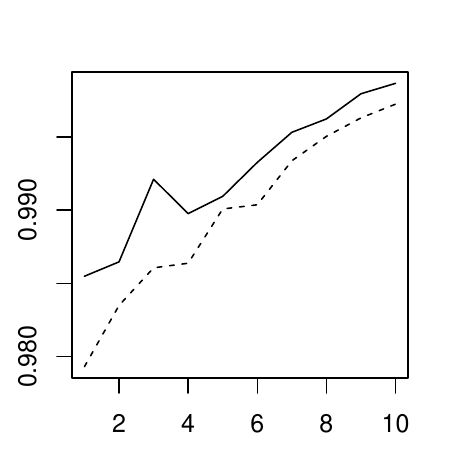}~ 
\includegraphics[width=0.24\textwidth]{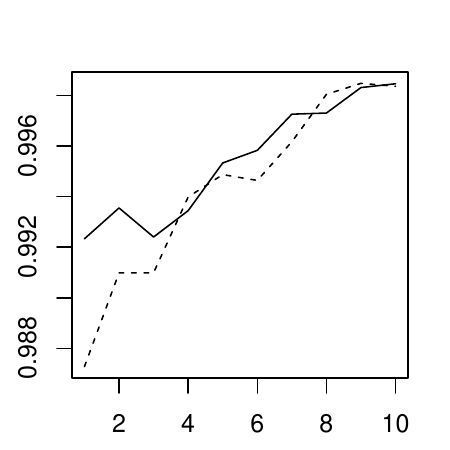} \\
\includegraphics[width=0.24\textwidth]{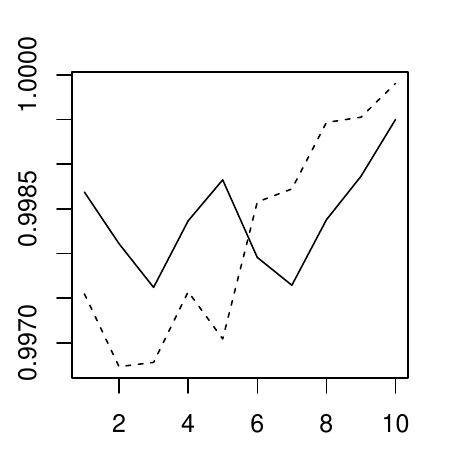}~
\includegraphics[width=0.24\textwidth]{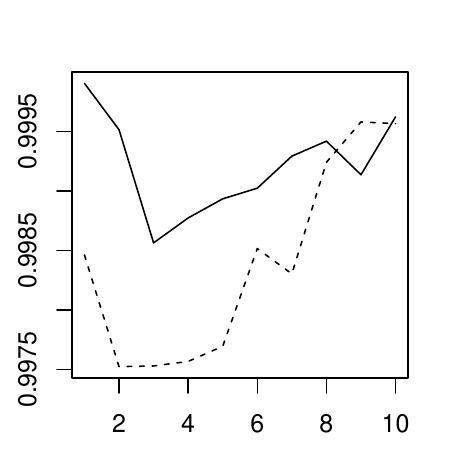}~ 
\includegraphics[width=0.24\textwidth]{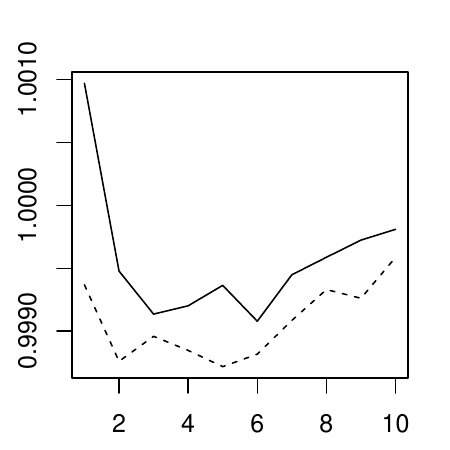}~ 
\includegraphics[width=0.24\textwidth]{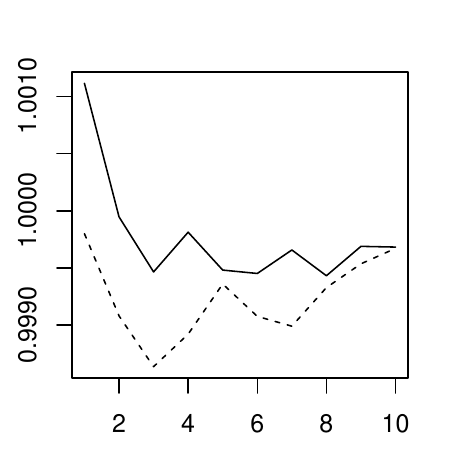}
   \caption{\it Plot of $h \mapsto R_{T,i}(h)$ for model \eqref{increasing4} and different values of $n$ (from left to right: n=100, n=200, n=500, n=1000 [first row], n=4000, n=6000, n=8000, n=10000 [second row]). Solid line: $i=3$ (test set), dashed line: $i=2$ (validation set 2).} \label{MSPEincreasing4}
\end{figure}

\begin{table} \scriptsize
\begin{center}
\begin{tabular}{|c||c|c|c||c|c|c|}
\hline
 \multirow{2}{*}{$n$} & & \multicolumn{2}{|c||}{$R^{{\rm s}, {\rm ls}}_{T,2}(1)$}  & & \multicolumn{2}{|c|}{$R^{{\rm s}, {\rm ls}}_{T,2}(5)$} \\ 
 & & $\geq 1.01$  & $< 1.01$  & & $\geq 1.01$  & $< 1.01$ \\ 
 \hline 
 \multirow{2}{*}{$100$} & $R^{{\rm s}, {\rm ls}}_{T,3}(1) \geq 1.01$  & 0.1155 & 0.1957  & $R^{{\rm s}, {\rm ls}}_{T,3}(5) \geq 1.01$  & 0.0802 & 0.1499 \\ 
& $R^{{\rm s}, {\rm ls}}_{T,3}(1) < 1.01$  & 0.1763 & 0.5125 & $R^{{\rm s}, {\rm ls}}_{T,3}(5) < 1.01$  & 0.1335 & 0.6364 \\

\hline\hline
 \multirow{2}{*}{$n$} & & \multicolumn{2}{|c||}{$R^{{\rm s}, {\rm ls}}_{T,2}(1)$}  & & \multicolumn{2}{|c|}{$R^{{\rm s}, {\rm ls}}_{T,2}(5)$} \\ 
 & & $\geq 1.2$  & $< 1.2$  & & $\geq 1.2$  & $< 1.2$ \\ 
 \hline 
 \multirow{2}{*}{$1000$} & $R^{{\rm s}, {\rm ls}}_{T,3}(1) \geq 1.2$  & 5e-04 & 0  & $R^{{\rm s}, {\rm ls}}_{T,3}(5) \geq 1.2$  & 0 & 9e-04 \\ 
& $R^{{\rm s}, {\rm ls}}_{T,3}(1) < 1.2$  & 1e-04 & 0.9994 & $R^{{\rm s}, {\rm ls}}_{T,3}(5) < 1.2$  & 5e-04 & 0.9986 \\

\hline\hline
 \multirow{2}{*}{$n$} & & \multicolumn{2}{|c||}{$R^{{\rm s}, {\rm ls}}_{T,2}(1)$}  & & \multicolumn{2}{|c|}{$R^{{\rm s}, {\rm ls}}_{T,2}(5)$} \\ 
 & & $\geq 1$  & $< 1$  & & $\geq 1$  & $< 1$ \\ 
 \hline 
 \multirow{2}{*}{$10000$} & $R^{{\rm s}, {\rm ls}}_{T,3}(1) \geq 1$  & 0.2792 & 0.2651  & $R^{{\rm s}, {\rm ls}}_{T,3}(5) \geq 1$  & 0.3615 & 0.1876 \\ 
& $R^{{\rm s}, {\rm ls}}_{T,3}(1) < 1$  & 0.2062 & 0.2495 & $R^{{\rm s}, {\rm ls}}_{T,3}(5) < 1$  & 0.1632 & 0.2877 \\

\hline\hline
 \multirow{2}{*}{$n$} & & \multicolumn{2}{|c||}{$R^{{\rm s}, {\rm ls}}_{T,2}(1)$}  & & \multicolumn{2}{|c|}{$R^{{\rm s}, {\rm ls}}_{T,2}(5)$} \\ 
 & & $\geq 1.01$  & $< 1.01$  & & $\geq 1.01$  & $< 1.01$ \\ 
 \hline 
 \multirow{2}{*}{$10000$} & $R^{{\rm s}, {\rm ls}}_{T,3}(1) \geq 1.01$  & 0.0037 & 0.0436  & $R^{{\rm s}, {\rm ls}}_{T,3}(5) \geq 1.01$  & 0.003 & 0.024 \\ 
& $R^{{\rm s}, {\rm ls}}_{T,3}(1) < 1.01$  & 0.0303 & 0.9224 & $R^{{\rm s}, {\rm ls}}_{T,3}(5) < 1.01$  & 0.0209 & 0.9521 \\

\hline\hline
 \multirow{2}{*}{$n$} & & \multicolumn{2}{|c||}{$R^{{\rm s}, {\rm ls}}_{T,2}(1)$}  & & \multicolumn{2}{|c|}{$R^{{\rm s}, {\rm ls}}_{T,2}(5)$} \\ 
 & & $\geq 1.05$  & $< 1.05$  & & $\geq 1.05$  & $< 1.05$ \\ 
 \hline 
 \multirow{2}{*}{$10000$} & $R^{{\rm s}, {\rm ls}}_{T,3}(1) \geq 1.05$  & 0 & 0  & $R^{{\rm s}, {\rm ls}}_{T,3}(5) \geq 1.05$  & 0 & 5e-04 \\ 
& $R^{{\rm s}, {\rm ls}}_{T,3}(1) < 1.05$  & 0 & 1 & $R^{{\rm s}, {\rm ls}}_{T,3}(5) < 1.05$  & 1e-04 & 0.9994 \\

\hline\hline
 \multirow{2}{*}{$n$} & & \multicolumn{2}{|c||}{$R^{{\rm s}, {\rm ls}}_{T,2}(1)$}  & & \multicolumn{2}{|c|}{$R^{{\rm s}, {\rm ls}}_{T,2}(5)$} \\ 
 & & $\geq 1.1$  & $< 1.1$  & & $\geq 1.1$  & $< 1.1$ \\ 
 \hline 
 \multirow{2}{*}{$10000$} & $R^{{\rm s}, {\rm ls}}_{T,3}(1) \geq 1.1$  & 0 & 0  & $R^{{\rm s}, {\rm ls}}_{T,3}(5) \geq 1.1$  & 0 & 0 \\ 
& $R^{{\rm s}, {\rm ls}}_{T,3}(1) < 1.1$  & 0 & 1 & $R^{{\rm s}, {\rm ls}}_{T,3}(5) < 1.1$  & 0 & 1 \\

\hline
\end{tabular}
\caption{\textit{Proportions of the individual events in~\eqref{samedecision} for the process \eqref{increasing4} and selected combinations of $n$ and $\delta$.}}  \label{MSPEanalysisincreasing4c}
\end{center}
\end{table}

\begin{table} \scriptsize
\begin{center}
	\begin{tabular}{|c||c|c|c|c|c|c|c|c|c|}\hline 
 &\multicolumn{9}{|c|}{\eqref{samedecision} holds for $h=1$} \\ 
\hline 
$\delta$ & $n$=100 & $n$=200 & $n$=500 & $n$=1000 & $n$=2000 & $n$=4000 & $n$=6000 & $n$=8000 & $n$=10000 \\ 
\hline \hline 
0  & 0.5996 & 0.5953 & 0.5831 & 0.5891 & 0.5691 & 0.5473 & 0.5436 & 0.5348 & 0.5287 \\ 
0.01  & 0.628 & 0.6487 & 0.6854 & 0.7301 & 0.7867 & 0.8543 & 0.886 & 0.9 & 0.9261 \\ 
0.05  & 0.7209 & 0.7983 & 0.921 & 0.9754 & 0.9951 & 0.9999 & 0.9998 & 0.9999 & 1 \\ 
0.1  & 0.8061 & 0.8999 & 0.9859 & 0.9985 & 1 & 1 & 1 & 1 & 1 \\ 
0.15  & 0.8645 & 0.9503 & 0.9966 & 0.9999 & 1 & 1 & 1 & 1 & 1 \\ 
0.2  & 0.9043 & 0.9715 & 0.9978 & 0.9999 & 1 & 1 & 1 & 1 & 1 \\ 
0.4  & 0.9649 & 0.9881 & 0.9968 & 0.9995 & 0.9999 & 1 & 1 & 1 & 1 \\ 
0.6  & 0.9808 & 0.9912 & 0.9977 & 0.9997 & 1 & 1 & 1 & 1 & 1 \\ 
\hline 
\end{tabular} 
	\begin{tabular}{|c||c|c|c|c|c|c|c|c|c|}\hline 
 &\multicolumn{9}{|c|}{\eqref{samedecision} holds for $h=5$} \\ 
\hline 
$\delta$ & $n$=100 & $n$=200 & $n$=500 & $n$=1000 & $n$=2000 & $n$=4000 & $n$=6000 & $n$=8000 & $n$=10000 \\ 
\hline \hline 
0  & 0.6994 & 0.6939 & 0.6911 & 0.6881 & 0.6782 & 0.6803 & 0.6729 & 0.6637 & 0.6492 \\ 
0.01  & 0.7166 & 0.7265 & 0.7633 & 0.7902 & 0.8411 & 0.8989 & 0.928 & 0.9376 & 0.9551 \\ 
0.05  & 0.7886 & 0.8455 & 0.9125 & 0.9522 & 0.9798 & 0.996 & 0.9987 & 0.9994 & 0.9994 \\ 
0.1  & 0.8548 & 0.9114 & 0.9631 & 0.9874 & 0.9971 & 0.9999 & 1 & 1 & 1 \\ 
0.15  & 0.8962 & 0.9438 & 0.9829 & 0.9957 & 0.9996 & 1 & 1 & 1 & 1 \\ 
0.2  & 0.9224 & 0.9626 & 0.9906 & 0.9986 & 0.9999 & 1 & 1 & 1 & 1 \\ 
0.4  & 0.9682 & 0.9886 & 0.9989 & 0.9999 & 1 & 1 & 1 & 1 & 1 \\ 
0.6  & 0.9817 & 0.9945 & 0.9998 & 1 & 1 & 1 & 1 & 1 & 1 \\ 
\hline 
\end{tabular}   
\caption{\textit{Proportion of \eqref{samedecision}  being fulfilled for the process \eqref{increasing4} and different values of $h$, $\delta$ and $n$.}}  \label{MSPEanalysisincreasing4b}
\end{center}
\end{table}

\begin{table} \scriptsize
\begin{center}
	\begin{tabular}{|c||c|c|c|c|c|c|c|c|c|}\hline 
 &\multicolumn{9}{|p{12.5cm}|}{\centering Value of $q(\delta)$, defined in \eqref{cond:f}, for $h=1$} \\ 
\hline 
$\delta$ & $n$=100 & $n$=200 & $n$=500 & $n$=1000 & $n$=2000 & $n$=4000 & $n$=6000 & $n$=8000 & $n$=10000 \\ 
\hline \hline 
0  & 0 & 0 & 0 & 0 & 0 & 0 & 0 & 0 & 0 \\ 
0.01  & 0.0081 & 0.0079 & 0.0078 & 0.0077 & 0.0076 & 0.0075 & 0.0075 & 0.0075 & 0.0075 \\ 
0.05  & 0.048 & 0.048 & 0.048 & 0.048 & 0.048 & 0.048 & 0.047 & 0.047 & 0.047 \\ 
0.1  & 0.098 & 0.098 & 0.098 & 0.098 & 0.098 & 0.098 & 0.097 & 0.097 & 0.097 \\ 
0.15  & 0.15 & 0.15 & 0.15 & 0.15 & 0.15 & 0.15 & 0.15 & 0.15 & 0.15 \\ 
0.2  & 0.2 & 0.2 & 0.2 & 0.2 & 0.2 & 0.2 & 0.2 & 0.2 & 0.2 \\ 
0.4  & 0.11 & 0.12 & 0.12 & 0.12 & 0.12 & 0.12 & 0.12 & 0.13 & 0.13 \\ 
0.6  & 0.086 & 0.083 & 0.08 & 0.078 & 0.077 & 0.076 & 0.075 & 0.075 & 0.074 \\ 
\hline 
\end{tabular} 
	\begin{tabular}{|c||c|c|c|c|c|c|c|c|c|}\hline 
 &\multicolumn{9}{|p{12.5cm}|}{\centering Value of $q(\delta)$, defined in \eqref{cond:f}, for $h=5$} \\ 
\hline 
$\delta$ & $n$=100 & $n$=200 & $n$=500 & $n$=1000 & $n$=2000 & $n$=4000 & $n$=6000 & $n$=8000 & $n$=10000 \\ 
\hline \hline 
0  & 0 & 0 & 0 & 0 & 0 & 0 & 0 & 0 & 0 \\ 
0.01  & 0.0087 & 0.0084 & 0.0081 & 0.008 & 0.0079 & 0.0079 & 0.0078 & 0.0078 & 0.0078 \\ 
0.05  & 0.069 & 0.069 & 0.069 & 0.069 & 0.069 & 0.069 & 0.069 & 0.069 & 0.069 \\ 
0.1  & 0.14 & 0.14 & 0.14 & 0.14 & 0.14 & 0.14 & 0.14 & 0.14 & 0.14 \\ 
0.15  & 0.22 & 0.22 & 0.22 & 0.22 & 0.22 & 0.22 & 0.22 & 0.22 & 0.22 \\ 
0.2  & 0.29 & 0.29 & 0.3 & 0.3 & 0.3 & 0.3 & 0.3 & 0.3 & 0.3 \\ 
0.4  & 0.59 & 0.6 & 0.6 & 0.6 & 0.6 & 0.6 & 0.6 & 0.6 & 0.6 \\ 
0.6  & 0.89 & 0.9 & 0.9 & 0.9 & 0.9 & 0.9 & 0.9 & 0.9 & 0.9 \\ 
\hline 
\end{tabular} 
\caption{\textit{Values of $q(\delta)$, defined in \eqref{cond:f}, for the process \eqref{increasing4} and different values of $h$, $\delta$ and $n$.}} \label{MSPEanalysisincreasing4d}
\end{center}
\end{table}

\begin{table} \scriptsize
\begin{center}
	\begin{tabular}{|c||c|c|c|c|c|c|c|c|c|}\hline 
 &\multicolumn{9}{|c|}{\eqref{samedecision} holds for $h=1$, $i=2$} \\ 
\hline 
$\delta$ & $n$=100 & $n$=200 & $n$=500 & $n$=1000 & $n$=2000 & $n$=4000 & $n$=6000 & $n$=8000 & $n$=10000 \\ 
\hline \hline 
0  & 0.3261 & 0.3146 & 0.3074 & 0.3079 & 0.3326 & 0.3748 & 0.4138 & 0.4533 & 0.4854 \\ 
0.01  & 0.2918 & 0.2563 & 0.2094 & 0.1632 & 0.1225 & 0.0783 & 0.0543 & 0.0465 & 0.034 \\ 
0.05  & 0.1929 & 0.1319 & 0.0464 & 0.0128 & 0.0021 & 1e-04 & 0 & 0 & 0 \\ 
0.1  & 0.1216 & 0.0655 & 0.0119 & 0.0011 & 1e-04 & 0 & 0 & 0 & 0 \\ 
0.15  & 0.0823 & 0.037 & 0.0062 & 6e-04 & 1e-04 & 0 & 0 & 0 & 0 \\ 
0.2  & 0.0564 & 0.025 & 0.0057 & 6e-04 & 1e-04 & 0 & 0 & 0 & 0 \\ 
0.4  & 0.0198 & 0.0064 & 0.0034 & 3e-04 & 0 & 0 & 0 & 0 & 0 \\ 
0.6  & 0.0086 & 0.0032 & 0.0014 & 3e-04 & 0 & 0 & 0 & 0 & 0 \\ 
\hline 
\end{tabular} 
	\begin{tabular}{|c||c|c|c|c|c|c|c|c|c|}\hline 
 &\multicolumn{9}{|c|}{\eqref{samedecision} holds for $h=1$, $i=3$} \\ 
\hline 
$\delta$ & $n$=100 & $n$=200 & $n$=500 & $n$=1000 & $n$=2000 & $n$=4000 & $n$=6000 & $n$=8000 & $n$=10000 \\ 
\hline \hline 
0  & 0.3433 & 0.3245 & 0.3285 & 0.3468 & 0.3643 & 0.4243 & 0.4736 & 0.5183 & 0.5443 \\ 
0.01  & 0.3112 & 0.2662 & 0.2216 & 0.1879 & 0.1332 & 0.0942 & 0.0743 & 0.0607 & 0.0473 \\ 
0.05  & 0.2092 & 0.1416 & 0.05 & 0.0154 & 0.003 & 0 & 2e-04 & 1e-04 & 0 \\ 
0.1  & 0.1383 & 0.0706 & 0.0142 & 0.0016 & 1e-04 & 0 & 0 & 0 & 0 \\ 
0.15  & 0.0958 & 0.0385 & 0.0074 & 5e-04 & 1e-04 & 0 & 0 & 0 & 0 \\ 
0.2  & 0.0701 & 0.0245 & 0.0053 & 5e-04 & 1e-04 & 0 & 0 & 0 & 0 \\ 
0.4  & 0.0261 & 0.0121 & 0.0028 & 4e-04 & 1e-04 & 0 & 0 & 0 & 0 \\ 
0.6  & 0.0146 & 0.0082 & 0.0013 & 2e-04 & 0 & 0 & 0 & 0 & 0 \\ 
\hline 
\end{tabular} 
	\begin{tabular}{|c||c|c|c|c|c|c|c|c|c|}\hline 
 &\multicolumn{9}{|c|}{\eqref{samedecision} holds for $h=5$, $i=2$} \\ 
\hline 
$\delta$ & $n$=100 & $n$=200 & $n$=500 & $n$=1000 & $n$=2000 & $n$=4000 & $n$=6000 & $n$=8000 & $n$=10000 \\ 
\hline \hline 
0  & 0.2534 & 0.2406 & 0.2373 & 0.2411 & 0.2366 & 0.2442 & 0.2542 & 0.271 & 0.2722 \\ 
0.01  & 0.2137 & 0.1868 & 0.1503 & 0.1291 & 0.0894 & 0.0543 & 0.0399 & 0.0339 & 0.0239 \\ 
0.05  & 0.1311 & 0.0915 & 0.046 & 0.0243 & 0.0085 & 0.0017 & 5e-04 & 3e-04 & 1e-04 \\ 
0.1  & 0.0775 & 0.0458 & 0.017 & 0.0057 & 0.001 & 0 & 0 & 0 & 0 \\ 
0.15  & 0.0514 & 0.0274 & 0.0071 & 0.0017 & 0 & 0 & 0 & 0 & 0 \\ 
0.2  & 0.0368 & 0.0171 & 0.0037 & 5e-04 & 0 & 0 & 0 & 0 & 0 \\ 
0.4  & 0.0132 & 0.0045 & 5e-04 & 1e-04 & 0 & 0 & 0 & 0 & 0 \\ 
0.6  & 0.0076 & 0.0019 & 1e-04 & 0 & 0 & 0 & 0 & 0 & 0 \\ 
\hline 
\end{tabular} 
	\begin{tabular}{|c||c|c|c|c|c|c|c|c|c|}\hline 
 &\multicolumn{9}{|c|}{\eqref{samedecision} holds for $h=5$, $i=3$} \\ 
\hline 
$\delta$ & $n$=100 & $n$=200 & $n$=500 & $n$=1000 & $n$=2000 & $n$=4000 & $n$=6000 & $n$=8000 & $n$=10000 \\ 
\hline \hline 
0  & 0.2632 & 0.2509 & 0.2426 & 0.2458 & 0.2564 & 0.2645 & 0.2729 & 0.2913 & 0.2966 \\ 
0.01  & 0.2301 & 0.1995 & 0.1584 & 0.1343 & 0.1019 & 0.0614 & 0.0425 & 0.0393 & 0.027 \\ 
0.05  & 0.1477 & 0.0982 & 0.0525 & 0.0261 & 0.0119 & 0.0023 & 8e-04 & 3e-04 & 5e-04 \\ 
0.1  & 0.0969 & 0.0546 & 0.0217 & 0.0071 & 0.0019 & 1e-04 & 0 & 0 & 0 \\ 
0.15  & 0.068 & 0.0342 & 0.0106 & 0.0026 & 4e-04 & 0 & 0 & 0 & 0 \\ 
0.2  & 0.0508 & 0.0239 & 0.0057 & 9e-04 & 1e-04 & 0 & 0 & 0 & 0 \\ 
0.4  & 0.0208 & 0.0073 & 6e-04 & 0 & 0 & 0 & 0 & 0 & 0 \\ 
0.6  & 0.0109 & 0.0036 & 1e-04 & 0 & 0 & 0 & 0 & 0 & 0 \\ 
\hline 
\end{tabular} 
\caption{\textit{Proportion of \eqref{decisionrule}  being fulfilled for the process \eqref{increasing4} and different values of $h$, $\delta$ and $n$.}}  \label{MSPEanalysisincreasing4}
\end{center}
\end{table}

\clearpage

\begin{figure}
\centering 
\includegraphics[width=0.24\textwidth]{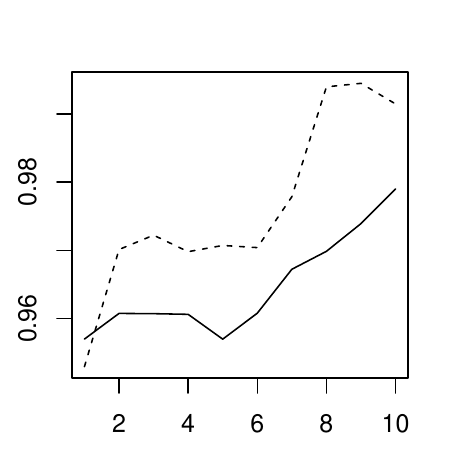}~
\includegraphics[width=0.24\textwidth]{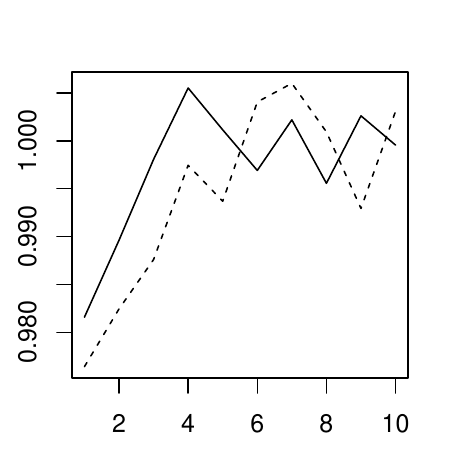}~
\includegraphics[width=0.24\textwidth]{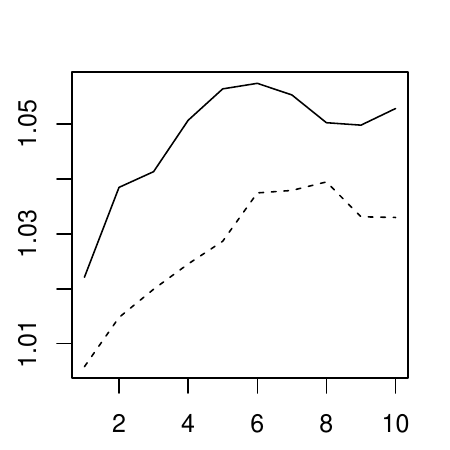}~ 
\includegraphics[width=0.24\textwidth]{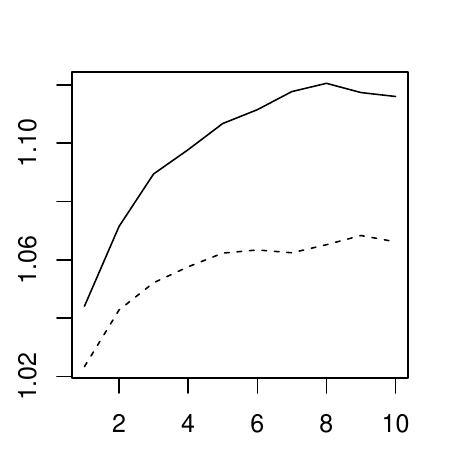}  \\
\includegraphics[width=0.24\textwidth]{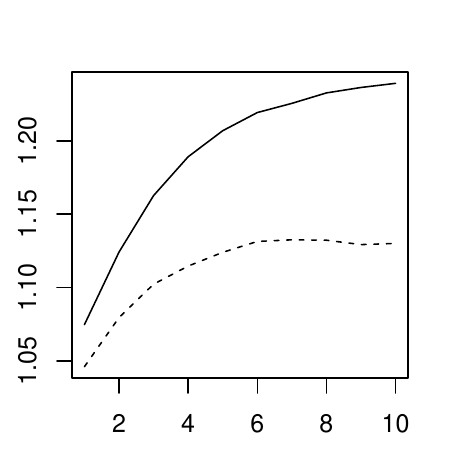}~
\includegraphics[width=0.24\textwidth]{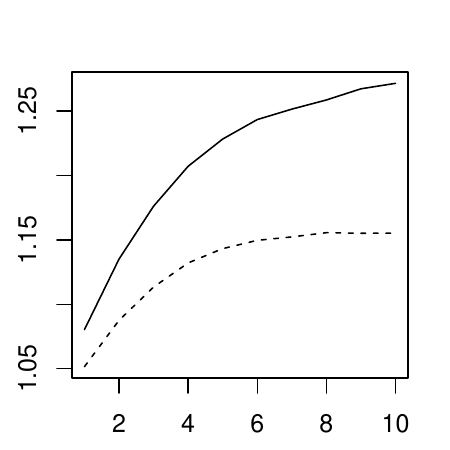}~ 
\includegraphics[width=0.24\textwidth]{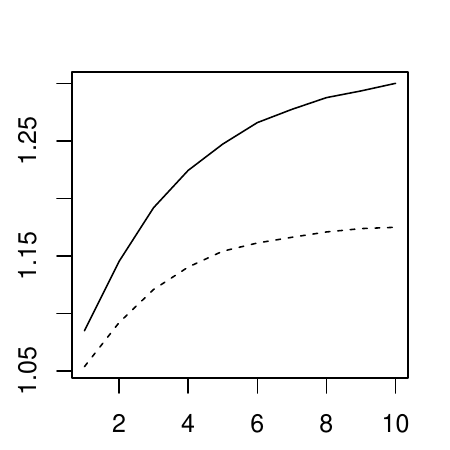}~ 
\includegraphics[width=0.24\textwidth]{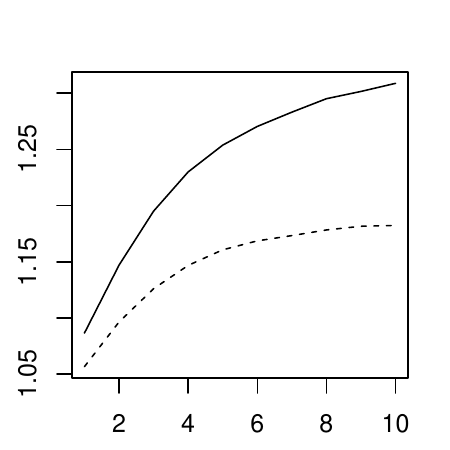}
   \caption{\it Plot of $h \mapsto R_{T,i}(h)$ for model \eqref{increasing1} and different values of $n$ (from left to right: n=100, n=200, n=500, n=1000 [first row], n=4000, n=6000, n=8000, n=10000 [second row]). Solid line: $i=3$ (test set), dashed line: $i=2$ (validation set 2). } \label{MSPEincreasing1}
\end{figure}

\begin{table} \scriptsize
\begin{center}
\begin{tabular}{|c||c|c|c||c|c|c|}
\hline
 \multirow{2}{*}{$n$} & & \multicolumn{2}{|c||}{$R^{{\rm s}, {\rm ls}}_{T,2}(1)$}  & & \multicolumn{2}{|c|}{$R^{{\rm s}, {\rm ls}}_{T,2}(5)$} \\ 
 & & $\geq 1.01$  & $< 1.01$  & & $\geq 1.01$  & $< 1.01$ \\ 
 \hline 
 \multirow{2}{*}{$100$} & $R^{{\rm s}, {\rm ls}}_{T,3}(1) \geq 1.01$  & 0.1612 & 0.1865  & $R^{{\rm s}, {\rm ls}}_{T,3}(5) \geq 1.01$  & 0.1174 & 0.1422 \\ 
& $R^{{\rm s}, {\rm ls}}_{T,3}(1) < 1.01$  & 0.1866 & 0.4657 & $R^{{\rm s}, {\rm ls}}_{T,3}(5) < 1.01$  & 0.1621 & 0.5783 \\

\hline\hline
 \multirow{2}{*}{$n$} & & \multicolumn{2}{|c||}{$R^{{\rm s}, {\rm ls}}_{T,2}(1)$}  & & \multicolumn{2}{|c|}{$R^{{\rm s}, {\rm ls}}_{T,2}(5)$} \\ 
 & & $\geq 1.2$  & $< 1.2$  & & $\geq 1.2$  & $< 1.2$ \\ 
 \hline 
 \multirow{2}{*}{$1000$} & $R^{{\rm s}, {\rm ls}}_{T,3}(1) \geq 1.2$  & 5e-04 & 0.0466  & $R^{{\rm s}, {\rm ls}}_{T,3}(5) \geq 1.2$  & 0.0538 & 0.216 \\ 
& $R^{{\rm s}, {\rm ls}}_{T,3}(1) < 1.2$  & 0.0097 & 0.9432 & $R^{{\rm s}, {\rm ls}}_{T,3}(5) < 1.2$  & 0.0855 & 0.6447 \\

\hline\hline
 \multirow{2}{*}{$n$} & & \multicolumn{2}{|c||}{$R^{{\rm s}, {\rm ls}}_{T,2}(1)$}  & & \multicolumn{2}{|c|}{$R^{{\rm s}, {\rm ls}}_{T,2}(5)$} \\ 
 & & $\geq 1$  & $< 1$  & & $\geq 1$  & $< 1$ \\ 
 \hline 
 \multirow{2}{*}{$10000$} & $R^{{\rm s}, {\rm ls}}_{T,3}(1) \geq 1$  & 0.991 & 0.0063  & $R^{{\rm s}, {\rm ls}}_{T,3}(5) \geq 1$  & 0.9856 & 0.01 \\ 
& $R^{{\rm s}, {\rm ls}}_{T,3}(1) < 1$  & 0.0027 & 0 & $R^{{\rm s}, {\rm ls}}_{T,3}(5) < 1$  & 0.0044 & 0 \\

\hline\hline
 \multirow{2}{*}{$n$} & & \multicolumn{2}{|c||}{$R^{{\rm s}, {\rm ls}}_{T,2}(1)$}  & & \multicolumn{2}{|c|}{$R^{{\rm s}, {\rm ls}}_{T,2}(5)$} \\ 
 & & $\geq 1.05$  & $< 1.05$  & & $\geq 1.05$  & $< 1.05$ \\ 
 \hline 
 \multirow{2}{*}{$10000$} & $R^{{\rm s}, {\rm ls}}_{T,3}(1) \geq 1.05$  & 0.4203 & 0.3675  & $R^{{\rm s}, {\rm ls}}_{T,3}(5) \geq 1.05$  & 0.8866 & 0.0836 \\ 
& $R^{{\rm s}, {\rm ls}}_{T,3}(1) < 1.05$  & 0.1236 & 0.0886 & $R^{{\rm s}, {\rm ls}}_{T,3}(5) < 1.05$  & 0.0285 & 0.0013 \\

\hline\hline
 \multirow{2}{*}{$n$} & & \multicolumn{2}{|c||}{$R^{{\rm s}, {\rm ls}}_{T,2}(1)$}  & & \multicolumn{2}{|c|}{$R^{{\rm s}, {\rm ls}}_{T,2}(5)$} \\ 
 & & $\geq 1.2$  & $< 1.2$  & & $\geq 1.2$  & $< 1.2$ \\ 
 \hline 
 \multirow{2}{*}{$10000$} & $R^{{\rm s}, {\rm ls}}_{T,3}(1) \geq 1.2$  & 1e-04 & 0.0295  & $R^{{\rm s}, {\rm ls}}_{T,3}(5) \geq 1.2$  & 0.1961 & 0.4468 \\ 
& $R^{{\rm s}, {\rm ls}}_{T,3}(1) < 1.2$  & 3e-04 & 0.9701 & $R^{{\rm s}, {\rm ls}}_{T,3}(5) < 1.2$  & 0.1251 & 0.232 \\

\hline\hline
 \multirow{2}{*}{$n$} & & \multicolumn{2}{|c||}{$R^{{\rm s}, {\rm ls}}_{T,2}(1)$}  & & \multicolumn{2}{|c|}{$R^{{\rm s}, {\rm ls}}_{T,2}(5)$} \\ 
 & & $\geq 1.4$  & $< 1.4$  & & $\geq 1.4$  & $< 1.4$ \\ 
 \hline 
 \multirow{2}{*}{$10000$} & $R^{{\rm s}, {\rm ls}}_{T,3}(1) \geq 1.4$  & 0 & 0  & $R^{{\rm s}, {\rm ls}}_{T,3}(5) \geq 1.4$  & 0.003 & 0.1818 \\ 
& $R^{{\rm s}, {\rm ls}}_{T,3}(1) < 1.4$  & 0 & 1 & $R^{{\rm s}, {\rm ls}}_{T,3}(5) < 1.4$  & 0.0158 & 0.7994 \\

\hline
\end{tabular}
\caption{\textit{Proportions of the individual events in~\eqref{samedecision} for the process \eqref{increasing1} and selected combinations of $n$ and $\delta$.}}  \label{MSPEanalysisincreasing1c}
\end{center}
\end{table}

\begin{table} \scriptsize
\begin{center}
	\begin{tabular}{|c||c|c|c|c|c|c|c|c|c|}\hline 
 &\multicolumn{9}{|c|}{\eqref{samedecision} holds for $h=1$} \\ 
\hline 
$\delta$ & $n$=100 & $n$=200 & $n$=500 & $n$=1000 & $n$=2000 & $n$=4000 & $n$=6000 & $n$=8000 & $n$=10000 \\ 
\hline \hline 
0  & 0.6103 & 0.5906 & 0.6028 & 0.6474 & 0.7636 & 0.9114 & 0.9647 & 0.9828 & 0.991 \\ 
0.01  & 0.6269 & 0.6132 & 0.5908 & 0.5882 & 0.6612 & 0.8098 & 0.8901 & 0.9322 & 0.9591 \\ 
0.05  & 0.7148 & 0.7239 & 0.6838 & 0.5815 & 0.5033 & 0.4719 & 0.4684 & 0.481 & 0.5089 \\ 
0.1  & 0.8014 & 0.8292 & 0.81 & 0.7449 & 0.6924 & 0.6393 & 0.617 & 0.5947 & 0.6006 \\ 
0.15  & 0.8638 & 0.9006 & 0.8996 & 0.8727 & 0.856 & 0.8489 & 0.8616 & 0.8593 & 0.8755 \\ 
0.2  & 0.9099 & 0.9415 & 0.9469 & 0.9437 & 0.9386 & 0.9471 & 0.9604 & 0.9615 & 0.9702 \\ 
0.4  & 0.9784 & 0.9947 & 0.9968 & 0.9986 & 0.9991 & 0.9996 & 1 & 0.9999 & 1 \\ 
0.6  & 0.9927 & 0.9995 & 0.9999 & 1 & 1 & 1 & 1 & 1 & 1 \\ 
\hline 
\end{tabular} 
	\begin{tabular}{|c||c|c|c|c|c|c|c|c|c|}\hline 
 &\multicolumn{9}{|c|}{\eqref{samedecision} holds for $h=5$} \\ 
\hline 
$\delta$ & $n$=100 & $n$=200 & $n$=500 & $n$=1000 & $n$=2000 & $n$=4000 & $n$=6000 & $n$=8000 & $n$=10000 \\ 
\hline \hline 
0  & 0.6906 & 0.6455 & 0.6048 & 0.639 & 0.7651 & 0.9016 & 0.9561 & 0.9744 & 0.9856 \\ 
0.01  & 0.6957 & 0.6537 & 0.6008 & 0.6138 & 0.7218 & 0.868 & 0.9328 & 0.9611 & 0.9766 \\ 
0.05  & 0.734 & 0.7067 & 0.626 & 0.5559 & 0.5839 & 0.6997 & 0.7902 & 0.8502 & 0.8879 \\ 
0.1  & 0.7807 & 0.7697 & 0.6889 & 0.5718 & 0.5168 & 0.5357 & 0.5843 & 0.628 & 0.685 \\ 
0.15  & 0.8217 & 0.8279 & 0.751 & 0.6288 & 0.53 & 0.4655 & 0.4703 & 0.4724 & 0.4989 \\ 
0.2  & 0.8564 & 0.8634 & 0.8064 & 0.6985 & 0.5832 & 0.4974 & 0.4573 & 0.4374 & 0.4281 \\ 
0.4  & 0.9295 & 0.9412 & 0.9294 & 0.8979 & 0.8512 & 0.818 & 0.8109 & 0.7998 & 0.8024 \\ 
0.6  & 0.9568 & 0.9691 & 0.973 & 0.9708 & 0.9619 & 0.9561 & 0.9628 & 0.9624 & 0.9655 \\ 
\hline 
\end{tabular}  
\caption{\textit{Proportion of \eqref{samedecision}  being fulfilled for the process \eqref{increasing1} and different values of $h$, $\delta$ and $n$.}}  \label{MSPEanalysisincreasing1b}
\end{center}
\end{table}

\begin{table} \scriptsize
\begin{center}
	\begin{tabular}{|c||c|c|c|c|c|c|c|c|c|}\hline 
 &\multicolumn{9}{|p{12.5cm}|}{\centering Value of $q(\delta)$, defined in \eqref{cond:f}, for $h=1$} \\ 
\hline 
$\delta$ & $n$=100 & $n$=200 & $n$=500 & $n$=1000 & $n$=2000 & $n$=4000 & $n$=6000 & $n$=8000 & $n$=10000 \\ 
\hline \hline 
0  & 0 & 0 & 0 & 0 & 0 & 0 & 0 & 0 & 0 \\ 
0.01  & 0.028 & 0.038 & 0.048 & 0.055 & 0.061 & 0.067 & 0.069 & 0.071 & 0.073 \\ 
0.05  & 1.4e-06 & 2.3e-05 & 0.0078 & 0.015 & 0.021 & 0.026 & 0.029 & 0.031 & 0.032 \\ 
0.1  & 0.041 & 0.033 & 0.024 & 0.017 & 0.012 & 0.0069 & 0.0041 & 0.0022 & 8e-04 \\ 
0.15  & 0.091 & 0.083 & 0.074 & 0.068 & 0.062 & 0.057 & 0.054 & 0.052 & 0.051 \\ 
0.2  & 0.14 & 0.13 & 0.12 & 0.12 & 0.11 & 0.11 & 0.1 & 0.1 & 0.1 \\ 
0.4  & 0.34 & 0.33 & 0.33 & 0.32 & 0.31 & 0.31 & 0.3 & 0.3 & 0.3 \\ 
0.6  & 0.55 & 0.54 & 0.53 & 0.52 & 0.51 & 0.51 & 0.51 & 0.5 & 0.5 \\ 
\hline 
\end{tabular} 
	\begin{tabular}{|c||c|c|c|c|c|c|c|c|c|}\hline 
 &\multicolumn{9}{|p{12.5cm}|}{\centering Value of $q(\delta)$, defined in \eqref{cond:f}, for $h=5$} \\ 
\hline 
$\delta$ & $n$=100 & $n$=200 & $n$=500 & $n$=1000 & $n$=2000 & $n$=4000 & $n$=6000 & $n$=8000 & $n$=10000 \\ 
\hline \hline 
0  & 0 & 0 & 0 & 0 & 0 & 0 & 0 & 0 & 0 \\ 
0.01  & 0.16 & 0.2 & 0.24 & 0.27 & 0.29 & 0.3 & 0.31 & 0.32 & 0.32 \\ 
0.05  & 0.097 & 0.28 & 0.52 & 0.67 & 0.79 & 0.89 & 0.95 & 0.98 & 1 \\ 
0.1  & 0.00045 & 0.045 & 0.28 & 0.42 & 0.55 & 0.65 & 0.71 & 0.74 & 0.77 \\ 
0.15  & 0.023 & 8.6e-05 & 0.035 & 0.18 & 0.3 & 0.41 & 0.47 & 0.5 & 0.53 \\ 
0.2  & 0.25 & 0.013 & 0.00031 & 5e-05 & 0.063 & 0.17 & 0.23 & 0.26 & 0.29 \\ 
0.4  & 1.1 & 0.93 & 0.69 & 0.53 & 0.4 & 0.29 & 0.23 & 0.19 & 0.16 \\ 
0.6  & 2 & 1.9 & 1.6 & 1.5 & 1.3 & 1.2 & 1.2 & 1.1 & 1.1 \\ 
\hline 
\end{tabular} 
\caption{\textit{Values of $q(\delta)$, defined in \eqref{cond:f}, for the process \eqref{increasing1} and different values of $h$, $\delta$ and $n$.}} \label{MSPEanalysisincreasing1d}
\end{center}
\end{table}

\begin{table} \scriptsize
\begin{center}
	\begin{tabular}{|c||c|c|c|c|c|c|c|c|c|}\hline 
 &\multicolumn{9}{|c|}{\eqref{samedecision} holds for $h=1$, $i=2$} \\ 
\hline 
$\delta$ & $n$=100 & $n$=200 & $n$=500 & $n$=1000 & $n$=2000 & $n$=4000 & $n$=6000 & $n$=8000 & $n$=10000 \\ 
\hline \hline 
0  & 0.3857 & 0.3985 & 0.5044 & 0.6621 & 0.8162 & 0.9405 & 0.9752 & 0.9879 & 0.9937 \\ 
0.01  & 0.3478 & 0.3413 & 0.4187 & 0.5562 & 0.714 & 0.8642 & 0.9216 & 0.9491 & 0.9701 \\ 
0.05  & 0.2269 & 0.1939 & 0.2046 & 0.2576 & 0.3306 & 0.4176 & 0.4771 & 0.5101 & 0.5439 \\ 
0.1  & 0.1479 & 0.1037 & 0.0887 & 0.1001 & 0.1057 & 0.108 & 0.1111 & 0.1054 & 0.0976 \\ 
0.15  & 0.0953 & 0.0553 & 0.0352 & 0.0336 & 0.0298 & 0.0236 & 0.0176 & 0.0157 & 0.0095 \\ 
0.2  & 0.0616 & 0.0318 & 0.0146 & 0.0102 & 0.0086 & 0.0038 & 0.0017 & 0.0013 & 4e-04 \\ 
0.4  & 0.0185 & 0.0049 & 9e-04 & 1e-04 & 0 & 0 & 0 & 0 & 0 \\ 
0.6  & 0.01 & 0.0026 & 1e-04 & 0 & 0 & 0 & 0 & 0 & 0 \\ 
\hline 
\end{tabular} 
	\begin{tabular}{|c||c|c|c|c|c|c|c|c|c|}\hline 
 &\multicolumn{9}{|c|}{\eqref{samedecision} holds for $h=1$, $i=3$} \\ 
\hline 
$\delta$ & $n$=100 & $n$=200 & $n$=500 & $n$=1000 & $n$=2000 & $n$=4000 & $n$=6000 & $n$=8000 & $n$=10000 \\ 
\hline \hline 
0  & 0.3902 & 0.4357 & 0.577 & 0.7529 & 0.8978 & 0.9651 & 0.9887 & 0.9949 & 0.9973 \\ 
0.01  & 0.3477 & 0.3799 & 0.4993 & 0.6656 & 0.8238 & 0.9214 & 0.9633 & 0.9823 & 0.989 \\ 
0.05  & 0.2357 & 0.23 & 0.2884 & 0.4085 & 0.5153 & 0.6419 & 0.7087 & 0.7563 & 0.7878 \\ 
0.1  & 0.1531 & 0.1301 & 0.1485 & 0.209 & 0.2563 & 0.3073 & 0.3301 & 0.3597 & 0.3608 \\ 
0.15  & 0.1003 & 0.0719 & 0.0786 & 0.1029 & 0.118 & 0.1303 & 0.1232 & 0.1268 & 0.1162 \\ 
0.2  & 0.0703 & 0.0409 & 0.0425 & 0.0471 & 0.053 & 0.0491 & 0.0379 & 0.0372 & 0.0296 \\ 
0.4  & 0.0225 & 0.0058 & 0.0023 & 0.0013 & 9e-04 & 4e-04 & 0 & 1e-04 & 0 \\ 
0.6  & 0.0129 & 0.0025 & 0 & 0 & 0 & 0 & 0 & 0 & 0 \\ 
\hline 
\end{tabular} 
	\begin{tabular}{|c||c|c|c|c|c|c|c|c|c|}\hline 
 &\multicolumn{9}{|c|}{\eqref{samedecision} holds for $h=5$, $i=2$} \\ 
\hline 
$\delta$ & $n$=100 & $n$=200 & $n$=500 & $n$=1000 & $n$=2000 & $n$=4000 & $n$=6000 & $n$=8000 & $n$=10000 \\ 
\hline \hline 
0  & 0.2969 & 0.3562 & 0.5013 & 0.6562 & 0.8177 & 0.9382 & 0.9726 & 0.9826 & 0.99 \\ 
0.01  & 0.2795 & 0.3299 & 0.4569 & 0.6096 & 0.7778 & 0.9147 & 0.9549 & 0.9727 & 0.9841 \\ 
0.05  & 0.2133 & 0.2407 & 0.3265 & 0.4452 & 0.613 & 0.7689 & 0.8485 & 0.888 & 0.9151 \\ 
0.1  & 0.1573 & 0.169 & 0.2154 & 0.3025 & 0.4285 & 0.5654 & 0.6474 & 0.693 & 0.7385 \\ 
0.15  & 0.1191 & 0.1198 & 0.1467 & 0.2048 & 0.285 & 0.3779 & 0.4449 & 0.486 & 0.5217 \\ 
0.2  & 0.0907 & 0.0867 & 0.1049 & 0.1393 & 0.1903 & 0.2454 & 0.2882 & 0.3035 & 0.3212 \\ 
0.4  & 0.0383 & 0.0334 & 0.0284 & 0.0275 & 0.0326 & 0.0308 & 0.0256 & 0.0216 & 0.0188 \\ 
0.6  & 0.0228 & 0.0187 & 0.0128 & 0.0048 & 0.0037 & 0.0031 & 0.0018 & 9e-04 & 2e-04 \\ 
\hline 
\end{tabular} 
	\begin{tabular}{|c||c|c|c|c|c|c|c|c|c|}\hline 
 &\multicolumn{9}{|c|}{\eqref{samedecision} holds for $h=5$, $i=3$} \\ 
\hline 
$\delta$ & $n$=100 & $n$=200 & $n$=500 & $n$=1000 & $n$=2000 & $n$=4000 & $n$=6000 & $n$=8000 & $n$=10000 \\ 
\hline \hline 
0  & 0.2763 & 0.3671 & 0.5549 & 0.7348 & 0.8884 & 0.958 & 0.9831 & 0.9916 & 0.9956 \\ 
0.01  & 0.2596 & 0.3426 & 0.5229 & 0.7018 & 0.8606 & 0.9449 & 0.9765 & 0.9878 & 0.9925 \\ 
0.05  & 0.2097 & 0.2662 & 0.4083 & 0.5761 & 0.7483 & 0.8722 & 0.9267 & 0.9558 & 0.9702 \\ 
0.1  & 0.1646 & 0.1959 & 0.3017 & 0.4549 & 0.6091 & 0.7535 & 0.8261 & 0.8702 & 0.8975 \\ 
0.15  & 0.1266 & 0.1445 & 0.2289 & 0.3514 & 0.4856 & 0.6266 & 0.691 & 0.746 & 0.7852 \\ 
0.2  & 0.1017 & 0.1127 & 0.1729 & 0.2698 & 0.3841 & 0.4982 & 0.5581 & 0.6111 & 0.6429 \\ 
0.4  & 0.0502 & 0.0496 & 0.0618 & 0.0844 & 0.1274 & 0.1596 & 0.1707 & 0.1848 & 0.1848 \\ 
0.6  & 0.033 & 0.0292 & 0.0248 & 0.0272 & 0.0348 & 0.0408 & 0.0358 & 0.0367 & 0.0343 \\ 
\hline 
\end{tabular} 
\caption{\textit{Proportion of \eqref{decisionrule}  being fulfilled for the process \eqref{increasing1} and different values of $h$, $\delta$ and $n$.}}  \label{MSPEanalysisincreasing1}
\end{center}
\end{table}

\clearpage

\begin{figure}
\centering 
\includegraphics[width=0.24\textwidth]{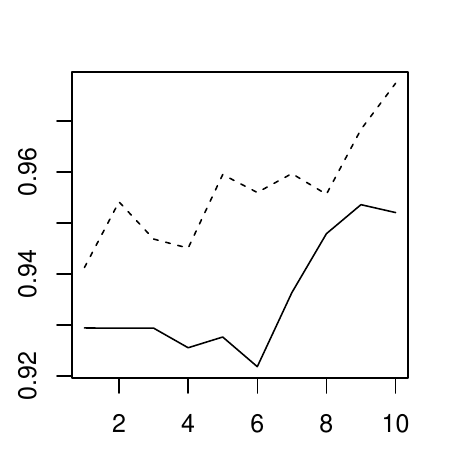}~
\includegraphics[width=0.24\textwidth]{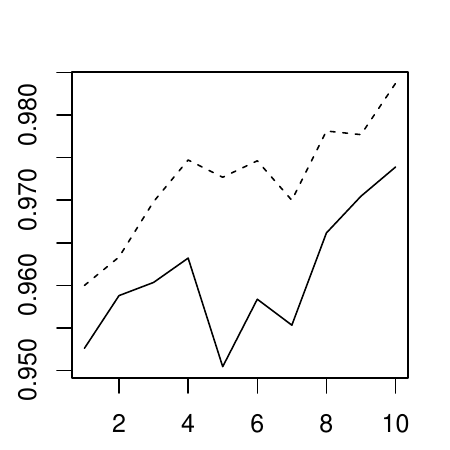}~
\includegraphics[width=0.24\textwidth]{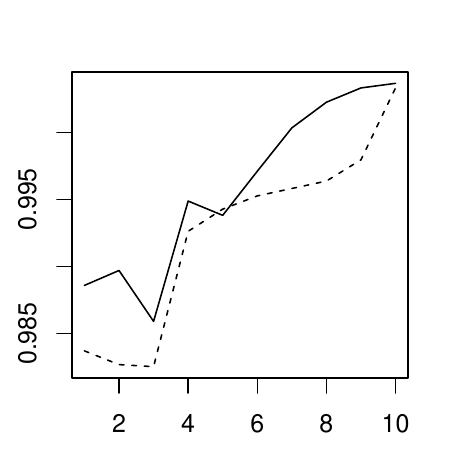}~ 
\includegraphics[width=0.24\textwidth]{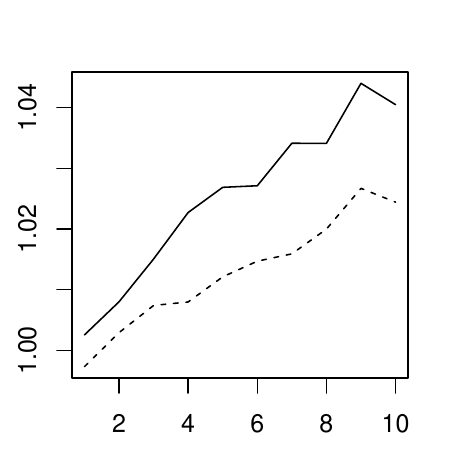} \\
\includegraphics[width=0.24\textwidth]{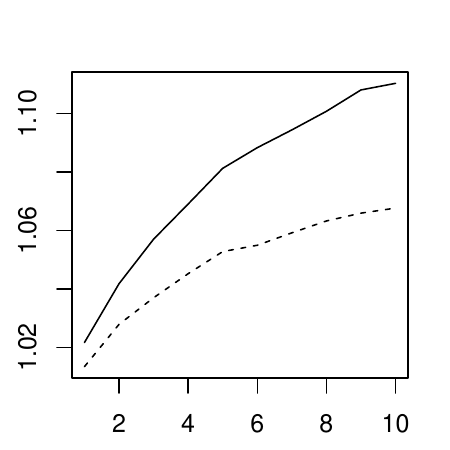}~
\includegraphics[width=0.24\textwidth]{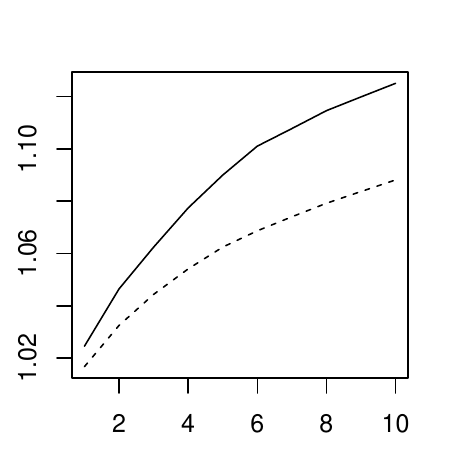}~ 
\includegraphics[width=0.24\textwidth]{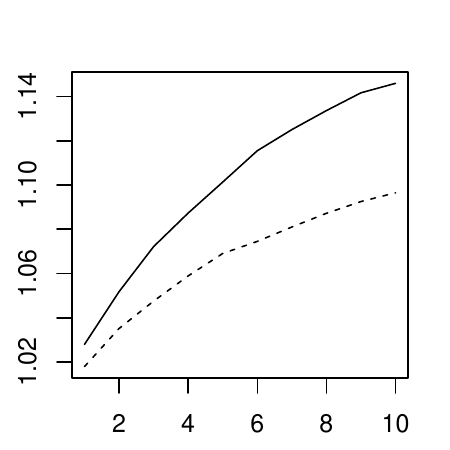}~ 
\includegraphics[width=0.24\textwidth]{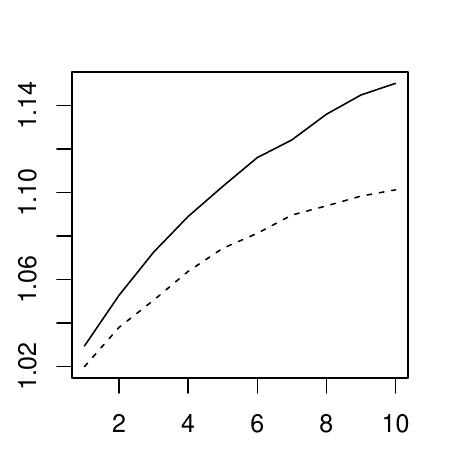}
   \caption{\it Plot of $h \mapsto R_{T,i}(h)$ for model \eqref{increasing3} and different values of $n$ (from left to right: n=100, n=200, n=500, n=1000 [first row], n=4000, n=6000, n=8000, n=10000 [second row]). Solid line: $i=3$ (test set), dashed line: $i=2$ (validation set 2). } \label{MSPEincreasing3}
\end{figure}

\begin{table} \scriptsize
\begin{center}
\begin{tabular}{|c||c|c|c||c|c|c|}
\hline
 \multirow{2}{*}{$n$} & & \multicolumn{2}{|c||}{$R^{{\rm s}, {\rm ls}}_{T,2}(1)$}  & & \multicolumn{2}{|c|}{$R^{{\rm s}, {\rm ls}}_{T,2}(5)$} \\ 
 & & $\geq 1.01$  & $< 1.01$  & & $\geq 1.01$  & $< 1.01$ \\ 
 \hline 
 \multirow{2}{*}{$100$} & $R^{{\rm s}, {\rm ls}}_{T,3}(1) \geq 1.01$  & 0.1309 & 0.1852  & $R^{{\rm s}, {\rm ls}}_{T,3}(5) \geq 1.01$  & 0.126 & 0.1506 \\ 
& $R^{{\rm s}, {\rm ls}}_{T,3}(1) < 1.01$  & 0.1886 & 0.4953 & $R^{{\rm s}, {\rm ls}}_{T,3}(5) < 1.01$  & 0.1741 & 0.5493 \\

\hline\hline
 \multirow{2}{*}{$n$} & & \multicolumn{2}{|c||}{$R^{{\rm s}, {\rm ls}}_{T,2}(1)$}  & & \multicolumn{2}{|c|}{$R^{{\rm s}, {\rm ls}}_{T,2}(5)$} \\ 
 & & $\geq 1.2$  & $< 1.2$  & & $\geq 1.2$  & $< 1.2$ \\ 
 \hline 
 \multirow{2}{*}{$1000$} & $R^{{\rm s}, {\rm ls}}_{T,3}(1) \geq 1.2$  & 0 & 0.0027  & $R^{{\rm s}, {\rm ls}}_{T,3}(5) \geq 1.2$  & 0.0131 & 0.0819 \\ 
& $R^{{\rm s}, {\rm ls}}_{T,3}(1) < 1.2$  & 0.0015 & 0.9958 & $R^{{\rm s}, {\rm ls}}_{T,3}(5) < 1.2$  & 0.0412 & 0.8638 \\

\hline\hline
 \multirow{2}{*}{$n$} & & \multicolumn{2}{|c||}{$R^{{\rm s}, {\rm ls}}_{T,2}(1)$}  & & \multicolumn{2}{|c|}{$R^{{\rm s}, {\rm ls}}_{T,2}(5)$} \\ 
 & & $\geq 1$  & $< 1$  & & $\geq 1$  & $< 1$ \\ 
 \hline 
 \multirow{2}{*}{$10000$} & $R^{{\rm s}, {\rm ls}}_{T,3}(1) \geq 1$  & 0.8792 & 0.0782  & $R^{{\rm s}, {\rm ls}}_{T,3}(5) \geq 1$  & 0.8927 & 0.0652 \\ 
& $R^{{\rm s}, {\rm ls}}_{T,3}(1) < 1$  & 0.0394 & 0.0032 & $R^{{\rm s}, {\rm ls}}_{T,3}(5) < 1$  & 0.0398 & 0.0023 \\

\hline\hline
 \multirow{2}{*}{$n$} & & \multicolumn{2}{|c||}{$R^{{\rm s}, {\rm ls}}_{T,2}(1)$}  & & \multicolumn{2}{|c|}{$R^{{\rm s}, {\rm ls}}_{T,2}(5)$} \\ 
 & & $\geq 1.05$  & $< 1.05$  & & $\geq 1.05$  & $< 1.05$ \\ 
 \hline 
 \multirow{2}{*}{$10000$} & $R^{{\rm s}, {\rm ls}}_{T,3}(1) \geq 1.05$  & 0.0094 & 0.1685  & $R^{{\rm s}, {\rm ls}}_{T,3}(5) \geq 1.05$  & 0.4644 & 0.2838 \\ 
& $R^{{\rm s}, {\rm ls}}_{T,3}(1) < 1.05$  & 0.0511 & 0.771 & $R^{{\rm s}, {\rm ls}}_{T,3}(5) < 1.05$  & 0.1504 & 0.1014 \\

\hline\hline
 \multirow{2}{*}{$n$} & & \multicolumn{2}{|c||}{$R^{{\rm s}, {\rm ls}}_{T,2}(1)$}  & & \multicolumn{2}{|c|}{$R^{{\rm s}, {\rm ls}}_{T,2}(5)$} \\ 
 & & $\geq 1.1$  & $< 1.1$  & & $\geq 1.1$  & $< 1.1$ \\ 
 \hline 
 \multirow{2}{*}{$10000$} & $R^{{\rm s}, {\rm ls}}_{T,3}(1) \geq 1.1$  & 0 & 0.0124  & $R^{{\rm s}, {\rm ls}}_{T,3}(5) \geq 1.1$  & 0.1515 & 0.337 \\ 
& $R^{{\rm s}, {\rm ls}}_{T,3}(1) < 1.1$  & 8e-04 & 0.9868 & $R^{{\rm s}, {\rm ls}}_{T,3}(5) < 1.1$  & 0.1561 & 0.3554 \\

\hline\hline
 \multirow{2}{*}{$n$} & & \multicolumn{2}{|c||}{$R^{{\rm s}, {\rm ls}}_{T,2}(1)$}  & & \multicolumn{2}{|c|}{$R^{{\rm s}, {\rm ls}}_{T,2}(5)$} \\ 
 & & $\geq 1.2$  & $< 1.2$  & & $\geq 1.2$  & $< 1.2$ \\ 
 \hline 
 \multirow{2}{*}{$10000$} & $R^{{\rm s}, {\rm ls}}_{T,3}(1) \geq 1.2$  & 0 & 0  & $R^{{\rm s}, {\rm ls}}_{T,3}(5) \geq 1.2$  & 0.0068 & 0.1493 \\ 
& $R^{{\rm s}, {\rm ls}}_{T,3}(1) < 1.2$  & 0 & 1 & $R^{{\rm s}, {\rm ls}}_{T,3}(5) < 1.2$  & 0.0414 & 0.8025 \\

\hline
\end{tabular}
\caption{\textit{Proportions of the individual events in~\eqref{samedecision} for the process \eqref{increasing3} and selected combinations of $n$ and $\delta$.}}  \label{MSPEanalysisincreasing3c}
\end{center}
\end{table}

\begin{table} \scriptsize
\begin{center}
	\begin{tabular}{|c||c|c|c|c|c|c|c|c|c|}\hline 
 &\multicolumn{9}{|c|}{\eqref{samedecision} holds for $h=1$} \\ 
\hline 
$\delta$ & $n$=100 & $n$=200 & $n$=500 & $n$=1000 & $n$=2000 & $n$=4000 & $n$=6000 & $n$=8000 & $n$=10000 \\ 
\hline \hline 
0  & 0.6055 & 0.5989 & 0.6087 & 0.5793 & 0.6071 & 0.6909 & 0.7765 & 0.8344 & 0.8824 \\ 
0.01  & 0.6262 & 0.6343 & 0.6616 & 0.6224 & 0.5773 & 0.5383 & 0.5659 & 0.5743 & 0.6261 \\ 
0.05  & 0.7433 & 0.7965 & 0.8423 & 0.8384 & 0.8143 & 0.7865 & 0.7837 & 0.7717 & 0.7804 \\ 
0.1  & 0.8363 & 0.9007 & 0.9465 & 0.9531 & 0.9568 & 0.9666 & 0.9809 & 0.9817 & 0.9868 \\ 
0.15  & 0.8964 & 0.952 & 0.9795 & 0.9849 & 0.9917 & 0.997 & 0.9989 & 0.9986 & 0.9995 \\ 
0.2  & 0.9322 & 0.9768 & 0.9918 & 0.9958 & 0.9987 & 0.9999 & 0.9998 & 1 & 1 \\ 
0.4  & 0.9867 & 0.9984 & 0.9999 & 1 & 1 & 1 & 1 & 1 & 1 \\ 
0.6  & 0.9948 & 1 & 1 & 1 & 1 & 1 & 1 & 1 & 1 \\ 
\hline 
\end{tabular} 
	\begin{tabular}{|c||c|c|c|c|c|c|c|c|c|}\hline 
 &\multicolumn{9}{|c|}{\eqref{samedecision} holds for $h=5$} \\ 
\hline 
$\delta$ & $n$=100 & $n$=200 & $n$=500 & $n$=1000 & $n$=2000 & $n$=4000 & $n$=6000 & $n$=8000 & $n$=10000 \\ 
\hline \hline 
0  & 0.672 & 0.6272 & 0.5775 & 0.5537 & 0.5976 & 0.7113 & 0.7987 & 0.8569 & 0.895 \\ 
0.01  & 0.6753 & 0.6333 & 0.5932 & 0.5594 & 0.5672 & 0.6393 & 0.7195 & 0.7776 & 0.8268 \\ 
0.05  & 0.7132 & 0.6955 & 0.6936 & 0.6496 & 0.5712 & 0.5199 & 0.5253 & 0.5265 & 0.5658 \\ 
0.1  & 0.7572 & 0.7733 & 0.7924 & 0.7437 & 0.6651 & 0.5881 & 0.5406 & 0.5087 & 0.5069 \\ 
0.15  & 0.7982 & 0.834 & 0.8545 & 0.8204 & 0.7681 & 0.7138 & 0.6788 & 0.6568 & 0.6603 \\ 
0.2  & 0.8322 & 0.8742 & 0.8983 & 0.8769 & 0.8428 & 0.8156 & 0.813 & 0.7988 & 0.8093 \\ 
0.4  & 0.92 & 0.9515 & 0.9731 & 0.9769 & 0.9786 & 0.9802 & 0.9878 & 0.9888 & 0.9906 \\ 
0.6  & 0.956 & 0.9732 & 0.9895 & 0.9953 & 0.9975 & 0.9988 & 0.9994 & 0.9995 & 0.9998 \\ 
\hline 
\end{tabular} 
\caption{\textit{Proportion of \eqref{samedecision}  being fulfilled for the process \eqref{increasing3} and different values of $h$, $\delta$ and $n$.}}  \label{MSPEanalysisincreasing3b}
\end{center}
\end{table}

\begin{table} \scriptsize
\begin{center}
	\begin{tabular}{|c||c|c|c|c|c|c|c|c|c|}\hline 
 &\multicolumn{9}{|p{12.5cm}|}{\centering Value of $q(\delta)$, defined in \eqref{cond:f}, for $h=1$} \\ 
\hline 
$\delta$ & $n$=100 & $n$=200 & $n$=500 & $n$=1000 & $n$=2000 & $n$=4000 & $n$=6000 & $n$=8000 & $n$=10000 \\ 
\hline \hline 
0  & 0 & 0 & 0 & 0 & 0 & 0 & 0 & 0 & 0 \\ 
0.01  & 0.0049 & 0.0081 & 0.012 & 0.014 & 0.016 & 0.018 & 0.018 & 0.019 & 0.019 \\ 
0.05  & 0.029 & 0.027 & 0.024 & 0.022 & 0.021 & 0.019 & 0.019 & 0.018 & 0.018 \\ 
0.1  & 0.079 & 0.077 & 0.074 & 0.072 & 0.071 & 0.069 & 0.069 & 0.068 & 0.068 \\ 
0.15  & 0.13 & 0.13 & 0.12 & 0.12 & 0.12 & 0.12 & 0.12 & 0.12 & 0.12 \\ 
0.2  & 0.18 & 0.18 & 0.17 & 0.17 & 0.17 & 0.17 & 0.17 & 0.17 & 0.17 \\ 
0.4  & 0.38 & 0.38 & 0.37 & 0.37 & 0.37 & 0.37 & 0.37 & 0.37 & 0.37 \\ 
0.6  & 0.58 & 0.58 & 0.57 & 0.57 & 0.57 & 0.57 & 0.57 & 0.57 & 0.57 \\ 
\hline 
\end{tabular} 
	\begin{tabular}{|c||c|c|c|c|c|c|c|c|c|}\hline 
 &\multicolumn{9}{|p{12.5cm}|}{\centering Value of $q(\delta)$, defined in \eqref{cond:f}, for $h=5$} \\ 
\hline 
$\delta$ & $n$=100 & $n$=200 & $n$=500 & $n$=1000 & $n$=2000 & $n$=4000 & $n$=6000 & $n$=8000 & $n$=10000 \\ 
\hline \hline 
0  & 0 & 0 & 0 & 0 & 0 & 0 & 0 & 0 & 0 \\ 
0.01  & 0.15 & 0.23 & 0.33 & 0.35 & 0.37 & 0.38 & 0.38 & 0.39 & 0.39 \\ 
0.05  & 8.7e-05 & 0.04 & 0.13 & 0.19 & 0.23 & 0.27 & 0.29 & 0.3 & 0.31 \\ 
0.1  & 0.16 & 0.073 & 3.8e-05 & 1.4e-06 & 0.00011 & 0.031 & 0.05 & 0.062 & 0.071 \\ 
0.15  & 0.39 & 0.31 & 0.23 & 0.18 & 0.14 & 0.1 & 0.088 & 0.077 & 0.069 \\ 
0.2  & 0.62 & 0.54 & 0.46 & 0.41 & 0.37 & 0.34 & 0.33 & 0.31 & 0.31 \\ 
0.4  & 1.6 & 1.5 & 1.4 & 1.4 & 1.3 & 1.3 & 1.3 & 1.3 & 1.3 \\ 
0.6  & 2.5 & 2.4 & 2.4 & 2.3 & 2.3 & 2.2 & 2.2 & 2.2 & 2.2 \\ 
\hline 
\end{tabular} 
\caption{\textit{Values of $q(\delta)$, defined in \eqref{cond:f}, for the process \eqref{increasing3} and different values of $h$, $\delta$ and $n$.}} \label{MSPEanalysisincreasing3d}
\end{center}
\end{table}

\begin{table} \scriptsize
\begin{center}
	\begin{tabular}{|c||c|c|c|c|c|c|c|c|c|}\hline 
 &\multicolumn{9}{|c|}{\eqref{samedecision} holds for $h=1$, $i=2$} \\ 
\hline 
$\delta$ & $n$=100 & $n$=200 & $n$=500 & $n$=1000 & $n$=2000 & $n$=4000 & $n$=6000 & $n$=8000 & $n$=10000 \\ 
\hline \hline 
0  & 0.367 & 0.3623 & 0.3821 & 0.4568 & 0.5783 & 0.7451 & 0.8346 & 0.8822 & 0.9186 \\ 
0.01  & 0.3195 & 0.2901 & 0.2728 & 0.312 & 0.3885 & 0.5035 & 0.5942 & 0.6439 & 0.69 \\ 
0.05  & 0.1965 & 0.1386 & 0.0932 & 0.0811 & 0.0806 & 0.0752 & 0.0724 & 0.0712 & 0.0605 \\ 
0.1  & 0.1164 & 0.062 & 0.0243 & 0.0172 & 0.0118 & 0.0068 & 0.0028 & 0.0026 & 8e-04 \\ 
0.15  & 0.072 & 0.0307 & 0.0073 & 0.0048 & 0.0018 & 1e-04 & 0 & 0 & 0 \\ 
0.2  & 0.0474 & 0.015 & 0.003 & 0.0015 & 3e-04 & 0 & 0 & 0 & 0 \\ 
0.4  & 0.0135 & 0.0028 & 1e-04 & 0 & 0 & 0 & 0 & 0 & 0 \\ 
0.6  & 0.0078 & 0.0021 & 0 & 0 & 0 & 0 & 0 & 0 & 0 \\ 
\hline 
\end{tabular} 
	\begin{tabular}{|c||c|c|c|c|c|c|c|c|c|}\hline 
 &\multicolumn{9}{|c|}{\eqref{samedecision} holds for $h=1$, $i=3$} \\ 
\hline 
$\delta$ & $n$=100 & $n$=200 & $n$=500 & $n$=1000 & $n$=2000 & $n$=4000 & $n$=6000 & $n$=8000 & $n$=10000 \\ 
\hline \hline 
0  & 0.3639 & 0.3674 & 0.3998 & 0.5157 & 0.658 & 0.8198 & 0.8963 & 0.937 & 0.9574 \\ 
0.01  & 0.3161 & 0.3004 & 0.2974 & 0.3778 & 0.4804 & 0.6284 & 0.7111 & 0.7684 & 0.8035 \\ 
0.05  & 0.187 & 0.1481 & 0.1135 & 0.1215 & 0.1373 & 0.1629 & 0.1671 & 0.1789 & 0.1779 \\ 
0.1  & 0.1131 & 0.0673 & 0.0388 & 0.0337 & 0.0324 & 0.0266 & 0.0165 & 0.0157 & 0.0124 \\ 
0.15  & 0.07 & 0.0311 & 0.0144 & 0.0105 & 0.0065 & 0.0029 & 0.0011 & 0.0014 & 5e-04 \\ 
0.2  & 0.0466 & 0.0158 & 0.0056 & 0.0027 & 0.001 & 1e-04 & 2e-04 & 0 & 0 \\ 
0.4  & 0.0156 & 0.0032 & 0 & 0 & 0 & 0 & 0 & 0 & 0 \\ 
0.6  & 0.0098 & 0.0021 & 0 & 0 & 0 & 0 & 0 & 0 & 0 \\ 
\hline 
\end{tabular} 
	\begin{tabular}{|c||c|c|c|c|c|c|c|c|c|}\hline 
 &\multicolumn{9}{|c|}{\eqref{samedecision} holds for $h=5$, $i=2$} \\ 
\hline 
$\delta$ & $n$=100 & $n$=200 & $n$=500 & $n$=1000 & $n$=2000 & $n$=4000 & $n$=6000 & $n$=8000 & $n$=10000 \\ 
\hline \hline 
0  & 0.3169 & 0.3634 & 0.4147 & 0.5017 & 0.6288 & 0.791 & 0.8659 & 0.9062 & 0.9325 \\ 
0.01  & 0.3001 & 0.3353 & 0.365 & 0.4364 & 0.5504 & 0.7089 & 0.7968 & 0.8458 & 0.8805 \\ 
0.05  & 0.2349 & 0.2381 & 0.226 & 0.2622 & 0.3364 & 0.4442 & 0.5257 & 0.5643 & 0.6148 \\ 
0.1  & 0.1787 & 0.1603 & 0.1386 & 0.1471 & 0.1859 & 0.2333 & 0.2737 & 0.2856 & 0.3076 \\ 
0.15  & 0.133 & 0.1093 & 0.0895 & 0.0875 & 0.1019 & 0.1174 & 0.132 & 0.135 & 0.1301 \\ 
0.2  & 0.1049 & 0.0761 & 0.059 & 0.0543 & 0.0586 & 0.0591 & 0.0564 & 0.0574 & 0.0482 \\ 
0.4  & 0.043 & 0.0257 & 0.0174 & 0.0089 & 0.0053 & 0.0041 & 0.0016 & 5e-04 & 5e-04 \\ 
0.6  & 0.0221 & 0.0118 & 0.0087 & 0.0032 & 9e-04 & 0 & 0 & 0 & 0 \\ 
\hline 
\end{tabular} 
	\begin{tabular}{|c||c|c|c|c|c|c|c|c|c|}\hline 
 &\multicolumn{9}{|c|}{\eqref{samedecision} holds for $h=5$, $i=3$} \\ 
\hline 
$\delta$ & $n$=100 & $n$=200 & $n$=500 & $n$=1000 & $n$=2000 & $n$=4000 & $n$=6000 & $n$=8000 & $n$=10000 \\ 
\hline \hline 
0  & 0.2923 & 0.3538 & 0.438 & 0.5548 & 0.702 & 0.8419 & 0.9024 & 0.9401 & 0.9579 \\ 
0.01  & 0.2766 & 0.3312 & 0.394 & 0.4992 & 0.639 & 0.7876 & 0.8593 & 0.903 & 0.9291 \\ 
0.05  & 0.2221 & 0.2456 & 0.2624 & 0.3324 & 0.4476 & 0.5747 & 0.6518 & 0.7082 & 0.7482 \\ 
0.1  & 0.1729 & 0.1722 & 0.1696 & 0.2146 & 0.2906 & 0.375 & 0.4165 & 0.4651 & 0.4885 \\ 
0.15  & 0.1372 & 0.1213 & 0.1166 & 0.1397 & 0.184 & 0.2298 & 0.2556 & 0.2784 & 0.2832 \\ 
0.2  & 0.1097 & 0.0941 & 0.0819 & 0.095 & 0.1192 & 0.1429 & 0.1464 & 0.1592 & 0.1561 \\ 
0.4  & 0.0542 & 0.0408 & 0.0263 & 0.0202 & 0.0175 & 0.0157 & 0.0108 & 0.0107 & 0.0089 \\ 
0.6  & 0.0331 & 0.0258 & 0.0146 & 0.0055 & 0.0022 & 0.0012 & 6e-04 & 5e-04 & 2e-04 \\ 
\hline 
\end{tabular} 
\caption{\textit{Proportion of \eqref{decisionrule}  being fulfilled for the process \eqref{increasing3} and different values of $h$, $\delta$ and $n$.}}  \label{MSPEanalysisincreasing3}
\end{center}
\end{table}

\clearpage

\begin{figure}
\centering 
\includegraphics[width=0.24\textwidth]{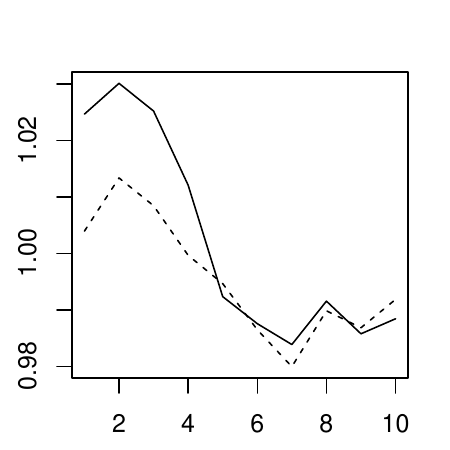}~
\includegraphics[width=0.24\textwidth]{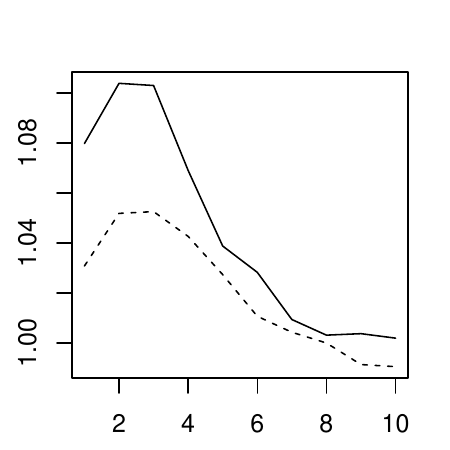}~
\includegraphics[width=0.24\textwidth]{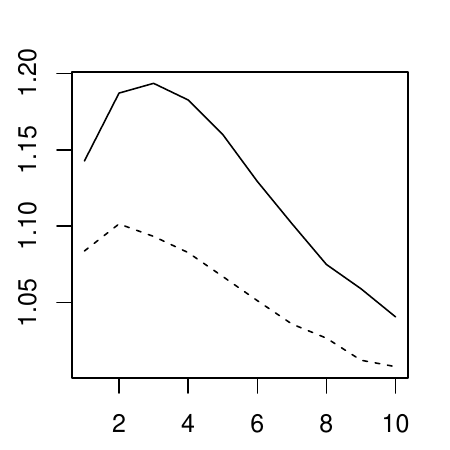}~ 
\includegraphics[width=0.24\textwidth]{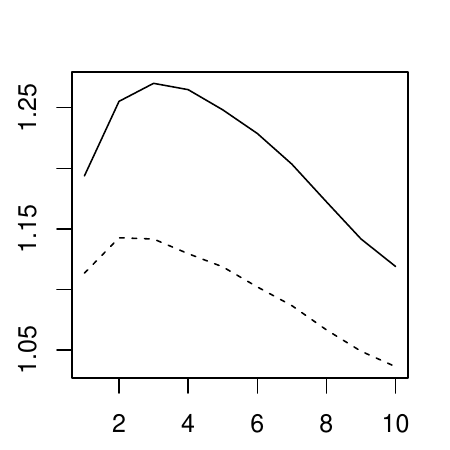} \\
\includegraphics[width=0.24\textwidth]{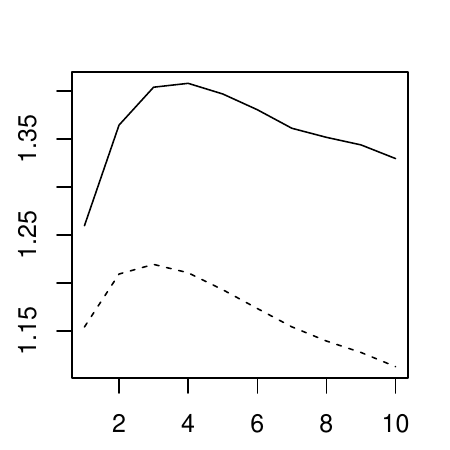}~
\includegraphics[width=0.24\textwidth]{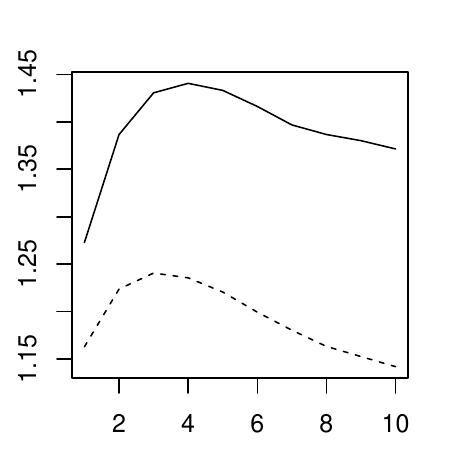}~ 
\includegraphics[width=0.24\textwidth]{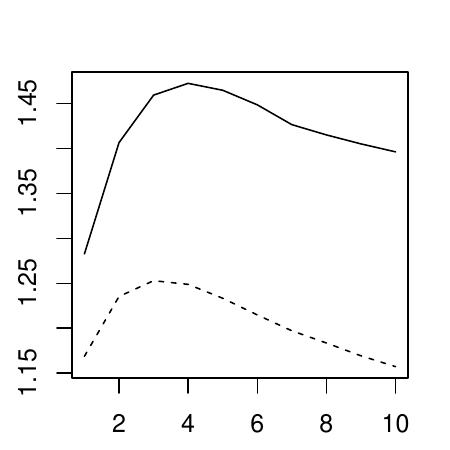}~ 
\includegraphics[width=0.24\textwidth]{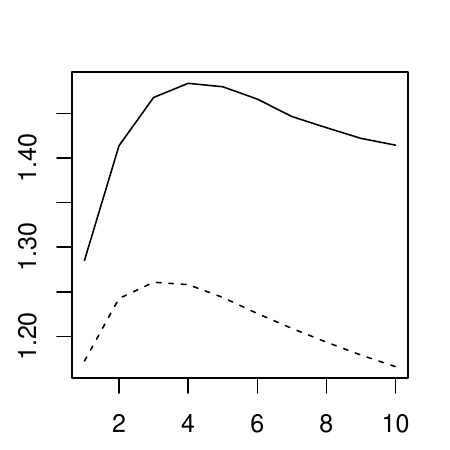}
   \caption{\it Plot of $h \mapsto R_{T,i}(h)$ for model \eqref{increasing5} and different values of $n$ (from left to right: n=100, n=200, n=500, n=1000 [first row], n=4000, n=6000, n=8000, n=10000 [second row]). Solid line: $i=3$ (test set), dashed line: $i=2$ (validation set 2). } \label{MSPEincreasing5}
\end{figure}

\begin{table} \scriptsize
\begin{center}
\begin{tabular}{|c||c|c|c||c|c|c|}
\hline
 \multirow{2}{*}{$n$} & & \multicolumn{2}{|c||}{$R^{{\rm s}, {\rm ls}}_{T,2}(1)$}  & & \multicolumn{2}{|c|}{$R^{{\rm s}, {\rm ls}}_{T,2}(5)$} \\ 
 & & $\geq 1.01$  & $< 1.01$  & & $\geq 1.01$  & $< 1.01$ \\ 
 \hline 
 \multirow{2}{*}{$100$} & $R^{{\rm s}, {\rm ls}}_{T,3}(1) \geq 1.01$  & 0.2666 & 0.2282  & $R^{{\rm s}, {\rm ls}}_{T,3}(5) \geq 1.01$  & 0.1158 & 0.1452 \\ 
& $R^{{\rm s}, {\rm ls}}_{T,3}(1) < 1.01$  & 0.1759 & 0.3293 & $R^{{\rm s}, {\rm ls}}_{T,3}(5) < 1.01$  & 0.1321 & 0.6069 \\

\hline\hline
 \multirow{2}{*}{$n$} & & \multicolumn{2}{|c||}{$R^{{\rm s}, {\rm ls}}_{T,2}(1)$}  & & \multicolumn{2}{|c|}{$R^{{\rm s}, {\rm ls}}_{T,2}(5)$} \\ 
 & & $\geq 1.2$  & $< 1.2$  & & $\geq 1.2$  & $< 1.2$ \\ 
 \hline 
 \multirow{2}{*}{$1000$} & $R^{{\rm s}, {\rm ls}}_{T,3}(1) \geq 1.2$  & 0.0882 & 0.3766  & $R^{{\rm s}, {\rm ls}}_{T,3}(5) \geq 1.2$  & 0.1447 & 0.3549 \\ 
& $R^{{\rm s}, {\rm ls}}_{T,3}(1) < 1.2$  & 0.1042 & 0.431 & $R^{{\rm s}, {\rm ls}}_{T,3}(5) < 1.2$  & 0.0834 & 0.417 \\

\hline\hline
 \multirow{2}{*}{$n$} & & \multicolumn{2}{|c||}{$R^{{\rm s}, {\rm ls}}_{T,2}(1)$}  & & \multicolumn{2}{|c|}{$R^{{\rm s}, {\rm ls}}_{T,2}(5)$} \\ 
 & & $\geq 1$  & $< 1$  & & $\geq 1$  & $< 1$ \\ 
 \hline 
 \multirow{2}{*}{$10000$} & $R^{{\rm s}, {\rm ls}}_{T,3}(1) \geq 1$  & 1 & 0  & $R^{{\rm s}, {\rm ls}}_{T,3}(5) \geq 1$  & 0.9999 & 1e-04 \\ 
& $R^{{\rm s}, {\rm ls}}_{T,3}(1) < 1$  & 0 & 0 & $R^{{\rm s}, {\rm ls}}_{T,3}(5) < 1$  & 0 & 0 \\

\hline\hline
 \multirow{2}{*}{$n$} & & \multicolumn{2}{|c||}{$R^{{\rm s}, {\rm ls}}_{T,2}(1)$}  & & \multicolumn{2}{|c|}{$R^{{\rm s}, {\rm ls}}_{T,2}(5)$} \\ 
 & & $\geq 1.2$  & $< 1.2$  & & $\geq 1.2$  & $< 1.2$ \\ 
 \hline 
 \multirow{2}{*}{$10000$} & $R^{{\rm s}, {\rm ls}}_{T,3}(1) \geq 1.2$  & 0.2372 & 0.58  & $R^{{\rm s}, {\rm ls}}_{T,3}(5) \geq 1.2$  & 0.6407 & 0.3255 \\ 
& $R^{{\rm s}, {\rm ls}}_{T,3}(1) < 1.2$  & 0.0699 & 0.1129 & $R^{{\rm s}, {\rm ls}}_{T,3}(5) < 1.2$  & 0.0253 & 0.0085 \\

\hline\hline
 \multirow{2}{*}{$n$} & & \multicolumn{2}{|c||}{$R^{{\rm s}, {\rm ls}}_{T,2}(1)$}  & & \multicolumn{2}{|c|}{$R^{{\rm s}, {\rm ls}}_{T,2}(5)$} \\ 
 & & $\geq 1.4$  & $< 1.4$  & & $\geq 1.4$  & $< 1.4$ \\ 
 \hline 
 \multirow{2}{*}{$10000$} & $R^{{\rm s}, {\rm ls}}_{T,3}(1) \geq 1.4$  & 1e-04 & 0.1515  & $R^{{\rm s}, {\rm ls}}_{T,3}(5) \geq 1.4$  & 0.0405 & 0.6207 \\ 
& $R^{{\rm s}, {\rm ls}}_{T,3}(1) < 1.4$  & 7e-04 & 0.8477 & $R^{{\rm s}, {\rm ls}}_{T,3}(5) < 1.4$  & 0.0276 & 0.3112 \\

\hline\hline
 \multirow{2}{*}{$n$} & & \multicolumn{2}{|c||}{$R^{{\rm s}, {\rm ls}}_{T,2}(1)$}  & & \multicolumn{2}{|c|}{$R^{{\rm s}, {\rm ls}}_{T,2}(5)$} \\ 
 & & $\geq 1.6$  & $< 1.6$  & & $\geq 1.6$  & $< 1.6$ \\ 
 \hline 
 \multirow{2}{*}{$10000$} & $R^{{\rm s}, {\rm ls}}_{T,3}(1) \geq 1.6$  & 0 & 0.0082  & $R^{{\rm s}, {\rm ls}}_{T,3}(5) \geq 1.6$  & 3e-04 & 0.2607 \\ 
& $R^{{\rm s}, {\rm ls}}_{T,3}(1) < 1.6$  & 0 & 0.9918 & $R^{{\rm s}, {\rm ls}}_{T,3}(5) < 1.6$  & 0.0011 & 0.7379 \\

\hline
\end{tabular}
\caption{\textit{Proportions of the individual events in~\eqref{samedecision} for the process \eqref{increasing5} and selected combinations of $n$ and $\delta$.}}  \label{MSPEanalysisincreasing5c}
\end{center}
\end{table}

\begin{table} \scriptsize
\begin{center}
	\begin{tabular}{|c||c|c|c|c|c|c|c|c|c|}\hline 
 &\multicolumn{9}{|c|}{\eqref{samedecision} holds for $h=1$} \\ 
\hline 
$\delta$ & $n$=100 & $n$=200 & $n$=500 & $n$=1000 & $n$=2000 & $n$=4000 & $n$=6000 & $n$=8000 & $n$=10000 \\ 
\hline \hline 
0  & 0.5959 & 0.6094 & 0.7471 & 0.8998 & 0.9789 & 0.9986 & 0.9999 & 1 & 1 \\ 
0.01  & 0.5959 & 0.5933 & 0.7039 & 0.8652 & 0.9643 & 0.9965 & 0.9995 & 0.9999 & 1 \\ 
0.05  & 0.6272 & 0.5771 & 0.5654 & 0.6863 & 0.834 & 0.9529 & 0.983 & 0.9935 & 0.9965 \\ 
0.1  & 0.6811 & 0.617 & 0.506 & 0.5113 & 0.6035 & 0.7443 & 0.8317 & 0.8793 & 0.9186 \\ 
0.15  & 0.7395 & 0.6773 & 0.5398 & 0.4615 & 0.438 & 0.4774 & 0.5348 & 0.5836 & 0.6215 \\ 
0.2  & 0.7989 & 0.7444 & 0.6148 & 0.5192 & 0.439 & 0.3831 & 0.3621 & 0.3461 & 0.3501 \\ 
0.4  & 0.9326 & 0.9303 & 0.8887 & 0.8513 & 0.8331 & 0.8301 & 0.8382 & 0.8348 & 0.8478 \\ 
0.6  & 0.9738 & 0.9828 & 0.9784 & 0.9702 & 0.9692 & 0.9778 & 0.9859 & 0.9881 & 0.9918 \\ 
\hline 
\end{tabular} 
	\begin{tabular}{|c||c|c|c|c|c|c|c|c|c|}\hline 
 &\multicolumn{9}{|c|}{\eqref{samedecision} holds for $h=5$} \\ 
\hline 
$\delta$ & $n$=100 & $n$=200 & $n$=500 & $n$=1000 & $n$=2000 & $n$=4000 & $n$=6000 & $n$=8000 & $n$=10000 \\ 
\hline \hline 
0  & 0.7186 & 0.6978 & 0.7147 & 0.807 & 0.9243 & 0.9841 & 0.9967 & 0.9979 & 0.9999 \\ 
0.01  & 0.7227 & 0.6912 & 0.6938 & 0.7741 & 0.9026 & 0.9768 & 0.9947 & 0.9975 & 0.9993 \\ 
0.05  & 0.7572 & 0.7175 & 0.6455 & 0.657 & 0.7931 & 0.9291 & 0.9709 & 0.9863 & 0.9937 \\ 
0.1  & 0.8011 & 0.7606 & 0.6375 & 0.5639 & 0.6302 & 0.7904 & 0.887 & 0.9268 & 0.9555 \\ 
0.15  & 0.8375 & 0.8066 & 0.6628 & 0.536 & 0.5055 & 0.6153 & 0.7199 & 0.7842 & 0.8376 \\ 
0.2  & 0.867 & 0.8402 & 0.7107 & 0.5617 & 0.453 & 0.4575 & 0.5281 & 0.588 & 0.6492 \\ 
0.4  & 0.9274 & 0.9188 & 0.864 & 0.7781 & 0.649 & 0.5198 & 0.4491 & 0.3891 & 0.3517 \\ 
0.6  & 0.9541 & 0.952 & 0.9302 & 0.9057 & 0.8606 & 0.8159 & 0.7848 & 0.7557 & 0.7382 \\ 
\hline 
\end{tabular} 
\caption{\textit{Proportion of \eqref{samedecision}  being fulfilled for the process \eqref{increasing5} and different values of $h$, $\delta$ and $n$.}}  \label{MSPEanalysisincreasing5b}
\end{center}
\end{table}

\begin{table} \scriptsize
\begin{center}
	\begin{tabular}{|c||c|c|c|c|c|c|c|c|c|}\hline 
 &\multicolumn{9}{|p{12.5cm}|}{\centering Value of $q(\delta)$, defined in \eqref{cond:f}, for $h=1$} \\ 
\hline 
$\delta$ & $n$=100 & $n$=200 & $n$=500 & $n$=1000 & $n$=2000 & $n$=4000 & $n$=6000 & $n$=8000 & $n$=10000 \\ 
\hline \hline 
0  & 0 & 0 & 0 & 0 & 0 & 0 & 0 & 0 & 0 \\ 
0.01  & 0.1 & 0.13 & 0.17 & 0.18 & 0.2 & 0.21 & 0.21 & 0.22 & 0.22 \\ 
0.05  & 0.058 & 0.089 & 0.13 & 0.15 & 0.17 & 0.19 & 0.2 & 0.21 & 0.22 \\ 
0.1  & 0.0046 & 0.036 & 0.075 & 0.1 & 0.12 & 0.14 & 0.15 & 0.16 & 0.17 \\ 
0.15  & 2.5e-05 & 7.8e-05 & 0.022 & 0.048 & 0.071 & 0.091 & 0.1 & 0.11 & 0.12 \\ 
0.2  & 8.1e-05 & 7.5e-05 & 2.1e-06 & 6.4e-06 & 0.019 & 0.04 & 0.051 & 0.058 & 0.064 \\ 
0.4  & 0.2 & 0.16 & 0.12 & 0.089 & 0.062 & 0.038 & 0.023 & 0.014 & 0.0063 \\ 
0.6  & 0.4 & 0.37 & 0.32 & 0.29 & 0.26 & 0.24 & 0.23 & 0.22 & 0.21 \\ 
\hline 
\end{tabular} 
	\begin{tabular}{|c||c|c|c|c|c|c|c|c|c|}\hline 
 &\multicolumn{9}{|p{12.5cm}|}{\centering Value of $q(\delta)$, defined in \eqref{cond:f}, for $h=5$} \\ 
\hline 
$\delta$ & $n$=100 & $n$=200 & $n$=500 & $n$=1000 & $n$=2000 & $n$=4000 & $n$=6000 & $n$=8000 & $n$=10000 \\ 
\hline \hline 
0  & 0 & 0 & 0 & 0 & 0 & 0 & 0 & 0 & 0 \\ 
0.01  & 0.079 & 0.11 & 0.14 & 0.17 & 0.19 & 0.2 & 0.21 & 0.22 & 0.22 \\ 
0.05  & 0.036 & 0.3 & 0.71 & 0.83 & 0.93 & 1 & 1.1 & 1.1 & 1.1 \\ 
0.1  & 0.00017 & 0.064 & 0.49 & 0.81 & 1.1 & 1.4 & 1.5 & 1.6 & 1.7 \\ 
0.15  & 0.00014 & 0.00092 & 0.25 & 0.56 & 0.85 & 1.1 & 1.3 & 1.4 & 1.4 \\ 
0.2  & 2e-04 & 0.0011 & 0.00046 & 0.31 & 0.61 & 0.88 & 1 & 1.1 & 1.2 \\ 
0.4  & 0.55 & 0.0037 & 0.0026 & 0.00099 & 8e-04 & 0.00016 & 0.047 & 0.14 & 0.22 \\ 
0.6  & 0.86 & 0.85 & 0.11 & 0.001 & 0.0018 & 0.0013 & 0.0026 & 0.00063 & 0.00066 \\ 
\hline 
\end{tabular} 
\caption{\textit{Values of $q(\delta)$, defined in \eqref{cond:f}, for the process \eqref{increasing5} and different values of $h$, $\delta$ and $n$.}} \label{MSPEanalysisincreasing5d}
\end{center}
\end{table}

\begin{table} \scriptsize
\begin{center}
	\begin{tabular}{|c||c|c|c|c|c|c|c|c|c|}\hline 
 &\multicolumn{9}{|c|}{\eqref{samedecision} holds for $h=1$, $i=2$} \\ 
\hline 
$\delta$ & $n$=100 & $n$=200 & $n$=500 & $n$=1000 & $n$=2000 & $n$=4000 & $n$=6000 & $n$=8000 & $n$=10000 \\ 
\hline \hline 
0  & 0.4784 & 0.5793 & 0.7938 & 0.9246 & 0.9847 & 0.9988 & 0.9999 & 1 & 1 \\ 
0.01  & 0.4425 & 0.5384 & 0.7544 & 0.8958 & 0.9734 & 0.997 & 0.9997 & 0.9999 & 1 \\ 
0.05  & 0.3283 & 0.3956 & 0.5743 & 0.7413 & 0.8703 & 0.9624 & 0.9854 & 0.9943 & 0.9966 \\ 
0.1  & 0.2391 & 0.2705 & 0.3853 & 0.5182 & 0.6534 & 0.7845 & 0.8529 & 0.892 & 0.9264 \\ 
0.15  & 0.176 & 0.186 & 0.2487 & 0.327 & 0.4087 & 0.5012 & 0.5676 & 0.6152 & 0.6512 \\ 
0.2  & 0.1276 & 0.1238 & 0.1524 & 0.1924 & 0.2263 & 0.2582 & 0.2831 & 0.292 & 0.3071 \\ 
0.4  & 0.039 & 0.0226 & 0.0175 & 0.0137 & 0.0112 & 0.006 & 0.0031 & 0.0019 & 8e-04 \\ 
0.6  & 0.017 & 0.006 & 0.0012 & 8e-04 & 4e-04 & 1e-04 & 0 & 0 & 0 \\ 
\hline 
\end{tabular} 
	\begin{tabular}{|c||c|c|c|c|c|c|c|c|c|}\hline 
 &\multicolumn{9}{|c|}{\eqref{samedecision} holds for $h=1$, $i=3$} \\ 
\hline 
$\delta$ & $n$=100 & $n$=200 & $n$=500 & $n$=1000 & $n$=2000 & $n$=4000 & $n$=6000 & $n$=8000 & $n$=10000 \\ 
\hline \hline 
0  & 0.5185 & 0.6707 & 0.8793 & 0.9682 & 0.994 & 0.9998 & 1 & 1 & 1 \\ 
0.01  & 0.4948 & 0.6379 & 0.8497 & 0.9556 & 0.9903 & 0.9995 & 0.9998 & 1 & 1 \\ 
0.05  & 0.3973 & 0.5165 & 0.7325 & 0.8742 & 0.9543 & 0.9903 & 0.9976 & 0.9992 & 0.9999 \\ 
0.1  & 0.3084 & 0.3911 & 0.5877 & 0.7427 & 0.8597 & 0.9402 & 0.9728 & 0.9859 & 0.9922 \\ 
0.15  & 0.2361 & 0.2915 & 0.4585 & 0.6027 & 0.7257 & 0.8338 & 0.8906 & 0.9294 & 0.9455 \\ 
0.2  & 0.1765 & 0.2206 & 0.3518 & 0.4648 & 0.5825 & 0.6819 & 0.7442 & 0.7949 & 0.8172 \\ 
0.4  & 0.0606 & 0.0583 & 0.0982 & 0.1388 & 0.1579 & 0.1659 & 0.1595 & 0.1637 & 0.1516 \\ 
0.6  & 0.0264 & 0.016 & 0.0206 & 0.029 & 0.0304 & 0.0221 & 0.0141 & 0.0119 & 0.0082 \\ 
\hline 
\end{tabular} 
	\begin{tabular}{|c||c|c|c|c|c|c|c|c|c|}\hline 
 &\multicolumn{9}{|c|}{\eqref{samedecision} holds for $h=5$, $i=2$} \\ 
\hline 
$\delta$ & $n$=100 & $n$=200 & $n$=500 & $n$=1000 & $n$=2000 & $n$=4000 & $n$=6000 & $n$=8000 & $n$=10000 \\ 
\hline \hline 
0  & 0.2714 & 0.3441 & 0.5552 & 0.7688 & 0.9229 & 0.9874 & 0.9973 & 0.9984 & 0.9999 \\ 
0.01  & 0.2479 & 0.3189 & 0.5217 & 0.7372 & 0.9055 & 0.9813 & 0.9958 & 0.998 & 0.9994 \\ 
0.05  & 0.1752 & 0.2324 & 0.4021 & 0.6074 & 0.8071 & 0.9416 & 0.9743 & 0.9878 & 0.9941 \\ 
0.1  & 0.1209 & 0.1607 & 0.2883 & 0.4485 & 0.6434 & 0.8181 & 0.8991 & 0.9325 & 0.9584 \\ 
0.15  & 0.0867 & 0.1155 & 0.2068 & 0.3256 & 0.4743 & 0.6434 & 0.7467 & 0.801 & 0.847 \\ 
0.2  & 0.0663 & 0.0846 & 0.1476 & 0.2281 & 0.3307 & 0.4529 & 0.5521 & 0.6099 & 0.666 \\ 
0.4  & 0.0305 & 0.0348 & 0.054 & 0.056 & 0.0548 & 0.0608 & 0.0626 & 0.0666 & 0.0681 \\ 
0.6  & 0.0182 & 0.0221 & 0.028 & 0.0214 & 0.0113 & 0.0068 & 0.0026 & 0.0029 & 0.0014 \\ 
\hline 
\end{tabular} 
	\input{tab1-m7-h5-i3.tex}
\caption{\textit{Proportion of \eqref{decisionrule}  being fulfilled for the process \eqref{increasing5} and different values of $h$, $\delta$ and $n$.}}  \label{MSPEanalysisincreasing5}
\end{center}
\end{table}

\clearpage

\begin{figure}
\centering 
\includegraphics[width=0.24\textwidth]{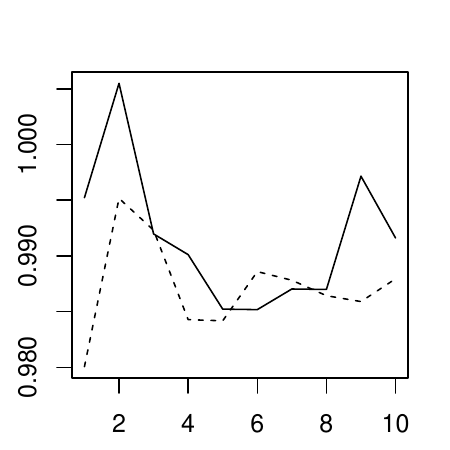}~
\includegraphics[width=0.24\textwidth]{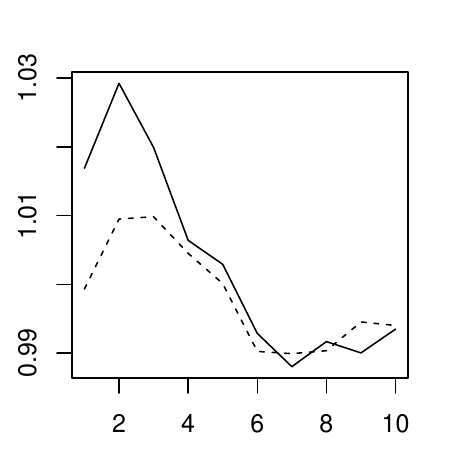}~
\includegraphics[width=0.24\textwidth]{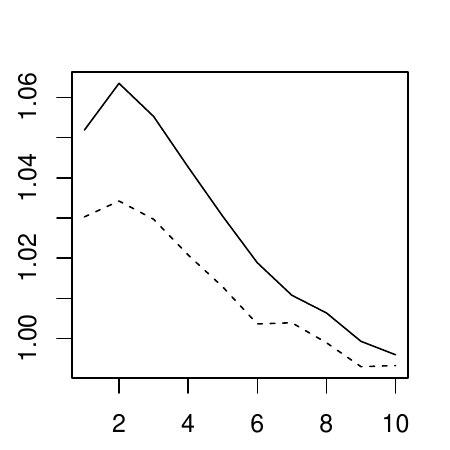}~ 
\includegraphics[width=0.24\textwidth]{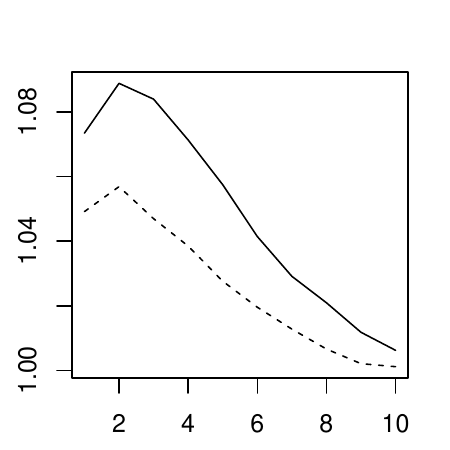} \\
\includegraphics[width=0.24\textwidth]{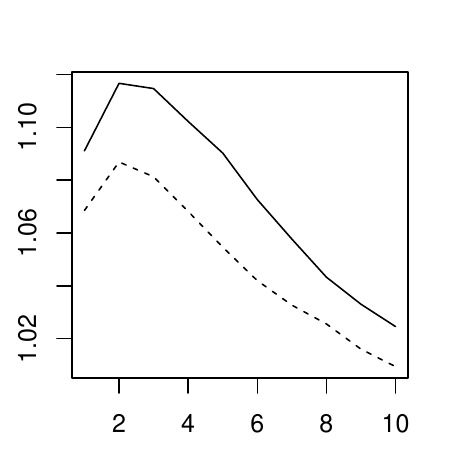}~
\includegraphics[width=0.24\textwidth]{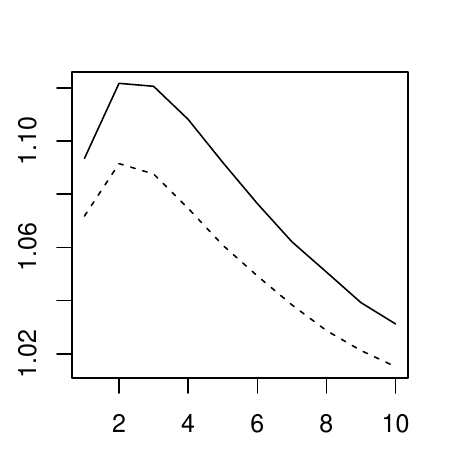}~ 
\includegraphics[width=0.24\textwidth]{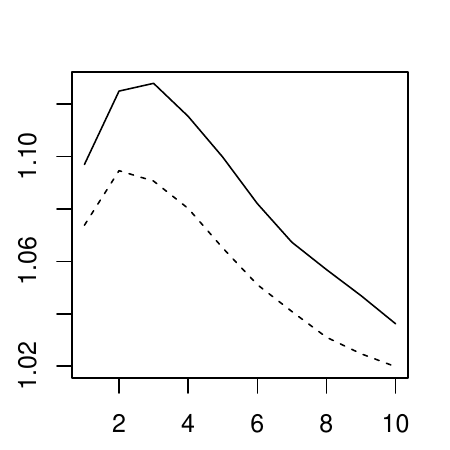}~ 
\includegraphics[width=0.24\textwidth]{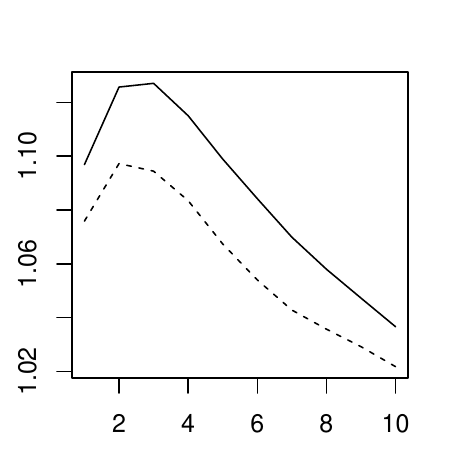}
   \caption{\it Plot of $h \mapsto R_{T,i}(h)$ for model \eqref{increasing6} and different values of $n$ (from left to right: n=100, n=200, n=500, n=1000 [first row], n=4000, n=6000, n=8000, n=10000 [second row]). Solid line: $i=3$ (test set), dashed line: $i=2$ (validation set 2). } \label{MSPEincreasing6}
\end{figure}

\begin{table} \scriptsize
\begin{center}
\begin{tabular}{|c||c|c|c||c|c|c|}
\hline
 \multirow{2}{*}{$n$} & & \multicolumn{2}{|c||}{$R^{{\rm s}, {\rm ls}}_{T,2}(1)$}  & & \multicolumn{2}{|c|}{$R^{{\rm s}, {\rm ls}}_{T,2}(5)$} \\ 
 & & $\geq 1.01$  & $< 1.01$  & & $\geq 1.01$  & $< 1.01$ \\ 
 \hline 
 \multirow{2}{*}{$100$} & $R^{{\rm s}, {\rm ls}}_{T,3}(1) \geq 1.01$  & 0.191 & 0.2226  & $R^{{\rm s}, {\rm ls}}_{T,3}(5) \geq 1.01$  & 0.0944 & 0.1508 \\ 
& $R^{{\rm s}, {\rm ls}}_{T,3}(1) < 1.01$  & 0.1888 & 0.3976 & $R^{{\rm s}, {\rm ls}}_{T,3}(5) < 1.01$  & 0.1354 & 0.6194 \\

\hline\hline
 \multirow{2}{*}{$n$} & & \multicolumn{2}{|c||}{$R^{{\rm s}, {\rm ls}}_{T,2}(1)$}  & & \multicolumn{2}{|c|}{$R^{{\rm s}, {\rm ls}}_{T,2}(5)$} \\ 
 & & $\geq 1.05$  & $< 1.05$  & & $\geq 1.05$  & $< 1.05$ \\ 
 \hline 
 \multirow{2}{*}{$1000$} & $R^{{\rm s}, {\rm ls}}_{T,3}(1) \geq 1.05$  & 0.2796 & 0.3124  & $R^{{\rm s}, {\rm ls}}_{T,3}(5) \geq 1.05$  & 0.1577 & 0.247 \\ 
& $R^{{\rm s}, {\rm ls}}_{T,3}(1) < 1.05$  & 0.1669 & 0.2411 & $R^{{\rm s}, {\rm ls}}_{T,3}(5) < 1.05$  & 0.1253 & 0.47 \\

\hline\hline
 \multirow{2}{*}{$n$} & & \multicolumn{2}{|c||}{$R^{{\rm s}, {\rm ls}}_{T,2}(1)$}  & & \multicolumn{2}{|c|}{$R^{{\rm s}, {\rm ls}}_{T,2}(5)$} \\ 
 & & $\geq 1$  & $< 1$  & & $\geq 1$  & $< 1$ \\ 
 \hline 
 \multirow{2}{*}{$10000$} & $R^{{\rm s}, {\rm ls}}_{T,3}(1) \geq 1$  & 0.9997 & 3e-04  & $R^{{\rm s}, {\rm ls}}_{T,3}(5) \geq 1$  & 0.975 & 0.018 \\ 
& $R^{{\rm s}, {\rm ls}}_{T,3}(1) < 1$  & 0 & 0 & $R^{{\rm s}, {\rm ls}}_{T,3}(5) < 1$  & 0.0067 & 3e-04 \\

\hline\hline
 \multirow{2}{*}{$n$} & & \multicolumn{2}{|c||}{$R^{{\rm s}, {\rm ls}}_{T,2}(1)$}  & & \multicolumn{2}{|c|}{$R^{{\rm s}, {\rm ls}}_{T,2}(5)$} \\ 
 & & $\geq 1.05$  & $< 1.05$  & & $\geq 1.05$  & $< 1.05$ \\ 
 \hline 
 \multirow{2}{*}{$10000$} & $R^{{\rm s}, {\rm ls}}_{T,3}(1) \geq 1.05$  & 0.7698 & 0.1674  & $R^{{\rm s}, {\rm ls}}_{T,3}(5) \geq 1.05$  & 0.5746 & 0.2863 \\ 
& $R^{{\rm s}, {\rm ls}}_{T,3}(1) < 1.05$  & 0.0536 & 0.0092 & $R^{{\rm s}, {\rm ls}}_{T,3}(5) < 1.05$  & 0.0939 & 0.0452 \\

\hline\hline
 \multirow{2}{*}{$n$} & & \multicolumn{2}{|c||}{$R^{{\rm s}, {\rm ls}}_{T,2}(1)$}  & & \multicolumn{2}{|c|}{$R^{{\rm s}, {\rm ls}}_{T,2}(5)$} \\ 
 & & $\geq 1.1$  & $< 1.1$  & & $\geq 1.1$  & $< 1.1$ \\ 
 \hline 
 \multirow{2}{*}{$10000$} & $R^{{\rm s}, {\rm ls}}_{T,3}(1) \geq 1.1$  & 0.0831 & 0.3572  & $R^{{\rm s}, {\rm ls}}_{T,3}(5) \geq 1.1$  & 0.0974 & 0.3738 \\ 
& $R^{{\rm s}, {\rm ls}}_{T,3}(1) < 1.1$  & 0.1107 & 0.449 & $R^{{\rm s}, {\rm ls}}_{T,3}(5) < 1.1$  & 0.1005 & 0.4283 \\

\hline\hline
 \multirow{2}{*}{$n$} & & \multicolumn{2}{|c||}{$R^{{\rm s}, {\rm ls}}_{T,2}(1)$}  & & \multicolumn{2}{|c|}{$R^{{\rm s}, {\rm ls}}_{T,2}(5)$} \\ 
 & & $\geq 1.15$  & $< 1.15$  & & $\geq 1.15$  & $< 1.15$ \\ 
 \hline 
 \multirow{2}{*}{$10000$} & $R^{{\rm s}, {\rm ls}}_{T,3}(1) \geq 1.15$  & 9e-04 & 0.0731  & $R^{{\rm s}, {\rm ls}}_{T,3}(5) \geq 1.15$  & 0.0055 & 0.1572 \\ 
& $R^{{\rm s}, {\rm ls}}_{T,3}(1) < 1.15$  & 0.0121 & 0.9139 & $R^{{\rm s}, {\rm ls}}_{T,3}(5) < 1.15$  & 0.0244 & 0.8129 \\

\hline
\end{tabular}
\caption{\textit{Proportions of the individual events in~\eqref{samedecision} for the process \eqref{increasing6} and selected combinations of $n$ and $\delta$.}}  \label{MSPEanalysisincreasing6c}
\end{center}
\end{table}

\begin{table} \scriptsize
\begin{center}
	\input{tab2-m15-h1.tex}
	\input{tab2-m15-h5.tex}
\caption{\textit{Proportion of \eqref{samedecision}  being fulfilled for the process \eqref{increasing6} and different values of $h$, $\delta$ and $n$.}}  \label{MSPEanalysisincreasing6b}
\end{center}
\end{table}

\begin{table} \scriptsize
\begin{center}
	\input{tab5-m15-h1.tex}
	\input{tab5-m15-h5.tex}
\caption{\textit{Values of $q(\delta)$, defined in \eqref{cond:f}, for the process \eqref{increasing6} and different values of $h$, $\delta$ and $n$.}} \label{MSPEanalysisincreasing6d}
\end{center}
\end{table}

\begin{table} \scriptsize
\begin{center}
	\input{tab1-m15-h1-i2.tex}
	\input{tab1-m15-h1-i3.tex}
	\input{tab1-m15-h5-i2.tex}
	\input{tab1-m15-h5-i3.tex}
\caption{\textit{Proportion of \eqref{decisionrule}  being fulfilled for the process \eqref{increasing6} and different values of $h$, $\delta$ and $n$.}}  \label{MSPEanalysisincreasing6}
\end{center}
\end{table}

\clearpage

\begin{figure}
\centering 
\includegraphics[width=0.24\textwidth]{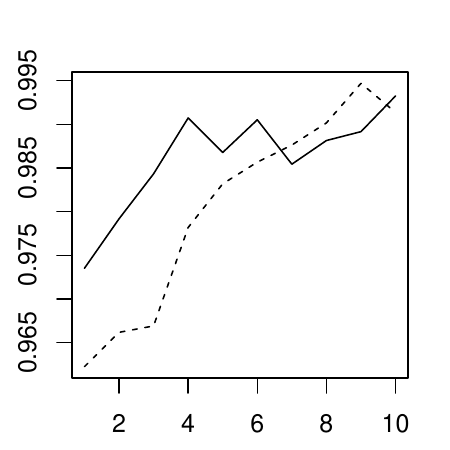}~
\includegraphics[width=0.24\textwidth]{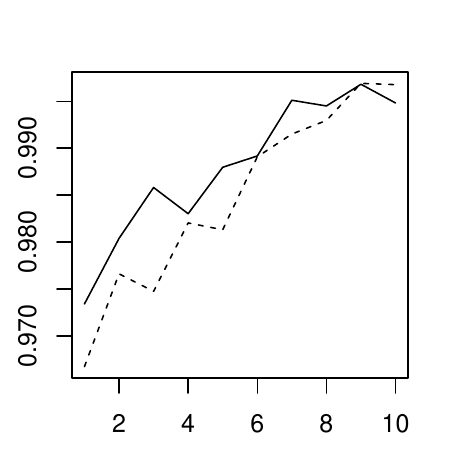}~
\includegraphics[width=0.24\textwidth]{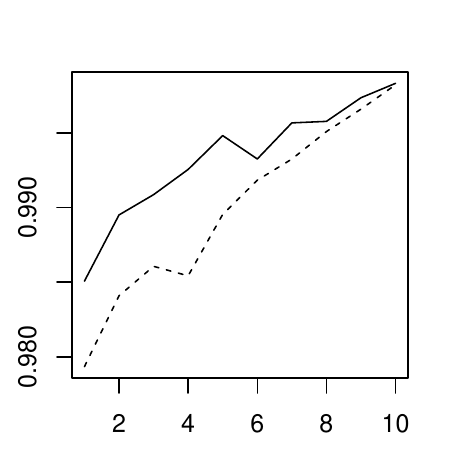}~ 
\includegraphics[width=0.24\textwidth]{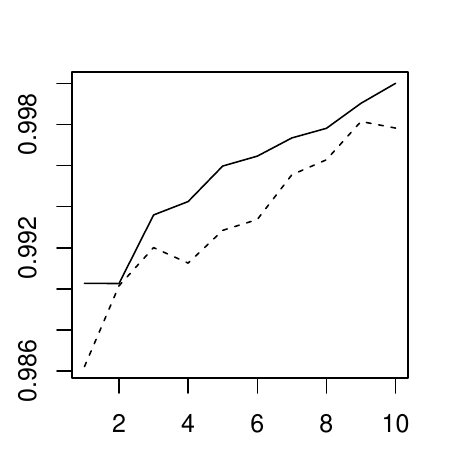} \\
\includegraphics[width=0.24\textwidth]{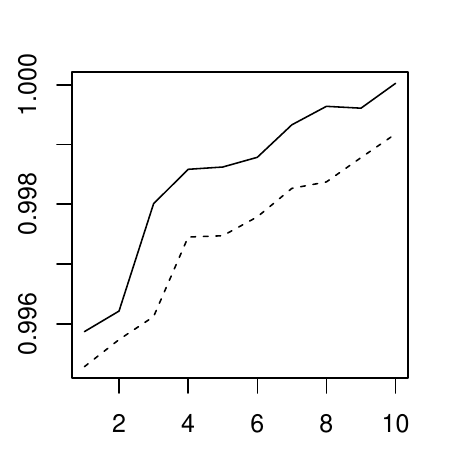}~
\includegraphics[width=0.24\textwidth]{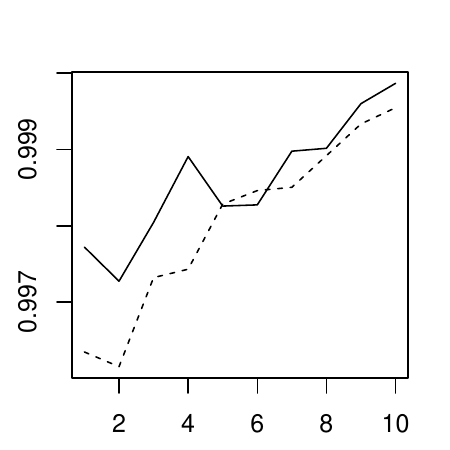}~ 
\includegraphics[width=0.24\textwidth]{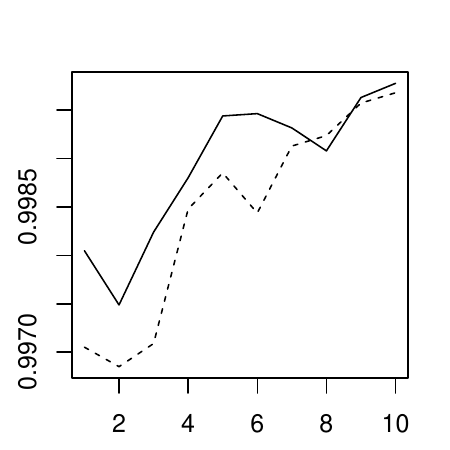}~ 
\includegraphics[width=0.24\textwidth]{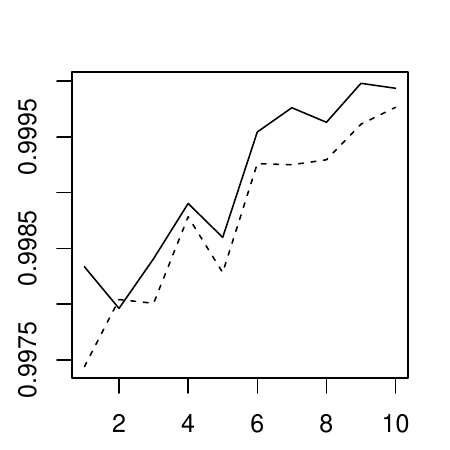}
   \caption{\it Plot of $h \mapsto R_{T,i}(h)$ for model \eqref{stationaryAR} and different values of $n$ (from left to right: n=100, n=200, n=500, n=1000 [first row], n=4000, n=6000, n=8000, n=10000 [second row]). Solid line: $i=3$ (test set), dashed line: $i=2$ (validation set 2).} \label{MSPEstationaryAR}
\end{figure}

\begin{table} \scriptsize
\begin{center}
\begin{tabular}{|c||c|c|c||c|c|c|}
\hline
 \multirow{2}{*}{$n$} & & \multicolumn{2}{|c||}{$R^{{\rm s}, {\rm ls}}_{T,2}(1)$}  & & \multicolumn{2}{|c|}{$R^{{\rm s}, {\rm ls}}_{T,2}(5)$} \\ 
 & & $\geq 1$  & $< 1$  & & $\geq 1$  & $< 1$ \\ 
 \hline 
 \multirow{2}{*}{$100$} & $R^{{\rm s}, {\rm ls}}_{T,3}(1) \geq 1$  & 0.2256 & 0.2131  & $R^{{\rm s}, {\rm ls}}_{T,3}(5) \geq 1$  & 0.4882 & 0.1574 \\ 
& $R^{{\rm s}, {\rm ls}}_{T,3}(1) < 1$  & 0.2002 & 0.3611 & $R^{{\rm s}, {\rm ls}}_{T,3}(5) < 1$  & 0.1549 & 0.1995 \\

\hline\hline
 \multirow{2}{*}{$n$} & & \multicolumn{2}{|c||}{$R^{{\rm s}, {\rm ls}}_{T,2}(1)$}  & & \multicolumn{2}{|c|}{$R^{{\rm s}, {\rm ls}}_{T,2}(5)$} \\ 
 & & $\geq 1$  & $< 1$  & & $\geq 1$  & $< 1$ \\ 
 \hline 
 \multirow{2}{*}{$1000$} & $R^{{\rm s}, {\rm ls}}_{T,3}(1) \geq 1$  & 0.0872 & 0.2186  & $R^{{\rm s}, {\rm ls}}_{T,3}(5) \geq 1$  & 0.4075 & 0.1793 \\ 
& $R^{{\rm s}, {\rm ls}}_{T,3}(1) < 1$  & 0.1895 & 0.5047 & $R^{{\rm s}, {\rm ls}}_{T,3}(5) < 1$  & 0.1516 & 0.2616 \\

\hline\hline
 \multirow{2}{*}{$n$} & & \multicolumn{2}{|c||}{$R^{{\rm s}, {\rm ls}}_{T,2}(1)$}  & & \multicolumn{2}{|c|}{$R^{{\rm s}, {\rm ls}}_{T,2}(5)$} \\ 
 & & $\geq 1$  & $< 1$  & & $\geq 1$  & $< 1$ \\ 
 \hline 
 \multirow{2}{*}{$10000$} & $R^{{\rm s}, {\rm ls}}_{T,3}(1) \geq 1$  & 0.0715 & 0.2003  & $R^{{\rm s}, {\rm ls}}_{T,3}(5) \geq 1$  & 0.313 & 0.185 \\ 
& $R^{{\rm s}, {\rm ls}}_{T,3}(1) < 1$  & 0.1819 & 0.5463 & $R^{{\rm s}, {\rm ls}}_{T,3}(5) < 1$  & 0.1731 & 0.3289 \\

\hline\hline
 \multirow{2}{*}{$n$} & & \multicolumn{2}{|c||}{$R^{{\rm s}, {\rm ls}}_{T,2}(1)$}  & & \multicolumn{2}{|c|}{$R^{{\rm s}, {\rm ls}}_{T,2}(5)$} \\ 
 & & $\geq 1.01$  & $< 1.01$  & & $\geq 1.01$  & $< 1.01$ \\ 
 \hline 
 \multirow{2}{*}{$10000$} & $R^{{\rm s}, {\rm ls}}_{T,3}(1) \geq 1.01$  & 0 & 0.0021  & $R^{{\rm s}, {\rm ls}}_{T,3}(5) \geq 1.01$  & 0.0029 & 0.021 \\ 
& $R^{{\rm s}, {\rm ls}}_{T,3}(1) < 1.01$  & 0.0012 & 0.9967 & $R^{{\rm s}, {\rm ls}}_{T,3}(5) < 1.01$  & 0.0144 & 0.9617 \\

\hline\hline
 \multirow{2}{*}{$n$} & & \multicolumn{2}{|c||}{$R^{{\rm s}, {\rm ls}}_{T,2}(1)$}  & & \multicolumn{2}{|c|}{$R^{{\rm s}, {\rm ls}}_{T,2}(5)$} \\ 
 & & $\geq 1.05$  & $< 1.05$  & & $\geq 1.05$  & $< 1.05$ \\ 
 \hline 
 \multirow{2}{*}{$10000$} & $R^{{\rm s}, {\rm ls}}_{T,3}(1) \geq 1.05$  & 0 & 0  & $R^{{\rm s}, {\rm ls}}_{T,3}(5) \geq 1.05$  & 0 & 2e-04 \\ 
& $R^{{\rm s}, {\rm ls}}_{T,3}(1) < 1.05$  & 0 & 1 & $R^{{\rm s}, {\rm ls}}_{T,3}(5) < 1.05$  & 1e-04 & 0.9997 \\

\hline\hline
 \multirow{2}{*}{$n$} & & \multicolumn{2}{|c||}{$R^{{\rm s}, {\rm ls}}_{T,2}(1)$}  & & \multicolumn{2}{|c|}{$R^{{\rm s}, {\rm ls}}_{T,2}(5)$} \\ 
 & & $\geq 1.1$  & $< 1.1$  & & $\geq 1.1$  & $< 1.1$ \\ 
 \hline 
 \multirow{2}{*}{$10000$} & $R^{{\rm s}, {\rm ls}}_{T,3}(1) \geq 1.1$  & 0 & 0  & $R^{{\rm s}, {\rm ls}}_{T,3}(5) \geq 1.1$  & 0 & 0 \\ 
& $R^{{\rm s}, {\rm ls}}_{T,3}(1) < 1.1$  & 0 & 1 & $R^{{\rm s}, {\rm ls}}_{T,3}(5) < 1.1$  & 0 & 1 \\

\hline
\end{tabular}
\caption{\textit{Proportions of the individual events in~\eqref{samedecision} for the process \eqref{stationaryAR} and selected combinations of $n$ and $\delta$.}}  \label{MSPEanalysisstationaryARc}
\end{center}
\end{table}

\begin{table} \scriptsize
\begin{center}
	\input{tab2-m8-h1.tex}
	\input{tab2-m8-h5.tex}
\caption{\textit{Proportion of \eqref{samedecision}  being fulfilled for the process \eqref{stationaryAR} and different values of $h$, $\delta$ and $n$.}}  \label{MSPEanalysisstationaryARb}
\end{center}
\end{table}

\begin{table} \scriptsize
\begin{center}
	\input{tab5-m8-h1.tex}
	\input{tab5-m8-h5.tex}
\caption{\textit{Values of $q(\delta)$, defined in \eqref{cond:f}, for the process \eqref{stationaryAR} and different values of $h$, $\delta$ and $n$.}} \label{MSPEanalysisstationaryARd}
\end{center}
\end{table}

\begin{table} \scriptsize
\begin{center}
	\input{tab1-m8-h1-i2.tex}
	\input{tab1-m8-h1-i3.tex}
	\input{tab1-m8-h5-i2.tex}
	\input{tab1-m8-h5-i3.tex}
\caption{\textit{Proportion of \eqref{decisionrule}  being fulfilled for the process \eqref{stationaryAR} and different values of $h$, $\delta$ and $n$.}}  \label{MSPEanalysisstationaryAR}
\end{center}
\end{table}

\clearpage

\begin{figure}
\centering 
\includegraphics[width=0.24\textwidth]{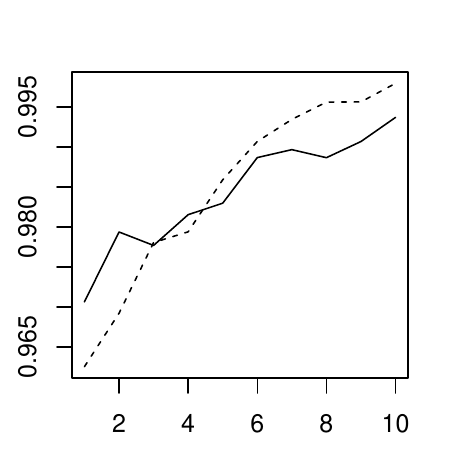}~
\includegraphics[width=0.24\textwidth]{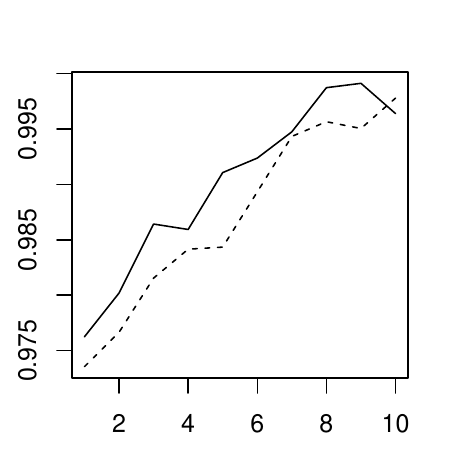}~
\includegraphics[width=0.24\textwidth]{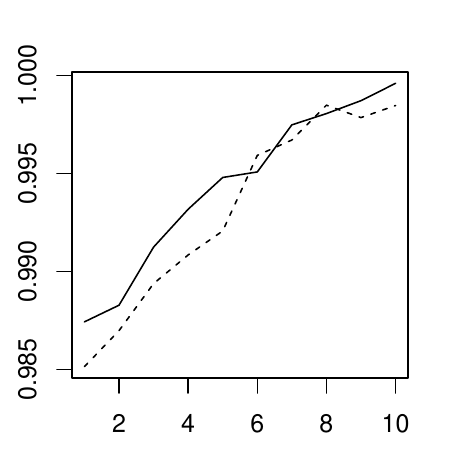}~ 
\includegraphics[width=0.24\textwidth]{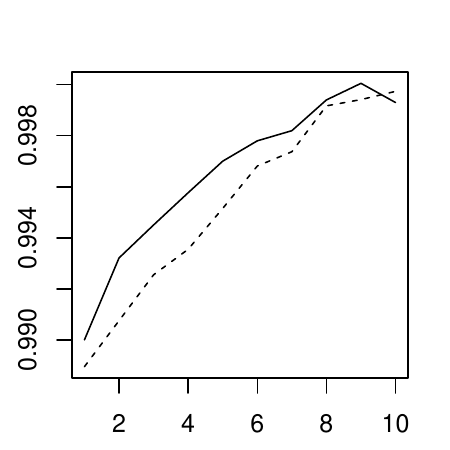} \\
\includegraphics[width=0.24\textwidth]{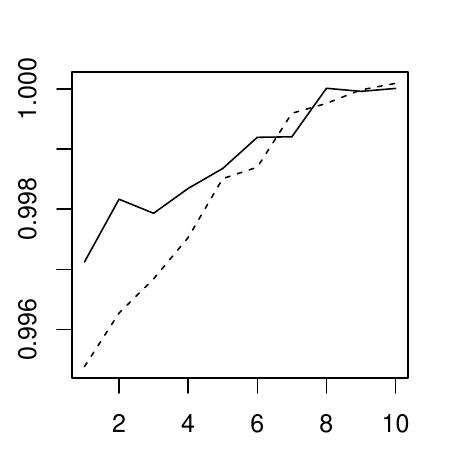}~
\includegraphics[width=0.24\textwidth]{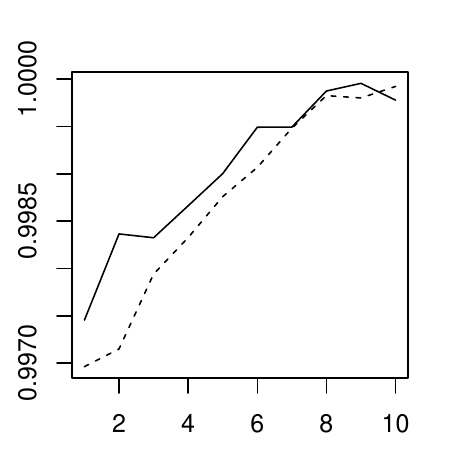}~ 
\includegraphics[width=0.24\textwidth]{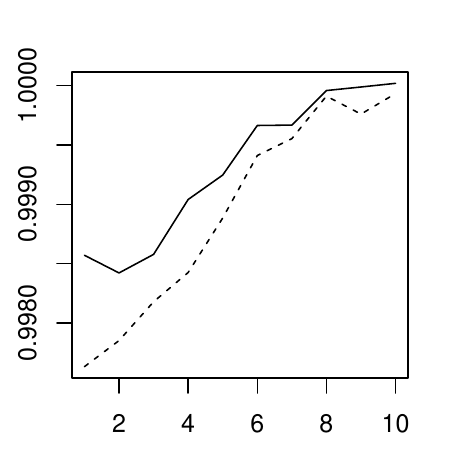}~ 
\includegraphics[width=0.24\textwidth]{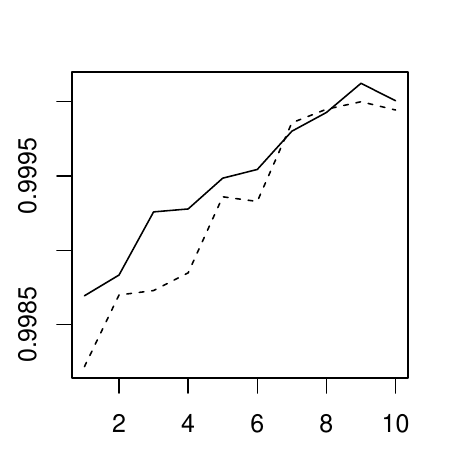}
   \caption{\it Plot of $h \mapsto R_{T,i}(h)$ for model \eqref{indepNonHetero} and different values of $n$ (from left to right: n=100, n=200, n=500, n=1000 [first row], n=4000, n=6000, n=8000, n=10000 [second row]). Solid line: $i=3$ (test set), dashed line: $i=2$ (validation set 2). }  \label{MSPEindepNonHetero}
\end{figure}

\begin{table} \scriptsize
\begin{center}
\begin{tabular}{|c||c|c|c||c|c|c|}
\hline
 \multirow{2}{*}{$n$} & & \multicolumn{2}{|c||}{$R^{{\rm s}, {\rm ls}}_{T,2}(1)$}  & & \multicolumn{2}{|c|}{$R^{{\rm s}, {\rm ls}}_{T,2}(5)$} \\ 
 & & $\geq 1$  & $< 1$  & & $\geq 1$  & $< 1$ \\ 
 \hline 
 \multirow{2}{*}{$100$} & $R^{{\rm s}, {\rm ls}}_{T,3}(1) \geq 1$  & 0.3454 & 0.1845  & $R^{{\rm s}, {\rm ls}}_{T,3}(5) \geq 1$  & 0.4974 & 0.1496 \\ 
& $R^{{\rm s}, {\rm ls}}_{T,3}(1) < 1$  & 0.159 & 0.3111 & $R^{{\rm s}, {\rm ls}}_{T,3}(5) < 1$  & 0.1513 & 0.2017 \\

\hline\hline
 \multirow{2}{*}{$n$} & & \multicolumn{2}{|c||}{$R^{{\rm s}, {\rm ls}}_{T,2}(1)$}  & & \multicolumn{2}{|c|}{$R^{{\rm s}, {\rm ls}}_{T,2}(5)$} \\ 
 & & $\geq 1$  & $< 1$  & & $\geq 1$  & $< 1$ \\ 
 \hline 
 \multirow{2}{*}{$1000$} & $R^{{\rm s}, {\rm ls}}_{T,3}(1) \geq 1$  & 0.3361 & 0.1589  & $R^{{\rm s}, {\rm ls}}_{T,3}(5) \geq 1$  & 0.4901 & 0.143 \\ 
& $R^{{\rm s}, {\rm ls}}_{T,3}(1) < 1$  & 0.1361 & 0.3689 & $R^{{\rm s}, {\rm ls}}_{T,3}(5) < 1$  & 0.131 & 0.2359 \\

\hline\hline
 \multirow{2}{*}{$n$} & & \multicolumn{2}{|c||}{$R^{{\rm s}, {\rm ls}}_{T,2}(1)$}  & & \multicolumn{2}{|c|}{$R^{{\rm s}, {\rm ls}}_{T,2}(5)$} \\ 
 & & $\geq 1$  & $< 1$  & & $\geq 1$  & $< 1$ \\ 
 \hline 
 \multirow{2}{*}{$10000$} & $R^{{\rm s}, {\rm ls}}_{T,3}(1) \geq 1$  & 0.322 & 0.151  & $R^{{\rm s}, {\rm ls}}_{T,3}(5) \geq 1$  & 0.4879 & 0.1411 \\ 
& $R^{{\rm s}, {\rm ls}}_{T,3}(1) < 1$  & 0.1408 & 0.3862 & $R^{{\rm s}, {\rm ls}}_{T,3}(5) < 1$  & 0.1215 & 0.2495 \\

\hline\hline
 \multirow{2}{*}{$n$} & & \multicolumn{2}{|c||}{$R^{{\rm s}, {\rm ls}}_{T,2}(1)$}  & & \multicolumn{2}{|c|}{$R^{{\rm s}, {\rm ls}}_{T,2}(5)$} \\ 
 & & $\geq 1.01$  & $< 1.01$  & & $\geq 1.01$  & $< 1.01$ \\ 
 \hline 
 \multirow{2}{*}{$10000$} & $R^{{\rm s}, {\rm ls}}_{T,3}(1) \geq 1.01$  & 0 & 0.0038  & $R^{{\rm s}, {\rm ls}}_{T,3}(5) \geq 1.01$  & 0 & 0.0053 \\ 
& $R^{{\rm s}, {\rm ls}}_{T,3}(1) < 1.01$  & 0.0031 & 0.9931 & $R^{{\rm s}, {\rm ls}}_{T,3}(5) < 1.01$  & 9e-04 & 0.9938 \\

\hline\hline
 \multirow{2}{*}{$n$} & & \multicolumn{2}{|c||}{$R^{{\rm s}, {\rm ls}}_{T,2}(1)$}  & & \multicolumn{2}{|c|}{$R^{{\rm s}, {\rm ls}}_{T,2}(5)$} \\ 
 & & $\geq 1.05$  & $< 1.05$  & & $\geq 1.05$  & $< 1.05$ \\ 
 \hline 
 \multirow{2}{*}{$10000$} & $R^{{\rm s}, {\rm ls}}_{T,3}(1) \geq 1.05$  & 0 & 0  & $R^{{\rm s}, {\rm ls}}_{T,3}(5) \geq 1.05$  & 0 & 0 \\ 
& $R^{{\rm s}, {\rm ls}}_{T,3}(1) < 1.05$  & 0 & 1 & $R^{{\rm s}, {\rm ls}}_{T,3}(5) < 1.05$  & 0 & 1 \\

\hline\hline
 \multirow{2}{*}{$n$} & & \multicolumn{2}{|c||}{$R^{{\rm s}, {\rm ls}}_{T,2}(1)$}  & & \multicolumn{2}{|c|}{$R^{{\rm s}, {\rm ls}}_{T,2}(5)$} \\ 
 & & $\geq 1.1$  & $< 1.1$  & & $\geq 1.1$  & $< 1.1$ \\ 
 \hline 
 \multirow{2}{*}{$10000$} & $R^{{\rm s}, {\rm ls}}_{T,3}(1) \geq 1.1$  & 0 & 0  & $R^{{\rm s}, {\rm ls}}_{T,3}(5) \geq 1.1$  & 0 & 0 \\ 
& $R^{{\rm s}, {\rm ls}}_{T,3}(1) < 1.1$  & 0 & 1 & $R^{{\rm s}, {\rm ls}}_{T,3}(5) < 1.1$  & 0 & 1 \\

\hline
\end{tabular}
\caption{\textit{Proportions of the individual events in~\eqref{samedecision} for the process \eqref{indepNonHetero} and selected combinations of $n$ and $\delta$.}}  \label{MSPEanalysisindepNonHeteroc}
\end{center}
\end{table}

\begin{table} \scriptsize
\begin{center}
	\input{tab2-m9-h1.tex}
	\input{tab2-m9-h5.tex}
\caption{\textit{Proportion of \eqref{samedecision}  being fulfilled for the process \eqref{indepNonHetero} and different values of $h$, $\delta$ and $n$.}}  \label{MSPEanalysisindepNonHeterob}
\end{center}
\end{table}

\begin{table} \scriptsize
\begin{center}
	\input{tab5-m9-h1.tex}
	\input{tab5-m9-h5.tex}
\caption{\textit{Values of $q(\delta)$, defined in \eqref{cond:f}, for the process \eqref{indepNonHetero} and different values of $h$, $\delta$ and $n$.}} \label{MSPEanalysisindepNonHeterod}
\end{center}
\end{table}

\begin{table} \scriptsize
\begin{center}
	\input{tab1-m9-h1-i2.tex}
	\input{tab1-m9-h1-i3.tex}
	\input{tab1-m9-h5-i2.tex}
	\input{tab1-m9-h5-i3.tex}
\caption{\textit{Proportion of \eqref{decisionrule}  being fulfilled for the process \eqref{indepNonHetero} and different values of $h$, $\delta$ and $n$.}}  \label{MSPEanalysisindepNonHetero}
\end{center}
\end{table}

\clearpage

\begin{figure}
\centering 
\includegraphics[width=0.24\textwidth]{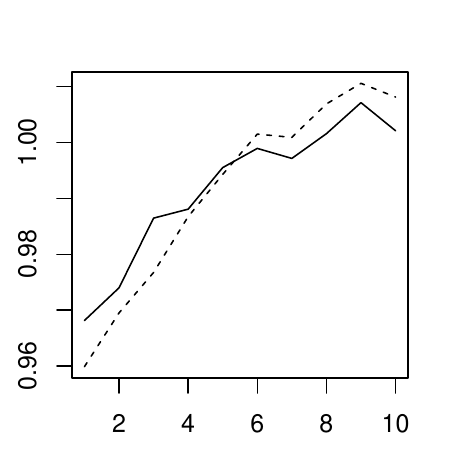}~
\includegraphics[width=0.24\textwidth]{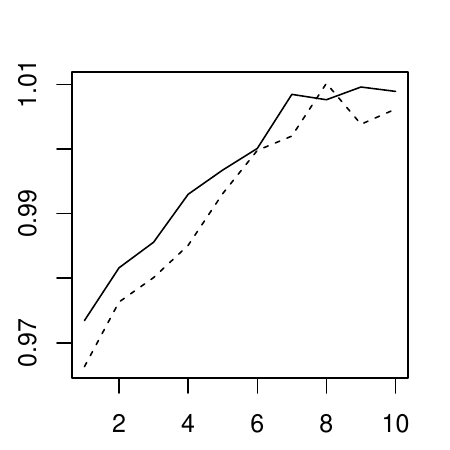}~
\includegraphics[width=0.24\textwidth]{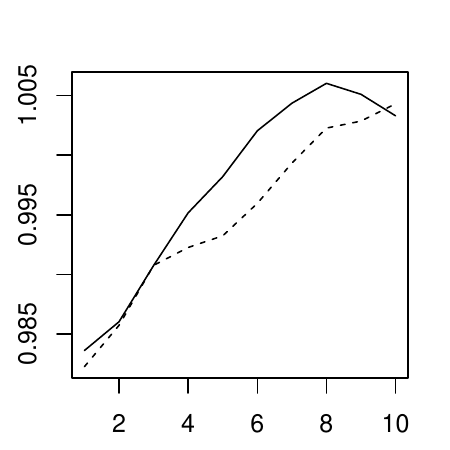}~ 
\includegraphics[width=0.24\textwidth]{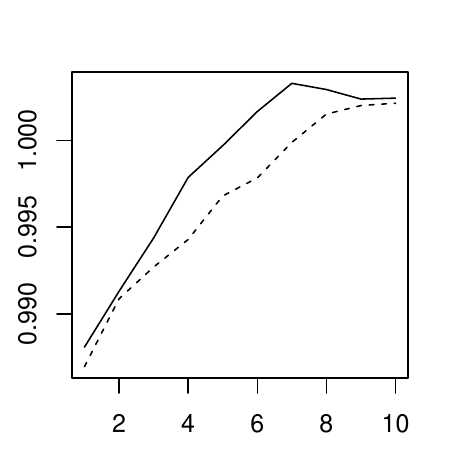} \\
\includegraphics[width=0.24\textwidth]{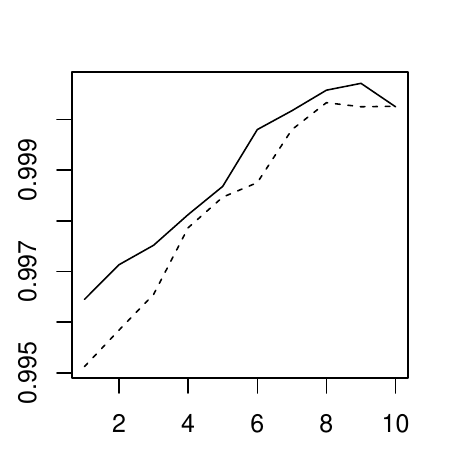}~
\includegraphics[width=0.24\textwidth]{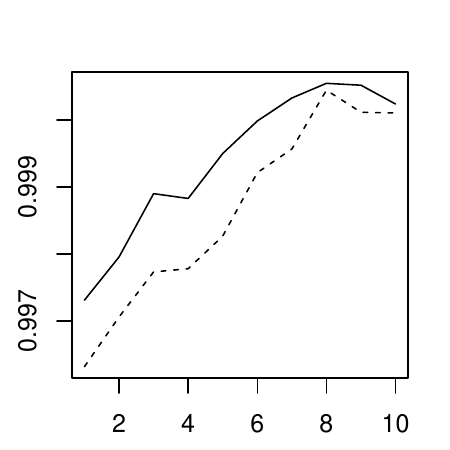}~ 
\includegraphics[width=0.24\textwidth]{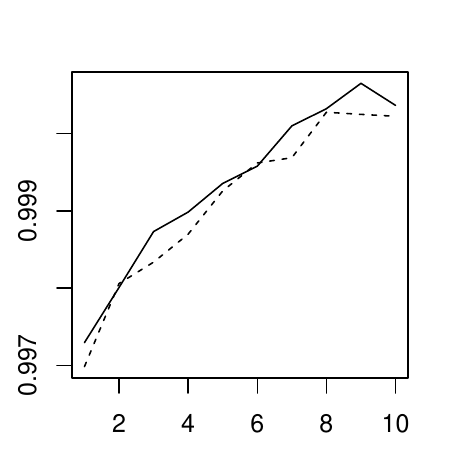}~ 
\includegraphics[width=0.24\textwidth]{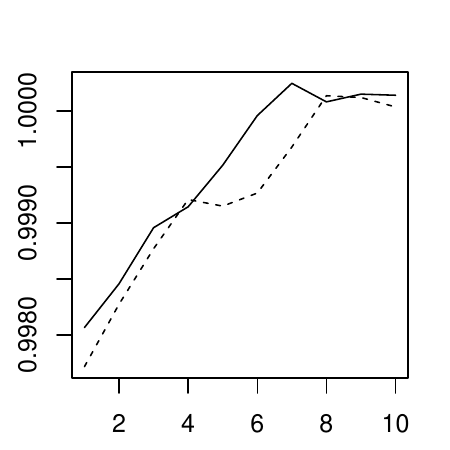}
   \caption{\it Plot of $h \mapsto R_{T,i}(h)$ for model \eqref{indepHetero} and different values of $n$ (from left to right: n=100, n=200, n=500, n=1000 [first row], n=4000, n=6000, n=8000, n=10000 [second row]). Solid line: $i=3$ (test set), dashed line: $i=2$ (validation set 2).} \label{MSPEindepHetero}
\end{figure}

\begin{table} \scriptsize
\begin{center}
\begin{tabular}{|c||c|c|c||c|c|c|}
\hline
 \multirow{2}{*}{$n$} & & \multicolumn{2}{|c||}{$R^{{\rm s}, {\rm ls}}_{T,2}(1)$}  & & \multicolumn{2}{|c|}{$R^{{\rm s}, {\rm ls}}_{T,2}(5)$} \\ 
 & & $\geq 1$  & $< 1$  & & $\geq 1$  & $< 1$ \\ 
 \hline 
 \multirow{2}{*}{$100$} & $R^{{\rm s}, {\rm ls}}_{T,3}(1) \geq 1$  & 0.3295 & 0.1982  & $R^{{\rm s}, {\rm ls}}_{T,3}(5) \geq 1$  & 0.4244 & 0.1675 \\ 
& $R^{{\rm s}, {\rm ls}}_{T,3}(1) < 1$  & 0.1677 & 0.3046 & $R^{{\rm s}, {\rm ls}}_{T,3}(5) < 1$  & 0.162 & 0.2461 \\

\hline\hline
 \multirow{2}{*}{$n$} & & \multicolumn{2}{|c||}{$R^{{\rm s}, {\rm ls}}_{T,2}(1)$}  & & \multicolumn{2}{|c|}{$R^{{\rm s}, {\rm ls}}_{T,2}(5)$} \\ 
 & & $\geq 1$  & $< 1$  & & $\geq 1$  & $< 1$ \\ 
 \hline 
 \multirow{2}{*}{$1000$} & $R^{{\rm s}, {\rm ls}}_{T,3}(1) \geq 1$  & 0.2815 & 0.1912  & $R^{{\rm s}, {\rm ls}}_{T,3}(5) \geq 1$  & 0.3784 & 0.1792 \\ 
& $R^{{\rm s}, {\rm ls}}_{T,3}(1) < 1$  & 0.1624 & 0.3649 & $R^{{\rm s}, {\rm ls}}_{T,3}(5) < 1$  & 0.1625 & 0.2799 \\

\hline\hline
 \multirow{2}{*}{$n$} & & \multicolumn{2}{|c||}{$R^{{\rm s}, {\rm ls}}_{T,2}(1)$}  & & \multicolumn{2}{|c|}{$R^{{\rm s}, {\rm ls}}_{T,2}(5)$} \\ 
 & & $\geq 1$  & $< 1$  & & $\geq 1$  & $< 1$ \\ 
 \hline 
 \multirow{2}{*}{$10000$} & $R^{{\rm s}, {\rm ls}}_{T,3}(1) \geq 1$  & 0.249 & 0.1898  & $R^{{\rm s}, {\rm ls}}_{T,3}(5) \geq 1$  & 0.3535 & 0.18 \\ 
& $R^{{\rm s}, {\rm ls}}_{T,3}(1) < 1$  & 0.172 & 0.3892 & $R^{{\rm s}, {\rm ls}}_{T,3}(5) < 1$  & 0.1588 & 0.3077 \\

\hline\hline
 \multirow{2}{*}{$n$} & & \multicolumn{2}{|c||}{$R^{{\rm s}, {\rm ls}}_{T,2}(1)$}  & & \multicolumn{2}{|c|}{$R^{{\rm s}, {\rm ls}}_{T,2}(5)$} \\ 
 & & $\geq 1.01$  & $< 1.01$  & & $\geq 1.01$  & $< 1.01$ \\ 
 \hline 
 \multirow{2}{*}{$10000$} & $R^{{\rm s}, {\rm ls}}_{T,3}(1) \geq 1.01$  & 0 & 0.009  & $R^{{\rm s}, {\rm ls}}_{T,3}(5) \geq 1.01$  & 0.0025 & 0.0212 \\ 
& $R^{{\rm s}, {\rm ls}}_{T,3}(1) < 1.01$  & 0.0048 & 0.9862 & $R^{{\rm s}, {\rm ls}}_{T,3}(5) < 1.01$  & 0.0116 & 0.9647 \\

\hline\hline
 \multirow{2}{*}{$n$} & & \multicolumn{2}{|c||}{$R^{{\rm s}, {\rm ls}}_{T,2}(1)$}  & & \multicolumn{2}{|c|}{$R^{{\rm s}, {\rm ls}}_{T,2}(5)$} \\ 
 & & $\geq 1.05$  & $< 1.05$  & & $\geq 1.05$  & $< 1.05$ \\ 
 \hline 
 \multirow{2}{*}{$10000$} & $R^{{\rm s}, {\rm ls}}_{T,3}(1) \geq 1.05$  & 0 & 0  & $R^{{\rm s}, {\rm ls}}_{T,3}(5) \geq 1.05$  & 0 & 0 \\ 
& $R^{{\rm s}, {\rm ls}}_{T,3}(1) < 1.05$  & 0 & 1 & $R^{{\rm s}, {\rm ls}}_{T,3}(5) < 1.05$  & 0 & 1 \\

\hline\hline
 \multirow{2}{*}{$n$} & & \multicolumn{2}{|c||}{$R^{{\rm s}, {\rm ls}}_{T,2}(1)$}  & & \multicolumn{2}{|c|}{$R^{{\rm s}, {\rm ls}}_{T,2}(5)$} \\ 
 & & $\geq 1.1$  & $< 1.1$  & & $\geq 1.1$  & $< 1.1$ \\ 
 \hline 
 \multirow{2}{*}{$10000$} & $R^{{\rm s}, {\rm ls}}_{T,3}(1) \geq 1.1$  & 0 & 0  & $R^{{\rm s}, {\rm ls}}_{T,3}(5) \geq 1.1$  & 0 & 0 \\ 
& $R^{{\rm s}, {\rm ls}}_{T,3}(1) < 1.1$  & 0 & 1 & $R^{{\rm s}, {\rm ls}}_{T,3}(5) < 1.1$  & 0 & 1 \\

\hline
\end{tabular}
\caption{\textit{Proportions of the individual events in~\eqref{samedecision} for the process \eqref{indepHetero} and selected combinations of $n$ and $\delta$.}}  \label{MSPEanalysisindepHeteroc}
\end{center}
\end{table}

\begin{table} \scriptsize
\begin{center}
	\input{tab2-m10-h1.tex}
	\input{tab2-m10-h5.tex}
\caption{\textit{Proportion of \eqref{samedecision}  being fulfilled for the process \eqref{indepHetero} and different values of $h$, $\delta$ and $n$.}}  \label{MSPEanalysisindepHeterob}
\end{center}
\end{table}

\begin{table} \scriptsize
\begin{center}
	\input{tab5-m10-h1.tex}
	\input{tab5-m10-h5.tex}
\caption{\textit{Values of $q(\delta)$, defined in \eqref{cond:f}, for the process \eqref{indepHetero} and different values of $h$, $\delta$ and $n$.}} \label{MSPEanalysisindepHeterod}
\end{center}
\end{table}

\begin{table} \scriptsize
\begin{center}
	\input{tab1-m10-h1-i2.tex}
	\input{tab1-m10-h1-i3.tex}
	\input{tab1-m10-h5-i2.tex}
	\input{tab1-m10-h5-i3.tex}
\caption{\textit{Proportion of \eqref{decisionrule}  being fulfilled for the process \eqref{indepHetero} and different values of $h$, $\delta$ and $n$.}}  \label{MSPEanalysisindepHetero}
\end{center}
\end{table}

\clearpage

\begin{figure}
\centering 
\includegraphics[width=0.24\textwidth]{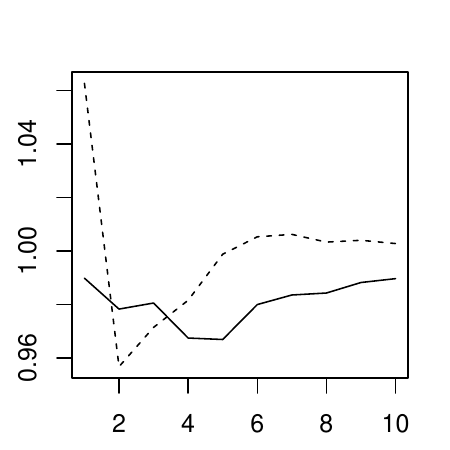}~
\includegraphics[width=0.24\textwidth]{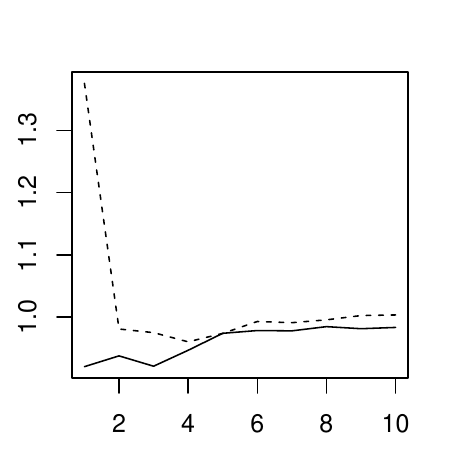}~
\includegraphics[width=0.24\textwidth]{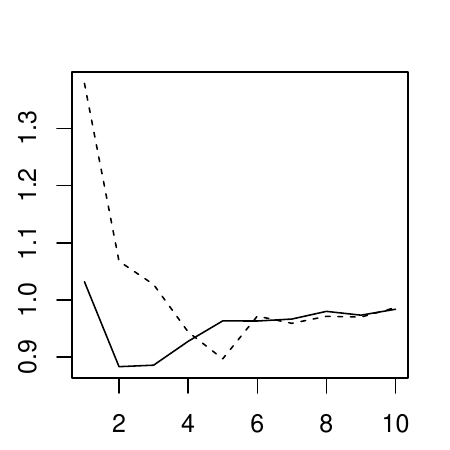}~ 
\includegraphics[width=0.24\textwidth]{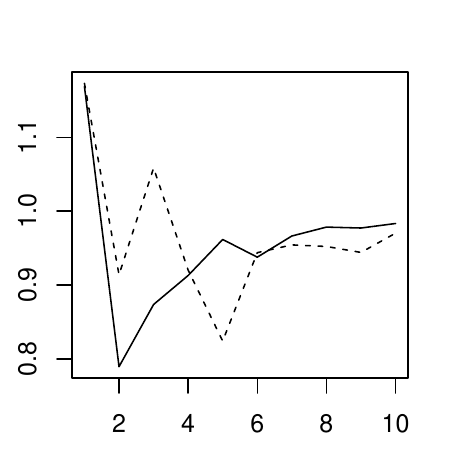} \\
\includegraphics[width=0.24\textwidth]{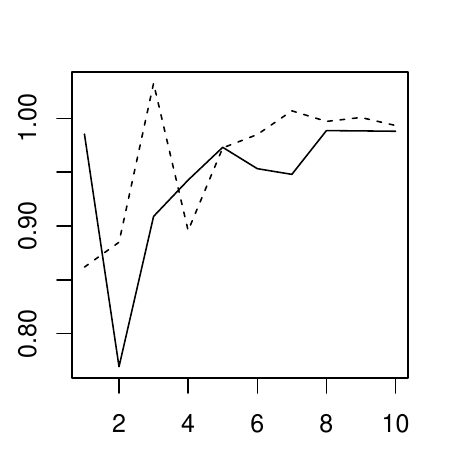}~
\includegraphics[width=0.24\textwidth]{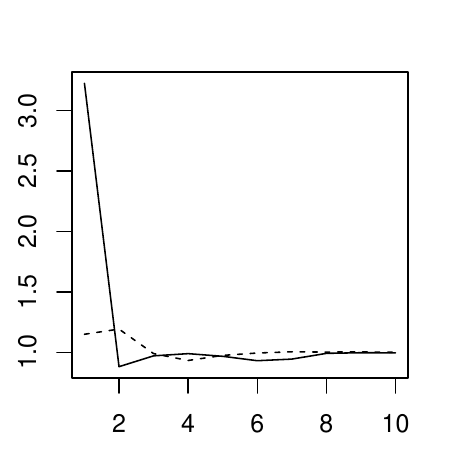}~ 
\includegraphics[width=0.24\textwidth]{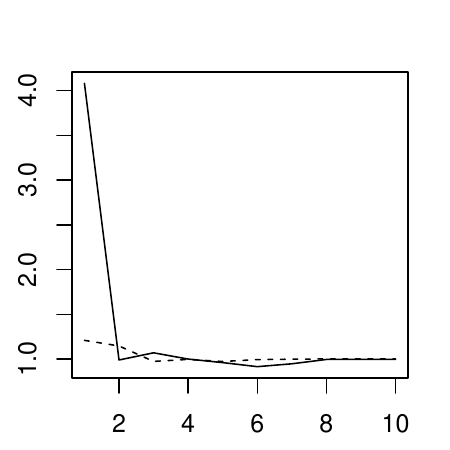}~ 
\includegraphics[width=0.24\textwidth]{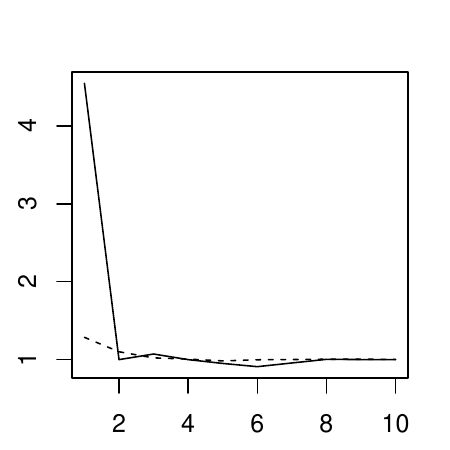}
   \caption{\it Plot of $h \mapsto R_{T,i}(h)$ for model \eqref{mdl11} and different values of $n$ (from left to right: n=100, n=200, n=500, n=1000 [first row], n=4000, n=6000, n=8000, n=10000 [second row]). Solid line: $i=3$ (test set), dashed line: $i=2$ (validation set 2). } \label{MSPEmdl11}
\end{figure}

\begin{table} \scriptsize
\begin{center}
\begin{tabular}{|c||c|c|c||c|c|c|}
\hline
 \multirow{2}{*}{$n$} & & \multicolumn{2}{|c||}{$R^{{\rm s}, {\rm ls}}_{T,2}(1)$}  & & \multicolumn{2}{|c|}{$R^{{\rm s}, {\rm ls}}_{T,2}(5)$} \\ 
 & & $\geq 1$  & $< 1$  & & $\geq 1$  & $< 1$ \\ 
 \hline 
 \multirow{2}{*}{$100$} & $R^{{\rm s}, {\rm ls}}_{T,3}(1) \geq 1$  & 0.6049 & 0.1541  & $R^{{\rm s}, {\rm ls}}_{T,3}(5) \geq 1$  & 0.6493 & 0.095 \\ 
& $R^{{\rm s}, {\rm ls}}_{T,3}(1) < 1$  & 0.1673 & 0.0737 & $R^{{\rm s}, {\rm ls}}_{T,3}(5) < 1$  & 0.1329 & 0.1228 \\

\hline\hline
 \multirow{2}{*}{$n$} & & \multicolumn{2}{|c||}{$R^{{\rm s}, {\rm ls}}_{T,2}(1)$}  & & \multicolumn{2}{|c|}{$R^{{\rm s}, {\rm ls}}_{T,2}(5)$} \\ 
 & & $\geq 1$  & $< 1$  & & $\geq 1$  & $< 1$ \\ 
 \hline 
 \multirow{2}{*}{$1000$} & $R^{{\rm s}, {\rm ls}}_{T,3}(1) \geq 1$  & 0.4839 & 0.2381  & $R^{{\rm s}, {\rm ls}}_{T,3}(5) \geq 1$  & 0.0445 & 0.3321 \\ 
& $R^{{\rm s}, {\rm ls}}_{T,3}(1) < 1$  & 0.2341 & 0.0439 & $R^{{\rm s}, {\rm ls}}_{T,3}(5) < 1$  & 0.0907 & 0.5327 \\

\hline\hline
 \multirow{2}{*}{$n$} & & \multicolumn{2}{|c||}{$R^{{\rm s}, {\rm ls}}_{T,2}(1)$}  & & \multicolumn{2}{|c|}{$R^{{\rm s}, {\rm ls}}_{T,2}(5)$} \\ 
 & & $\geq 1$  & $< 1$  & & $\geq 1$  & $< 1$ \\ 
 \hline 
 \multirow{2}{*}{$10000$} & $R^{{\rm s}, {\rm ls}}_{T,3}(1) \geq 1$  & 0.9967 & 0.0033  & $R^{{\rm s}, {\rm ls}}_{T,3}(5) \geq 1$  & 4e-04 & 5e-04 \\ 
& $R^{{\rm s}, {\rm ls}}_{T,3}(1) < 1$  & 0 & 0 & $R^{{\rm s}, {\rm ls}}_{T,3}(5) < 1$  & 0.0677 & 0.9314 \\

\hline\hline
 \multirow{2}{*}{$n$} & & \multicolumn{2}{|c||}{$R^{{\rm s}, {\rm ls}}_{T,2}(1)$}  & & \multicolumn{2}{|c|}{$R^{{\rm s}, {\rm ls}}_{T,2}(5)$} \\ 
 & & $\geq 1.01$  & $< 1.01$  & & $\geq 1.01$  & $< 1.01$ \\ 
 \hline 
 \multirow{2}{*}{$10000$} & $R^{{\rm s}, {\rm ls}}_{T,3}(1) \geq 1.01$  & 0.9957 & 0.0043  & $R^{{\rm s}, {\rm ls}}_{T,3}(5) \geq 1.01$  & 0 & 2e-04 \\ 
& $R^{{\rm s}, {\rm ls}}_{T,3}(1) < 1.01$  & 0 & 0 & $R^{{\rm s}, {\rm ls}}_{T,3}(5) < 1.01$  & 0.0476 & 0.9522 \\

\hline\hline
 \multirow{2}{*}{$n$} & & \multicolumn{2}{|c||}{$R^{{\rm s}, {\rm ls}}_{T,2}(1)$}  & & \multicolumn{2}{|c|}{$R^{{\rm s}, {\rm ls}}_{T,2}(5)$} \\ 
 & & $\geq 1.05$  & $< 1.05$  & & $\geq 1.05$  & $< 1.05$ \\ 
 \hline 
 \multirow{2}{*}{$10000$} & $R^{{\rm s}, {\rm ls}}_{T,3}(1) \geq 1.05$  & 0.9859 & 0.0141  & $R^{{\rm s}, {\rm ls}}_{T,3}(5) \geq 1.05$  & 0 & 0 \\ 
& $R^{{\rm s}, {\rm ls}}_{T,3}(1) < 1.05$  & 0 & 0 & $R^{{\rm s}, {\rm ls}}_{T,3}(5) < 1.05$  & 0.0162 & 0.9838 \\

\hline\hline
 \multirow{2}{*}{$n$} & & \multicolumn{2}{|c||}{$R^{{\rm s}, {\rm ls}}_{T,2}(1)$}  & & \multicolumn{2}{|c|}{$R^{{\rm s}, {\rm ls}}_{T,2}(5)$} \\ 
 & & $\geq 1.1$  & $< 1.1$  & & $\geq 1.1$  & $< 1.1$ \\ 
 \hline 
 \multirow{2}{*}{$10000$} & $R^{{\rm s}, {\rm ls}}_{T,3}(1) \geq 1.1$  & 0.9492 & 0.0508  & $R^{{\rm s}, {\rm ls}}_{T,3}(5) \geq 1.1$  & 0 & 0 \\ 
& $R^{{\rm s}, {\rm ls}}_{T,3}(1) < 1.1$  & 0 & 0 & $R^{{\rm s}, {\rm ls}}_{T,3}(5) < 1.1$  & 0.0025 & 0.9975 \\

\hline
\end{tabular}
\caption{\textit{Proportions of the individual events in~\eqref{samedecision} for the process \eqref{mdl11} and selected combinations of $n$ and $\delta$.}}  \label{MSPEanalysismdl11c}
\end{center}
\end{table}

\begin{table} \scriptsize
\begin{center}
	\input{tab2-m11-h1.tex}
	\input{tab2-m11-h5.tex}
\caption{\textit{Proportion of \eqref{samedecision}  being fulfilled for the process \eqref{mdl11} and different values of $h$, $\delta$ and $n$.}}  \label{MSPEanalysismdl11b}
\end{center}
\end{table}

\begin{table} \scriptsize
\begin{center}
	\input{tab5-m11-h1.tex}
	\input{tab5-m11-h5.tex}
\caption{\textit{Values of $q(\delta)$, defined in \eqref{cond:f}, for the process \eqref{mdl11} and different values of $h$, $\delta$ and $n$.}} \label{MSPEanalysismdl11d}
\end{center}
\end{table}

\begin{table} \scriptsize
\begin{center}
	\input{tab1-m11-h1-i2.tex}
	\input{tab1-m11-h1-i3.tex}
	\input{tab1-m11-h5-i2.tex}
	\input{tab1-m11-h5-i3.tex}
\caption{\textit{Proportion of \eqref{decisionrule}  being fulfilled for the process \eqref{mdl11} and different values of $h$, $\delta$ and $n$.}}  \label{MSPEanalysismdl11}
\end{center}
\end{table}

\clearpage

\begin{figure}
\centering 
\includegraphics[width=0.24\textwidth]{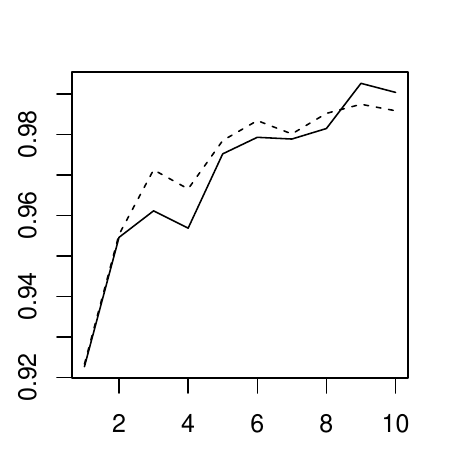}~
\includegraphics[width=0.24\textwidth]{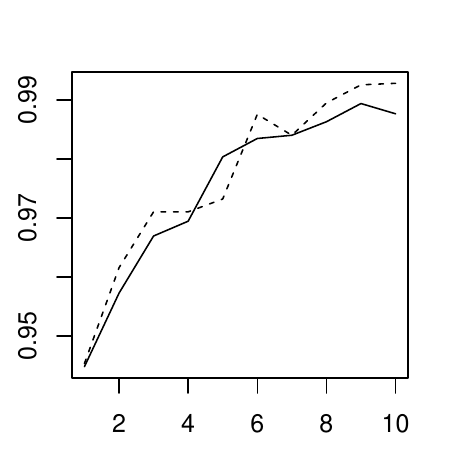}~
\includegraphics[width=0.24\textwidth]{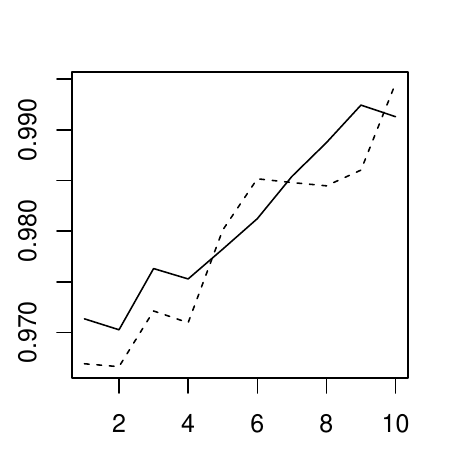}~ 
\includegraphics[width=0.24\textwidth]{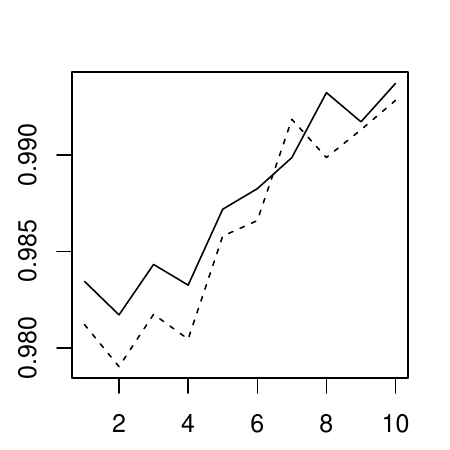} \\
\includegraphics[width=0.24\textwidth]{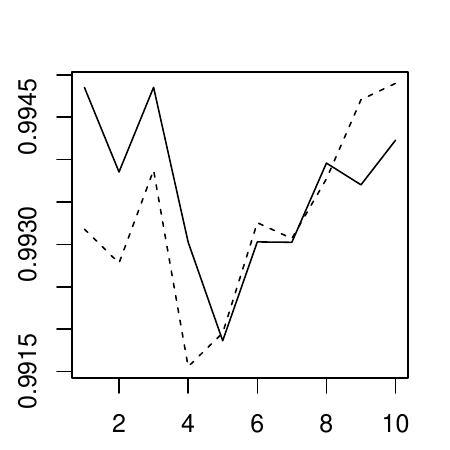}~
\includegraphics[width=0.24\textwidth]{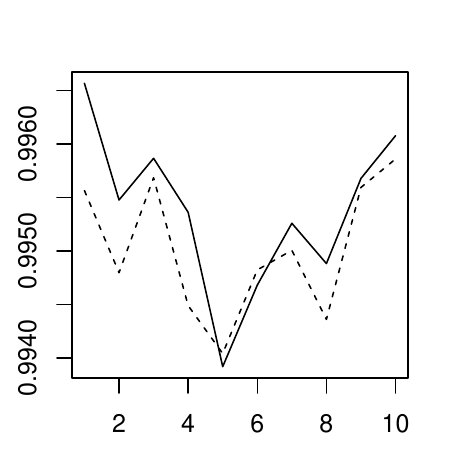}~ 
\includegraphics[width=0.24\textwidth]{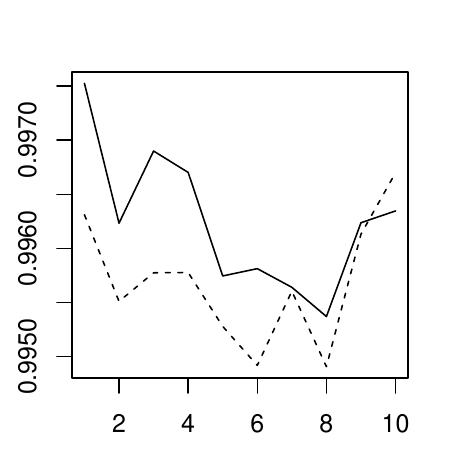}~ 
\includegraphics[width=0.24\textwidth]{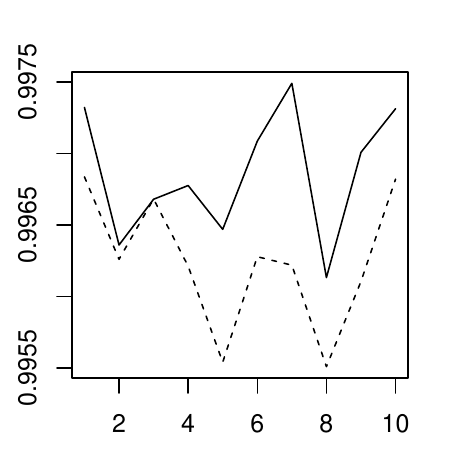}
   \caption{\it Plot of $h \mapsto R_{T,i}(h)$ for model \eqref{mdl12} and different values of $n$ (from left to right: n=100, n=200, n=500, n=1000 [first row], n=4000, n=6000, n=8000, n=10000 [second row]). Solid line: $i=3$ (test set), dashed line: $i=2$ (validation set 2). } \label{MSPEmdl12}
\end{figure}

\begin{table} \scriptsize
\begin{center}
\begin{tabular}{|c||c|c|c||c|c|c|}
\hline
 \multirow{2}{*}{$n$} & & \multicolumn{2}{|c||}{$R^{{\rm s}, {\rm ls}}_{T,2}(1)$}  & & \multicolumn{2}{|c|}{$R^{{\rm s}, {\rm ls}}_{T,2}(5)$} \\ 
 & & $\geq 1$  & $< 1$  & & $\geq 1$  & $< 1$ \\ 
 \hline 
 \multirow{2}{*}{$100$} & $R^{{\rm s}, {\rm ls}}_{T,3}(1) \geq 1$  & 0.1122 & 0.2114  & $R^{{\rm s}, {\rm ls}}_{T,3}(5) \geq 1$  & 0.4153 & 0.1682 \\ 
& $R^{{\rm s}, {\rm ls}}_{T,3}(1) < 1$  & 0.2042 & 0.4722 & $R^{{\rm s}, {\rm ls}}_{T,3}(5) < 1$  & 0.1728 & 0.2437 \\

\hline\hline
 \multirow{2}{*}{$n$} & & \multicolumn{2}{|c||}{$R^{{\rm s}, {\rm ls}}_{T,2}(1)$}  & & \multicolumn{2}{|c|}{$R^{{\rm s}, {\rm ls}}_{T,2}(5)$} \\ 
 & & $\geq 1$  & $< 1$  & & $\geq 1$  & $< 1$ \\ 
 \hline 
 \multirow{2}{*}{$1000$} & $R^{{\rm s}, {\rm ls}}_{T,3}(1) \geq 1$  & 0.0714 & 0.2008  & $R^{{\rm s}, {\rm ls}}_{T,3}(5) \geq 1$  & 0.1572 & 0.2043 \\ 
& $R^{{\rm s}, {\rm ls}}_{T,3}(1) < 1$  & 0.1824 & 0.5454 & $R^{{\rm s}, {\rm ls}}_{T,3}(5) < 1$  & 0.1987 & 0.4398 \\

\hline\hline
 \multirow{2}{*}{$n$} & & \multicolumn{2}{|c||}{$R^{{\rm s}, {\rm ls}}_{T,2}(1)$}  & & \multicolumn{2}{|c|}{$R^{{\rm s}, {\rm ls}}_{T,2}(5)$} \\ 
 & & $\geq 1$  & $< 1$  & & $\geq 1$  & $< 1$ \\ 
 \hline 
 \multirow{2}{*}{$10000$} & $R^{{\rm s}, {\rm ls}}_{T,3}(1) \geq 1$  & 0.0507 & 0.1897  & $R^{{\rm s}, {\rm ls}}_{T,3}(5) \geq 1$  & 0.0639 & 0.2063 \\ 
& $R^{{\rm s}, {\rm ls}}_{T,3}(1) < 1$  & 0.1756 & 0.584 & $R^{{\rm s}, {\rm ls}}_{T,3}(5) < 1$  & 0.1903 & 0.5395 \\

\hline\hline
 \multirow{2}{*}{$n$} & & \multicolumn{2}{|c||}{$R^{{\rm s}, {\rm ls}}_{T,2}(1)$}  & & \multicolumn{2}{|c|}{$R^{{\rm s}, {\rm ls}}_{T,2}(5)$} \\ 
 & & $\geq 1.01$  & $< 1.01$  & & $\geq 1.01$  & $< 1.01$ \\ 
 \hline 
 \multirow{2}{*}{$10000$} & $R^{{\rm s}, {\rm ls}}_{T,3}(1) \geq 1.01$  & 0 & 0.0028  & $R^{{\rm s}, {\rm ls}}_{T,3}(5) \geq 1.01$  & 0.0011 & 0.0279 \\ 
& $R^{{\rm s}, {\rm ls}}_{T,3}(1) < 1.01$  & 0.0022 & 0.995 & $R^{{\rm s}, {\rm ls}}_{T,3}(5) < 1.01$  & 0.0255 & 0.9455 \\

\hline\hline
 \multirow{2}{*}{$n$} & & \multicolumn{2}{|c||}{$R^{{\rm s}, {\rm ls}}_{T,2}(1)$}  & & \multicolumn{2}{|c|}{$R^{{\rm s}, {\rm ls}}_{T,2}(5)$} \\ 
 & & $\geq 1.05$  & $< 1.05$  & & $\geq 1.05$  & $< 1.05$ \\ 
 \hline 
 \multirow{2}{*}{$10000$} & $R^{{\rm s}, {\rm ls}}_{T,3}(1) \geq 1.05$  & 0 & 0  & $R^{{\rm s}, {\rm ls}}_{T,3}(5) \geq 1.05$  & 0 & 0 \\ 
& $R^{{\rm s}, {\rm ls}}_{T,3}(1) < 1.05$  & 0 & 1 & $R^{{\rm s}, {\rm ls}}_{T,3}(5) < 1.05$  & 0 & 1 \\

\hline\hline
 \multirow{2}{*}{$n$} & & \multicolumn{2}{|c||}{$R^{{\rm s}, {\rm ls}}_{T,2}(1)$}  & & \multicolumn{2}{|c|}{$R^{{\rm s}, {\rm ls}}_{T,2}(5)$} \\ 
 & & $\geq 1.1$  & $< 1.1$  & & $\geq 1.1$  & $< 1.1$ \\ 
 \hline 
 \multirow{2}{*}{$10000$} & $R^{{\rm s}, {\rm ls}}_{T,3}(1) \geq 1.1$  & 0 & 0  & $R^{{\rm s}, {\rm ls}}_{T,3}(5) \geq 1.1$  & 0 & 0 \\ 
& $R^{{\rm s}, {\rm ls}}_{T,3}(1) < 1.1$  & 0 & 1 & $R^{{\rm s}, {\rm ls}}_{T,3}(5) < 1.1$  & 0 & 1 \\

\hline
\end{tabular}
\caption{\textit{Proportions of the individual events in~\eqref{samedecision} for the process \eqref{mdl12} and selected combinations of $n$ and $\delta$.}}  \label{MSPEanalysismdl12c}
\end{center}
\end{table}

\begin{table} \scriptsize
\begin{center}
	\input{tab2-m12-h1.tex}
	\input{tab2-m12-h5.tex} 
\caption{\textit{Proportion of \eqref{samedecision}  being fulfilled for the process \eqref{mdl12} and different values of $h$, $\delta$ and $n$.}}  \label{MSPEanalysismdl12b}
\end{center}
\end{table}

\begin{table} \scriptsize
\begin{center}
	\input{tab5-m12-h1.tex}
	\input{tab5-m12-h5.tex}
\caption{\textit{Values of $q(\delta)$, defined in \eqref{cond:f}, for the process \eqref{mdl12} and different values of $h$, $\delta$ and $n$.}} \label{MSPEanalysismdl12d}
\end{center}
\end{table}

\begin{table} \scriptsize
\begin{center}
	\input{tab1-m12-h1-i2.tex}
	\input{tab1-m12-h1-i3.tex}
	\input{tab1-m12-h5-i2.tex}
	\input{tab1-m12-h5-i3.tex}
\caption{\textit{Proportion of \eqref{decisionrule}  being fulfilled for the process \eqref{mdl12} and different values of $h$, $\delta$ and $n$.}}  \label{MSPEanalysismdl12}
\end{center}
\end{table}

\clearpage

\begin{figure}
\centering 
\includegraphics[width=0.24\textwidth]{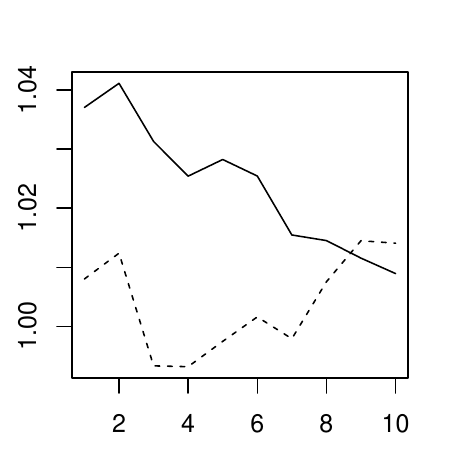}~
\includegraphics[width=0.24\textwidth]{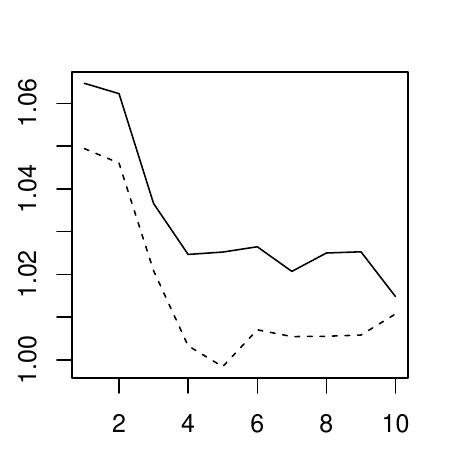}~
\includegraphics[width=0.24\textwidth]{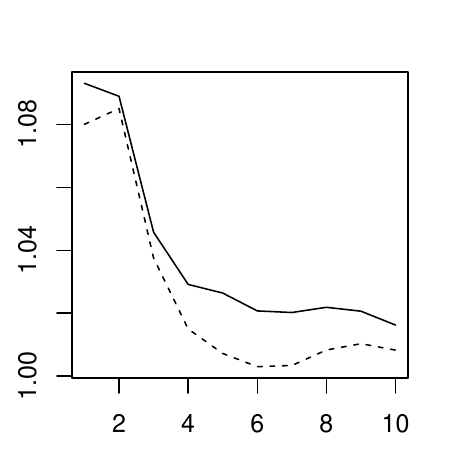}~ 
\includegraphics[width=0.24\textwidth]{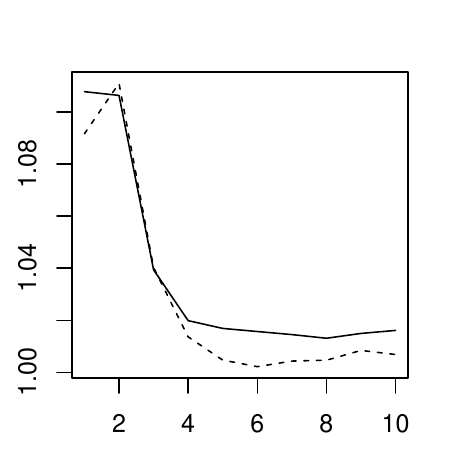} \\
\includegraphics[width=0.24\textwidth]{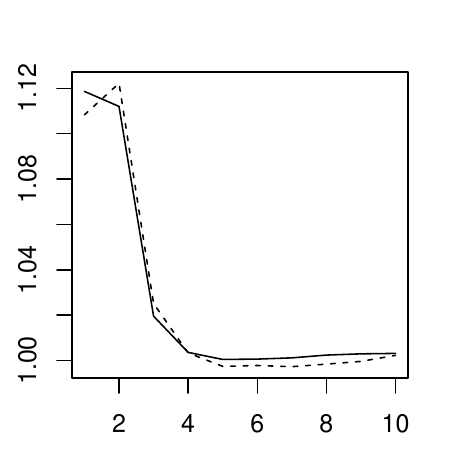}~
\includegraphics[width=0.24\textwidth]{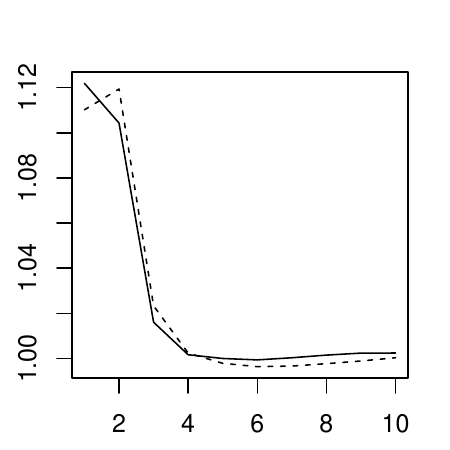}~ 
\includegraphics[width=0.24\textwidth]{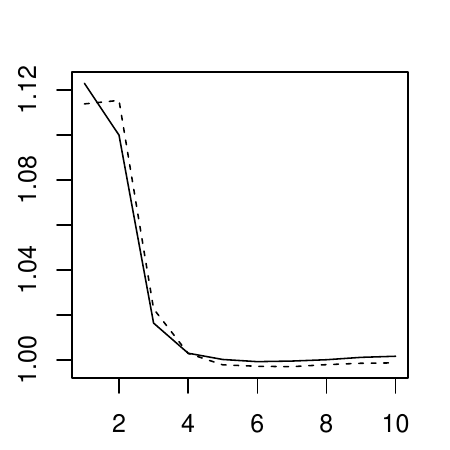}~ 
\includegraphics[width=0.24\textwidth]{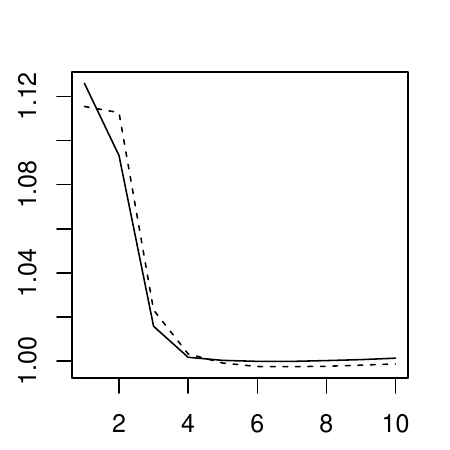}
   \caption{\it Plot of $h \mapsto R_{T,i}(h)$ for model \eqref{decreasing1} and different values of $n$ (from left to right: n=100, n=200, n=500, n=1000 [first row], n=4000, n=6000, n=8000, n=10000 [second row]). Solid line: $i=3$ (test set), dashed line: $i=2$ (validation set 2). } \label{MSPEdecreasing1}
\end{figure}

\begin{table} \scriptsize
\begin{center}
\begin{tabular}{|c||c|c|c||c|c|c|}
\hline
 \multirow{2}{*}{$n$} & & \multicolumn{2}{|c||}{$R^{{\rm s}, {\rm ls}}_{T,2}(1)$}  & & \multicolumn{2}{|c|}{$R^{{\rm s}, {\rm ls}}_{T,2}(5)$} \\ 
 & & $\geq 1.01$  & $< 1.01$  & & $\geq 1.01$  & $< 1.01$ \\ 
 \hline 
 \multirow{2}{*}{$100$} & $R^{{\rm s}, {\rm ls}}_{T,3}(1) \geq 1.01$  & 0.2793 & 0.2405  & $R^{{\rm s}, {\rm ls}}_{T,3}(5) \geq 1.01$  & 0.1792 & 0.1822 \\ 
& $R^{{\rm s}, {\rm ls}}_{T,3}(1) < 1.01$  & 0.1822 & 0.298 & $R^{{\rm s}, {\rm ls}}_{T,3}(5) < 1.01$  & 0.1483 & 0.4903 \\

\hline\hline
 \multirow{2}{*}{$n$} & & \multicolumn{2}{|c||}{$R^{{\rm s}, {\rm ls}}_{T,2}(1)$}  & & \multicolumn{2}{|c|}{$R^{{\rm s}, {\rm ls}}_{T,2}(5)$} \\ 
 & & $\geq 1.1$  & $< 1.1$  & & $\geq 1.1$  & $< 1.1$ \\ 
 \hline 
 \multirow{2}{*}{$1000$} & $R^{{\rm s}, {\rm ls}}_{T,3}(1) \geq 1.1$  & 0.263 & 0.2781  & $R^{{\rm s}, {\rm ls}}_{T,3}(5) \geq 1.1$  & 0.0635 & 0.0688 \\ 
& $R^{{\rm s}, {\rm ls}}_{T,3}(1) < 1.1$  & 0.1587 & 0.3002 & $R^{{\rm s}, {\rm ls}}_{T,3}(5) < 1.1$  & 0.0471 & 0.8206 \\

\hline\hline
 \multirow{2}{*}{$n$} & & \multicolumn{2}{|c||}{$R^{{\rm s}, {\rm ls}}_{T,2}(1)$}  & & \multicolumn{2}{|c|}{$R^{{\rm s}, {\rm ls}}_{T,2}(5)$} \\ 
 & & $\geq 1$  & $< 1$  & & $\geq 1$  & $< 1$ \\ 
 \hline 
 \multirow{2}{*}{$10000$} & $R^{{\rm s}, {\rm ls}}_{T,3}(1) \geq 1$  & 1 & 0  & $R^{{\rm s}, {\rm ls}}_{T,3}(5) \geq 1$  & 0.3935 & 0.235 \\ 
& $R^{{\rm s}, {\rm ls}}_{T,3}(1) < 1$  & 0 & 0 & $R^{{\rm s}, {\rm ls}}_{T,3}(5) < 1$  & 0.2006 & 0.1709 \\

\hline\hline
 \multirow{2}{*}{$n$} & & \multicolumn{2}{|c||}{$R^{{\rm s}, {\rm ls}}_{T,2}(1)$}  & & \multicolumn{2}{|c|}{$R^{{\rm s}, {\rm ls}}_{T,2}(5)$} \\ 
 & & $\geq 1.01$  & $< 1.01$  & & $\geq 1.01$  & $< 1.01$ \\ 
 \hline 
 \multirow{2}{*}{$10000$} & $R^{{\rm s}, {\rm ls}}_{T,3}(1) \geq 1.01$  & 1 & 0  & $R^{{\rm s}, {\rm ls}}_{T,3}(5) \geq 1.01$  & 0.0068 & 0.1104 \\ 
& $R^{{\rm s}, {\rm ls}}_{T,3}(1) < 1.01$  & 0 & 0 & $R^{{\rm s}, {\rm ls}}_{T,3}(5) < 1.01$  & 0.09 & 0.7928 \\

\hline\hline
 \multirow{2}{*}{$n$} & & \multicolumn{2}{|c||}{$R^{{\rm s}, {\rm ls}}_{T,2}(1)$}  & & \multicolumn{2}{|c|}{$R^{{\rm s}, {\rm ls}}_{T,2}(5)$} \\ 
 & & $\geq 1.05$  & $< 1.05$  & & $\geq 1.05$  & $< 1.05$ \\ 
 \hline 
 \multirow{2}{*}{$10000$} & $R^{{\rm s}, {\rm ls}}_{T,3}(1) \geq 1.05$  & 0.9995 & 3e-04  & $R^{{\rm s}, {\rm ls}}_{T,3}(5) \geq 1.05$  & 7e-04 & 8e-04 \\ 
& $R^{{\rm s}, {\rm ls}}_{T,3}(1) < 1.05$  & 2e-04 & 0 & $R^{{\rm s}, {\rm ls}}_{T,3}(5) < 1.05$  & 0 & 0.9985 \\

\hline\hline
 \multirow{2}{*}{$n$} & & \multicolumn{2}{|c||}{$R^{{\rm s}, {\rm ls}}_{T,2}(1)$}  & & \multicolumn{2}{|c|}{$R^{{\rm s}, {\rm ls}}_{T,2}(5)$} \\ 
 & & $\geq 1.1$  & $< 1.1$  & & $\geq 1.1$  & $< 1.1$ \\ 
 \hline 
 \multirow{2}{*}{$10000$} & $R^{{\rm s}, {\rm ls}}_{T,3}(1) \geq 1.1$  & 0.7133 & 0.1876  & $R^{{\rm s}, {\rm ls}}_{T,3}(5) \geq 1.1$  & 5e-04 & 1e-04 \\ 
& $R^{{\rm s}, {\rm ls}}_{T,3}(1) < 1.1$  & 0.0678 & 0.0313 & $R^{{\rm s}, {\rm ls}}_{T,3}(5) < 1.1$  & 1e-04 & 0.9993 \\

\hline
\end{tabular}
\caption{\textit{Proportions of the individual events in~\eqref{samedecision} for the process \eqref{decreasing1} and selected combinations of $n$ and $\delta$.}}  \label{MSPEanalysisdecreasing1c}
\end{center}
\end{table}

\begin{table} \scriptsize
\begin{center}
	\input{tab2-m13-h1.tex}
	\input{tab2-m13-h5.tex}  
\caption{\textit{Proportion of \eqref{samedecision}  being fulfilled for the process \eqref{decreasing1} and different values of $h$, $\delta$ and $n$.}}  \label{MSPEanalysisdecreasing1b}
\end{center}
\end{table}

\begin{table} \scriptsize
\begin{center}
	\input{tab5-m13-h1.tex}
	\input{tab5-m13-h5.tex}
\caption{\textit{Values of $q(\delta)$, defined in \eqref{cond:f}, for the process \eqref{decreasing1} and different values of $h$, $\delta$ and $n$.}} \label{MSPEanalysisdecreasing1d}
\end{center}
\end{table}

\begin{table} \scriptsize
\begin{center}
	\input{tab1-m13-h1-i2.tex}
	\input{tab1-m13-h1-i3.tex}
	\input{tab1-m13-h5-i2.tex}
	\input{tab1-m13-h5-i3.tex} 
\caption{\textit{Proportion of \eqref{decisionrule}  being fulfilled for the process \eqref{decreasing1} and different values of $h$, $\delta$ and $n$.}}  \label{MSPEanalysisdecreasing1}
\end{center}
\end{table}

\clearpage

\begin{figure} 
\centering 
\includegraphics[width=0.24\textwidth]{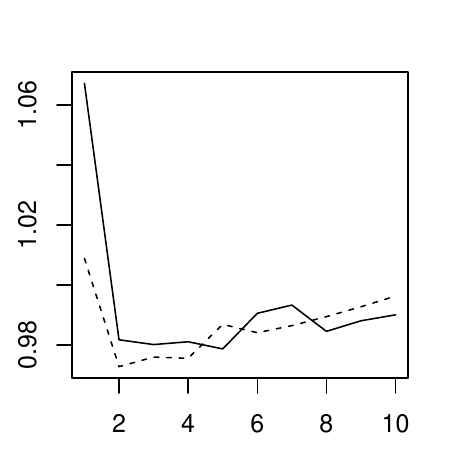}~
\includegraphics[width=0.24\textwidth]{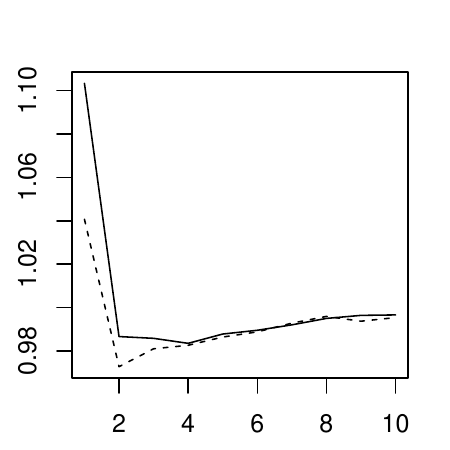}~
\includegraphics[width=0.24\textwidth]{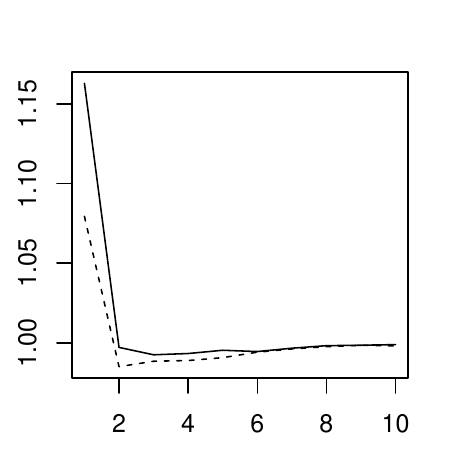}~ 
\includegraphics[width=0.24\textwidth]{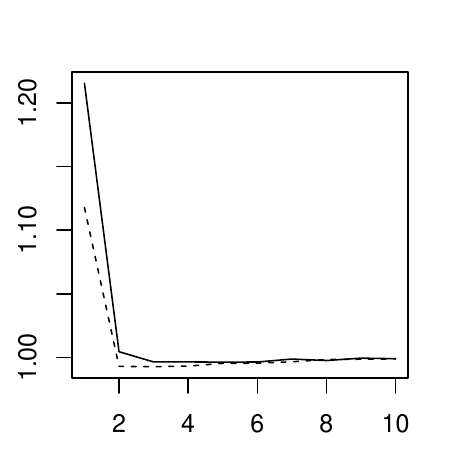} \\
\includegraphics[width=0.24\textwidth]{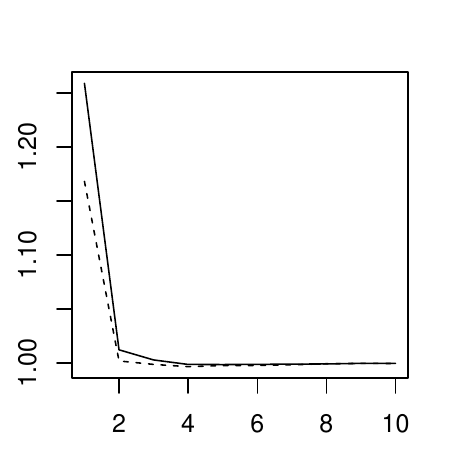}~
\includegraphics[width=0.24\textwidth]{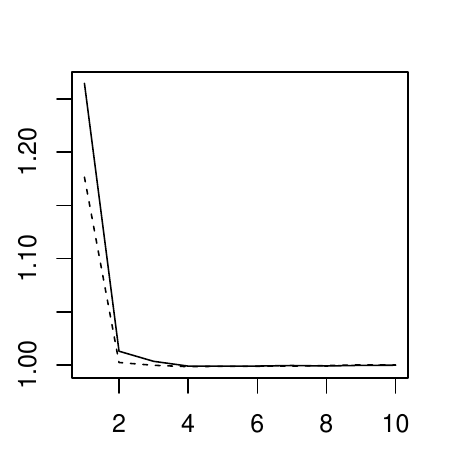}~ 
\includegraphics[width=0.24\textwidth]{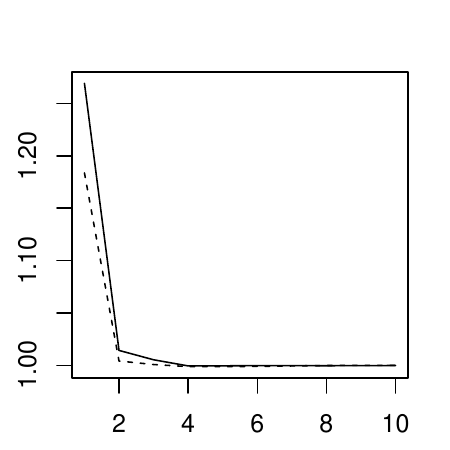}~ 
\includegraphics[width=0.24\textwidth]{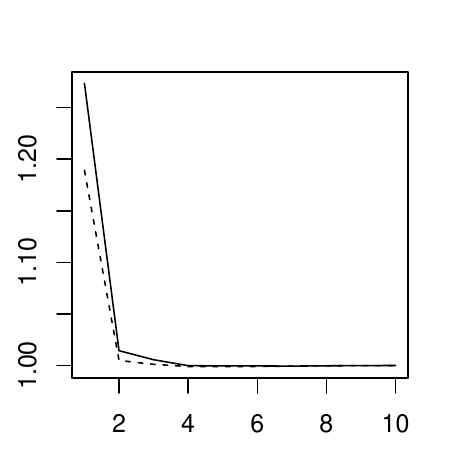}
   \caption{\it Plot of $h \mapsto R_{T,i}(h)$ for model \eqref{decreasing2} and different values of $n$ (from left to right: n=100, n=200, n=500, n=1000 [first row], n=4000, n=6000, n=8000, n=10000 [second row]). Solid line: $i=3$ (test set), dashed line: $i=2$ (validation set 2).} \label{MSPEdecreasing2}
\end{figure}

\begin{table} \scriptsize
\begin{center}
\begin{tabular}{|c||c|c|c||c|c|c|}
\hline
 \multirow{2}{*}{$n$} & & \multicolumn{2}{|c||}{$R^{{\rm s}, {\rm ls}}_{T,2}(1)$}  & & \multicolumn{2}{|c|}{$R^{{\rm s}, {\rm ls}}_{T,2}(5)$} \\ 
 & & $\geq 1.01$  & $< 1.01$  & & $\geq 1.01$  & $< 1.01$ \\ 
 \hline 
 \multirow{2}{*}{$100$} & $R^{{\rm s}, {\rm ls}}_{T,3}(1) \geq 1.01$  & 0.2399 & 0.2282  & $R^{{\rm s}, {\rm ls}}_{T,3}(5) \geq 1.01$  & 0.0657 & 0.1335 \\ 
& $R^{{\rm s}, {\rm ls}}_{T,3}(1) < 1.01$  & 0.1416 & 0.3903 & $R^{{\rm s}, {\rm ls}}_{T,3}(5) < 1.01$  & 0.1381 & 0.6627 \\

\hline\hline
 \multirow{2}{*}{$n$} & & \multicolumn{2}{|c||}{$R^{{\rm s}, {\rm ls}}_{T,2}(1)$}  & & \multicolumn{2}{|c|}{$R^{{\rm s}, {\rm ls}}_{T,2}(5)$} \\ 
 & & $\geq 1.15$  & $< 1.15$  & & $\geq 1.15$  & $< 1.15$ \\ 
 \hline 
 \multirow{2}{*}{$1000$} & $R^{{\rm s}, {\rm ls}}_{T,3}(1) \geq 1.15$  & 0.2464 & 0.4449  & $R^{{\rm s}, {\rm ls}}_{T,3}(5) \geq 1.15$  & 0 & 8e-04 \\ 
& $R^{{\rm s}, {\rm ls}}_{T,3}(1) < 1.15$  & 0.0872 & 0.2215 & $R^{{\rm s}, {\rm ls}}_{T,3}(5) < 1.15$  & 1e-04 & 0.9991 \\

\hline\hline
 \multirow{2}{*}{$n$} & & \multicolumn{2}{|c||}{$R^{{\rm s}, {\rm ls}}_{T,2}(1)$}  & & \multicolumn{2}{|c|}{$R^{{\rm s}, {\rm ls}}_{T,2}(5)$} \\ 
 & & $\geq 1$  & $< 1$  & & $\geq 1$  & $< 1$ \\ 
 \hline 
 \multirow{2}{*}{$10000$} & $R^{{\rm s}, {\rm ls}}_{T,3}(1) \geq 1$  & 1 & 0  & $R^{{\rm s}, {\rm ls}}_{T,3}(5) \geq 1$  & 0.4839 & 0.1483 \\ 
& $R^{{\rm s}, {\rm ls}}_{T,3}(1) < 1$  & 0 & 0 & $R^{{\rm s}, {\rm ls}}_{T,3}(5) < 1$  & 0.1228 & 0.245 \\

\hline\hline
 \multirow{2}{*}{$n$} & & \multicolumn{2}{|c||}{$R^{{\rm s}, {\rm ls}}_{T,2}(1)$}  & & \multicolumn{2}{|c|}{$R^{{\rm s}, {\rm ls}}_{T,2}(5)$} \\ 
 & & $\geq 1.01$  & $< 1.01$  & & $\geq 1.01$  & $< 1.01$ \\ 
 \hline 
 \multirow{2}{*}{$10000$} & $R^{{\rm s}, {\rm ls}}_{T,3}(1) \geq 1.01$  & 1 & 0  & $R^{{\rm s}, {\rm ls}}_{T,3}(5) \geq 1.01$  & 1e-04 & 0.0113 \\ 
& $R^{{\rm s}, {\rm ls}}_{T,3}(1) < 1.01$  & 0 & 0 & $R^{{\rm s}, {\rm ls}}_{T,3}(5) < 1.01$  & 0.0054 & 0.9832 \\

\hline\hline
 \multirow{2}{*}{$n$} & & \multicolumn{2}{|c||}{$R^{{\rm s}, {\rm ls}}_{T,2}(1)$}  & & \multicolumn{2}{|c|}{$R^{{\rm s}, {\rm ls}}_{T,2}(5)$} \\ 
 & & $\geq 1.1$  & $< 1.1$  & & $\geq 1.1$  & $< 1.1$ \\ 
 \hline 
 \multirow{2}{*}{$10000$} & $R^{{\rm s}, {\rm ls}}_{T,3}(1) \geq 1.1$  & 0.994 & 0.0059  & $R^{{\rm s}, {\rm ls}}_{T,3}(5) \geq 1.1$  & 0 & 0 \\ 
& $R^{{\rm s}, {\rm ls}}_{T,3}(1) < 1.1$  & 1e-04 & 0 & $R^{{\rm s}, {\rm ls}}_{T,3}(5) < 1.1$  & 0 & 1 \\

\hline\hline
 \multirow{2}{*}{$n$} & & \multicolumn{2}{|c||}{$R^{{\rm s}, {\rm ls}}_{T,2}(1)$}  & & \multicolumn{2}{|c|}{$R^{{\rm s}, {\rm ls}}_{T,2}(5)$} \\ 
 & & $\geq 1.2$  & $< 1.2$  & & $\geq 1.2$  & $< 1.2$ \\ 
 \hline 
 \multirow{2}{*}{$10000$} & $R^{{\rm s}, {\rm ls}}_{T,3}(1) \geq 1.2$  & 0.3563 & 0.5784  & $R^{{\rm s}, {\rm ls}}_{T,3}(5) \geq 1.2$  & 0 & 0 \\ 
& $R^{{\rm s}, {\rm ls}}_{T,3}(1) < 1.2$  & 0.0234 & 0.0419 & $R^{{\rm s}, {\rm ls}}_{T,3}(5) < 1.2$  & 0 & 1 \\

\hline
\end{tabular}
\caption{\textit{Proportions of the individual events in~\eqref{samedecision} for the process \eqref{decreasing2} and selected combinations of $n$ and $\delta$.}}  \label{MSPEanalysisdecreasing2c}
\end{center}
\end{table}

\begin{table} \scriptsize
\begin{center}
	\input{tab2-m14-h1.tex}
	\input{tab2-m14-h5.tex}  
\caption{\textit{Proportion of \eqref{samedecision}  being fulfilled for the process \eqref{decreasing2} and different values of $h$, $\delta$ and $n$.}}  \label{MSPEanalysisdecreasing2b}
\end{center}
\end{table}

\begin{table} \scriptsize
\begin{center}
	\input{tab5-m14-h1.tex}
	\input{tab5-m14-h5.tex}
\caption{\textit{Values of $q(\delta)$, defined in \eqref{cond:f}, for the process \eqref{decreasing2} and different values of $h$, $\delta$ and $n$.}} \label{MSPEanalysisdecreasing2d}
\end{center}
\end{table}

\begin{table} \scriptsize
\begin{center}
	\input{tab1-m14-h1-i2.tex}
	\input{tab1-m14-h1-i3.tex}
	\input{tab1-m14-h5-i2.tex}
	\input{tab1-m14-h5-i3.tex}
\caption{\textit{Proportion of \eqref{decisionrule}  being fulfilled for the process \eqref{decreasing2} and different values of $h$, $\delta$ and $n$.}}  \label{MSPEanalysisdecreasing2}
\end{center}
\end{table}

\end{document}